\documentclass[11pt]{amsart}
\usepackage{amscd}
\usepackage{amssymb}

\newtheorem{thm}{Theorem}[section]
\newtheorem{lem}[thm]{Lemma}

\newtheorem{cor}[thm]{Corollary}

\newtheorem{prop}[thm]{Proposition}

\theoremstyle{definition}

\renewcommand{\thecase}{}

\newtheorem{conj}[thm]{Conjecture}

\newtheorem{defn}[thm]{Definition}

\newtheorem{rmk}[thm]{Remark} 
\renewcommand{\thestep}{}

\theoremstyle{remark}

\makeatletter
\def\alphenumi{
  \def\theenumi{\alph{enumi}}
  \def\p@enumi{\theenumi}
  \def\labelenumi{(\@alph\c@enumi)}}
\makeatother

\makeatletter
\def\thecase{\@arabic\c@case}
\makeatother

\numberwithin{equation}{section}

\makeatletter
\def\thestep{\@arabic\c@step}
\makeatother

\renewcommand\emptyset{\varnothing}

\newcommand\embed{\hookrightarrow}

\newcommand\barM{{\bar{M}}}

\newcommand\barmu{{\bar\mu}}

\newcommand\ubarRR{{\underline{\mathbb{R}}}}

\newcommand\AAA{\mathbb{A}}

\newcommand\CC{\mathbb{C}}

\newcommand\EE{\mathbb{E}}
\newcommand\FF{\mathbb{F}}

\newcommand\HH{\mathbb{H}}

\newcommand\LL{\mathbb{L}}

\newcommand\NN{\mathbb{N}}

\newcommand\PP{\mathbb{P}}
\newcommand\QQ{\mathbb{Q}}
\newcommand\RR{\mathbb{R}}

\newcommand\VV{\mathbb{V}}
\newcommand\WW{\mathbb{W}}

\newcommand\ZZ{\mathbb{Z}}

\newcommand\bga{{\boldsymbol{\gamma}}}
\newcommand\bgamma{{\boldsymbol{\gamma}}}
\newcommand\bchi{{\boldsymbol{\chi}}}

\newcommand\bvarphi{{\boldsymbol{\varphi}}}

\newcommand\bA{{\mathbf{A}}}

\newcommand\bB{{\mathbf{B}}}
\newcommand\bc{{\mathbf{c}}}

\newcommand\bD{{\mathbf{D}}}

\newcommand\bE{{\mathbf{E}}}

\newcommand\bg{{\mathbf{g}}}

\newcommand\bL{{\mathbf{L}}}

\newcommand\bS{{\mathbf{S}}}

\newcommand\bV{{\mathbf{V}}}

\newcommand\bW{{\mathbf{W}}}
\newcommand\bx{{\mathbf{x}}}

\newcommand{\cov}{\nabla}

\newcommand{\rd}{\partial}

\newcommand\HG{{}_2 F_1}

\newcommand\thalf{{\textstyle{\frac{1}{2}}}}
\newcommand\tquarter{{\textstyle{\frac{1}{4}}}}

\newcommand\tthreequarter{{\textstyle{\frac{3}{4}}}}
\newcommand\tfivehalf{{\textstyle{\frac{5}{2}}}}

\newcommand\half{{{\frac{1}{2}}}}

\newcommand\quarter{{{\frac{1}{4}}}}
\newcommand\threehalf{{{\frac{3}{2}}}}

\newcommand\fb{{\mathfrak{b}}}

\newcommand\fg{{\mathfrak{g}}}

\newcommand\fs{{\mathfrak{s}}}

\newcommand\ft{{\mathfrak{t}}}

\newcommand\fV{{\mathfrak{V}}}

\newcommand\eps{\varepsilon}
\newcommand\ga{\gamma}

\newcommand\la{\lambda}
\newcommand\La{\Lambda}
\newcommand\ka{\kappa}
\newcommand\om{\omega}
\newcommand\Om{\Omega}
\newcommand\si{\sigma}
\newcommand\Si{\Sigma}

\newcommand\so{{\mathfrak{s}\mathfrak{o}}}
\newcommand\su{{\mathfrak{s}\mathfrak{u}}}
\newcommand\fu{{\mathfrak{u}}}

\newcommand\BS{\operatorname{BS}}

\newcommand\ES{\operatorname{ES}}

\newcommand\GL{\operatorname{GL}}

\newcommand\SO{\operatorname{SO}}

\newcommand\Spin{\operatorname{Spin}}
\newcommand\SU{\operatorname{SU}}
\newcommand\U{\operatorname{U}}

\newcommand\less{\setminus}

\newcommand{\8}{\infty}

\newcommand\ad{{\operatorname{ad}}}
\newcommand\Ad{{\operatorname{Ad}}}

\newcommand\asd{{\operatorname{asd}}}

\newcommand\Cl{\operatorname{C\ell}}
\newcommand\CCl{\operatorname{{\mathbb{C}\ell}}}

\newcommand\Coker{\operatorname{Coker}}

\newcommand\dist{\operatorname{dist}}

\newcommand\End{\operatorname{End}}

\newcommand\Fr{\operatorname{Fr}}

\newcommand\Gl{\operatorname{Gl}}

\newcommand\Hom{\operatorname{Hom}}

\newcommand\Ind{\operatorname{Index}}

\newcommand\Imag{\operatorname{Im}}
\newcommand\Isom{\operatorname{Isom}}

\newcommand\Ker{\operatorname{Ker}}

\newcommand\PD{\operatorname{PD}}

\newcommand\red{\operatorname{red}}

\newcommand\SW{SW}
\newcommand\Sym{\operatorname{Sym}}

\newcommand\Th{\operatorname{Th}}

\newcommand\Tr{\operatorname{Tr}}

\newcommand\vol{\operatorname{vol}}

\newcommand\cl{{\mathrm{cl}}}

\newcommand\even{{\mathrm{even}}}

\newcommand\id{{\mathrm{id}}}

\newcommand\sing{\mathrm{low}}
\newcommand\stab{\mathrm{vir}}

\newcommand\spinc{\text{$\text{spin}^c$ }}
\newcommand\spinu{\text{$\text{spin}^u$ }}
\newcommand\Spinc{\text{$\text{Spin}^c$}}
\newcommand\Spinu{\text{$\text{Spin}^u$}}

\newcommand\sA{{\mathcal{A}}}
\newcommand\sB{{\mathcal{B}}}
\newcommand\sC{{\mathcal{C}}}
\newcommand\sD{{\mathcal{D}}}

\newcommand\sG{{\mathcal{G}}}
\newcommand\sH{{\mathcal{H}}}
\newcommand\sI{{\mathcal{I}}}

\newcommand\sK{{\mathcal{K}}}

\newcommand\sM{{\mathcal{M}}}

\newcommand\sO{{\mathcal{O}}}
\newcommand\sP{{\mathcal{P}}}

\newcommand\sS{{\mathcal{S}}}

\newcommand\sU{{\mathcal{U}}}
\newcommand\sV{{\mathcal{V}}}
\newcommand\sW{{\mathcal{W}}}

\newcommand\tsC{{\tilde\sC}}

\newcommand\tM{{\tilde M}}
\newcommand\tN{{\tilde N}}

\begin{document}
\title[SO(3) Monopoles and Level-One Seiberg-Witten Moduli Spaces]
{SO(3) Monopoles, Level-One Seiberg-Witten Moduli Spaces, and Witten's
Conjecture in Low Degrees}
\author[Paul M. N. Feehan]{Paul M. N. Feehan}
\address{Department of Mathematics\\
Rutgers University\\
Piscataway, NJ 08854-8019}
\email{feehan@math.rutgers.edu}
\curraddr{Department of Mathematics\\
Trinity College\\
University of Dublin\\
Dublin 2\\
Ireland}
\email{feehan@maths.tcd.ie}
\urladdr{http://www.maths.tcd.ie/$\sim$feehan}
\author[Thomas G. Leness]{Thomas G. Leness}
\address{Department of Mathematics\\
Florida International University\\
Miami, FL 33199}
\email{lenesst@fiu.edu}
\urladdr{http://www.fiu.edu/$\sim$lenesst}
\dedicatory{}
\thanks{PMNF was supported in part by NSF grants
DMS 9704174 and, through the Institute for Advanced Study, by DMS 9729992.
TGL was supported in part by a Florida International University
Provost's Office Summer Research Grant.}
\date{This version: August 6, 2001. math.DG/0106238. Topology and its
Applications, to appear.}
\begin{abstract}
  We prove Witten's formula relating the Donaldson and Seiberg-Witten series
  modulo powers of degree $c+2$, with
  $c=-\frac{1}{4}(7\chi+11\sigma)$, for four-manifolds obeying some mild
  conditions, where $\chi$ and $\sigma$ are their Euler characteristic and
  signature. We use the moduli space of $\SO(3)$ monopoles as a cobordism
  between a link of the Donaldson moduli space of anti-self-dual $\SO(3)$
  connections and links of the moduli spaces of Seiberg-Witten
  monopoles. Gluing techniques allow us to compute contributions from
  Seiberg-Witten moduli spaces lying in the first (or `one-bubble') level
  of the Uhlenbeck compactification of the moduli space of $\SO(3)$
  monopoles.
\end{abstract}
\maketitle


\section{Introduction}
\label{sec:Introduction}
\subsection{Main results}
\label{subsec:MainResults}
In the present article we extend our results in \cite{FL2a, FL2b},
showing that Witten's conjecture \cite{Witten} relating the Donaldson and
Seiberg-Witten series holds in `low degrees' for a broad class of
four-manifolds. We apply our work on gluing $\SO(3)$ monopoles \cite{FL3,
FL5} --- restricting to the case of one `instanton
bubble' in this article --- to prove that these two series are equivalent
through a higher range of degrees than was possible in
\cite{FL2a, FL2b}. We assume throughout that $X$
is a closed, connected, smooth four-manifold with an orientation for which
$b_2^+(X)>0$. The Seiberg-Witten (SW) invariants (see
\S \ref{subsubsec:SWMonopoles}) comprise a function, $SW_X:\Spinc(X)\to\ZZ$,
where $\Spinc(X)$ is the set of isomorphism classes of \spinc structures on
$X$. For $w\in H^2(X;\ZZ)$, define the {\em Seiberg-Witten series\/} by
\begin{equation}
\label{eq:SWSeries}
\bS\bW_X^{w}(h)
=
\sum_{\fs \in \Spinc(X)}(-1)^{\half(w^{2}+c_{1}(\fs)\cdot w)}
SW_X(\fs)e^{\langle c_{1}(\fs),h\rangle},
\quad
h\in H_2(X;\RR),
\end{equation}
and let $\bD_X^w(h)$ denote the Donaldson series 
(see \cite[Theorem 1.7]{KMStructure} or \S \ref{subsec:Donaldsonseries} here). 
There is a map $c_1:\Spinc(X)\to
H^2(X;\ZZ)$ and the image of the support of $SW_X$ is the set $B$ of
SW-basic classes \cite{Witten}. A four-manifold $X$ has SW-simple type, when
$b_1(X)=0$, if $c_1(\fs)^2=2\chi+3\sigma$ for all $c_1(\fs)\in B$,
where $\chi$ and $\sigma$ are the Euler characteristic and signature of
$X$.  Let $B^\perp\subset H^2(X;\ZZ)$ denote the orthogonal complement
of $B$ with respect to the intersection form $Q_X$ on $H^2(X;\ZZ)$. Denote
$c(X)=-\frac{1}{4}(7\chi+11\sigma)$. Our main result is

\begin{thm}
\label{thm:WCL1}
Let $X$ be four-manifold with $b_1(X)=0$ and odd $b_2^+(X)\geq 3$. Assume
$X$ is abundant, SW-simple type, and effective. Then there exist $\La\in
B^\perp$ and $w\in H^2(X;\ZZ)$ for which $\La^2=4-(\chi+\si)$ and
$w-\Lambda\equiv w_2(X)\pmod{2}$.  For any such $\La$ and $w$, and any
$h\in H_2(X;\RR)$, one has
\begin{equation}
\label{eq:WCL1Equation}
\begin{aligned}
\bD^{w}_X(h) 
&\equiv 0\equiv \bS\bW^{w}_X(h)
\pmod{h^{c(X)-2}},
\\
\bD^{w}_X(h) 
&\equiv 
2^{2-c(X)}e^{\half h\cdot h}\bS\bW^{w}_X(h)
\pmod{h^{c(X)+2}}.
\end{aligned}
\end{equation}
\end{thm}

Noting that $c(X)=\chi_h(X)-c_1^2(X)$ (see \S \ref{subsec:Applications}),
Theorem \ref{thm:WCL1} can be used to compute Donaldson invariants for
four-manifolds in the region $c_1^2\leq \chi_h$ of the $(\chi_h,c_1^2)$
plane; all simply connected, compact, complex algebraic surfaces with odd
$b_2^+\geq 3$ are abundant (see below) and SW-simple type.

We shall explain below the terminology and notation in the statement of
Theorem \ref{thm:WCL1}.  Witten's conjecture \cite{Witten}
asserts that a four-manifold $X$ with $b_1(X)=0$ and odd $b_2^+(X)\geq 3$
has SW-simple type if and only if it has KM-simple type, that is,
simple type in the sense of Kronheimer and Mrowka (see Definition 1.4 in
\cite{KMStructure}), and that the SW-basic and KM-basic classes (see
Theorem 1.7 in \cite{KMStructure}) coincide; if $X$ has simple type, then
\begin{equation}
\label{eq:WConjecture}
\bD^{w}_X(h) 
=
2^{2-c(X)}e^{\frac{1}{2}Q_X(h,h)}\bS\bW^{w}_X(h),
\quad
h\in H_2(X;\RR).
\end{equation}
The quantum field theory argument giving equation
\eqref{eq:WConjecture} when $b_2^+(X)\ge 3$ has been extended by Moore and
Witten \cite{MooreWitten} to allow $b_2^+(X) \ge 1$, $b_1(X)\ge 0$, and
four-manifolds $X$ of non-simple type. Recall that $b_2^+(X)$ is the
dimension of a maximal positive-definite linear subspace $H^{2,+}(X;\RR)$
for the intersection pairing $Q_X$ on $H^2(X;\RR)$. 

In \cite{FL2a,FL2b} we proved that the second equation
\eqref{eq:WCL1Equation} holds modulo $h^{c(X)}$, with $\La^2=2-(\chi+\si)$
but otherwise identical hypotheses. If one desires a mod $h^\delta$
relation such as \eqref{eq:WCL1Equation} for larger values of $\delta$
relative to $c(X)$, one must allow more bubbles and the difficulty of
the calculations rapidly increases: see \cite[\S 1]{FL2b} for a more
detailed discussion. In particular, to prove Witten's conjecture
\eqref{eq:WConjecture} in full, we would need to prove that equation
\eqref{eq:WCL1Equation} holds modulo $h^{\delta}$ for all
$\delta\in\ZZ_{\geq 0}$: this is the goal of the remaining papers
\cite{FL3,FL4,FL5} in our series.

If the hypotheses concerning a class $\Lambda$ are omitted, one still
obtains a formula for Donaldson invariants in terms of Seiberg-Witten
invariants (see Theorem
\ref{thm:SumFormula}), but it has a more complicated structure and one
cannot immediately see if it has the shape \eqref{eq:WCL1Equation}
consistent with Witten's conjecture.

\begin{defn}
\label{defn:Abundant}
  \cite[p. 169]{FKLM} We say that a closed, oriented four-manifold $X$ is
  {\em abundant\/} if the restriction of the intersection form to
  $B^{\bot}$ contains a hyperbolic sublattice.
\end{defn}

The abundance condition is just a convenient way of formulating the weaker,
but more technical condition that one can find (for example) classes
$\Lambda_j\in B^\perp$ such that $\Lambda_j^2=2j-(\chi+\sigma)$, for
$j=1,2$: this is the only property of $Q_X|B^\perp$ which we use to prove
Theorem \ref{thm:WCL1}.  All compact, complex algebraic, simply-connected
surfaces with $b_2^+\geq 3$ are abundant \cite[Theorem A.1]{FL2a}.

We expect the second equation \eqref{eq:WCL1Equation} to hold modulo
$h^{\delta}$, for any $\delta$, without reference to whether a class
$\Lambda$ exists with the properties stated and without the constraint on
$w$. For an explanation of why the limitations of the present article
prevent us from relaxing the constraints on $w$ in Theorem
\ref{thm:WCL1}, see Remark \ref{rmk:ProblemWithW}.

In the present article and its companion \cite{FL2b} we prove Theorem
\ref{thm:WCL1} using the moduli space $\sM_{\ft}$ of $\SO(3)$
monopoles \cite{PTLocal} to provide a cobordism between the link of the
moduli space $M_\kappa^w$ of anti-self-dual connections and links of moduli
spaces of Seiberg-Witten moduli spaces, $M_{\fs}$, these moduli spaces
being (topologically) embedded in $\sM_{\ft}$. Pairing certain cohomology
classes with the links of these moduli spaces gives multiples of the
Donaldson and Seiberg-Witten invariants, respectively, and thus provides a
relation between the two types of invariants. Let $\bar\sM_{\ft}$ denote
the Uhlenbeck compactification (see Theorem \ref{thm:Compactness}) of
$\sM_{\ft}$ in the space of ideal $\SO(3)$ monopoles,
$\cup_{\ell=0}^\8(\sM_{\ft_\ell}\times\Sym^\ell(X))$.

\begin{defn}
\label{defn:Effective}
\cite[Definition 1.3]{FL2a}
  We say that a closed, oriented, smooth four-manifold $X$ with $b_1(X)=0$
  and $b_2^+(X)\geq 1$ is {\em effective\/} if $X$ satisfies Conjecture
  3.1 in \cite{FKLM}.  This conjecture asserts that for a moduli space of
  ideal Seiberg-Witten monopoles, $M_{\fs}\times\Sym^\ell(X)$,
  appearing in level $\ell\geq 0$ of $\bar\sM_{\ft}$, the intersection number
\begin{equation}
\label{eq:LinkIntersection0}
{}\#  \left(
\bar\sV(z)\cap \bar\sW^{\eta}\cap\bar\bL_{\ft,\fs}
\right)
\end{equation}
  is a multiple of the Seiberg-Witten invariant $SW_X(\fs)$.  Here,
  $\bar\sV(z)$ and $\bar\sW$ are the closures in $\bar\sM_{\ft}/S^1$ of
  geometric representatives of Donaldson-type cohomology classes on the top
  stratum of $\sM_{\ft}/S^1$ and $\bar\bL_{\ft,\fs}$ is the link of
  $(M_{\fs}\times\Sym^\ell(X))\cap\bar\sM_{\ft}/S^1$ in
  $\bar\sM_{\ft}/S^1$, and $\eta\geq 0$ is an integer for which
  $\deg(z)+2(\eta+1)=\dim\sM_{\ft}$. In particular, the intersection number
  \eqref{eq:LinkIntersection0} is zero when the Seiberg-Witten invariant
  $SW_X(\fs)$ is zero.
\end{defn}

When $b_1(X)>0$, one replaces the Seiberg-Witten invariants, $SW_X(\fs)$,
mentioned in Definition \ref{defn:Effective} with values of the
Seiberg-Witten function $SW_{X,\fs}$ --- see \cite[Equation (1.7) and \S
4.1]{FL2b} for definitions.
The motivation for Conjecture 3.1 of \cite{FKLM} and a more detailed 
explanation appears in \cite[\S 3.1]{FKLM} together with an explanation of
its role in the proofs of the main results of that paper; see also
\cite{FLGeorgia}. It is almost certainly true that
this conjecture holds for all four-manifolds, based on our work in
\cite{FL3, FL4, FL5}, and it is a simple
consequence of our Conjecture \ref{conj:PTConjecture} for the general form
of the pairings of Donaldson-type cohomology classes with links of ideal
Seiberg-Witten moduli spaces in $\bar\sM_{\ft}$. 

We verify Conjecture 3.1 in
\cite{FKLM} by direct calculation in the present article
(see Theorem \ref{thm:LevelOne} and Proposition \ref{prop:LevelOneBlowUp})
for Donaldson invariants defined by $M_\kappa^w\embed \sM_{\ft}$
and Seiberg-Witten moduli spaces $M_{\fs}$ embedded in the first level,
$\sM_{\ft_1}\times X$, while
in \cite{FL2a, FL2b} we verified the conjecture for Seiberg-Witten
moduli spaces $M_{\fs}$ contained in the top level, $\sM_{\ft}$, of the
Uhlenbeck compactification $\bar\sM_{\ft}$
(see Theorem 4.13 and Proposition 4.22 in \cite{FL2b}).
However, we strongly
expect the conjecture to hold for Seiberg-Witten moduli spaces
$M_{\fs}$ contained in any level of the compactification
$\bar\sM_{\ft}$.

Equation \eqref{eq:WCL1Equation} is a special case of a more general
formula for Donaldson invariants which we now describe; the hypotheses
still include an important restriction which guarantees that Seiberg-Witten
moduli spaces with non-trivial invariants do not lie in the second or lower
levels of the compactified $\SO(3)$-monopole moduli space.  For $\Lambda\in
H^{2}(X;\ZZ)$, define
\begin{equation}
\label{eq:PositiveDiracIndexFunction}
i(\La)=\Lambda^2 + c(X) + \chi + \sigma.
\end{equation}
If $S(X)\subset\Spinc(X)$ denotes the support of 
the function $\SW_X:\Spinc(X)\to\ZZ$, let
\begin{equation}
\label{eq:SWInTopLevelFunction}
r(\Lambda,c_1(\fs))
=
-(c_1(\fs)-\Lambda)^2 - \textstyle{\frac{3}{4}}(\chi+\sigma)
\quad\text{and}\quad
r(\Lambda)
=
\mathop{\min}\limits_{\fs\in S(X)} r(\Lambda,c_1(\fs)).
\end{equation}
See Remark 3.36 in \cite{FL2b} for a discussion of the significance of
$r(\Lambda,c_1(\fs))$ and $r(\Lambda)$, while the significance of
$i(\Lambda)$ is explained in \cite[\S 1.2 \& \S 4.6]{FL2b}; see also
\cite{FKLM}. For the statement of Theorem \ref{thm:Main} below, we refer
the reader to \S \ref{subsec:Donaldsonseries} for a definition of the
Donaldson invariants, $D^w_X(h^{\delta-2m}x^m)$. Recall that
$P^{a,b}_d(\zeta)$, the {\em Jacobi polynomial\/} \cite[\S
8.960]{GradshteynRyzhik6}, is defined by
\begin{equation}
\label{eq:DefineJacobiPolynomial}
P^{a,b}_d(\zeta) 
=
\frac{1}{2^d}
\sum_{v=0}^d \binom{d+a}{v}\binom{d+b}{d-v}(\zeta-1)^{d-v}(\zeta+1)^{v},
\quad \zeta \in \CC,
\end{equation}
just as in \cite[Equation (4.27)]{FL2b}; if $d=0$, then $P^{a,b}_{d}(0)=1$.

\begin{thm}
\label{thm:Main}
Let $X$ be a four-manifold with $b_1(X)=0$ and odd $b_2^+(X) \geq 1$. 
Assume $X$ is effective. Suppose
$\Lambda, w\in H^{2}(X;\ZZ)$ are classes such that $w-\Lambda\equiv
w_{2}(X) \pmod 2$ and, if $b_2^+(X)=1$, the class $w\pmod 2$ admits no
torsion integral lifts.  Suppose $h\in H_2(X;\RR)$, $x\in H_0(X;\ZZ)$
is the positive generator, and
$\delta$ is a non-negative integer. If $\delta < i(\Lambda)$ and
$\delta=r(\Lambda)+4$, then 
\begin{equation}
\label{eq:Main}
\begin{aligned}
D^w_X(h^{\delta-2m}x^m) 
&= 
2^{1-\frac{1}{4}i(\La)-\frac{3}{4}\delta}
(-1)^{m+\half(\sigma-w^2)}
\sum_{\substack{\fs\in\text{\em Spin}^c(X)\\  r(\Lambda,c_1(\fs))=\delta}}
(-1)^{\frac{1}{2}(w^{2}+c_{1}(\fs)\cdot(w-\Lambda))}
\\
&\quad
\times(-2)^{d_s(\fs)/2}P^{a-1,b}_{d_s(\fs)/2}(0)SW_X(\fs)
\langle c_1(\fs)-\La,h\rangle^{\delta-2m}
\\
&\quad+
2^{1-\frac{1}{4}i(\La)-\frac{3}{4}\delta}
(-1)^{m+\half(\sigma-w^2)}
\sum_{\substack{\fs\in\text{\em Spin}^c(X)\\  r(\Lambda,c_1(\fs))=\delta-4}}
(-1)^{\frac{1}{2}(w^{2}+c_{1}(\fs)\cdot(w-\Lambda))}
\\
\\
&\quad
\times(-2)^{d_s(\fs)/2}P^{a,b}_{d_s(\fs)/2}(0)SW_X(\fs)
\\
&\quad
\times\left(a_0\langle c_1(\fs)-\La,h\rangle^{\delta-2m}
+ b_0\langle c_1(\fs)-\La,h\rangle^{\delta-2m-1}\langle\La,h\rangle
\right.
\\
&\qquad
\left.
+ a_1\langle c_1(\fs)-\La,h\rangle^{\delta-2m-2}Q_X(h,h)
\right),
\end{aligned}
\end{equation}
where all terms on the right which would have a negative exponent are
omitted and the coefficients $a_0(c_1(\fs),\Lambda,\delta,m)$, $a_1(\delta,m)$,
and $b_0(\delta,m)$ are given by
\begin{equation}
\label{eq:a0a1b1}
\begin{aligned}
a_0&=
3(c_1(\fs)-\Lambda)^2+c_1^2(X)+2(c_1(\fs)-\Lambda)\cdot\Lambda+4\delta-12m,
\\
b_0&= 2(\delta-2m)\frac{P^{a-1,b+1}_{d_s(\fs)/2}(0)}{P^{a,b}_{d_s(\fs)/2}(0)},
\\
a_1&= 4\binom{\delta-2m}{2},
\end{aligned}
\end{equation}
with 
$$
a=- \textstyle{\frac{1}{2}}d_s(\fs)
+ \textstyle{\frac{1}{4}}(i(\Lambda)-\delta) 
\quad\text{and}\quad
b = -\textstyle{\frac{1}{2}}d_s(\fs)
-\textstyle{\frac{1}{4}}(\chi+\sigma), 
$$
and $P^{a,b}_{d_s(\fs)/2}(0)$ given by definition
\eqref{eq:DefineJacobiPolynomial}.  If $d_s(\fs)=0$, then
$P^{a,b}_{d_s(\fs)/2}(0)=1$.  If $b_2^+(X)=1$, then all invariants in
equation \eqref{eq:Main} are evaluated with respect to the chambers
determined by the same period point in the positive cone of $H^2(X;\RR)$.
\end{thm}

Theorem 1.2 in \cite{FL2b} provides a formula analogous to \eqref{eq:Main}
under similar hypotheses for $D^w_X(h^{\delta-2m}x^m)$, but for
$\delta=r(\Lambda)$, and a vanishing result for $\delta<r(\Lambda)$.

Specializing Theorem \ref{thm:Main} to case where $X$ has SW-simple
type and $\Lambda$ is in $B^\perp$ and recalling that in this situation one
has \cite[\S 4.6]{FL2b}
\begin{equation}
\label{eq:SimpleTypeBperpformR}
r(\Lambda)=-\Lambda^2+c(X)-(\chi+\sigma),
\end{equation}
yields

\begin{thm}
\label{thm:SumFormula}
Continue the hypotheses of Theorem \ref{thm:Main}. We further assume that
$X$ has $b_2^+(X)>1$, Seiberg-Witten simple type and that $\La\in
B^\perp$. Then  
\begin{equation}
\label{eq:LevelOneSum}
\begin{aligned}
D^w_X(h^{\delta-2m}x^m)
&=2^{-\half(c(X)+\delta)}(-1)^{m+\half(\sigma-w^2)}
\\
&\quad\times
\sum_{\fs\in\text{\em Spin}^c(X)}
(-1)^{\half(w^2+w\cdot c_1(\fs))}SW_X(\fs)
\\
&\quad\times\left(
a_0\langle c_1(\fs)-\La,h\rangle^{\delta-2m}
+
b_0\langle c_1(\fs)-\La,h\rangle^{\delta-2m-1}\langle \La,h\rangle
\right.
\\
&\qquad\left.
+a_1\langle c_1(\fs)-\La,h\rangle^{\delta-2m-2}Q_X(h,h)
\right),
\end{aligned}
\end{equation}
where the coefficients $a_0(\Lambda,\delta,m)$, $b_0(\delta,m)$, and
$a_1(\delta,m)$ are given by 
\begin{equation}
\label{eq:Simplea0a1b1}
\begin{aligned}
a_0&=4c_1^2(X)+\La^2+4\delta -12m,
\\
b_0&=2(\delta-2m),
\\
a_1&=4\binom{\delta-2m}{2}.
\end{aligned}
\end{equation}
\end{thm}

Theorem 1.4 in \cite{FL2b} provides a formula analogous to
\eqref{eq:LevelOneSum} under similar hypotheses for
$D^w_X(h^{\delta-2m}x^m)$, but for $\delta=r(\Lambda)$, and a vanishing
result for $\delta<r(\Lambda)$.

Our main result, Theorem \ref{thm:WCL1}, then follows from Theorem
\ref{thm:SumFormula} and \cite[Theorem 1.1]{FL2b}.

\subsection{Some applications to minimal surfaces of general type}
\label{subsec:Applications}
If $X$ is any closed, oriented, smooth four-manifold we may define (by
analogy with their values when $X$ is a complex surface)
\begin{equation}
\label{eq:Definec1Squared}
c_1^2(X) = 2\chi+3\sigma,
\end{equation}
and
\begin{equation}
\label{eq:DefineHolcEulerChar}
\chi_h(X) = \textstyle{\frac{1}{4}}(\chi+\sigma).
\end{equation}
Then, the topological invariant $c(X)$ acquires a more familiar
interpretation:  
\begin{equation}
\label{eq:HolcInterpcX}
c(X)
=
-\textstyle{\frac{1}{4}}(7\chi+11\sigma)
=
\chi_h(X)-c_1^2(X).
\end{equation}
Theorem \ref{thm:WCL1} has some immediate applications to surfaces of
general type: it can be used to compute previously unknown Donaldson
invariants.  If $X$ is a simply-connected, minimal surface of general type,
then, according to \cite[Theorem VII.1.1($\text{\em iv}'$)]{BPV}, these
surfaces obey the Noether inequality:
$$
\chi_h(X)  
\le 
\begin{cases}
\thalf c_1^2(X) + 3, &\text{if $c_1^2(X)$ is even},
\\
\thalf c_1^2(X) + \tfivehalf, &\text{if $c_1^2(X)$ is odd}.
\end{cases}
$$
Hence, Theorem \ref{thm:WCL1} computes previously unknown Donaldson
invariants $D_X^w(x^mh^{\delta-2m})$ for minimal surfaces of general type
for $\delta$ obeying $\delta\leq c(X)$ and
$$
\delta 
\le 
\begin{cases}
3 - \thalf c_1^2(X), &\text{if $c_1^2(X)$ is even},
\\
\tfivehalf -\thalf c_1^2(X), &\text{if $c_1^2(X)$ is odd}.
\end{cases}
$$
This gives non-trivial results for Donaldson invariants when $c_1^2(X)\le 6$.

\subsection{Discussion of the hypotheses of Theorems \protect{\ref{thm:WCL1}},
\protect{\ref{thm:Main}}, and \protect{\ref{thm:SumFormula}}}
To prove Theorem \ref{thm:Main} (and thus Theorems
\ref{thm:SumFormula} and \ref{thm:WCL1}), we employ the
compactified moduli space of
$\SO(3)$ monopoles, $\bar\sM_{\ft}/S^1$, as a cobordism between a link of the
moduli space of anti-self-dual connections, $M^w_{\ka}$, and the links of
moduli spaces of ideal Seiberg-Witten monopoles, $M_{\fs}\times\Sym^\ell(X)$.
Our application of
the cobordism method in this article requires that
\begin{enumerate}
\item
The codimension of $M^w_{\ka}$ in $\sM_{\ft}$,
given by twice the complex index of a Dirac operator, is positive
(used in Proposition 3.29 in \cite{FL2b}), and
\item
Only the top or first levels of the Uhlenbeck compactification
$\bar\sM_{\ft}$ can contain Seiberg-Witten moduli spaces
$M_{\fs}\times\Sym^\ell(X)$ (where $\ell=0,1$) with
non-trivial invariants
(used in \S \ref{subsec:LinkCobordism} to eliminate the more intractable
terms in the sum \eqref{eq:CobordismSum}).
\end{enumerate}
In the proof of Theorem \ref{thm:Main}, one has
$2\delta=\deg(z)=\dim M^w_{\ka}$ and the hypotheses $\delta< i(\La)$ and
$\delta=r(\La)+4$ ensure that Conditions (1) and (2) hold, respectively.

The hypotheses imply that the cobordism $\bar\sM_{\ft}/S^1$ yields an equality
between pairings with the link of $M_\kappa^w$ and a sum of pairings with
the links of $M_{\fs}\times\Sym^\ell(X)$ where $\ell=0,1$.
The same remarks apply to the hypotheses in
Theorem \ref{thm:SumFormula}.

The assumption that $b_1(X)=0$ could be relaxed by generalizing our
calculation of the Segre classes for the virtual normal bundle
$N_{\ft,\fs}\to M_{\fs}$ to the case $b_1(X)>0$. We computed those Segre
classes in \cite[Corollary 3.32]{FL2a}, \cite[Lemma 4.7]{FL2b} when
$b_1(X)>0$, under a technical assumption on $H^1(X;\RR)$, but we have not
used this generalization here as it would greatly complicate our main
formulae.

When $b_2^+(X)=1$, we assume that $w\pmod 2$ does not admit a torsion
integral lift in order to avoid complications in defining the chamber in
the positive cone of $H^2(X;\RR)$ with respect to which the Donaldson and
Seiberg-Witten invariants are computed.  See the comments at the end of \S
3.4.2 in \cite{FL2b} and before Lemma 4.1 in \cite{FL2b} for further
discussion.

The proof of Theorem \ref{thm:WCL1} requires one to choose
classes $\Lambda\in B^\perp$ with optimally prescribed even square in order
to obtain the indicated vanishing results for the Donaldson and
Seiberg-Witten series, as well as compute the first non-vanishing
terms. The hypothesis that $X$ is abundant guarantees that one can find
such classes, though such choices are also possible for some non-abundant
four-manifolds \cite{FKLM}. The constraints on the pair $w$, $\Lambda$ were
discussed in \S \ref{subsec:MainResults}.

\subsection{Role of the present article in the SO(3)-monopole program}
Since our series of articles on the $\SO(3)$-monopole approach to the proof
of Witten's conjecture is quite long, it is perhaps worth mentioning how
the present article fits into this program. We gave a broad (but now
slightly dated) review of the program in \cite{FLGeorgia}, while we provided
a more modern outline in the introductions to \cite{FL2a} and \cite{FL2b}.

We proved the basic transversality and compactness properties of the moduli
space of $\SO(3)$ monopoles in \cite{FL1, FeehanGenericMetric}. The local
gluing results for $\SO(3)$ 
monopoles --- focusing on analytical steps involved in construction of the
links of lower-level Seiberg-Witten moduli spaces --- are the subject of
\cite{FL3, FL4}, while the global gluing results --- focusing on
topological steps involved in construction of these links --- are
considered in \cite{FL5}. At the time of writing, work on \cite{FL4} and
\cite{FL5} is still in progress, though we believe we have surmounted most
of the difficult technicalities.

Turning to computational aspects of our work on Witten's conjecture, we
show in \cite{FL2a, FL2b} that, accepting some mild hypotheses,
the Donaldson series is equal to $2^{2-c(X)}e^{\frac{1}{2}Q_X}$ times the
Seiberg-Witten series, at least through terms $h^\delta$ of degree
$\delta\leq c(X)$ (compare Theorem \ref{thm:WCL1} here).  In our articles
\cite{FL2a, FL2b} we did not consider contributions from
Seiberg-Witten moduli spaces other than those in the top level, so we could
not compare higher-degree terms in those articles. In order to compare
terms of arbitrarily high degree we need to compute contributions
of the form \eqref{eq:LinkIntersection0}
from Seiberg-Witten moduli spaces in arbitrarily low levels: partial
computations in this vein are given in
\cite{FL5}. Though we only compute contributions from the first level in
the present article, thus extending the calculations of \cite{FL2a,
FL2b}, whenever possible we present the calculations in sufficient
generality that they apply to arbitrary levels: see \S \ref{sec:DualCohom} and
\S \ref{sec:comp} for details. Finally, the article \cite{FKLM} is an
application of the main conclusions of \cite{FL2b} and \cite{FL5}.

\subsection{Comparison with \protect{\cite{LenessWC}}}
We use our gluing theorems
\cite{FL3, FL4} to describe the topology of an open neighborhood in
$\bar\sM_{\ft}$ of the Seiberg-Witten `stratum'
$(M_{\fs}\times\Sym^\ell(X))\cap\bar\sM_{\ft}$, restricting our attention
in the present article to the case $\ell=1$. We then compute pairings with
(see Theorem \ref{thm:LevelOne}) the circle-quotient of the boundary of
this open neighborhood, the link $\bar\bL_{\ft,\fs}$, by generalizing
techniques developed in
\cite{Yang} and \cite{LenessWC} to compute wall-crossing contributions to
Donaldson invariants of four-manifolds with $b_2^+=1$ (see Lemma
\ref{lem:PairingsWithInstantonLink}) to compute contributions to
Donaldson invariants Seiberg-Witten strata in $\bar\sM_{\ft}$.
There are two novel features here:
\begin{enumerate}
\item
The moduli space of `reducibles', $M_\fs$, can be positive-dimensional.
\item
There is an obstruction to gluing arising from the cokernel of a
Dirac operator on a complex-rank eight Clifford module over $S^4$.
\end{enumerate}
In \cite{LenessWC, Yang} the strata of `reducibles' have the form
$([A]\times\Sym^\ell(X))\cap\bar M_\kappa^w(g_I)$, where $[A]$ is a single
point represented by a reducible anti-self-dual connection and
$M_\kappa^w(g_I)$ is the parametrized (by a path of metrics $g_I$
with $I=(-1,1)$) moduli
space of anti-self-dual connections on an $\SO(3)$ bundle over $X$ with
second Stiefel-Whitney class $w\pmod{2}$ and Pontrjagin number
$-\kappa/4$. However, as we have already computed
\cite{FL2a, FL2b} the Segre classes of a `virtual' or `stabilized' normal
bundle of the stratum $M_{\fs}\embed \sM_{\ft_\ell}$,
the problem that $M_\fs$ is positive-dimensional
can be addressed via an elementary formula,
Proposition \ref{prop:S1Localization}, relating the Thom class of an
equivariant normal bundle with the Segre classes.  (The authors are not
aware if this formula is already known --- the derivation is straightforward,
but it clarifies the relation between the equivariant localization
computations of Seiberg-Witten wall-crossing formula in \cite{CaoZhouWC}
and the Segre class computations of the same formula in \cite{LiLiu} and
\cite{OTWall}.)

The second novel feature is the presence of an obstruction to gluing,
described in \S \ref{subsubsec:InstantonObstruction},
arising from the cokernel of the twisted Dirac operator on $S^4$: this
produces a term in the formula \eqref{eq:LevelOne} with coefficient $b_0$,
which has no counterpart in the corresponding formulae of \cite{LenessWC}
and \cite{Yang}. 

\subsection{Outline of the computation of Seiberg-Witten link pairings}
\label{subsec:OutlineSWLinkPairing}
The central problem in the application of the cobordism method to a proof
of Witten's conjecture is to compute  the intersection number
\eqref{eq:LinkIntersection0}.
It seems worthwhile to outline the basic steps involved in this calculation
for the case $\ell=1$.
These steps comprise the proof of Theorem
\ref{thm:LevelOne}, the computational core of the current article.
In our outline we try to draw a distinction
between the technical points which are at the heart of the difficulty
underlying the cobordism approach to a proof of Witten's conjecture and
those which are more tractable.

\subsubsection{Defining the virtual link}
The domain of the gluing map (see Theorem \ref{thm:GluingThm})
is the space
\begin{equation}
\label{eq:GluingDomain0}
\bar\sM_{\ft,\fs}^{\stab}/S^1
=
N_{\ft_1,\fs}(\eps)\times_{S^1}\bar{\Gl}_{\ft_1}(\delta),
\end{equation}
where $N_{\ft_1,\fs}(\eps)\to M_{\fs}$ is a complex disk bundle (of radius
$\eps$) homeomorphic to a neighborhood of $M_{\fs}$ in a virtual moduli
space containing a neighborhood of $M_{\fs}$ in the moduli space
$\sM_{\ft_1}$ and $\bar{\Gl}_{\ft_1}(\delta)$ is a space of gluing
parameters, containing the moduli space of instantons on $S^4$ and the
bundle frames necessary to splice instantons onto $\SO(3)$ monopoles
represented by points in $N_{\ft_1,\fs}(\eps)$.  We define a {\em virtual
link\/}, $\bar\bL^{\stab}_{\ft,\fs}$, as the boundary of the domain
\eqref{eq:GluingDomain0}.

\subsubsection{Passage from an intersection product to a pairing of dual
cohomology classes with the fundamental class of the virtual link}
In \S \ref{subsec:Duality}, we prove that the intersection number
\eqref{eq:LinkIntersection0} is equal
to a pairing of a product of cohomology classes with a homology class
$[\bar\bL_{\ft,\fs}^{\stab}]$ which is, effectively, the fundamental class of
the virtual link.  We first show that the intersection number is equal to
the pairing of a relative cohomology class defined by the geometric
representatives and the relative Euler class of the obstruction bundle and
obstruction section with the relative fundamental class of a
codimension-zero subspace of the top stratum of the virtual link (see
\eqref{eq:Intersection1}).  The proofs given in \S
\ref{subsec:Duality} are specific to the topology of the case
$\ell=1$, but the results should hold for all $\ell>0$.

\subsubsection{Integration over the fiber}
We then consider a subspace,
$M_{\fs}\times\partial\bar{\Gl}_{\ft_1}(\delta)/S^1$,
of the virtual link, $\bar\bL_{\ft,\fs}^{\stab}$.
We use division by the Poincar\'e dual of the subspace 
$M_{\fs}\times\partial\bar{\Gl}_{\ft_1}(\delta)/S^1$ to relate pairings
with $[\bar\bL_{\ft,\fs}^{\stab}]$ to pairings with the
fundamental class of this subspace (see equation \eqref{eq:LinkPairingSegre}).
An analogous subspace appears in the virtual link for
when $\ell>1$ and the Poincar\'e dual division argument
will translate to the general case.

\subsubsection{Kunneth-type formula and Seiberg-Witten pairings}
Next we write pairings with
$[M_{\fs}\times\partial\bar{\Gl}_{\ft_1}(\delta)/S^1]$ as a sum of products
of pairings with $[M_{\fs}]$, yielding multiples of the Seiberg-Witten
invariant, $\SW_X(\fs)$, and pairings with
$[\partial\bar{\Gl}_{\ft_1}(\delta)/S^1]$; see equations
\eqref{eq:ProductPairing}, \eqref{eq:CohomClass2b}, and the identities that
follow. 

\subsubsection{Instanton link pairings}
The pairings with $[\partial\bar{\Gl}_{\ft_1}(\delta)/S^1]$ are almost
identical with those found in \cite{KotschickMorgan}. It is the calculation
of these pairings (and indeed the definition of the space
$\bar{\Gl}_{\ft_\ell}(\delta)$ for arbitrary levels $\ell>0$) which causes
by far the most difficulty.  In this article we use the computations of
\cite{LenessWC}, recorded here in equation \eqref{eq:InstLinkPairing}, to
compute these pairings; see equation \eqref{eq:Pairing3}.

\subsection{Extension of results from level-one case to higher levels} 
Although Theorem \ref{thm:WCL1} falls well short of a complete proof of
Witten's conjecture, it nonetheless provides further confirmation of the
promise of the $\SO(3)$-monopole approach. Our proof of Theorem
\ref{thm:WCL1} addresses many of the technical issues needed for our later
work on a proof of Conjecture \ref{conj:PTConjecture} and
hence Witten's conjecture. As explained in more detail in
\S \ref{subsec:KotschickMorgan}, our proof of Theorem
\ref{thm:WCL1} (in particular, Theorem \ref{thm:LevelOne}) 
illustrates that a crucial
step is to evaluate a certain `instanton link' pairing
(see \eqref{eq:KoMPairing2}), of the kind that
arise in the G\"ottsche-Kotschick-Morgan wall-crossing formula for
Donaldson invariants \cite{Goettsche, KotschickMorgan}.

\subsubsection{Level two}
The most immediate extension would involve level-two Seiberg-Witten moduli
spaces. The first main change would be to replace Lemma
\ref{lem:PairingsWithInstantonLink}, which records formulae from
\cite{LenessWC} for pairings with the boundary of the level-one gluing-data
bundle $\rd\bar{\Gl}_{\ft_1}(\delta)/S^1$, with the analogous results from
\cite{LenessWC} for the boundary of the level-two gluing-data
$\rd\bar{\Gl}_{\ft_2}(\delta)/S^1$. The second main change would be to
compute the Euler class of the instanton obstruction bundle for the
level-two case, generalizing the current Lemma \ref{lem:InstantonEuler}
which addresses the level-one case. For the Euler class calculation, the
description of the cokernel of the twisted Dirac operator
\cite[Lemma 3.3.28]{DK} might provide a useful starting point.
(Note that the description of this cokernel bundle in
\cite{AtiyahJones} does not suffice as we must work with an equivariant
extension over the gluing-data bundle.)
We would then need to replace in the proof of Theorem
\ref{thm:LevelOne} the resulting formulae for pairings with
$[\rd\bar{\Gl}_{\ft_1}(\delta)/S^1]$ by formulae for level-two pairings and
replace the Euler class of the instanton obstruction bundle with its
level-two counterpart. The end result would be versions of Theorems
\ref{thm:Main} and \ref{thm:SumFormula} valid for $\delta\leq r(\Lambda)+8$ and
an extension of Theorem \ref{thm:WCL1} to a mod $h^{c(X)+4}$ equivalence in
equation \eqref{eq:WCL1Equation} (using $\delta=c(X)+2$ when $\Lambda\in
B^\perp$ obeys $\Lambda^2=6-(\chi+\sigma)$, noting that $\delta$ must obey
the inequality $\delta<i(\Lambda)$, with $i(\Lambda)$ as given in
definition \eqref{eq:PositiveDiracIndexFunction}, and obey $\delta\equiv
c(X)+\Lambda^2 \pmod{4}$ to give non-zero Donaldson invariants.)  {}From
\cite[Theorem 1.3]{FKLM} we know that the Donaldson invariants 
$D^w_X(z)$ vanish, for suitable $w$ (for example, $w$ characteristic), when
$\deg(z)<c(X)-2$ (at least for abundant four-manifolds of
Seiberg-Witten simple type), so this calculation should give a direct
verification that Kronheimer-Mrowka simple type implies Seiberg-Witten
simple type for a non-trivial case with $\deg(z)\leq c(X)+2$:
$$
D^w_X(h^{c(X)-2}x^2)=4D^w_X(h^{c(X)-2}).
$$
Hence, this computation could be very useful for future attempts to
understand the relationship between Kronheimer-Mrowka and Seiberg-Witten
simple type. 

\subsubsection{Higher levels and comparison with the Kotschick-Morgan
conjecture} 
\label{subsec:KotschickMorgan}
When $b_2^+(X)=1$, the work of Kotschick and Morgan (see
\cite[Theorem 3.0.1]{KotschickMorgan} or \cite[Theorem 2.3]{Goettsche})
expresses the difference in Donaldson invariants of $X$, defined by two
different chambers $\sC_{\pm}$ in the same connected 
component of the positive cone of $H^2(X;\RR)$, as an alternating sum
\begin{align*}
{}&D_{X,\sC_+}^w(x^mh^{\delta-2m}) - D_{X,\sC_-}^w(x^mh^{\delta-2m})
\\
&\quad=
\sum_{\xi\in S(\sC_+,\sC_-,w,\delta)}
(-1)^{\eps(w,\xi,\delta)}
\delta_{X,\xi}^w(x^mh^{\delta-2m}),
\end{align*}
where the precise definition of the index set 
$S(\sC_+,\sC_-,w,\delta)\subset H^2(X;\ZZ)$
and sign $\eps(w,\xi,\delta)$ are not relevant to this discussion.
The classes $\xi\in S(\sC_+,\sC_-,w,\delta)$ correspond to points
$[A_\xi]\in M_{\kappa-\ell}^w(g_0)$ represented by reducible connections.
The terms $\delta_{X,\xi}^w(x^mh^{\delta-2m})$ in the preceding sum are
given by 
\cite[Definition 5.1.1]{KotschickMorgan}
\begin{equation}
\label{eq:KoMPairing}
\delta_{X,\xi}^w(x^mh^{\delta-2m})
=
\left\langle \bar\mu(x^mh^{\delta-2m}),
[\rd\left(B^n_{\CC}(\eps)\times_{S^1}
\bar{\Gl}_{\xi,\ell}(\delta)\right)]\right\rangle,
\end{equation}
where $B^n_{\CC}(\eps)\subset\CC^n$ is the ball 
of radius $\eps$ centered at the origin.  The space
$B^n_{\CC}(\eps)\times_{S^1}\bar{\Gl}_{\xi,\ell}(\delta)$ is homeomorphic
to a neighborhood of the reducible ideal anti-self-dual connections
$[A_\xi]\times\Sym^\ell(X)$ in a 
compactified moduli space $\bar M^w_\ka(g_I)$ of anti-self-dual ideal
connections parametrized by a path of metrics $g_I$.
A neighborhood of the reducible connection $[A_\xi]$ in the
parametrized moduli space $M^w_{\ka-\ell}(g_I)$ is homeomorphic to
$B^n_{\CC}(\eps)/S^1$.  The space of global gluing data,
$\bar{\Gl}_{\xi,\ell}(\delta)$, is (roughly) the union of the spaces of
stratum-wise gluing data, $\bar{\Gl}_{\xi,\ell}(\Si,\delta)$, as $\Sigma$
ranges over the smooth strata of $\Sym^\ell(X)$, patched together via
(non-canonical) transition maps obeying a cocycle condition. Naturally, the
case $\ell=1$ is simplest, as no transition map is needed, while if
$\ell=2$ the transition map does not need to satisfy a cocycle condition:
the real problem of constructing $\bar{\Gl}_{\xi,\ell}(\delta)$ arises when
$\ell\geq 3$.  An approach to constructing the space of global gluing data
$\bar{\Gl}_{\xi,\ell}(\delta)$ is indicated in \cite[Theorem
4.4.2]{KotschickMorgan}, although the
construction of the space of global gluing data
and the global gluing map --- in a form suitable for our purposes --- is a
difficult analytical and topological problem.

The terms $\delta_{X,\xi}^w(x^mh^{\delta-2m})$ are computed directly in
\cite{LenessWC} for the cases $\ell=1$ and $\ell=2$; we exploit the
computation when $\ell=1$ in the present article, where the result is
recorded as Lemma \ref{lem:PairingsWithInstantonLink}.
It is natural then to ask if one can use the
computations of $\delta_{X,\xi}^w(x^mh^{\delta-2m})$ in
\cite{Goettsche} to compute the intersection numbers
\eqref{eq:LinkIntersection0} with links in the moduli space
of $\SO(3)$ monopoles.

There are some differences worth noting
between the link of $[A_\xi]\times\Sym^\ell(X)$ in
$\barM_\kappa^w(g_I)$, 
$$
\bar\bL_{\kappa,\xi}^w
:=
\rd(B^n_{\CC}(\eps)\times_{S^1}\bar{\Gl}_{\xi,\ell}(\delta)),
$$ 
and the links $\bar\bL^{\stab}_{\ft,\fs}$ of $M_{\fs}\times\Sym^\ell(X)$ in
$\bar\sM_{\ft,\fs}^{\stab}/S^1$ or the links $\bar\bL_{\ft,\fs}$ of
$M_{\fs}\times\Sym^\ell(X)$ in $\bar\sM_{\ft}/S^1$.
A neighborhood of the point $[A_\xi]$
in the background moduli space $M_{\kappa-\ell}^w(g_0)$
is homeomorphic to $B^n_{\CC}(\eps)/S^1$,
while a neighborhood of the Seiberg-Witten moduli space 
(which can have positive dimension) in the background
virtual moduli space, $\sM_{\ft_1,\fs}^{\stab}/S^1$,
is homeomorphic to the $S^1$-quotient of the
complex disk bundle, $N_{\ft_1,\fs}(\eps)/S^1\to M_{\fs}$.  However, the
Poincar\'e-dual division argument mentioned in \S
\ref{subsec:OutlineSWLinkPairing} reduces both the pairing
\eqref{eq:KoMPairing} and the intersection number
\eqref{eq:LinkIntersection0} to a sum of terms of the form
\begin{equation}
\label{eq:KoMPairing2}
\left\langle \nu_{\ft_\ell}^n\smile \pi_X^*\alpha,[\rd\bar{\Gl}_{\ft_\ell}/S^1]
\right\rangle,
\end{equation}
where $\nu_{\ft_\ell}$ is the first Chern class of the circle bundle
$\rd\bar{\Gl}_{\ft_\ell}\to \rd\bar{\Gl}_{\ft_\ell}/S^1$, the map $\pi_X$
is a projection $\rd\bar{\Gl}_{\ft_\ell}/S^1\to\Sym^\ell(X)$, and $\alpha$
is a cohomology class on $\Sym^\ell(X)$.  (One has $\xi=\Lambda-c_1(\fs)$
in the proof of Theorem \ref{thm:LevelOne}; see equation
\eqref{eq:ProductPairing}.)  Thus, when comparing the Kotschick-Morgan
conjecture and conjectures about the intersection number
\eqref{eq:LinkIntersection0}, we will discuss pairings of the form
\eqref{eq:KoMPairing2}.

We note that the terms $\delta_{X,\xi}^w(x^mh^{\delta-2m})$ in
\cite{Goettsche} are computed under the assumption that the four-manifold
$X$ has $b^+_2(X)=1$, so the computational methods of \cite{Goettsche} do
not immediately apply to the pairings \eqref{eq:KoMPairing2}, where we
allow $b_2^+(X)\geq 1$.

While pairings with links of reducible $\SO(3)$ monopoles involve some
cohomology classes differing from those appearing in \eqref{eq:KoMPairing}
(specifically, a cohomology class \eqref{eq:DefineMuC1}
associated with the circle action on $\sM_{\ft}^{*,0}$), 
these pairings can still be expressed in
terms of pairings of the form \eqref{eq:KoMPairing2}.  Thus, no new
difficulties are posed by those additional cohomology classes.

A final difference between computation of the pairings in
\cite{KotschickMorgan} and the computation of pairings with
Seiberg-Witten links is due to the presence of obstructions to gluing
$\SO(3)$ monopoles.  The link of a family of $M_{\fs}\times\Sym^\ell(X)$ in
$\bar\sM_{\ft}/S^1$ is given not by the space $\bar\bL^{\stab}_{\ft,\fs}$,
the $\SO(3)$-monopole analogue of $\bar\bL_{\kappa,\xi}^w$, but rather by
the zero-locus of a section of an obstruction bundle over
$\bar\bL^{\stab}_{\ft,\fs}$.  However, in Lemmas
\ref{lem:EulerOfBackgroundObstruction} and
\ref{lem:InstantonEuler} we prove that for $\ell=1$,
the Euler class of this obstruction bundle can be
expressed in terms of $\nu_{\ft_\ell}$ and $\pi_X^*\La$.  We expect
a similar result will hold for $\ell>1$.  In addition,
we expect the arguments of \S \ref{subsec:Duality} representing
the zero-locus of the obstruction section as dual to extensions
of the Euler class of the obstruction bundle
to hold for $\ell>1$.  Thus, in spite of  the
obstruction to gluing, the pairings with the link of reducible $\SO(3)$
monopoles can still be understood in terms of the pairings
\eqref{eq:KoMPairing2}.

The crux of the proof of Witten's conjecture, via $\SO(3)$ monopoles, is
then to define the link components $\rd\bar{\Gl}_{\ft_\ell}(\delta)/S^1$
and compute the pairings \eqref{eq:KoMPairing2}: one
can see how they arise in the present article in the course of the proof of
Theorem \ref{thm:LevelOne}. 
Short of proving Witten's conjecture outright, one can aim
instead to prove that (by analogy with the conjecture of Kotschick and
Morgan, as phrased by G\"ottsche \cite[Remark 4.3]{Goettsche}):

\begin{conj}
\label{conj:PTConjecture}
Let $X$ be a closed, oriented, smooth four-manifold with
$b_1(X)=0$ and odd $b_2^+(X)> 1$. Suppose $w, \Lambda\in H^2(X;\ZZ)$ obey
$w-\Lambda \equiv w_2(X)\pmod{2}$. Let $\delta,m$ be non-negative integers
for which $m\leq [\delta/2]$ and $\delta\equiv
-w^2-\frac{3}{4}(\chi+\sigma)\pmod{4}$. 
Then let $\ft$ be the \spinu structure over $(X,g)$ with $c_1(\ft)=\Lambda$,
$p_1(\ft)$ determined by $\delta=-p_1(\ft)-\frac{3}{4}(\chi+\sigma)$,
and $w_2(\ft)\equiv w\pmod{2}$. (See \S \ref{subsec:ModuliSpace} here
or \S 2.1 and \S 2.2 in \cite{FL2a}.)
Suppose $\fs$ is a \spinc structure over $(X,g)$ with $\dim M_{\fs}=0$ and
$\ell\equiv\frac{1}{4}(\delta-r(\Lambda,c_1(\fs)))\geq 0$. 
Let $\bar\bL_{\ft,\fs}$ denote the link of $M_{\fs}\times\Sym^\ell(X)$ in
$\bar\sM_\ft/S^1$ and denote
$\eta=\frac{1}{4}(p_1(\ft)+\Lambda^2-\sigma)-1$. 
Let $x\in H_0(X;\ZZ)$ be the positive generator.
Then for any $h\in H_2(X;\RR)$, and $z=h^{\delta-2m}x^m$, one has
\begin{equation}
\label{eq:PTConjecture}
\begin{aligned}
{}&\#\left( \bar\sV(z)\cap\bar\sW^{\eta}\cap \bar\bL_{\ft,\fs}\right)
\\
&\quad=
\pm 2^{a(\chi,\sigma,\delta,m,\ell)}
SW_X(\fs)
\sum_{j=0}^{q}
p_{\delta,m,,\ell,j}
\left(\langle c_1(\fs),h\rangle,\langle\Lambda,h\rangle\right)
Q_X^{j}(h,h),
\end{aligned}
\end{equation}
where $q = \min(\ell,[\delta/2]-m)$, $a(\cdot)$ is a linear polynomial, and
$p_{\delta,m,\ell,j}$ 
is a degree-$(\delta-2m-2j)$ homogeneous polynomial in two
variables with coefficients which are degree-$(\ell-j)$ polynomials in
$2\chi\pm 3\sigma$, $(c_1(\fs)-\Lambda)^2$, $\Lambda^2$,
$(c_1(\fs)-\Lambda)\cdot c_1(\fs)$, and also depend on
$\delta,m,\ell$. 
\end{conj}

The assumption that $\dim M_{\fs}=0$ can be relaxed; the coefficients on
the right-hand side of equation \eqref{eq:PTConjecture} will then 
additionally depend on
$d_s(\fs)=\dim M_{\fs}$.  Conjecture \ref{conj:PTConjecture} is motivated
by our explicit calculations of the pairings \eqref{eq:PTConjecture} when
$\ell=0$ in \cite{FL2b} and $\ell=1$ here, as well as partial calculations when
$\ell \geq 2$. Our development of a proof of Conjecture
\ref{conj:PTConjecture} is work in progress \cite{FL5}.
Conjecture \ref{conj:PTConjecture} also implies the `multiplicity
conjecture',  namely Conjecture \ref{conj:Multiplicity},
on which our proof of the Mari\~no-Moore-Peradze conjecture
\cite{FKLM} relies.
Conjecture \ref{conj:PTConjecture} is of interest to us because of

\begin{conj}
\label{conj:PUKotschickMorganImpliesWitten}
Witten's formula \eqref{eq:WConjecture} is implied by Conjecture
\ref{conj:PTConjecture}.
\end{conj}

G\"ottsche's calculation of the wall-crossing formula for Donaldson invariants
\cite{Goettsche}, assuming the Kotschick-Morgan conjecture
\cite{KotschickMorgan}  and Fintushel and Stern's proof of the general
blow-up formula for Donaldson invariants lends weight to our expectation
that Conjecture \ref{conj:PUKotschickMorganImpliesWitten} holds.  Even though
there are many unknown universal coefficients in our conjectured formula
\eqref{eq:PTConjecture} for the pairings for arbitrary $\ell(\ft,\fs)\geq
0$, that formula still yields a qualitative version of Witten's
conjecture --- that Donaldson invariants are determined by Seiberg-Witten
invariants --- and implies that a formula \eqref{eq:WConjecture} of the
type predicted by Witten should exist.

\subsection{Guide to the article}
We now outline the contents of the remainder of this article and indicate
the principal steps in the proofs of our principal results.

In \S \ref{sec:prelim} we recall the definition of the $\SO(3)$ monopole
moduli space and its basic properties, describe the strata given by the
moduli spaces of anti-self-dual $\SO(3)$ connections and Seiberg-Witten
monopoles, and review the definitions of the Donaldson and Seiberg-Witten
invariants. We briefly describe the links of these strata 
(deferring a detailed account of the Seiberg-Witten links until 
\S \ref{sec:gluing}) and the resulting
$\SO(3)$-monopole cobordism formula relating the invariants.

In \S \ref{sec:gluing} we define the domain of the gluing map
\cite{FL3, FL4}, the topological model space which parametrizes
a neighborhood of the stratum $(M_{\fs}\times\Sym^\ell(X))\cap\bar\sM_{\ft}$
when $\ell=1$; the detailed definition occupies \S
\ref{subsec:CliffordModuleBases} through \S \ref{subsec:SplicingMap}.
We identify the obstruction bundle in \S \ref{subsec:ObstructionBundle}
and discuss the gluing theorem in \S \ref{subsec:GluingMap}.
We define the virtual link in \S \ref{subsec:DefnOfLink} as the boundary
of the domain of the gluing map; the actual link is given by the image of
the zero-locus of a section of the obstruction bundle under the gluing map.
In \S \ref{subsec:Orient}, we define an orientation for the virtual link
and compare this orientation to orientations previously defined in
\cite{FL2b}.

In \S \ref{sec:cohom} we identify and extend the universal cohomology
classes. Specifically,  in \S
\ref{subsec:ExtensionCohomClass} we compute the pullback of the cohomology
classes on $\sM_{\ft}$ (both Donaldson-type and those specific to the
$\SO(3)$-monopole moduli space) to the domain of the gluing map. In \S
\ref{subsec:EulerClassObstruction} we compute the Euler classes of the
instanton and Seiberg-Witten components of the obstruction bundle.

The proof that the intersection number in \eqref{eq:LinkIntersection0} can be
expressed cohomologically appears in \S \ref{subsec:Duality}.
The topological computations described in \S \ref{subsec:OutlineSWLinkPairing}
for the cohomology of the virtual link  appear in \S
\ref{subsec:PrelimComp}. 

We give the computations leading to Theorem \ref{thm:WCL1} in \S
\ref{sec:comp}.  The main technical result of the article is
Theorem \ref{thm:LevelOne} in which we compute the intersection
number appearing on the left-hand-side of \eqref{eq:PTConjecture}
for a level-one link.  Theorem \ref{thm:LevelOne}
follows from the results of
\S \ref{sec:DualCohom} and some algebraic
computations performed 
in \S \ref{subsec:Algebraic}.  Because the Donaldson and Seiberg-Witten
invariants are defined with the aid of the blow-up formula (to avoid
technical difficulties associated with flat (reducible) $\SO(3)$ or $\U(1)$
connections), in \S \ref{subsec:BlowUp} we prove a blow-up formula for
level-one links. The proofs of Theorems \ref{thm:Main},
\ref{thm:SumFormula}, and \ref{thm:WCL1} are then completed in \S
\ref{subsec:WCL1}.

\subsubsection*{Acknowledgments}
We are grateful for the generous support of the Institute for Advanced
 Study (Princeton), the Max Planck Institut f\"ur Mathematik (Bonn), and
 the Department of Mathematics at Columbia University, and Trinity College
 Dublin (half-year visit by Leness) where parts of this article were
 completed during visits to these institutions.


\section{Preliminaries}
\label{sec:prelim}
We begin in \S \ref{subsec:ModuliSpace} by recalling the definition of the
moduli space of $\SO(3)$ monopoles and its basic properties 
\cite{FL1, FeehanGenericMetric}. In \S \ref{subsec:ASDsingularities} we
describe 
the stratum of zero-section monopoles, or anti-self-dual connections. In \S
\ref{subsec:Reducibles} we discuss the strata of reducible, or
Seiberg-Witten monopoles, together with their `virtual' neighborhoods and
normal bundles.  In \S
\ref{subsec:Cohomology} we define the cohomology classes which will be
paired with the links of the anti-self-dual and Seiberg-Witten moduli
spaces.  In \S \ref{subsec:Donaldsonseries} we review the definition of the
Donaldson series.  Lastly, in \S \ref{subsec:LinkCobordism} we describe the
basic relation between the pairings with links of the anti-self-dual and
Seiberg-Witten moduli spaces provided by the $\SO(3)$-monopole cobordism.

\subsection{The moduli space of \boldmath{$\SO(3)$} monopoles}
\label{subsec:ModuliSpace}
Throughout this article, $(X,g)$ will denote a closed, connected, oriented,
smooth, Riemannian four-manifold.

\subsubsection{Clifford modules}
\label{subsubsec:SpincuStr}
Let $V$ be a Hermitian vector bundle over $(X,g)$ and let $\rho:T^*X\to
\End_\CC(V)$ be a real-linear map satisfying
\begin{equation}
\label{eq:CliffordMapDefn}
\rho(\alpha)^2 = -g(\alpha,\alpha)\id_{V}
\quad\text{and}\quad
\rho(\alpha)^\dagger = -\rho(\alpha),
\quad \alpha \in C^\8(T^*X).
\end{equation}
The map $\rho$ uniquely extends to an algebra homomorphism,
$\rho:\Lambda^{\bullet}(T^*X)\otimes_\RR\CC\to\End_\CC(V)$, and gives $V$
the structure of a Hermitian Clifford module for the complex Clifford
algebra $\CCl(T^*X)$.  There is a splitting $V=V^+\oplus V^-$, where
$V^\pm$ are the $\mp 1$ eigenspaces of $\rho(\vol)$. A unitary connection
$A$ on $V$ is {\em spin \/} if
\begin{equation}
\label{eq:SpinConnection}
[\nabla_A,\rho(\alpha)] =\rho(\nabla\alpha)
\quad\text{on }C^\8(V),
\end{equation}
for any $\alpha\in C^\8(T^*X)$, where $\nabla$ is the Levi-Civita connection.

A Hermitian Clifford module $\fs=(\rho,W)$ is a \spinc structure when $W$
has complex rank four; it defines a class
\begin{equation}
\label{eq:DefineChernClassOfSpinc}
c_1(\fs)=c_1(W^+),
\end{equation}
and every class in $H^2(X;\ZZ)$ lifting the second Stiefel-Whitney class,
$w_2(X)\in H^2(X;\ZZ/2\ZZ)$, arises this way.

We call a Hermitian Clifford module $\ft=(\rho,V)$ a \spinu structure when
$V$ has complex rank eight. Recall that $\fg_{V}\subset\su(V)$ is the
$\SO(3)$ subbundle given by the span of the sections of the bundle $\su(V)$
which commute with the action of $\CCl(T^*X)$ on $V$. We obtain a splitting
\begin{equation}
\label{eq:EndSplitting}
\su(V^+)
\cong 
\rho(\Lambda^+)\oplus i\rho(\Lambda^+)\otimes_\RR\fg_{V}
\oplus \fg_{V},
\end{equation}
and similarly for $\su(V^-)$. The fibers $V_x^+$ define complex lines whose
tensor-product square is $\det(V^+_x)$ and thus a complex line bundle over
$X$, 
\begin{equation}
\label{eq:CliffordDeterminantBundle}
{\det}^{\frac{1}{2}}(V^+).
\end{equation}
A \spinu structure $\ft$ thus defines classes,
\begin{equation}
\label{eq:SpinUCharacteristics}
c_1(\ft)=\textstyle{\frac{1}{2}} c_1(V^+),
\quad
p_1(\ft) = p_1(\fg_{V}),
\quad
\text{and}\quad 
w_2(\ft)=w_2(\fg_{V}).
\end{equation}
Given $W$, one has an isomorphism $V\cong W\otimes_\CC E$ of Hermitian
Clifford modules, where $E$ is a rank-two Hermitian vector bundle
\cite[Lemma 2.3]{FL2a}; then 
$$
\fg_V
=
\su(E) 
\quad\text{and}\quad
{\det}^{\frac{1}{2}}(V^+) 
=
\det(W^+)\otimes_\CC\det(E). 
$$

\subsubsection{$\SO(3)$ monopoles}
\label{subsubsec:PU2Monopoles}
We fix a smooth unitary connection $A_\La$ on the line bundle
$\det^{\frac{1}{2}}(V^+)$, let $k\geq 2$ be an integer, and let $\sA_{\ft}$
be the affine space of $L^2_k$ spin connections on $V$ which induce the
connection $2A_\La$ on $\det(V^+)$. If $A$ is a spin connection on $V$ then
it defines an $\SO(3)$ connection $\hat A$ on the subbundle
$\fg_V\subset\su(V)$
\cite[Lemma 2.5]{FL2a}; conversely, every $\SO(3)$ connection on $\fg_V$
lifts to a unique spin connection on $V$ inducing the connection
$2A_\Lambda$ on $\det(V^+)$ \cite[Lemma 2.11]{FL2a}.

Let $\sG_{\ft}$ denote the group of $L^2_{k+1}$ unitary automorphisms of
$V$ which commute with $\CCl(T^*X)$ and which have Clifford-determinant one
(see \cite[Definition 2.6]{FL2a}). Define
\begin{equation}
\label{eq:SpinUPreConfiguration}
\tsC_\ft(X) = \sA_\ft(X)\times L^2_k(X,V^+)
\quad\text{and}\quad
\sC_\ft = \tsC_\ft/\sG_\ft.
\end{equation}
The action of $\sG_{\ft}$ on $V$ induces an adjoint action on
$\End_\CC(V)$, acting as the identity on
$\rho(\Lambda^\bullet_\CC)\subset\End_\CC(V)$ and inducing an adjoint
action on $\fg_{V}\subset\End_\CC(V)$ (see
\cite[Lemma 2.7]{FL2a}). The space $\tsC_\ft$ and hence $\sC_\ft$ carry
circle actions induced by scalar multiplication on $V$:
\begin{equation}
\label{eq:S1ZAction}
S^1\times V \to V, 
\quad (e^{i\theta},\Phi)\mapsto e^{i\theta}\Phi.
\end{equation}
Because this action commutes with that of $\sG_{\ft}$, the
action \eqref{eq:S1ZAction} also defines an action on $\sC_{\ft}$.
Note that $-1\in S^1$ acts trivially on $\sC_{\ft}$.  Let
$\sC^0_\ft\subset\sC_\ft$ be the subspace represented by pairs whose spinor
components are not identically zero, let $\sC^*_\ft\subset\sC_\ft$ be the
subspace represented by pairs where the induced $\SO(3)$ connections on
$\fg_{V}$ are irreducible, and let $\sC^{*,0}_\ft=\sC^*_\ft\cap
\sC^0_\ft$.

We call a pair $(A,\Phi)$ in $\tsC_\ft$ a {\em $\SO(3)$ monopole\/} if 
\begin{equation}
\label{eq:UnperturbedMonopoles}
\begin{aligned}
\ad^{-1}(F^+_{\hat A}) - \tau\rho^{-1}(\Phi\otimes\Phi^*)_{00} 
&=0,
\\
D_A\Phi + \rho(\vartheta)
&=0.
\end{aligned}
\end{equation}
Here, $D_A=\rho\circ \cov_A:C^\8(X,V^+)\to C^\8(X,V^-)$ is the Dirac
operator; the isomorphism $\ad:\fg_{V}\to \so(\fg_{V})$ identifies the
self-dual curvature $F_{\hat A}^+$, a section of
$\La^+\otimes\so(\fg_{V})$, with $\ad^{-1}(F^+_{\hat A})$, a section of
$\La^+\otimes\fg_{V}$; the section $\tau$ of $\GL(\Lambda^+)$ is a
perturbation close to the identity; the perturbation $\vartheta$ is a
complex one-form close to zero; $\Phi^* \in
\Hom(V^+,\CC)$ is the pointwise Hermitian dual $\langle\cdot,\Phi\rangle$
of $\Phi$; and $(\Phi\otimes\Phi^*)_{00}$ is the component of the section
$\Phi\otimes\Phi^*$ of $i\fu(V^+)$ lying in $\rho(\Lambda^+)\otimes\fg_{V}$
with respect to the splitting $\fu(V^+)=i\underline{\RR}\oplus\su(V^+)$ and 
decomposition \eqref{eq:EndSplitting} of $\su(V^+)$.

Equation \eqref{eq:UnperturbedMonopoles} is invariant under the action of
$\sG_\ft$. We let
$\sM_{\ft}\subset \sC_\ft$ be the subspace represented by pairs satisfying
equation
\eqref{eq:UnperturbedMonopoles} and write
\begin{equation}
\label{eq:PU2MonopoleSubspaces}
\sM^*_{\ft} = \sM_{\ft}\cap\sC^*_\ft, \quad
\sM^0_{\ft}= \sM_{\ft}\cap\sC^0_\ft, \quad\text{and}\quad
\sM^{*,0}_\ft = \sM_{\ft}\cap\sC^{*,0}_\ft.
\end{equation}
Since equation \eqref{eq:UnperturbedMonopoles} is invariant under the
circle action induced by scalar multiplication on $V$, the subspaces
\eqref{eq:PU2MonopoleSubspaces} of $\sC_{\ft}$ are also invariant under
this action.

\begin{thm}
\label{thm:Transv}
\cite[Theorem 1.1]{FeehanGenericMetric}, \cite{TelemanGenericMetric}
Let $X$ be a closed, oriented, smooth four-manifold and let $V$ be a
complex rank-eight, Hermitian vector bundle over $X$. Then for parameters
$(\rho,g,\tau,\vartheta)$, which are generic in the sense of
\cite{FeehanGenericMetric}, and $\ft=(\rho,V)$, the space
$\sM^{*,0}_{\ft}$ is a smooth manifold of the expected dimension,
\begin{equation}
\label{eq:Transv}
\begin{aligned}
\dim \sM^{*,0}_{\ft}
= d_a(\ft)+2n_a(\ft),
\quad\text{where }
d_a(\ft)
&=-2p_1(\ft)- \textstyle{\frac{3}{2}}(\chi+\sigma),
\\
n_a(\ft)
&=\textstyle{\frac{1}{4}}(p_1(\ft)+c_1(\ft)^2-\sigma).
\end{aligned}
\end{equation}
\end{thm}

For the remainder of the article, we assume that the perturbation
parameters in \eqref{eq:UnperturbedMonopoles} are chosen as indicated in
Theorem \ref{thm:Transv}. For each non-negative integer $\ell$, let
$\ft_\ell=(\rho,V_{\ell})$, where
$$
c_1(V_\ell)=c_1(V),
\quad
p_1(\fg_{V_\ell})=p_1(\fg_{V})+4\ell,
\quad\text{and}\quad
w_2(\fg_{V_\ell})=w_2(\fg_{V}).
$$
We let $\bar\sM_{\ft}$ denote the closure of $\sM_{\ft}$ in the space of
ideal monopoles, 
\begin{equation}
\label{eq:idealmonopoles}
I\sM_{\ft} = \bigsqcup_{\ell=0}^\8(\sM_{\ft_{\ell}}\times\Sym^\ell(X)),
\end{equation}
with respect to an Uhlenbeck topology
\cite[Definition 4.19]{FL1} and call the intersection of $\bar\sM_{\ft}$
with $\sM_{\ft_{\ell}}\times \Sym^\ell(X)$ its {\em $\ell$-th level\/}.

\begin{thm}
\label{thm:Compactness}
\cite[Theorem 1.1]{FL1}
Let $X$ be a Riemannian four-manifold with \spinu structure $\ft$.  Then
there is a positive integer $N$, depending at most on the curvature of the
chosen unitary connection on $\det(V^+)$ together with $p_1(\ft)$, such
that the Uhlenbeck closure $\bar\sM_{\ft}$ of $\sM_{\ft}$ in 
$\sqcup_{\ell=0}^N(\sM_{\ft_{\ell}}\times\Sym^\ell(X))$
is a second countable, compact, Hausdorff space. The space $\bar\sM_{\ft}$
carries a continuous circle action,
which restricts to the circle action defined on $\sM_{\ft_\ell}$ on each level.
\end{thm}

\subsection{Stratum of anti-self-dual or zero-section solutions}
\label{subsec:ASDsingularities}
{}From equation \eqref{eq:UnperturbedMonopoles}, we see that the stratum of
$\sM_{\ft}$ represented by pairs with zero spinor is identified with 
$$
\{A\in\sA_{\ft}: F_{\hat A}^+ = 0\}/\sG_{\ft} \cong M_\kappa^w(X,g),
$$
the moduli space of $g$-anti-self-dual connections on the $\SO(3)$ bundle
$\fg_{V}$, where $\ka=-\quarter p_1(\ft)$ and $w\equiv w_2(\ft)\pmod 2$.
For a generic Riemannian metric $g$, the space $M_\kappa^w(X,g)$ is a
smooth manifold of the expected dimension $-2p_1(\ft) -
\threehalf(\chi+\sigma)=d_a(\ft)$.

As explained in \cite[\S 3.4.1]{FL2a}, it is desirable to choose
$w\pmod{2}$ so as to exclude points in $\bar\sM_{\ft}$ with associated flat
$\SO(3)$ connections, so we have a {\em disjoint\/} union,
\begin{equation}
\label{eq:StratificationCptPU(2)Space}
\bar\sM_{\ft} 
\cong 
\bar\sM_{\ft}^{*,0} \sqcup \bar M_\kappa^w \sqcup \bar\sM_{\ft}^{\red},
\end{equation}
where $\bar\sM_{\ft}^*\subset\bar\sM_{\ft}$ is the subspace represented by
triples whose associated $\SO(3)$ connections are irreducible,
$\bar\sM_{\ft}^0\subset\bar\sM_{\ft}$ is the subspace represented by
triples whose spinors are not identically zero,
$\bar\sM_{\ft}^{*,0} = \bar\sM_{\ft}^{*}\cap\bar\sM_{\ft}^{0}$, 
while $\bar\sM_{\ft}^{\red}\subset\bar\sM_{\ft}$
is the subspace $\bar\sM_{\ft}-\bar\sM_{\ft}^*$
represented by triples whose associated $\SO(3)$
connections are reducible. We recall the 

\begin{defn}
\label{defn:Good}
\cite[Definition 3.20]{FL2b} 
A class $v\in H^2(X;\ZZ/2)$ is {\em good\/} if no integral lift of $v$ is
torsion.
\end{defn}

If $w\pmod{2}$ is good, then the union
\eqref{eq:StratificationCptPU(2)Space} is disjoint, as desired. In
practice, rather than constraining $w\pmod{2}$ itself, we replace $X$ by
the blow-up $X\#\overline{\CC\PP}^2$ and $w$ by $w+\PD[e]$ (where $e\in
H_2(X;\ZZ)$ is the exceptional class), noting that $w+\PD[e]\pmod{2}$ is
always good, and define gauge-theoretic invariants of $X$ in terms of
moduli spaces on $X\#\overline{\CC\PP}^2$. When $w\pmod{2}$ is good, we
define \cite[Definition 3.7]{FL2a} the link of $\barM^w_{\ka}$ in
$\bar\sM_{\ft}/S^1$ by
\begin{equation}
\label{eq:DefineASDLink}
\bL^{w}_{\ft,\ka}
=
\{[A,\Phi,\bx]\in\bar\sM_{\ft}/S^1: \|\Phi\|_{L^2}^2=\eps\},
\end{equation}
where $\eps$ is a small positive constant; for generic $\eps$, the link
$\bL^w_{\ft,\ka}$ is a smoothly-stratified, codimension-one subspace of
$\bar\sM_{\ft}/S^1$.

\subsection{Strata of Seiberg-Witten or reducible solutions}
\label{subsec:Reducibles}
We call a pair $(A,\Phi)\in\tsC_\ft$ {\em reducible\/} if the connection
$A$ on $V$ respects a splitting,
\begin{equation}
\label{eq:BasicSplitting}
V = W\oplus W\otimes L = W\otimes (\underline{\CC}\oplus L),
\end{equation}
for some \spinc structure $\fs=(\rho,W)$ and complex line bundle $L$, in
which case $c_1(L)=c_1(\ft)-c_1(\fs)$.  
A spin connection $A$ on $V$ is reducible with respect to the splitting
\eqref{eq:BasicSplitting} if and only if $\hat A$ is reducible
with respect to the splitting $\fg_{V}\cong \underline{\RR}\oplus L$,
\cite[Lemma 2.9]{FL2a}. If $A$ is reducible, we can write $A=B\oplus
B\otimes A_L$, where $B$ is a spin connection on $W$ and $A_L$ is a unitary
connection on $L$; then $\hat A= d_{\RR}\oplus A_L$ and $A_L = A_\La\otimes
(B^{\det})^*$, where $B^{\det}$ is the connection on $\det(W^+)$ induced by
$B$ on $W$.

\subsubsection{Seiberg-Witten monopoles} 
\label{subsubsec:SWMonopoles}
Given a \spinc structure $\fs=(\rho,W)$ on $X$, let $\sA_{\fs}$ denote the
affine space of $L^2_k$ spin connections on $W$. Let $\sG_{\fs}$
denote the group of $L^2_{k+1}$ unitary automorphisms of $W$, commuting with
$\CCl(T^*X)$, which we identify with $L^2_{k+1}(X,S^1)$. We then define
\begin{equation}
\label{eq:SpincPreConfig}
\tsC_{\fs}= \sA_{\fs}\times L^2_k(W^+)
\quad\text{and}\quad
\sC_{\fs}=\tsC_{\fs}/\sG_{\fs},
\end{equation}
where $\sG_{\fs}$ acts on $\tsC_{\fs}$ by
\begin{equation}
\label{eq:SWGaugeGroupAction}
(s,(B,\Psi))
\mapsto
s(B,\Psi)
=
(B- (s^{-1}ds)\id_W, s\Psi).
\end{equation}
We call a pair $(B,\Psi)\in\tsC_{\fs}$ a Seiberg-Witten monopole if
\begin{equation}
\label{eq:SeibergWitten}
\begin{aligned}
\Tr(F^+_B) - \tau\rho^{-1}(\Psi\otimes\Psi^*)_{0} - F^+(A_{\Lambda}) 
&=0,
\\
D_B\Psi + \rho(\vartheta)\Psi
&=0,
\end{aligned}
\end{equation}
where $\Tr:\fu(W^+)\to i\RR$ is defined by the trace on $2\times 2$ complex
matrices, $(\Psi\otimes\Psi^*)_0$ is the component of the section
$\Psi\otimes\Psi^*$ of $i\fu(W^+)$ contained in $i\su(W^+)$, 
$D_B:C^\8(W^+)\to C^\8(W^-)$ is the Dirac operator, and $A_\Lambda$ is a
unitary connection on a line bundle with first Chern class $\Lambda\in
H^2(X;\ZZ)$.  The perturbations are chosen so that solutions to equation
\eqref{eq:SeibergWitten} are identified with reducible solutions to
\eqref{eq:UnperturbedMonopoles} when $c_1(\ft)=\Lambda$.
Let $\tM_{\fs}\subset\tsC_{\fs}$ be the subspace cut out by equation
\eqref{eq:SeibergWitten} and denote the moduli space of Seiberg-Witten
monopoles by $M_{\fs}=\tM_{\fs}/\sG_{\fs}$.

\subsubsection{Seiberg-Witten invariants}
\label{subsubsec:SWInvariants}
We let $\sC_{\fs}^0\subset\sC_{\fs}$ be the open subspace represented by
pairs whose spinor components which are not identically zero and
define a complex line bundle over $\sC_\fs^0\times X$ by
\begin{equation}
\label{eq:DefineSWUniversal}
\LL_{\fs}= \tsC_\fs^0\times_{\sG_\fs}\underline{\CC},
\end{equation}
where $\underline{\CC}=X\times\CC$ and
$s\in \sG_{\fs}$ acts on $(B,\Psi)\in \tsC_{\fs}$ and $(x,\zeta)\in
\underline{\CC}$ by
\begin{equation}
\label{eq:DefineSWUniversalS1Action}
((B,\Psi),(x,\zeta))
\mapsto
(s(B,\Psi),(x,s(x)^{-1}\zeta)).
\end{equation}
If $x\in H_0(X;\ZZ)$ denotes the positive generator, we set
\begin{equation}
\label{eq:SWClass}
\mu_{\fs} = c_1(\LL_{\fs})/x \in H^2(\sC_\fs^0;\ZZ).
\end{equation}
Equivalently, $\mu_{\fs}$ is the first Chern class of the $S^1$ base-point
fibration over $\sC_\fs^0$. If $b_1(X)=0$, then $c_1(\LL_{\fs})=\mu_{\fs}\times
1$ by \cite[Lemma 2.14]{FL2a}.

For $b_2^+(X)>0$ and generic Riemannian metrics on $X$, the space $M_{\fs}$
contains zero-section pairs if and only if $c_1(\fs)-\Lambda$ is a torsion
class by \cite[Proposition 6.3.1]{MorganSWNotes}. If
$M_{\fs}$ contains no zero-section pairs then, for generic perturbations,
it is a compact, oriented, smooth manifold of dimension
\begin{equation}
\label{eq:DimSW}
d_s(\fs)
= 
\dim M_{\fs}
=
\tquarter(c_1(\fs)^2 -2\chi -3\sigma).
\end{equation}
Let $\tilde X=X\#\overline{\CC\PP}^2$ denote the blow-up of $X$ with
exceptional class $e\in H_2(\tilde X;\ZZ)$ and denote its Poincar\'e dual
by $\PD[e]\in H^2(\tilde X;\ZZ)$. Let $\fs^\pm=(\tilde\rho,\tilde W)$
denote the \spinc structure on $\tilde X$ with $c_1(\fs^\pm)=c_1(\fs)\pm
\PD[e]$ obtained by splicing the \spinc structure $\fs=(\rho,W)$ on $X$
with the
\spinc structure on $\overline{\CC\PP}^2$ with first Chern class $\pm \PD[e]$.
(See \cite[\S 4.5]{FL2b} for an explanation of the relation between \spinc
structures on $X$ and $\tilde X$.) Now
$$
c_1(\fs)\pm \PD[e]-\Lambda \in H^2(\tilde X;\ZZ)
$$
is not a torsion class and so --- for $b_2^+(X)>0$, generic Riemannian
metrics on $X$ and related metrics on the connected sum $\tilde X$ --- the
moduli spaces $M_{\fs^\pm}(\tilde X)$ contain no zero-section pairs. Thus,
for our choice of generic perturbations, the moduli spaces
$M_{\fs^\pm}(\tilde X)$ are compact, oriented, smooth manifolds, both of
dimension $\dim M_{\fs}(X)$.

For $b_1(X)=0$ and odd $b_2^+(X)>1$, we define the {\em Seiberg-Witten
invariant\/} by \cite[\S 4.1]{FL2b}
\begin{equation}
\label{eq:DefSW}
  SW_{X}(\fs) 
= 
\langle\mu_{\fs^+}^{d},[M_{\fs^+}(\tilde X)]\rangle
= 
\langle\mu_{\fs^-}^{d},[M_{\fs^-}(\tilde X)]\rangle,
\end{equation}
where $2d=d_s(\fs)=d_s(\fs^\pm)$.
When $b^+_2(X)=1$ the pairing on the right-hand side of definition
\eqref{eq:DefSW} depends on the chamber in the positive cone of
$H^2(\tilde X;\RR)$ determined by the period point of the Riemannian metric
on $\tilde X$. The definition of the Seiberg-Witten invariant for this
case is also given in \cite[\S 4.1]{FL2b}: we assume that the class
$w_2(X)-\Lambda\pmod{2}$ is good to avoid technical difficulties involved
in chamber specification. Since $w\equiv w_2(X)-\Lambda \pmod{2}$, this
coincides with the constraint we use to define the Donaldson invariants in
\S \ref{subsec:Donaldsonseries} when $b_2^+(X)=1$. We refer to \cite[Lemma
4.1 \& Remark 4.2]{FL2b} for a comparison of the chamber structures
required for the definition of Donaldson and Seiberg-Witten invariants when
$b_2^+(X)=1$.

\subsubsection{Reducible $\SO(3)$ monopoles}
\label{subsubsec:RedPU2Monopole}
If $\ft=(\rho,V)$ and $\fs=(\rho,W)$ with $V=W\oplus W\otimes L$,
then there is an embedding
\begin{equation}
\label{eq:DefnOfIota}
\iota: \tsC_{\fs} \embed \tsC_{\ft},
\quad
(B,\Psi)\mapsto (B\oplus B\otimes A_\La\otimes B^{\det,*},\Psi\oplus 0),
\end{equation}
which is gauge-equivariant with respect to the homomorphism
\begin{equation}
\label{eq:GaugeGroupInclusion}
\varrho:\sG_{\fs}\embed \sG_{\ft},
\quad
s\mapsto \id_W\otimes\begin{pmatrix}s & 0 \\ 0 & s^{-1}\end{pmatrix}.
\end{equation}
According to \cite[Lemma 3.13]{FL2a}, the map \eqref{eq:DefnOfIota} defines
a topological embedding $\iota:M_{\fs}^0\embed\sM_{\ft}$, where
$M_{\fs}^0=M_{\fs}\cap\sC_{\fs}^0$ and an embedding of $M_{\fs}$ if
$w_2(\ft)\neq 0$ or $b_1(X)=0$; its image in $\sM_{\ft}$ is represented
by pairs which are reducible with respect to the splitting $V=W\oplus
W\otimes L$. Henceforth, we shall not distinguish between $M_\fs$ and its
image in $\sM_{\ft}$ under this embedding.

\subsubsection{Circle actions}
\label{subsubsec:S1Actions}
When $V=W\oplus W\otimes L$ and $\ft=(\rho,V)$, the space $\tsC_{\ft}$
inherits a circle action defined by
\begin{equation}
\label{eq:DefineS1LAction}
S^1\times V \to V, 
\quad (e^{i\theta},\Psi\oplus\Psi')\mapsto \Psi\oplus e^{i\theta}\Psi',
\end{equation}
where $\Psi\in C^\8(W)$ and $\Psi'\in C^\8(W\otimes L)$. With respect to
the splitting $V=W\oplus W\otimes L$, the actions
\eqref{eq:DefineS1LAction} and \eqref{eq:S1ZAction} are related by
\begin{equation}
\label{eq:RelateS1Actions}
\begin{pmatrix}1 & 0 \\ 0 & e^{i2\theta}\end{pmatrix}
=
e^{i\theta}u,
\text{ where }
u = \begin{pmatrix}e^{-i\theta} & 0 \\ 0 & e^{i\theta}\end{pmatrix}
\in\sG_{\ft},
\end{equation}
and so, when we pass to the induced circle actions on the quotient
$\sC_{\ft}=\tsC_{\ft}/\sG_{\ft}$, the actions
\eqref{eq:DefineS1LAction} and \eqref{eq:S1ZAction} on
$\sC_{\ft}$ differ only in their multiplicity. Recall \cite[Lemma
3.11]{FL2a} that the image  in $\tsC_{\ft}$ of the map \eqref{eq:DefnOfIota}
contains all pairs which are fixed by the circle action
\eqref{eq:DefineS1LAction}.

\subsubsection{The virtual normal bundle of the Seiberg-Witten moduli space}
\label{subsubsec:ThickenedNeighborhood}
Suppose $\ft=(\rho,V)$ and $\fs=(\rho,W)$, with $V=W\oplus W\otimes L$, so
we have a topological embedding $M_{\fs}\embed \sM_{\ft}$; we assume
$M_{\fs}$ contains no zero-section monopoles.  Recall from \cite[\S
3.5]{FL2a} that there exist finite-rank, complex vector bundles,
\begin{equation}
\label{eq:VirtualNormalandObstBundles}
\pi_{\Xi}:\Xi_{\ft,\fs}\to M_{\fs}
\quad\text{and}\quad
\pi_{N}:N_{\ft,\fs}\to M_{\fs},
\end{equation}
with $\Xi_{\ft,\fs}\cong M_{\fs}\times\CC^{r_\Xi}$, called the {\em
obstruction bundles\/} and {\em virtual normal bundles\/} of
$M_{\fs}\embed\sM_{\ft}$, respectively.  For a small enough positive radius
$\eps$, there are a topological embedding \cite[Theorem 3.21]{FL2a} of an
open tubular neighborhood,
\begin{equation}
\label{eq:BackgroundConfigEmbedding}
\bga_{\fs}:N_{\ft,\fs}(\eps)
\embed
\sC_{\ft},
\end{equation}
and a smooth section $\bchi_{\fs}$ of the pulled-back
complex vector bundle, 
\begin{equation}
\label{eq:PulledBackObstructionBundle}
\pi_{N}^*\Xi_{\ft,\fs} \to N_{\ft,\fs}(\eps),
\end{equation}
such that the tubular map yields a homeomorphism
\begin{equation}
\label{eq:HomeoSWEmbedding}
\bga_{\fs}: \bchi_{\fs}^{-1}(0)\cap N_{\ft,\fs}(\eps)
\cong
\sM_{\ft}\cap \bga_{\fs}(N_{\ft,\fs}(\eps)),
\end{equation}
restricting to a diffeomorphism on the complement of $M_{\fs}$ and
identifying $M_{\fs}$ with its image $\iota(M_{\fs})\subset\sM_{\ft}$.  
The image $\bga_{\fs}(N_{\ft,\fs}(\eps))$ is an open subset of a {\em
virtual moduli space\/}, $\sM_{\ft,\fs}^{\stab}\subset \sC_{\ft}$.

Our terminology is loosely motivated by that of \cite{GraberPand} and
\cite{RuanSW}, where the goal (translated to our setting) would be to
construct a {\em virtual fundamental class\/} for $\sM_{\ft}$, given by
the cap product of the fundamental class of an ambient space containing
$\sM_{\ft}$ with the Euler class of a vector bundle over this ambient
space, where $\sM_{\ft}$ is the zero locus of a (possibly
non-transversally vanishing) section. Here, $\sM_{\ft,\fs}^{\stab}$ plays the
role of the ambient space and (the pushforward of) $\Xi_{\ft,\fs}$ the
vector bundle with zero section yielding (an open neighborhood in)
$\sM_{\ft}$.  Then, $N_{\ft,\fs}$ is the normal bundle of $M_{\fs}\embed
\sM_{\ft,\fs}^{\stab}$, while $[N_{\ft,\fs}]-[\Xi_{\ft,\fs}]$ would more
properly be called the `virtual normal bundle' of $M_{\fs}\embed
\sM_{\ft}$, in the language of $K$-theory.

Recall that minus the index of the $\SO(3)$-monopole elliptic deformation
complex \cite[Equations (2.47) \& (3.35)]{FL2a} at a reducible solution can
be written as 
\begin{equation}
\label{eq:SWDimRelations}
\dim\sM_{\ft}
= 
2n_s(\ft,\fs) + d_{s}(\fs),
\end{equation}
where $d_s(\fs)$ is the expected dimension of the Seiberg-Witten moduli
space $M_{\fs}$ (see equation \eqref{eq:DimSW}), while
$n_s(\ft,\fs)=n_s'(\ft,\fs)+n_s''(\ft,\fs)$ is minus the complex index of
the normal deformation operator \cite[Equations (3.71) \& (3.72)]{FL2a}),
with
\begin{equation}
\label{eq:NormalComponentDims}
\begin{aligned}
n_s'(\ft,\fs) 
&= -(c_1(\ft)-c_1(\fs))^2-\textstyle{\frac{1}{2}}(\chi+\sigma),
\\
n_s''(\ft,\fs)
&=
\textstyle{\frac{1}{8}}(c_1(\fs)-2c_1(\ft))^2-\sigma).
\end{aligned}
\end{equation}
If $r_\Xi$ is the complex rank of $\Xi_{\ft,\fs}\to M_{\fs}$, and
$r_N(\ft,\fs)$ is the complex rank of $N_{\ft,\fs}\to M_{\fs}$, then
\begin{equation}
\label{eq:SWBundleRankRelations}
r_N(\ft,\fs)
=
n_s(\ft,\fs) + r_\Xi,
\end{equation}
as we can see from the dimension relation \eqref{eq:SWDimRelations}
and the topological model \eqref{eq:HomeoSWEmbedding}.

The map \eqref{eq:BackgroundConfigEmbedding} is circle equivariant when the
circle acts trivially on $M_{\fs}$, by scalar multiplication on the fibers
of $N_{\ft,\fs}(\eps)$, and by the action
\eqref{eq:DefineS1LAction} on $\sC_{\ft}$. The bundle
\eqref{eq:PulledBackObstructionBundle} and section $\bchi_{\fs}$ are
circle equivariant if the circle acts on $N_{\ft,\fs}$ and the fibers
of $\bga_{\fs}^*\Xi_{\ft,\fs}$ by scalar multiplication.

Let $\tN_{\ft,\fs}\to \tM_{\fs}$ be the pullback of $N_{\ft,\fs}$
by the projection $\tM_{\fs}\to M_{\fs}=\tM_{\fs}/\sG_{\fs}$, so
$\tN_{\ft,\fs}$ is a $\sG_{\fs}$-equivariant bundle, 
where $\sG_{\fs}$ acts on the base $\tM_{\fs}$ by the usual gauge group
action \eqref{eq:SWGaugeGroupAction} and the induced action on the total space,
\begin{equation}
\label{eq:VirtualNormalBundleTotalSpace}
\tN_{\ft,\fs}
\subset
\tM_{\fs}\times
L^2_k(\Lambda^1\otimes_\RR L) \oplus L^2(W^+\otimes L)
\subset
\tM_{\fs}\times
L^2_k(\Lambda^1\otimes_\RR \fg_{V}) \oplus L^2(V^+),
\end{equation}
via the embedding \eqref{eq:GaugeGroupInclusion} of $\sG_{\fs}$ into
$\sG_{\ft}$ and the splittings $\fg_{V}\cong\underline{\RR}\oplus L$ 
\cite[Lemma 3.10]{FL2a} and
$V=W\oplus W\otimes L$.  Thus, $s\in\sG_{\fs}$ acts by scalar multiplication
by $s^{-2}$ on sections of $\Lambda^1\otimes_\RR L$ and by $s^{-1}$ on
sections of $W^+\otimes L$ \cite[\S 3.5.4]{FL2a}.

For a small enough positive $\eps$, there is a smooth embedding \cite[\S
3.5.4]{FL2a} of the open tubular neighborhood $\tN_{\ft,\fs}(\eps)$,
\begin{equation}
\label{eq:BackgroundPreConfigEmbedding}
\tilde\bga_{\fs}:\tN_{\ft,\fs}(\eps)
\to
\tsC_{\ft},
\end{equation}
which is gauge equivariant with respect to the preceding action of
$\sG_{\fs}$, and covers the topological embedding
\eqref{eq:BackgroundConfigEmbedding}.  The map
\eqref{eq:BackgroundPreConfigEmbedding} is circle equivariant, where the
circle acts trivially on $\tM_{\fs}$, by scalar multiplication on the
fibers of $\tN_{\ft,\fs}(\eps)$, and by the action
\eqref{eq:DefineS1LAction} on $\sC_{\ft}$.
We note that the map \eqref{eq:BackgroundPreConfigEmbedding} is also circle
equivariant with respect to the action \eqref{eq:S1ZAction} on $\sC_{\ft}$,
if the circle acts on $\tN_{\ft,\fs}$ by 
\begin{equation}
\label{eq:S1ZActionOnN}
(e^{i\theta},(B,\Psi,\beta,\psi))
\mapsto
\varrho(e^{i\theta})( B,\Psi,e^{2i\theta}\beta,e^{2i\theta}\psi)
=
(B,e^{i\theta}\Psi,\beta,e^{i\theta}\psi),
\end{equation}
where $(B,\Psi)\in \tM_{\fs}$,
$(\beta,\psi)\in L^2_k(\Lambda^1\otimes_\RR L) \oplus L^2(W^+\otimes L)$,
so $(B,\Psi,\beta,\psi)\in \tN_{\ft,\fs}$, and
$\varrho:\sG_{\fs}\to\sG_{\ft}$ is the homomorphism
\eqref{eq:GaugeGroupInclusion}. This equivariance follows from
the relation \eqref{eq:RelateS1Actions}
between the actions \eqref{eq:DefineS1LAction}
and \eqref{eq:S1ZAction}.

We recall the following calculation of the Segre classes of the virtual
normal bundle $N_{\ft,\fs}$:

\begin{lem}
\label{lem:TopOfReducibleNormal}
\cite[Lemma 4.11]{FL2b}
Assume $M_{\fs}\subset \sM_{\ft}$ contains no zero-section pairs and let
$\mu_{\fs}\in H^2(M_{\fs};\ZZ)$ be the cohomology class 
\eqref{eq:SWClass}.  If $b_1(X)=0$, then the Segre
classes of the complex vector bundle $N_{\ft,\fs}\to M_{\fs}$ are given by
\begin{equation}
\label{eq:TopOfReducibleNormal}
s_{p}(N_{\ft,\fs})
=
\mu_{\fs}^{p}\sum_{j=0}^{p} 2^j \binom{-n_s'}{j}\binom{-n_s''}{p-j},
\quad p=0,1,2,\dots
\end{equation}
\end{lem}

\subsection{Cohomology classes on the moduli space of SO(3) monopoles}
\label{subsec:Cohomology}
The identity \eqref{eq:WCL1Equation} arises as an equality of pairings of
suitable cohomology classes with a link in $\bar\sM_{\ft}^{*,0}/S^1$ of the
anti-self-dual moduli subspace and the links of the Seiberg-Witten moduli
subspaces in $\bar\sM_{\ft}/S^1$. We recall the definitions of these
cohomology classes and their dual geometric representatives given in
\cite[\S 3]{FL2b}.

The first kind of cohomology class is defined on $\sM^*_{\ft}/S^1$, via the
associated $\SO(3)$ bundle,
\begin{equation}
\label{eq:UniversalU2Bundle}
\FF_{\ft}
=
\tsC^{*}_\ft/S^1\times_{\sG_\ft}\fg_{V}
\to
\sC^{*}_\ft/S^1\times X.
\end{equation}
The group $\sG_{\ft}$ acts diagonally in \eqref{eq:UniversalU2Bundle}, 
with $\sG_{\ft}$ acting on the left on $\fg_{V}$. We define \cite[\S
3.1]{FL2b} 
\begin{equation}
\label{eq:DefnMuClasses}
\mu_p: H_\bullet(X;\RR)\to H^{4-\bullet}(\sC^*_\ft/S^1;\RR),
\quad
\beta\mapsto -\tquarter p_1(\FF_{\ft})/\beta.
\end{equation}
On restriction to $M^w_{\ka}\hookrightarrow \sM_{\ft}$, the cohomology classes
$\mu_p(\beta)$ coincide with those used in the definition of Donaldson
invariants \cite[Lemma 3.1]{FL2b}. We assume $b_1(X)=0$ 
and thus let
\begin{equation}
\label{eq:DefineAAA}
\AAA(X) = \Sym\left( H_{\even}(X;\RR)\right)
\end{equation}
be the graded algebra, with $z=\beta_1\beta_2\cdots\beta_r$ having total
degree $\deg(z) = \sum_p(4-i_p)$, when $\beta_p\in H_{i_p}(X;\RR)$.  A
point $x \in X$ gives a distinguished generator still called $x$ in
$\AAA(X)$ of degree four. Then $\mu_p$ extends in the usual way to
a homomorphism of graded real algebras,
$$
\mu_p:\AAA(X) \to \Sym(H^{\even}(\sC^*_\ft/S^1;\RR)),
$$
which preserves degrees. Next, we define a complex line bundle over
$\sC^{*,0}_{\ft}/S^1$, 
\begin{equation}
\label{eq:DefineDetLineBundle}
\LL_{\ft}= \sC^{*,0}_{\ft}\times_{(S^1,\times -2)}\CC,
\end{equation}
where the circle action is given, for $[A,\Phi]\in\sC^{*,0}_{\ft}$
and $\zeta\in\CC$, by
\begin{equation}
\label{eq:DeterminantS1Action}
([A,\Phi],\zeta) \mapsto ([A,e^{i\theta}\Phi], e^{2i\theta}\zeta).
\end{equation}
Then we define the second kind of cohomology class on $\sM_{\ft}^{*,0}/S^1$ by
\begin{equation}
\label{eq:DefineMuC1}
\mu_c = c_1(\LL_{\ft}) \in H^2(\sC_{\ft}^{*,0}/S^1;\RR).
\end{equation}
For monomials $z\in\AAA(X)$, we constructed
\cite[\S 3.2]{FL2b} geometric representatives $\sV(z)$
dual to $\mu_p(z)$ and $\sW$ dual to $\mu_c$, defined on $\sM_{\ft}^*/S^1$ and
$\sM_{\ft}^{*,0}/S^1$, respectively; their closures in $\bar\sM_{\ft}/S^1$ are
denoted by $\bar\sV(z)$ and $\bar\sW$ \cite[Definition 3.14]{FL2b}. When
$$
\deg(z) +2\eta
=
\dim(\sM^{*,0}_{\ft}/S^1) -1,
$$
and $\deg(z)\ge \dim M^w_\ka$
it follows that \cite[\S 3.3]{FL2b} the intersection
\begin{equation}
\label{eq:GeomReprIntersection}
\bar\sV(z) \cap \bar\sW^{\eta} \cap \bar\sM^{*,0}_{\ft}/S^1,
\end{equation}
is an oriented one-manifold (not necessarily connected) whose closure in
$\bar\sM_{\ft}/S^1$ can only intersect $(\bar\sM_{\ft}-\sM_{\ft})/S^1$ at
points in $\sM_{\ft}^{\red}\cong\cup(M_{\fs}\times\Sym^\ell(X))$, where the
union is over $\ell\geq 0$ and $\fs\in\Spinc(X)$ \cite[Corollary 3.18]{FL2b}.

\subsection{Donaldson invariants}
\label{subsec:Donaldsonseries}
We first recall the definition \cite[\S 2]{KMStructure} of the Donaldson series
when $b_1(X)=0$ and $b_2^+(X)>1$ is odd, so that, in this case,
$\chi+\sigma\equiv 0\pmod{4}$. See also \S 3.4.2 in \cite{FL2b}, especially
for a definition of the Donaldson invariants when $b_2^+(X)=1$.
For any choice of $w\in H^{2}(X;\ZZ)$, the
Donaldson invariant is a linear function
$$
D^{w}_{X}:\AAA(X) \to \RR.
$$
Let $\tilde X=X\#\overline{\CC\PP}^2$ be the blow-up of $X$ and let $e\in
H_2(\tilde X;\ZZ)$ be the exceptional class, with Poincar\'e dual
$\PD[e]\in H^2(\tilde X;\ZZ)$.  If $z\in\AAA(X)$ is a monomial, we define
$D_X^w(z)=0$ unless
\begin{equation}
    \mathrm{deg}(z)\equiv -2w^{2}-\textstyle{\frac{3}{2}}(\chi+\sigma)\pmod{8}.
    \label{mod8}
\end{equation}
If $\deg(z)$ obeys equation \eqref{mod8}, we let
$\kappa\in\frac{1}{4}\ZZ$ be defined by 
$$
\deg(z)=8\ka -\textstyle{\frac{3}{2}}\left(\chi+\si\right).
$$
There exists an $\SO(3)$ bundle over $\tilde X$ with first Pontrjagin
number $-4\ka -1$ and second Stiefel-Whitney class $w+\PD[e]\pmod{2}$.
One then defines the Donaldson invariant on monomials by
\begin{equation}
\label{eq:DefineDonaldson}
D^w_X(z)
= 
\#\left( \bar\sV(ze)\cap\barM^{w+\PD[e]}_{\ka+1/4}(\tilde X)\right),
\end{equation}
and extends to a real linear function on $\AAA(X)$. Note that
$w+\PD[e]\pmod{2}$ is good in the sense of Definition
\ref{defn:Good}. If $w'\equiv w \pmod{2}$, then \cite{DonOrient}
\begin{equation}
\label{eq:DonaldsonsSignChange}
D^{w'}_{X}=(-1)^{\frac{1}{4}(w'-w)^{2}}D^{w}_{X}.
\end{equation}
The Donaldson series is a formal power series, 
\begin{equation}
\label{eq:DefineDonaldsonSeries}
\bD^{w}_{X}(h) = D^{w}_{X}((1+\textstyle{\frac{1}{2}} x)e^{h}),
\quad h \in H_{2}(X;\RR).
\end{equation}
By equation~\eqref{mod8}, the series $\bD^{w}_{X}$ is even if
$$
-w^{2}-\tthreequarter(\chi+\sigma)\equiv  0 \pmod 2,
$$
and odd otherwise. A four-manifold has KM-simple type if for some $w$ and
all $z\in \AAA(X)$,  
$$
D^{w}_{X}(x^{2}z)=4D^{w}_{X}(z).
$$
According to \cite[Theorem 1.7]{KMStructure}, when $X$ has KM-simple type
the series $\bD^{w}_{X}(h)$ is an analytic function of $h$ and there are
finitely many characteristic cohomology classes $K_{1},\ldots,K_{m}$ (the
KM-basic classes) and constants $a_{1},\ldots,a_{m}$ (independent of $w$)
so that
$$
\bD^{w}_X(h)
= 
e^{\half h \cdot h}
\sum_{i=1}^{r}(-1)^{\half(w^{2}+K_{i}\cdot w)}a_{i}e^{\langle
K_{i},h\rangle}.
$$
Witten's conjecture \eqref{eq:WConjecture} then relates the Donaldson and
Seiberg-Witten series for four-manifolds of simple type.

When $b^+_2(X)=1$ the pairing on the right-hand side of definition
\eqref{eq:DefineDonaldson} depends on the chamber in the positive cone of
$H^2(\tilde X;\RR)$ determined by the period point of the Riemannian metric
on $\tilde X$, just as in the case of Seiberg-Witten invariants described in \S
\ref{subsubsec:SWInvariants}. We refer to \S 3.4.2 in \cite{FL2b} for a
detailed discussion of this case and, as in \cite{FL2b}, we assume that the
class $w\pmod{2}$ is good in order to avoid technical difficulties involved
in chamber specification.

\subsection{Links and the cobordism}
\label{subsec:LinkCobordism}
Since the endpoints of the components of the one-manifold
\eqref{eq:GeomReprIntersection} either lie in $M_\kappa^w$ or
$M_{\fs}\times\Sym^\ell(X)$, for some $\fs$ and $\ell(\ft,\fs)\geq 0$ for
which $\ft_{\ell} = \fs\oplus\fs\otimes L$, we have when $w_2(\ft)$ is good
in the sense of Definition \ref{defn:Good} that
\begin{equation}
\label{eq:RawCobordismSum}
\#\left(\bar\sV(z)\cap \bar\sW^{n_a-1}\cap \bar\bL^w_{\ft,\ka}\right)
=
-\sum_{\fs\in\Spinc(X)}
\# \left( \bar\sV(z)\cap \bar\sW^{n_a-1}\cap \bar\bL_{\ft,\fs}\right),
\end{equation}
where $\bar\bL^w_{\ft,\ka}$ is the link of $\bar M^w_\ka$
in $\bar\sM_{\ft}/S^1$ (see \cite[Definition 3.7]{FL2a})
and
where $\bar\bL_{\ft,\fs}$ is empty if $\ell(\ft,\fs)<0$
and is the boundary of an open neighborhood of the Seiberg-Witten stratum
$M_{\fs}\times\Sym^\ell(X)$ in $\bar\sM_{\ft}/S^1$
if $\ell=\ell(\ft,\fs)\ge 0$.  By construction, the intersection of
$\bar\bL_{\ft,\fs}$ with the top stratum of $\bar\sM_{\ft}/S^1$
is a smooth manifold and the intersection of $\bar\bL_{\ft,\fs}$ 
with the geometric representatives is in this top stratum.
The precise definition of $\bar\bL_{\ft,\fs}$ is given
in \cite[Definition 3.22]{FL2a} for $\ell=0$, in \eqref{eq:DefineLink} for
$\ell=1$, and will be given in \cite{FL5} for $\ell>1$.

When $\deg(z) = \dim M^w_{\ka}$, the intersection of the one-manifold
\eqref{eq:GeomReprIntersection} with the link $\bar\bL^w_{\ft,\ka}$ is given by
\cite[Lemma 3.30]{FL2b}
\begin{equation}
\label{eq:ASDPairing}
2^{1-n_a}\#\left( \bar\sV(z)\cap \bar\sW^{n_a-1}\cap \bar\bL^w_{\ft,\ka}\right)
=
\#\left(\bar\sV(z)\cap\barM^w_{\ka}\right).
\end{equation}
Applying this identity to the blow-up, $X\#\overline{\CC\PP}^2$, when
$n_a(\ft)>0$, we recover the Donaldson invariant $D_X^w(z)$ on the
right-hand side of
\eqref{eq:ASDPairing} via definition \eqref{eq:DefineDonaldson}.  

Our proof of Theorem
\ref{thm:WCL1} relies on the following conjecture whose motivation is
discussed in \cite[p. 179]{FKLM}:

\begin{conj}
\label{conj:Multiplicity}
\cite[Conjecture 3.1]{FKLM}
When $b_1(X)=0$ and $\eta\geq 0$ is an integer for which
$\deg(z)+2(\eta+1)=\dim\sM_{\ft}$, the intersection number 
$$
\#\left(\bar\sV(z)\cap\bar\sW^{\eta} \cap \bar\bL_{\ft,\fs}\right)
$$ 
is a multiple of $SW_X(\fs)$ and thus vanishes if $SW_X(\fs)=0$.
\end{conj}

We need only assume this conjecture holds for $\ell(\ft,\fs)>1$, as we
shall prove it here when $\ell(\ft,\fs)=1$; the relevant statement for
$b_1(X)\geq 0$ is given as Conjecture 3.34 in \cite{FL2b}, which we proved
when $\ell(\ft,\fs)=0$ \cite[Remark 4.15]{FL2b}. The identities
\eqref{eq:ASDPairing} and \eqref{eq:RawCobordismSum} yield
\begin{equation}
\label{eq:CobordismSum}
2^{n_a-1}\#\left(\bar\sV(z)\cap\barM^w_{\ka}\right)
=
-\sum_{\fs\in\Spinc(X)}
\# \left( \bar\sV(z)\cap \bar\sW^{n_a-1}\cap \bar\bL_{\ft,\fs}\right),
\end{equation}
where, assuming Conjecture \ref{conj:Multiplicity}, the intersection
numbers on the right-hand side are non-zero only when $SW_X(\fs)\neq 0$
(and $\ell(\ft,\fs)\geq 0$). The hypotheses of Theorem \ref{thm:WCL1}
ensure that we only have $\ell(\ft,\fs)=0,1$ in \eqref{eq:CobordismSum}
when $SW_X(\fs)\neq 0$ and, since we addressed the case $\ell=0$ in
\cite{FL2b}, the remainder of the present article concerns the case
$\ell=1$.


\section{Gluing SO(3) monopoles}
\label{sec:gluing}
Our goal in this section is to construct a topological model for an open
neighborhood in $\bar\sM_{\ft}$, and thus a link, of the Seiberg-Witten stratum
$$
(M_{\fs}\times\Sym^\ell(X))\cap\bar\sM_{\ft}
\subset
\bar\sM_{\ft},
$$
where $M_{\fs}\subset\sM_{\ft_\ell}$ and $p_1(\ft_\ell)=p_1(\ft)+4\ell$
with $\ell\geq 0$; in this article, we shall only carry out this
construction when $\ell=1$, whereas in \cite{FL5} we consider the general
case $\ell\geq 1$. For the remainder of this article, however, it will be
more convenient to denote the pair $(\ft_\ell,\ft)$ by $(\ft,\ft')$ in
order to minimize the number of subscripts and simplify notation. Hence,
for the remainder of the article, we shall denote the Seiberg-Witten
stratum of interest by
\begin{equation}
\label{eq:SWStratumInLevelEll}
(M_{\fs}\times\Sym^\ell(X))\cap\bar\sM_{\ft'}
\subset
\bar\sM_{\ft'},
\end{equation}
where $M_{\fs}\subset\sM_{\ft}$ and $p_1(\ft)=p_1(\ft')+4\ell$ with
$\ft=(\rho,V)$, $\ft'=(\rho',V')$, and $\ell\geq 0$.

The topological models for open neighborhoods in $\bar\sM_{\ft'}$ of the strata
\eqref{eq:SWStratumInLevelEll} are constructed 
by gluing $\SO(3)$ monopoles. In \S
\ref{subsec:SplicingBasics} we describe the cut-and-paste process of
splicing instantons from the four-sphere (corresponding to $\SO(3)$
monopoles with zero-spinor) onto a family of $\SO(3)$ monopoles over $X$
parametrized by a tubular neighborhood $N_{\ft,\fs}(\eps)$ of
$M_{\fs}$. The supporting technical details required to define this
splicing process are discussed in
\S \ref{subsec:CliffordModuleBases}, where we review the required Clifford
module and spin group representation theory, in \S
\ref{subsec:StructureGroups} where we describe the Clifford frame bundles
and their structure groups, and in \S
\ref{subsec:S4ModuliSp} where we describe the connections and
spinors on the four-sphere to be spliced onto the four-manifold $X$.
The space of gluing data and splicing map $\bgamma'$ are defined in \S
\ref{subsec:SplicingMap}, while in \S \ref{subsec:GluingMap}
we recall a gluing result (Theorem \ref{thm:GluingThm}) from \cite{FL3,
FL4} which asserts that the splicing map can be perturbed to a gluing map, thus
giving a topological model for an open neighborhood of the stratum
\eqref{eq:SWStratumInLevelEll} in 
$\bar\sM_{\ft'}$ when $\ell=1$; the extension to the general case $\ell\geq
1$ is discussed in \cite{FL5}.  The obstruction bundle
appearing in this model is described in \S
\ref{subsec:ObstructionBundle}. In
particular, the model allows us to define in \S \ref{subsec:DefnOfLink} a
link $\bar\bL_{\ft',\fs}$ of the stratum $M_{\fs}\times X$ in
$\bar\sM_{\ft'}$. The relationships among the orientations of the link and
the associated moduli spaces are described in \S \ref{subsec:Orient}.

\subsection{Clifford modules and representation theory}
\label{subsec:CliffordModuleBases}
We recall the representation theory we shall need for the
complex Clifford algebra $\CCl(\RR^4)$. (In Lawson-Michelsohn \cite{LM},
this algebra is denoted $\CCl_4 = \Cl_4\otimes_\RR\CC$, with
$\Cl_4=\Cl_{4,0}$.)  Recall from \cite[Theorem I.4.3]{LM} that
$\CCl(\RR^4)\cong M_4(\CC)$, where $M_d(\CC)$ denotes the algebra of
complex $d\times d$ matrices. Then the natural representation of $M_4(\CC)$
on $\Delta=\CC^4$ is the only irreducible complex representation of
$M_4(\CC)$ \cite[Theorems I.5.6--7]{LM}, up to equivalence in the sense of
\cite[Definition I.5.1]{LM}. In particular, $\Delta$ is the only
irreducible complex representation of $\CCl(\RR^4)$ \cite[Theorems
I.5.8]{LM}.

If $\VV$ is a complex vector space of dimension $4d$ and
$\CCl(\RR^4)$-module, then $\VV$ is isomorphic to
$\Delta\oplus\cdots\oplus\Delta$ ($d$ times), or $\Delta\otimes\CC^d$, as a
$\CCl(\RR^4)$-module by \cite[Proposition I.5.4]{LM}.
From the remarks surrounding
\cite[Equation (2.28)]{FL2a} for the case $d=2$, the map
\begin{equation}
\label{eq:MatrixIsCliffModEndo}
M_d(\CC)\to \End_{\CCl(\RR^4)}(\Delta\otimes\CC^d),
\quad
M\mapsto \id_{\Delta}\otimes M,
\end{equation}
is an isomorphism of complex algebras.

Recall that a representation $\rho:\CCl(\RR^4)\to \End_\CC(\VV)$ is
{\em unitary}, in the sense implied by \cite[Theorem I.5.16]{LM}, when
$\VV$ is equipped with a Hermitian metric such that
$$
\langle\rho(\alpha)\Phi,\rho(\alpha)\Psi\rangle
=
\langle\Phi,\Psi\rangle,
$$
for all $\alpha\in\RR^4$ with $|\alpha|=1$ and $\Phi,\Psi\in \VV$.

\begin{lem}
\label{lem:CliffModuleBasis}
Let $\VV$ be a complex vector space of dimension $4d$ and
$\CCl(\RR^4)$-module. If $F_1, F_2:\Delta\otimes\CC^d\to\VV$ are
isomorphisms of $\CCl(\RR^4)$-modules, then
$$
F_1=F_2\circ(\id_\Delta\otimes M),
$$
where $M\in\GL(d,\CC)$; if $F_1,F_2$ respect Hermitian metrics, then
$M\in\U(d)$.
\end{lem}

\begin{proof}
Given $\CCl(\RR^4)$-module isomorphisms $F_1,F_2$, then $F_2^{-1}F_1$ is a
$\CCl(\RR^4)$-automorphism of the $\CCl(\RR^4)$-module $\Delta\otimes\CC^d$
and so $F_2^{-1}F_1=\id_{\Delta}\otimes M\in\GL(d,\CC)$ by equation
\eqref{eq:MatrixIsCliffModEndo}. If $F_1,F_2$ respect Hermitian metrics, then
$M^\dagger M=\id_{\CC^d}$ and so $M\in\U(d)$.
\end{proof}

If $\VV$ is a complex rank-eight vector space and complex
$\CCl(\RR^4)$-module, then $\VV=\VV^+\oplus\VV^-$ and we have an associated
complex determinant line $\det(\VV^+)=\Lambda^4(\VV^+)$. Given a complex
rank-four vector space and $\CCl(\RR^4)$-module $\WW$, there is an
isomorphism $F:\WW\to\Delta$ of $\CCl(\RR^4)$-modules (determined up to
multiplication by some $e^{i\theta}\in S^1$). Setting
$\EE=\Hom_{\CCl(\RR^4)}(\WW,\VV)$ (now determined up to multiplication by
$e^{-i\theta}\in S^1$), we obtain an isomorphism of $\CCl(\RR^4)$-modules,
$$
\WW\otimes_\CC\EE \to \VV, 
$$
induced by the map $\Phi\otimes H \mapsto H(\Phi)$, where $\Phi\in\WW$
and $H\in E$. We then obtain a square root of $\det(\VV^+)$,
\begin{equation}
\label{eq:DefnDetVV}
{\det}^{\frac{1}{2}}(\VV^+)
=
\det(\WW^+)\otimes_\CC\det(\EE)
=
\Lambda^2(\WW^+)\otimes_\CC\Lambda^2(\EE),
\end{equation}
which is independent of the choice of $\WW$.

\subsection{Structure groups and associated bundles}
\label{subsec:StructureGroups}
Given a Hermitian Clifford module $(\rho,V)$ over $(X,g)$, where $V$ has
complex rank $4d$, the Hermitian subbundles $V^\pm \subset V$ are defined
as the $\mp 1$ eigenspaces of $\rho(\vol)$ on $V$. The determinant bundles
$\det(V^\pm)$ are canonically isomorphic via an isomorphism
$\det(\rho(\alpha)):\det(V^+)\to \det(V^-)$, $\alpha\in C^\8(U,T^*X)$ with
$|\alpha|=1$, which is independent of the choice of local one-forms
\cite[\S 5.2]{SalamonSWBook}.

We let $\Fr(T^*X)\to X$ be the principal $\SO(4)$ bundle with fiber
$\Isom_\RR^+(\RR^4,T^*X|_x)$ over $x\in X$ given by the space of
orientation-preserving isometries, with $\SO(4)$ action induced by the
standard action on $\RR^4$, and analogously define $\Fr(TX)$.

Next, we define the frame bundles for the Clifford modules.
Given $\Delta\otimes_\CC\CC^d$ with its $\CCl(\RR^4)$-module structure
induced by the Clifford map $\rho_0:\RR^4\to\End_\CC(\Delta\otimes_\CC\CC^d)$,
and $V_x$ with its $\CCl(T^*X|_x)$-module structure induced by the
Clifford map $\rho:T^*X|_x\to \End_\CC(V_x)$, for each $x\in X$ we let
$\Isom_{\CCl}(\Delta\otimes_\CC\CC^d,V_x)$ be the space of unitary isomorphisms
$\tilde F:\Delta\otimes\CC^d\to V_x$ which are Clifford-module isomorphisms
with respect to a Clifford-algebra isomorphism
$F:\CCl(\RR^4)\to\CCl(T^*X|_x)$ determined by a frame $F:\RR^4\to
T^*X|_x$. Thus each $\tilde F\in\Isom_{\CCl}(\Delta\otimes_\CC\CC^d,V_x)$
covers a frame $F\in\Isom_\RR^+(\RR^4,T^*X|_x)$ and corresponding
Clifford-algebra isomorphism yielding a commutative diagram, for all
$\alpha\in\CCl(\RR^4)$, 
$$
\begin{CD}
\Delta\otimes_\CC\CC^d @>{\tilde F}>> V_x
\\
@V{\rho_0(\alpha)}VV @VV{\rho(F\alpha)}V
\\
\Delta\otimes_\CC\CC^d @>{\tilde F}>> V_x
\end{CD}
$$
Hence, we define the Clifford-module frame bundle for $V$,
\begin{equation}
\label{eq:CliffordFrameBundlesWV}
\pi:\Fr_{\CCl(T^*X)}(V)\to \Fr(T^*X),
\end{equation}
to have fibers $\Isom_{\CCl}(\Delta\otimes_\CC\CC^d,V_x)$, where $\tilde
F\in\Isom_{\CCl}(\Delta\otimes_\CC\CC^d,V_x)$ and $\pi(\tilde F)=F$ obey
\begin{equation}
\label{eq:CliffordModuleHom}
\tilde F\rho_0(\alpha)\tilde F^{-1} = \rho(F\alpha),
\quad
\alpha\in\CCl(\RR^4).
\end{equation}
If $\tilde F_1,\tilde F_2\in\Fr_{\CCl(T^*X)}(V_x)$ both lie in the fiber of
the projection \eqref{eq:CliffordFrameBundlesWV} over $F\in\Fr(T^*X|_x)$, then
$\tilde F^{-1}_2\tilde F_1$ is a Hermitian $\CCl(\RR^4)$-module
automorphism of $\Delta\otimes_\CC\CC^d$ and hence $\tilde F^{-1}_2\tilde
F_1=\id_\Delta\otimes M$, for some $M\in\U(d)$ by Lemma
\ref{lem:CliffModuleBasis}. Therefore, we have

\begin{lem}
\label{lem:CliffordFrameBundleOverTX}
The space $\Fr_{\CCl(T^*X)}(V)$ is a principal $\U(d)$ bundle over
$\Fr(T^*X)$. 
\end{lem}

We now consider the case $d=1$ and let $(\rho,W)$ be a \spinc structure over
$X$. 

\begin{lem}
\label{lem:SpincFrameBundleOverX}
The space $\Fr_{\CCl(T^*X)}(W)$ is a principal {\em $\Spinc(4)$} bundle
over $X$. 
\end{lem}

\begin{proof}
If $\tilde F_1,\tilde F_2\in\Fr_{\CCl(T^*X)}(W_x)$ lie in the fibers
over $F_1,F_2\in\Fr(T^*X|_x)$ of the projection
\eqref{eq:CliffordFrameBundlesWV}, respectively, then $F_2^{-1}F_1=R\in\SO(4)$
is an isometric automorphism of $\RR^4$ while $\tilde F_2^{-1}\tilde
F_1=\tilde R\in\U(4)$ is a Hermitian automorphism of $\Delta$, yielding a
commutative diagram
$$
\begin{CD}
\Delta @>{\tilde R}>> \Delta
\\
@V{\rho_0(\alpha)}VV @VV{\rho_0(R\alpha)}V
\\
\Delta @>{\tilde R}>> \Delta
\end{CD}
$$
Thus, just as in \eqref{eq:CliffordModuleHom}, these automorphisms obey
\begin{equation}
\label{eq:CliffordR} 
\tilde R\rho_0(\alpha) \tilde R^{-1}
=
\rho_0(R\alpha),
\quad
\alpha\in\CCl(\RR^4).
\end{equation}
Because $R$ has determinant one, it preserves the volume form.  Thus
identity \eqref{eq:CliffordR} implies that $\tilde R$ preserves the
splitting $\Delta=\Delta^+\oplus \Delta^-$ and hence, we may write
$$
\tilde R=(\tilde R_+,\tilde R_-)\in\U(2)\times\U(2)\subset \U(4).
$$
Identity \eqref{eq:CliffordR} then implies that
$\tilde R_-\rho_0(\alpha)\tilde R_+^{-1}=\rho_0(R\alpha)$ on $\Delta^+$,
for $\alpha\in\RR^4$, and since the isomorphism
$\det(\rho_0(\alpha)):\Lambda^2(\Delta^+)\cong\Lambda^2(\Delta^-)$
of complex lines is independent of the choice of unit-norm
$\alpha\in\RR^4$, we must have
$$
\det(\tilde R_+)=\det(\tilde R_-).
$$
Recall that $\Spin(4)=\SU(2)\times\SU(2)$. Thus, 
\begin{equation}
\label{eq:ReductionFromU(4)toSpinc(4)}
\begin{aligned}
\tilde R &\in \{(M,N)\in \U(2)\times\U(2): \det(M)=\det(N)\}
\\
&= (\SU(2)\times\SU(2)\times S^1)/\{\pm 1\}
\\
&= (\Spin(4)\times S^1)/\{\pm 1\} 
\\
&= \Spinc(4).
\end{aligned}
\end{equation}
Therefore, $\Fr_{\CCl(T^*X)}(W)\to X$ is a principal $\Spinc(4)$ bundle, as
claimed. 
\end{proof}

For $A\in\Spinc(4)$, $\Ad(A)$ denotes an element of $\SO(\End(\Delta))$.
However, we note that $\Ad(A)$ preserves $\rho(\RR^4)\subset\End(\Delta)$
and thus defines a linear transformation of $\RR^4$ which is an element
of $\SO(4)$.  This defines a homomorphism
$$
\Ad^c:\Spinc(4) \to \SO(4).
$$
By \cite[Corollary I.5.19]{LM} the homomorphism $\Ad^c$ is
equal to the projection $\Spinc(4)\to\SO(4)$
implicit in the description
of $\Spinc(4)$ in \eqref{eq:ReductionFromU(4)toSpinc(4)},
recalling that $\Spin(4)=(\SU(2)\times\SU(2))/\{\pm 1\}$.
The $S^1$ factor in the
presentation \eqref{eq:ReductionFromU(4)toSpinc(4)} of
$\Spin^c(4)$ is the kernel of the homomorphism $\Ad^c$. Since
$\Fr_{\CCl(T^*X)}(W)/S^1\cong \Fr(T^*X)$ by Lemma
\ref{lem:CliffordFrameBundleOverTX}, we have
\begin{equation}
\label{eq:SpincAssocS04}
\Fr_{\CCl(T^*X)}(W)\times_{(\Spinc(4),\Ad^c)}\RR^4 \cong T^*X.
\end{equation}
{}From the description of $\Spinc(4)$ in
\eqref{eq:ReductionFromU(4)toSpinc(4)}, we see that there is a homomorphism,
$$
{\det}^c:\Spinc(4)\to S^1, 
\quad
(M,N)\mapsto \det(M)=\det(N),
$$
and thus,
\begin{equation}
\label{eq:DetcLineBundle}
\Fr_{\CCl(T^*X)}(W)\times_{(\Spinc(4),{\det}^c)}\CC \cong\det(W^+).
\end{equation}
This completes our discussion of the Clifford frame bundle for $W$ and its
associated vector bundles.

We now determine the structure group of the bundle $\Fr_{\CCl(T^*X)}(V)$
over $X$, where $(\rho,V)$ is a complex-rank eight Hermitian Clifford
module over $X$. Suppose $\tilde F_1,\tilde F_2\in\Fr_{\CCl(T^*X)}(V_x)$
lie in the fibers over $F_1,F_2\in\Fr(T^*X|_x)$ of the projection
\eqref{eq:CliffordFrameBundlesWV}, so that $F_2^{-1}F_1=R=\Ad^c(\tilde
R)\in\SO(4)$, for some $\tilde R\in\Spinc(4)$. Therefore, according to
\eqref{eq:CliffordModuleHom}, we have
$$
\tilde F_2\rho_0(\alpha)\tilde F_2^{-1} = \rho(F_2\alpha),
\quad \alpha\in\CCl(\RR^4),
$$
and hence, replacing $\alpha$ by $R\alpha$ and applying 
identity \eqref{eq:CliffordR}, we see that
\begin{align*}
\tilde F_2\tilde R\rho_0(\alpha)(\tilde F_2\tilde R)^{-1} 
&= 
\tilde F_2\rho_0(R\alpha)\tilde F_2^{-1} 
\\
&=
\rho(F_2R\alpha)
\\
&=
\rho(F_1\alpha)
\quad \alpha\in\CCl(\RR^4),
\end{align*}
so that $\tilde F_2\tilde R\in\Fr_{\CCl(T^*X)}(V_x)$ also lies in the
fiber of the projection \eqref{eq:CliffordFrameBundlesWV} over $F_1\in
\Fr(T^*X|_x)$. Thus, $(\tilde F_2\tilde R)^{-1}\tilde F_1
=\tilde R^{-1}\tilde F_2^{-1}\tilde F_1$
is a Hermitian $\CCl(\RR^4)$-module
automorphism of $\Delta\otimes_\CC\CC^2$ and hence, by Lemma
\ref{lem:CliffModuleBasis}, we have
$$
(\tilde R^{-1}\otimes\id_{\CC^2})\circ\tilde F_2^{-1}\tilde F_1
=
\id_\Delta\otimes M,
$$
for some $M\in\U(2)$. Therefore, 
$$
\tilde F_1
=
\tilde F_2\circ(\tilde R\otimes M).
$$
Hence, any two elements of the fiber of $\Fr_{\CCl(T^*X)}(V)$ over $x\in X$
differ by an element of 
\begin{equation}
\label{eq:DefnSpinu(4)}
\begin{aligned}
\Spinu(4) &= (\Spinc(4)\times \U(2))/S^1
\\
&= \Spin(4)\times_{\{\pm 1\}}\U(2),
\\
&= (\SU(2)\times\SU(2))\times_{\{\pm 1\}}(\SU(2)\times_{\{\pm
1\}}S^1),
\end{aligned}
\end{equation}
noting that $\U(2)=(\SU(2)\times S^1)/\{\pm 1\}$. Therefore, we have

\begin{lem}
\label{lem:SpinuFrameBundleOverX}
The space $\Fr_{\CCl(T^*X)}(V)$ is a principal 
{\em $\Spinu(4)$} bundle over $X$.
\end{lem}

The composition of projection onto the factor $S^1$ of $\Spinu(4)$ in the
presentation \eqref{eq:DefnSpinu(4)} and squaring defines a homomorphism
$$
{\det}^u:\Spinu(4)\to S^1.
$$
Equivalently, writing elements of $\Spinu(4)$ as $\tilde R\otimes M$, with
$\tilde R\in\Spinc(4)$ and $M\in \U(2)$ using the first isomorphism in
\eqref{eq:DefnSpinu(4)}, this homomorphism is given by
\begin{equation}
\label{eq:SpinuDeterminant}
{\det}^u:\Spinu(4)\to S^1,
\quad
\tilde R\otimes M
\mapsto
{\det}^c(\tilde R)\det(M).
\end{equation}
Definition \eqref{eq:SpinuDeterminant}, the identification 
\eqref{eq:DetcLineBundle} of the line
bundle associated to ${\det}^c$, and the definition \eqref{eq:DefnDetVV} of
the line bundle $\det^{\frac{1}{2}}(V^+)$ imply that
\begin{equation}
\label{eq:DetuBundle}
\Fr_{\CCl(T^*X)}(V)\times_{(\Spinu(4),{\det}^u)}\CC \cong
{\det}^{\frac{1}{2}}(V^+).
\end{equation}
The homomorphisms
\begin{equation}
\label{eq:DefineSpinuHom}
\begin{aligned}
\Ad^u_{\SO(4)}&:\Spinu(4)\to \SO(4),
\quad \tilde R\otimes M \mapsto \Ad^c(\tilde R),
\\
\Ad^u_{\SO(3)}&:\Spinu(4)\to \SO(3),
\quad \tilde R\otimes M \mapsto \Ad(M),
\end{aligned}
\end{equation}
induce bundle isomorphisms
\begin{equation}
\label{eq:AssociatedSpinuBundles}
\begin{aligned}
\Fr_{\CCl(T^*X)}(V)\times_{(\Spinu(4),\Ad^u_{\SO(4)})}\RR^4 &\cong T^*X,
\\
\Fr_{\CCl(T^*X)}(V)\times_{(\Spinu(4),\Ad^u_{\SO(3)})}\su(2) &\cong
\fg_V.
\end{aligned}
\end{equation}
The first isomorphism in \eqref{eq:AssociatedSpinuBundles}
follows from the argument giving
the isomorphism \eqref{eq:SpincAssocS04}.
To prove the second isomorphism in \eqref{eq:AssociatedSpinuBundles}, 
observe that one has a bundle isomorphism
\begin{equation}
\label{eq:AdjointBundle1}
\Fr_{\CCl(T^*X)}(V)\times_{(\Spinu(4),\Ad)}\su(\Delta\otimes_\CC\CC^2)
\cong \su(V),
\end{equation}
defined via the adjoint representation,
$$
\Ad:\Spinu(4)\to \SO(\su(\Delta\otimes_\CC\CC^2)),
\quad
U\mapsto U(\,\cdot\,)U^\dagger.
$$
Recall from \cite[Equation (2.18)]{FL2a} that the $\SO(3)$ subbundle
$\fg_V\subset\su(V)$ is characterized as the
span of the sections which commute with the
sections $\rho(\om)$ of $\End_\CC(V)$ for all $\om\in\CCl(T^*X)$. In the
decomposition,
\begin{equation}
\label{eq:Decompose} 
\su(\Delta\otimes_\CC\CC^2) 
\cong
\su(\Delta)\oplus \su(\Delta)\otimes\su(2)\oplus \su(2),
\end{equation}
one has a similar characterization of the subspace
\begin{equation}
\label{eq:su(2)subspace}
\su(2)\cong\id_\Delta\otimes\su(2)\subset\su(\Delta\otimes_\CC \CC^2).
\end{equation}
Because any frame $\tilde F\in
\Fr_{\Cl(T^*X)}(V_x)$ defines an isomorphism of Hermitian Clifford modules,
$\Delta\otimes_\CC\CC^2\cong V_x$,
the induced isomorphism 
$\su(\Delta\otimes_\CC\CC^2)\cong \su(V_x)$ preserves the splitting
\eqref{eq:Decompose} and restricts to an isomorphism $\su(2)\cong
\fg_V|_x$. Hence, the isomorphism \eqref{eq:AdjointBundle1} gives
an isomorphism of subbundles,
\begin{equation}
\label{eq:AdjointBundle2}
\Fr_{\CCl(T^*X)}(V)\times_{(\Spinu(4),\Ad)} \su(2) \cong \fg_V.
\end{equation}
For $\tilde R\in\Spinc(4)$ and $M\in\U(2)$, one can write
$$
\Ad(\tilde R\otimes M)
=
\Ad^c(\tilde R)\oplus(\Ad^c(\tilde R)\otimes \Ad(M))\oplus \Ad(M)
$$
with respect to the decomposition \eqref{eq:Decompose},
the restrictions of $\Ad(\tilde R\otimes M)$ and $\Ad^u_{\SO(3)}(\tilde
R\otimes M)$ to the subspace $\su(2)$ in \eqref{eq:su(2)subspace}
agree, so the
isomorphism \eqref{eq:AdjointBundle2} yields the second bundle isomorphism in
\eqref{eq:AssociatedSpinuBundles}.

\subsection{Anti-self-dual and spin connections over the four-sphere}
\label{subsec:S4ModuliSp}
Assume $S^4\subset\RR^5$ has the standard embedding, with round metric of
radius one. Up to equivalence, there is a unique \spinc structure
$(\rho,\bW)$ over $S^4$ and unitary connection on $\bW$ which is spin with
respect to the Levi-Civita connection on $T^*S^4$ and Clifford map
$\rho:T^*S^4\to\End_\CC(\bW)$. The sphere $S^4$ has north pole
$n=(0,0,0,0,1)$, south pole $s=(0,0,0,0,-1)$, and local parametrizations 
\begin{equation}
\label{eq:SphereParams}
\begin{aligned}
\varphi_n:\RR^4&\to S^4\less\{s\}, \quad x\mapsto \varphi_n(x),
\\
\varphi_s:\RR^4&\to S^4\less\{n\}, \quad y\mapsto \varphi_s(y),
\end{aligned}
\end{equation}
defined by inverse stereographic projection from the south and north poles,
respectively.

Let $\bE\to S^4$ be a Hermitian, rank-two bundle with $c_2(\bE)=k\geq 1$
and set $\bV=\bW\otimes\bE$, yielding a complex rank-eight, Hermitian
Clifford module $(\rho,\bV)$. We may fix isomorphisms of complex line
bundles $\det(\bW^\pm) \cong S^4\times\CC$ and $\det(\bE) \cong S^4\times\CC$.

Let $\sA_{k}(S^4)$ be the space of spin connections on $\bV$, which induce
the product connection on $\det(\bV^+)\cong S^4\times\CC$, let $\sG_k$ be
the group of determinant-one, unitary automorphisms of $\bE$, identify $\sG_k$ 
with the group of \spinu automorphisms of $\bV$ via $u\mapsto
\id_{\bW}\otimes u$, and set
\begin{equation}
\label{eq:DefineFramedQuotientSpace}
\sB^s_{k}(S^4) = \sA_{k}(S^4)\times_{\sG_k}\Fr(\fg_{\bV}|_s)
\end{equation}
where $\Fr(\fg_{\bV})$ denotes the frame bundle for the $\SO(3)$ bundle
$\fg_{\bV}$. This is a principal $\SO(3)$ bundle over the quotient
space $\sB_k(S^4)$ of unframed spin connections on $\bV$.

A point $[\bA]\in \sB_k(S^4)$ has an associated `center' and `scale'
defined by \cite[pp. 343--344]{TauFrame}
\begin{equation}
\label{eq:DefnCenterScale}
\begin{aligned}
z[\bA]
&=
\frac{1}{8\pi^2 {k}}\int_{\RR^4} x|\varphi_n^*F(\hat\bA)|^2\,d^4x,
\\
\lambda^2[\bA]
&=
\frac{1}{8\pi^2 {k}}\int_{\RR^4} |x-z[\bA]|^2|\varphi_n^*F(\hat\bA)|^2\,d^4x.
\end{aligned}
\end{equation}
Let $M_{k}^s(S^4)\subset\sB_{k}^s(S^4)$ be the moduli space of pairs consisting
of frames for the fiber
$\fg_{\bV}|_s$ and `instanton' connections on $\bV$ --- the spin connections
$\bA$ on $\bV$ corresponding to anti-self-dual $\SO(3)$ connections
$\hat\bA$ on $\fg_{\bV}=\su(\bE)$.  Let $M^{\natural}_{k}(S^4) \subset
M^{s}_{k}(S^4)$ denote the subspace of `mass-centered' instantons, so
$z[\bA]=0$ for $[\bA]\in M^{\natural}_{k}(S^4)$ and
similarly define $M_{k}^{s,\natural}(S^4)$. Let
\begin{equation}
\label{eq:CenteredConcInstantons}
M^{\natural}_{k}(S^4,\delta) \subset M_{k}^{\natural}(S^4)
\end{equation}
denote the subspace of points $[\bA]\in M_{k}^{\natural}(S^4)$ with
scales $\lambda[\bA]$ in $(0,\delta]$, where $\delta$ is a small positive
constant,
and similarly define $M^{s,\natural}_{k}(S^4,\delta)$. Hence,
points in $M^{s,\natural}_{\ka}(S^4,\delta)$ have connections with curvature
density with center of mass at and energy concentrated near the north pole.

Now suppose $k=1$.  Because there is, up to gauge transformation and
rescaling, a unique mass-centered anti-self-dual connection on $\fg_{\bV}$,
there is diffeomorphism \cite[\S 3.1]{FL3}, \cite[\S 3]{TauFrame},
\begin{equation}
\label{eq:1InstModIsScaleFrame}
M^{s,\natural}_1(S^4,\delta)
\cong
(0,\delta]\times \SO(3).
\end{equation}
The Uhlenbeck compactification $\barM^{s,\natural}_1(S^4,\delta)$ of
$M^{s,\natural}_1(S^4,\delta)$ is obtained by adjoining the ideal point
$\{n\}\times[\Theta]\in S^4\times M_0(S^4)$, where $\Theta$ is
the tensor product of the Levi-Civita connection on $\bW$ with
the product connection on $S^4\times\CC^2$, and forgetting the frame
for $\fg_{\bV}|_s$. Thus, we have a homeomorphism
\begin{equation}
\label{eq:UhlenbeckInstantonConeHomeo}
\barM^{s,\natural}_1(S^4,\delta) \cong c(\SO(3)),
\end{equation}
where $c(\SO(3))$ is the cone on $\SO(3)$.

An element $R\in\SO(4)$ acts by rotation on $\RR^4$ and hence on
$S^4\subset\RR^5$ via the fixed stereographic projection from the south pole,
\begin{equation}
\label{eq:SphereStereoProjFromSouth}
S^4\less \{s\}\cong\RR^4,
\end{equation}
where the north pole is identified with the origin in $\RR^4$. The group
$\SO(4)$ in turn acts by pullback via \eqref{eq:SphereStereoProjFromSouth}
on sections of and spin connections on $\bV\to S^4$.

\subsection{Splicing Clifford modules, connections, and spinors}
\label{subsec:SplicingBasics}
In this section we define the process of splicing instantons from the
four-sphere onto $\SO(3)$ monopoles over $X$, together with the required
gluing data choices.

\subsubsection{Splicing Clifford modules}
Let $B(x_0,r_0)$ be a geodesic ball in $X$ with center $x_0\in X$ and
radius $r_0$. Choose a frame $F\in\Fr(T^*X|_{x_0})$, that is, an isomorphism
of oriented real inner-product spaces, 
\begin{equation}
\label{eq:TanCotanBundleIsomFiber}
F:\RR^4\cong T^*X|_{x_0}.
\end{equation}
This determines a dual isomorphism $\RR^4\cong TX|_{x_0}$ and an isomorphism of
Clifford algebras,
\begin{equation}
\label{eq:CliffAlgIsomFiber}
F:\CCl(\RR^4)\cong\CCl(T^*X|_{x_0}).
\end{equation}
Extend the frame \eqref{eq:TanCotanBundleIsomFiber} by parallel translation
along radial geodesics for the Levi-Civita connection on $(X,g)$ to a local
frame for $T^*X$ (and dual frame for $TX$) over $B(x_0,r_0)$ and geodesic,
normal coordinates on this ball, yielding a geodesic, normal coordinate
chart
\begin{equation}
\label{eq:ChartOnX}
\varphi_{x_0}^{-1}:B(x_0,r_0)\subset X\to \RR^4,
\end{equation}
and trivializations 
\begin{equation}
\label{eq:FInducedTrivTX}
T^*X|_{B(x_0,r_0)}\cong B(x_0,r_0)\times\RR^4
\quad\text{and}\quad
\CCl(T^*X)|_{B(x_0,r_0)}\cong B(x_0,r_0)\times\CCl(\RR^4),
\end{equation}
induced by parallel translation with respect to the Levi-Civita connection
along radial geodesics of the frame $F\in \Fr(T^*X|_{x_0})$.

Let $(\rho,V)$ be a complex-rank eight, Hermitian-Clifford module over
$(X,g)$, equipped with unitary connection $A$ which is spin with respect to
the Levi-Civita connection on $T^*X$ and induces the fixed connection
$2A_\Lambda$ on the fixed determinant line bundle $\det(V^+)$, where
$\Lambda=\det^{\frac{1}{2}}(V^+)$.  Choose a frame $\tilde
F\in\Fr(V_{x_0})$ in the fiber of the bundle
\eqref{eq:CliffordFrameBundlesWV} over $F$,  
so $\tilde F$ is an isomorphism of Hermitian-Clifford modules with respect
to the isomorphism
\eqref{eq:CliffAlgIsomFiber},
\begin{equation}
\label{eq:RedCliffModIsomFiber}
\tilde F:\Delta\otimes_\CC\CC^2\cong V_{x_0}.
\end{equation}
It will be useful to note the

\begin{lem}
\label{lem:TrivIsCliffModIsom}
Let $(U,\rho)$ be a $\CCl(T^*Y)$-module, of arbitrary complex rank, over an
oriented, Riemannian manifold $(Y,g)$. Then parallel translation with
respect to a spin connection on $U$ gives $\CCl(T^*Y)$-module isomorphisms
of the fibers of $U\to Y$, with respect to the isomorphisms of the fibers
of $T^*Y\to Y$ given by parallel translation with respect to the
Levi-Civita connection on $T^*Y$.
\end{lem}

\begin{proof}
Suppose $\gamma$ is a smooth curve in $Y$ with initial point $y_0$. Let
$\alpha\in\Omega^1(Y,\CC)$, let $\Phi\in C^\8(Y,U)$, and let $\cov$ be a
spin connection on $U$. Because $\cov$ is a Clifford module derivation of
$C^\8(Y,U)$, we have (see equation \eqref{eq:SpinConnection})
$$
\cov_{\dot\gamma}(\rho(\alpha)\Phi)
=
\rho(\cov_{\dot\gamma}\alpha)\Phi
+
\rho(\alpha)\cov_{\dot\gamma}\Phi,
$$ 
where $\cov$ also denotes the Levi-Civita connection on $T^*X$. Thus, if
$\alpha$ and $\Phi$ are parallel along $\gamma$ with respect to $\cov$,
then $\rho(\alpha)\Phi$ must be parallel along $\gamma$ with respect to
$\cov$.  Let $P_y$ denote parallel translation, with respect to the
Levi-Civita connection on $T^*Y$ or spin connection on $U$, along $\gamma$
from $y_0$ to a point $y\in Y$. Therefore, $\alpha_y=P_y\alpha_0$ and
$\Phi_y=P_y\Phi_0$. The sections of $U$ over $\gamma\subset Y$ given by
$$
\rho_y(P_y\alpha_0)P_y\Phi_0 \quad\text{and}\quad
P_y(\rho_0(\alpha_0)\Phi_0), 
\quad y\in \gamma,
$$
agree when $y=y_0$ and are parallel along the curve $\gamma$ from
$y_0$. Hence, parallel translation along curves in $Y$ with respect to spin
connections commutes with Clifford multiplication and thus yields
Clifford-module isomorphisms of the fibers of $U\to Y$.
\end{proof}

Therefore, given a spin connection $A$ on $V$ and a choice of frame $\tilde
F:\Delta\otimes_\CC\CC^2\cong V_{x_0}$, parallel translation with respect
to $A$ along radial geodesics in $X$ from $x_0$ yields a Clifford-module
isomorphism,
\begin{equation}
\label{eq:TrivSpinRank8X}
V|_{B(x_0,r_0)}\cong B(x_0,r_0)\times\Delta\otimes_\CC\CC^2,
\end{equation}
with respect to the Clifford algebra isomorphism \eqref{eq:FInducedTrivTX}.
The trivialization \eqref{eq:TrivSpinRank8X} also induces a trivialization
\begin{equation}
\label{eq:TrivgV}
\fg_{V}|_{B(x_0,r_0)}\cong B(x_0,r_0)\times\su(2).
\end{equation}
This completes our discussion of the Clifford modules and bundles over $X$.

Over $S^4$, we have the standard \spinc structure $(\rho,\bW)$ with
Clifford multiplication $\rho:T^*S^4\to\End_\CC(\bW)$ and unitary
connection which is spin with respect to the Levi-Civita connection on
$T^*S^4$ for the standard round metric of radius one. We may fix, once and
for all, a unit-norm $\CCl(T^*S|_n)$-frame $\Delta\cong\bW|_{n}$ for
$\bW|_{n}$ and a trivialization defined by parallel radial translation
which is a Clifford-module isomorphism (see Lemma
\ref{lem:TrivIsCliffModIsom}),
\begin{equation}
\label{eq:TrivSpinRank4Sphere}
\bW|_{S^4\less\{s\}} \cong \RR^4\times \Delta,
\end{equation}
with respect to the Clifford algebra isomorphism
\begin{equation}
\label{eq:SphereLessSCliffAlgIsom}
\CCl(T^*S^4)|_{S^4\less\{s\}}\cong S^4\less\{s\}\times\CCl(\RR^4).
\end{equation}
A spin connection $\bA$ on $\bV=\bW\otimes_\CC\bE\to S^4$ and a choice of
frame $\su(2)\cong\fg_{\bV|_{s}}$, that is, an
isomorphism of oriented real inner product
spaces, define a trivialization of $\SO(3)$ bundles,
\begin{equation}
\label{eq:TrivSO(3)Sphere}
\fg_{\bV}|_{S^4\less\{n\}}\cong\RR^4\times\su(2),
\end{equation}
via parallel translation with respect to the $\SO(3)$ connection $\hat\bA$
on $\fg_{\bV}$.  Up to an ambiguity corresponding to the action of $\{\pm
1\}$, the frame $\su(2)\cong\su(\bE)|_{s}=\fg_{\bV|_{s}}$ lifts to a frame
$\CC^2\cong\bE_{s}$, an isomorphism of complex inner product spaces.
This ambiguity will vanish when we take the quotient by the gauge group.
Thus, up to this ambiguity, we obtain an isomorphism lifting
\eqref{eq:TrivSO(3)Sphere},
\begin{equation}
\label{eq:TrivHermRank2Sphere}
\bE|_{S^4\less\{n\}}\cong\RR^4\times\CC^2,
\end{equation}
coinciding with the trivialization defined by the isomorphism
$\CC^2\cong\bE_{s}$ and radial parallel translation with respect to the
unique lift of $\hat\bA$ to an $\SU(2)$ connection on $\bE$.

The isomorphisms \eqref{eq:TrivSpinRank4Sphere} and
\eqref{eq:TrivHermRank2Sphere} combine to give a Clifford-module
isomorphism
\begin{equation}
\label{eq:TrivSpinRank8Sphere}
\bV|_{S^4\less\{n,s\}} \cong \RR^4\less\{0\}\times \Delta\otimes_\CC\CC^2,
\end{equation}
with respect to the Clifford algebra isomorphism 
\eqref{eq:SphereLessSCliffAlgIsom}.

We then define a \spinu structure $\ft'=(\rho,V')$ over $(X,g)$ by setting
\begin{equation}
\label{eq:SplicedSpinu}
V'= (V|_{X\less\{x_0}\})\# (\bV|_{S^4\less\{s\}}),
\end{equation}
where the bundles $V$ and $\bV$ are identified
via the trivializations \eqref{eq:TrivSpinRank8X}
and \eqref{eq:TrivSpinRank8Sphere} and the embedding
\begin{equation}
\label{eq:BallInXToBallInSphere}
\varphi_n\circ\varphi_{x_0}^{-1}:B(x_0,r_0)\to S^4\less\{s\},
\end{equation}
identifying a ball in $X$ centered at $x_0$ with a ball in $S^4$ centered
at the north pole. Because these trivializations are Clifford module
isomorphisms, the bundle $V'$ has a $\CCl(T^*X)$-module structure. The
Clifford module \eqref{eq:SplicedSpinu} induces an associated
$\SO(3)$-bundle $\fg_{V'}$,
\begin{equation}
\label{eq:SplicedSO(3)Bundle}
\fg_{V'} = (\fg_V|_{X\less\{x_0\}})\# (\fg_{\bV}|_{S^4\less\{s\}}),
\end{equation}
where the bundles $\fg_V$ and $\fg_{\bV}$ are identified
over small annuli centered at $x\in X$ and $n\in S^4$
by the trivializations \eqref{eq:TrivgV} and
\eqref{eq:TrivSO(3)Sphere} and the embedding
\eqref{eq:BallInXToBallInSphere}. The characteristic classes of
$\ft=(\rho,V)$ and $\ft'=(\rho,V')$ are related by
\begin{equation}
\label{eq:HowTheCharClassesChangeUnderSplicing}
p_1(\fg_{V'})=p_1(\fg_{V})-4,
\quad
c_1(V')=c_1(V),
\quad\text{and}\quad
w_2(\fg_{V'})=w_2(\fg_{V}),
\end{equation}
as we shall later exploit, often without comment.

We give a more detailed description of where the
bundles $V'$ and $\fg_{V'}$ are spliced
below. Define a smooth cutoff function $\chi_{x_0,\eps}:X\to [0,1]$ by setting
\begin{equation}
\label{eq:ChiCutoffFunctionDefn}
\chi_{x_0,\eps}(x) 
:= 
\chi(\dist(x,x_0)/\eps),
\end{equation}
where $\chi:\RR\to [0,1]$ is a smooth function such that $\chi(t)=1$
for $t\ge 1$ and $\chi(t)=0$ for $t\le 1/2$.  Thus, we have
$$
\chi_{x_0,\eps}(x) 
=
\begin{cases}
1 &\text{for } x\in X - B(x_0,\eps),
\\
0 &\text{for } x\in B(x_0,\eps/2).
\end{cases}
$$
The orientation-preserving embedding \eqref{eq:BallInXToBallInSphere} therefore
identifies the annulus
$\varphi_n(\Omega(0,\thalf\sqrt{\lambda},2\sqrt{\lambda}))$ in $S^4$, where
$$
\Omega\left(0,\thalf\sqrt{\lambda},2\sqrt{\lambda}\right)
:=
\left\{x\in\RR^4:\thalf\sqrt{\lambda}<|x|<2\sqrt{\lambda}\right\}
\subset \RR^4,
$$
with the annulus in $X$,
$$
\Omega\left(x_0,\thalf\sqrt{\lambda},2\sqrt{\lambda}\right)
:=
\left\{x\in X: \thalf\sqrt{\lambda} < \dist_g(x,x_0) < 2\sqrt{\lambda} 
\right\}
\subset X.
$$
Hence, the spliced bundle \eqref{eq:SplicedSpinu}
is defined explicitly by setting
$$
V'
= 
\begin{cases}
V &\text{over $X - B(x_0,\half\sqrt{\lambda})$,} \\
\bV &\text{over $B(x_0,2\sqrt{\lambda})$.} 
\end{cases}
$$
The bundles $V$ and $\bV$ are identified over the annulus
$\Omega(x_0,\half\sqrt{\lambda},2\sqrt{\lambda})$ in $X$ via
the orientation-preserving diffeomorphism \eqref{eq:BallInXToBallInSphere}
identifying the annulus
$\Omega(x_0,\half\sqrt{\lambda},2\sqrt{\lambda})$ with the
corresponding annulus in $S^4$ and the bundle map defined by the
trivializations \eqref{eq:TrivSpinRank8X}
and \eqref{eq:TrivSpinRank8Sphere}.

Similarly, the spliced $\SO(3)$ bundle \eqref{eq:SplicedSO(3)Bundle} 
is defined explicitly by setting
$$
\fg_{V'}
= 
\begin{cases}
\fg_{V} &\text{over $X - B(x_0,\half\sqrt{\lambda})$,} \\
\fg_{\bV} &\text{over $B(x_0,2\sqrt{\lambda})$.} 
\end{cases}
\label{eq:DefnGluedUpSO(3)Bundle}
$$
The bundles $\fg_{V}$ and $\fg_{\bV}$ are identified over the annulus
$\Omega(x_0,\half\sqrt{\lambda},2\sqrt{\lambda}\})$ in $X$ via 
the diffeomorphism \eqref{eq:BallInXToBallInSphere} and
the $\SO(3)$ bundle map defined by the trivializations \eqref{eq:TrivgV} and
\eqref{eq:TrivSO(3)Sphere}.

\subsubsection{Splicing connections and cutting off background spinors}
To construct a spliced spin connection on the Clifford module $V'$, we
choose the following data:
\begin{itemize}
\item
An $\SO(3)$ connection $\hat A$ on $\fg_V\to X$,
\item
An $\SO(4)$ frame $F\in\Fr(T^*X|_{x_0})$, and induced chart
$B(x_0,r_0)\subset X\to\RR^4$ as in \eqref{eq:ChartOnX},
\item
An $\SO(3)$ frame for $\Fr(\fg_{V}|_{x_0})$, and
corresponding trivialization $\fg_{V}|_{B(x_0,r_0)}\cong
B(x_0,r_0)\times\su(2)$ as in \eqref{eq:TrivgV},
\item
An $\SO(3)$ connection $\hat\bA$ on $\fg_{\bV}\to S^4$, and associated
scale parameter $\lambda\in (0,\delta]$,
\item
An $\SO(3)$ frame for $\fg_{\bV}|_{s}$, and corresponding
trivialization $\fg_{\bV}|_{S^4\less\{n\}}\cong S^4\less\{n\}\times\su(2)$
as in \eqref{eq:TrivSO(3)Sphere},
together with the fixed chart $S^4\less\{s\}\cong\RR^4$
as in \eqref{eq:SphereParams}.
\end{itemize}
This yields a spliced connection $\hat A'$ on the spliced $\SO(3)$ bundle
$\fg_{V'}$ over $X$, as in \eqref{eq:SplicedSO(3)Bundle}. 

We give a more detailed description of where the $\SO(3)$ connections
are spliced below. 
We first define a cut-off $\SO(3)$ connection on the bundle $\fg_{V}$ over
$X$ by  setting
\begin{equation}
\label{eq:CutOffBackgroundConnection}
\chi_{x_0,4\sqrt{\lambda}}\hat A
:=
\begin{cases}
\hat A
&\text{over $X - B(x_0,4\sqrt{\lambda})$}, 
\\
\Gamma + \chi_{x_0,4\sqrt{\lambda}}\sigma_0^*\hat A
&\text{over $\Omega(x_0,2\sqrt{\lambda},4\sqrt{\lambda})$,} 
\\
\Gamma 
&\text{over $B(x_0,2\sqrt{\lambda_0})$,} 
\end{cases} 
\end{equation}
where $\Gamma$ denotes the product connection on $B(x_0,r_0)\times \SO(3)$
and $\sigma_0$ is the section of $\fg_{V}$ over $B(x_0,r_0)$
defined by the trivialization \eqref{eq:TrivgV}. 

We define a cut-off $\SO(3)$
connection on the bundle $\fg_{\bV}$ over $S^4$, with mass center at the
north pole, by setting
\begin{equation}
(1-\chi_{n,\sqrt{\lambda}/2})\hat\bA
:=
\begin{cases}
\hat\bA
&\text{over $\varphi_n(B(0,\quarter\sqrt{\lambda}))$}, 
\\
\Gamma + (1-\chi_{n,\sqrt{\lambda}/2})\tau_0^*\hat\bA
&\text{over $\varphi_n(\Omega(0,\quarter\sqrt{\lambda},
\half\sqrt{\lambda}))$,} 
\\
\Gamma 
&\text{over $S^4 - \varphi_n(B(0,\half\sqrt{\lambda}))$,} 
\end{cases} 
\end{equation}
where $\Gamma$ denotes the product connection on $(S^4-\{n\})\times \SO(3)$
and $\tau_0$ is the section of $\fg_{\bV}$ over $S^4-\{n\}$
defined by the trivialization \eqref{eq:TrivSO(3)Sphere}. 

Hence, we define a spliced connection $\hat A'$ on $\fg_{V'}$ by setting
\begin{equation}
\label{eq:SplicedConnection}
\hat A'
:= 
\begin{cases}
  \hat A &\text{over $X - B(x_0,4\sqrt{\lambda})$},
  \\
  \Gamma +
  \chi_{x_0,4\sqrt{\lambda}}\hat A +
  (1-\chi_{x_0,\sqrt{\lambda}/2})\hat\bA &\text{over
    $\Omega(x_0,\quarter\sqrt{\lambda},4\sqrt{\lambda})$},
  \\
  \hat\bA &\text{over $B(x_0,\quarter\sqrt{\lambda})$,}
\end{cases} 
\end{equation}
where the cut-off connections are defined as above; the bundle and annulus
identifications are understood.

According to Lemma 2.11 in \cite{FL2a},
we then obtain a unitary connection $A'$ on $V'$ defined by the $\SO(3)$ 
connection $\hat A'$ on $\fg_{V'}$ and by the
\begin{itemize}
\item
Frame $\tilde F\in\Fr_{\CCl(T^*X)}(V_{x_0})$ covering $F\in
\Fr(T^*X|_{x_0})$, and corresponding Clifford-module isomorphism
$V|_{B(x_0,r_0)}\cong B(x_0,r_0)\times\Delta\otimes_\CC\CC^2$ as in
\eqref{eq:TrivSpinRank8X},
\item
Fixed $\U(1)$ connection $2A_\Lambda$ on $\det({V'}^+)$, and the
\item
Requirement that $A'$ be spin with respect to the Levi-Civita connection on
$T^*X$ for the metric $g$ and Clifford map $\rho':T^*X\to\End_\CC(V')$.
\end{itemize}

Finally, given a section $\Phi$ of $V\to X$, we obtain a section $\Phi'$ of
$V'\to X$ by cutting off $\Phi$ in a ball around $x_0\in X$, with radius
determined by $\lambda$. 

More explicitly, we define the spinor $\Phi'$ on the
bundle $V'$ over $X$ by setting
\begin{equation}
\label{eq:SplicedSpinor}
\Phi' \equiv \chi_{x_0,8\lambda^{1/3}}\Phi
:=
\begin{cases}
\Phi
&\text{over $X - B(x_0,8\lambda^{1/3})$}, 
\\
\chi_{x_0,8\lambda^{1/3}}\Phi
&\text{over $\Omega(x_0,4\lambda^{1/3},8\lambda^{1/3})$,} 
\\
0
&\text{over $B(x_0,4\lambda^{1/3})$.} 
\end{cases} 
\end{equation}
The reason for the different choice of annulus
radii is explained in \cite{FL3}.

\subsubsection{Splicing in spinors from the four-sphere}
\label{subsubsec:SpliceSpinorsFromSphereToX}
To construct a spliced section of the Clifford module $V'$, given a spinor
over the four-sphere, we choose the following data:
\begin{itemize}
\item
A spin connection $A$ on $V\to X$,
\item
An $\SO(4)$ frame $F\in \Fr(T^*X|_{x_0})$ and induced chart
$B(x_0,r_0)\subset X\to\RR^4$ as in \eqref{eq:ChartOnX},
\item
A frame $\tilde F\in\Fr_{\CCl(T^*X)}(V_{x_0})$ covering $F\in
\Fr(T^*X|_{x_0})$, and corresponding Clifford-module isomorphism
$V|_{B(x_0,r_0)}\cong B(x_0,r_0)\times\Delta\otimes_\CC\CC^2$ as in
\eqref{eq:TrivSpinRank8X},
\item
A spin connection $\bA$ on $\bV\to S^4$, and associated
scale parameter $\lambda\in (0,\delta]$,
\item
An oriented, orthonormal frame for $\bE_{s}$, and corresponding
trivialization $\bE|_{S^4\less\{n\}}\cong S^4\less\{n\}\times\CC^2$,
together with the fixed chart $S^4\less\{s\}\cong\RR^4$
as in \eqref{eq:SphereParams}, fixed unit-norm
$\CCl(T^*S^4|_{n})$-frame for $\bW|_{n}$, and fixed Clifford-module
isomorphism $\bW|_{S^4\less\{s\}}\cong S^4\less\{s\}\times\Delta$. The
bundle isomorphisms combine to give a Clifford-module isomorphism
$\bV|_{S^4\less\{n,s\}}\cong S^4\less\{n,s\}\times\Delta\otimes\CC^2$,
which we employ in conjunction with the chart $S^4\less\{s\}\cong\RR^4$.
\end{itemize}
Spinors $\Psi$ on $\bV\to S^4$ are then spliced from $\RR^4\cong
S^4\less\{s\}$ onto 
the ball $B(x_0,r_0)\subset X$, via the preceding data, with the north pole
$n\in S^4$ being identified with $x_0\in X$.

More explicitly, we define the spinor $\Psi'$
on the bundle $V'$ over $X$ by setting
\begin{equation}
\label{eq:DefnSplicedSpinorFromS4}
\begin{aligned}
\Psi'
&:=
(1-\chi_{n,\sqrt{\lambda}/2})\Psi
\\
&:=
\begin{cases}
\Psi
&\text{over $\varphi_n(B(0,\quarter\sqrt{\lambda}))$}, 
\\
(1-\chi_{n,\sqrt{\lambda}/2})\Psi
&\text{over $\varphi_n(\Omega(0,\quarter\sqrt{\lambda},
\half\sqrt{\lambda}))$,} 
\\
0
&\text{over $X - B(x_0,\half\sqrt{\lambda})$,} 
\end{cases} 
\end{aligned}
\end{equation}
where the identification of the bundles $V$ over the annulus
$\Omega(x_0,\quarter\sqrt{\lambda},\half\sqrt{\lambda})$ in $X$
and the bundle $\bV$
over the annulus 
$\varphi_n(\Omega(0,\quarter\sqrt{\lambda},\half\sqrt{\lambda}))$ in $S^4$,
together with these annuli, are described following the definition 
\eqref{eq:SplicedSpinu} of $V'$.

\subsection{Space of gluing data and the definition of the splicing map}
\label{subsec:SplicingMap}
The patching constructions described in \S \ref{subsec:SplicingBasics}
yield a pre-splicing map
\begin{equation}
\label{eq:PreSplicingMap1}
\tilde\bgamma':
\tilde N_{\ft,\fs}(\eps)\times\Fr(\fg_V)\times_X\Fr(T^*X)
\times\Fr(\fg_{\bV}|_s)\times \tilde M_1^{\natural}(S^4,\delta)
\to
\tilde\sC_{\ft'},
\end{equation}
where
$$
\tilde M_1^{\natural}(S^4,\delta) 
\subset 
\sA_1(S^4)
$$
is the preimage of the subspace $M_1^{\natural}(S^4,\delta)$ under the
projection $\sA_1(S^4)\to \sB_1(S^4)$.

The map \eqref{eq:PreSplicingMap1} is equivariant with respect to the
\begin{enumerate}
\item
Diagonal action of
$\sG_{\fs}$ on $\tilde N_{\ft,\fs}(\eps)\times\Fr(\fg_{V})$,
where $\sG_{\fs}$ acts on $\tilde N_{\ft,\fs}$ as described
following equation \eqref{eq:VirtualNormalBundleTotalSpace} and
on $\Fr(\fg_V)$ by the inclusion \eqref{eq:GaugeGroupInclusion}
of $\sG_{\fs}$ in $\sG_{\ft}$,
\item
Diagonal action of $\sG_1$ on $\Fr(\fg_{\bV}|_s)\times \tilde
M_1^{\natural}(S^4,\delta)$,
\item
Diagonal action of $\SO(4)$ on
$\Fr(T^*X)\times \tilde M_1^{\natural}(S^4,\delta)$
(where the action of $\SO(4)$ is as described following
\eqref{eq:SphereStereoProjFromSouth}),
\item
Diagonal action of
$\SO(3)$ on $\Fr(\fg_{V})\times\Fr(\fg_{\bV}|_s)$.
\end{enumerate}
Therefore, if we define the space
\begin{equation}
\label{eq:PreGluingDataSpace}
\tilde{\Gl}_{\ft}(\delta)
=
\left(
\Fr(\fg_V)\times_X\Fr(T^*X)
\times
\Fr(\fg_{\bV}|_s)\times \tilde M_1^{\natural}(S^4,\delta)
\right)/(\SO(3)\times\SO(4)),
\end{equation}
where $\SO(3)\times\SO(4)$ act as described prior to
\eqref{eq:PreGluingDataSpace}, the pre-splicing map \eqref{eq:PreSplicingMap1}
descends to a map
\begin{equation}
\label{eq:PreSplicingMap}
\tilde\bga':
    \tilde N_{\ft,\fs}(\eps)\times\tilde{\Gl}_{\ft}(\delta)
    \to
    \tsC^*_{\ft'}.
\end{equation}
The map \eqref{eq:PreSplicingMap} is gauge equivariant and thus,
if we define
\begin{equation}
\label{eq:DefnXGluingParameters}
\begin{aligned}
\Gl_{\ft}(\delta)
&=
\tilde{\Gl}_{\ft}(\delta)/\sG_1
\\
&=
\left(\Fr(\fg_{V})\times_X\Fr(T^*X)\times M_1^{s,\natural}(S^4,\delta)\right)
/(\SO(3)\times\SO(4)),
\end{aligned}
\end{equation}
(recalling that $M_1^{s,\natural}(S^4,\delta)
=\Fr(\fg_{\bV}|_s)\times_{\sG_1} \tilde M_1^{\natural}(S^4,\delta)$),
the pre-splicing map descends to the gauge group quotients
and gives the {\em splicing map\/}
\begin{equation}
\label{eq:SplicingMap}
\bgamma':
\sM_{\ft',\fs}^{\stab}
\to
\sC_{\ft'},
\end{equation}
where we have defined
\begin{equation}
\label{eq:TopStratumOfSplicingDomain}
\sM_{\ft',\fs}^{\stab}
=
\tN_{\ft,\fs}(\eps)\times_{\sG_{\fs}}\Gl_{\ft}(\delta).
\end{equation}
We refer to the space \eqref{eq:TopStratumOfSplicingDomain} as
a {\em virtual moduli space\/}; it is defined by the
stabilizing bundle used to construct $N_{\ft,\fs}$.
It can be shown \cite{FL3, FL4} that $\bgamma'$ is a smooth
embedding, provided $\eps$ and $\delta$ are sufficiently small.

By analogy with the definition of the Uhlenbeck compactification in \S
\ref{subsubsec:PU2Monopoles}, we define 
\begin{equation}
\label{eq:CompactifiedGluingBundle}
\bar{\Gl}_{\ft}(\delta)
=
\left(\Fr(\fg_{V})\times_X\Fr(T^*X)\times
\barM_1^{s,\natural}(S^4,\delta)\right)
/(\SO(3)\times\SO(4)).
\end{equation}
Points in $\bar{\Gl}_{\ft}(\delta)- \Gl_{\ft}(\delta)$ then correspond to
the cone point in \eqref{eq:UhlenbeckInstantonConeHomeo}.  The
action of $\SO(4)\times\SO(3)$ in \eqref{eq:CompactifiedGluingBundle}
is trivial on these cone points, so
$\bar{\Gl}_{\ft}(\delta)- \Gl_{\ft}(\delta)=X$. The cone completion of
$\Gl_{\ft}(\delta)$ then gives an `Uhlenbeck compactification'
\begin{equation}
\label{eq:UhlenbeckCompactifiedSplicingDomain}
\bar\sM_{\ft',\fs}^{\stab}
=
\tN_{\ft,\fs}(\eps)\times_{\sG_{\fs}}\bar{\Gl}_{\ft}(\delta)
\end{equation}
of the space $\sM_{\ft',\fs}^{\stab}$, where $\sG_{\fs}$ acts trivially
on $\bar{\Gl}_{\ft}(\delta)-\Gl_{\ft}(\delta)$.

The splitting $\fg_V\cong\ubarRR\oplus L$, induced by the splitting
$V=W\oplus W\otimes L$, and the circle action given in
\eqref{eq:DefineS1LAction} induce a circle action on $\Fr(\fg_V)$ and thus
on $\bar{\Gl}_{\ft}(\delta)$ by scalar multiplication on $L$.  If we write this
action as $(e^{i\theta},\bg)\mapsto e^{i\theta}\bg$, where $e^{i\theta}\in
S^1$ and $\bg\in\bar{\Gl}_{\ft}(\delta)$, then the description of the action of
$\sG_{\fs}$ on $\fg_V$ following equation
\eqref{eq:VirtualNormalBundleTotalSpace} implies that
\begin{equation}
\label{eq:GaugeGrpS1OnGluing}
\varrho(e^{i\theta})\bg
=
e^{-2i\theta}\bg,
\quad\text{for $e^{i\theta}\in S^1$ and $\bg\in\bar{\Gl}_{\ft}(\delta)$}.
\end{equation}
If we set
$$
\bar{\sC}_{\ft'}
=
\bigsqcup_{\ell=0}^\8(\sC_{\ft'_{\ell}}\times\Sym^\ell(X)),
$$
and extend the splicing map $\bga'$ to a map
\begin{equation}
\label{eq:ExtendedSplicing}
\bgamma':
\bar\sM_{\ft',\fs}^{\stab}
\to
\bar{\sC}_{\ft'},
\end{equation}
by setting $\bga'$ equal to
$\bga_{\fs}\times\id_X: N_{\ft,\fs}(\eps)\times X\to \sC_{\ft}\times X$
on the cone points, where $\bga_{\fs}$ is defined in
\eqref{eq:BackgroundConfigEmbedding},
then this extension is continuous with respect to the Uhlenbeck topology.  
The extended splicing map
\eqref{eq:ExtendedSplicing} is $S^1$-equivariant if the circle acts
by the action \eqref{eq:S1ZAction} on $\bar{\sC}_{\ft'}$ and diagonally
with weight two on each factor in the product $\tilde
N_{\ft,\fs}\times\bar{\Gl}_{\ft}(\delta)$: 

\begin{lem}
\label{lem:S1EquivarianceSplicing}
The following two circle actions on
$\tilde N_{\ft,\fs}(\eps)\times_{\sG_{\fs}}\bar{\Gl}_{\ft}(\delta)$
are equivalent:
\begin{enumerate}
\item
The action \eqref{eq:S1ZActionOnN} on $\tN_{\ft,\fs}(\eps)$ and
the trivial action on $\bar{\Gl}_{\ft}(\delta)$, and
\item
The diagonal action
with weight two on the fibers of $\tilde N_{\ft,\fs}(\eps)$
and with weight two on $\bar{\Gl}_{\ft}(\delta)$.
\end{enumerate}
Furthermore, the extended splicing and gluing maps
$$
\bga,\bga':
\tilde N_{\ft,\fs}(\eps)\times_{\sG_{\fs}}\bar{\Gl}_{\ft}(\delta)
\to
\bar{\sC}_{\ft'}
$$
are circle-equivariant if the circle acts on $\bar{\sC}_{\ft'}$
by the action \eqref{eq:S1ZAction} and
on the domain by either of the above equivalent actions.
\end{lem}

\begin{proof}
The embedding $\tilde N_{\ft,\fs}(\eps)\to\tsC_{\ft}$ is $S^1$-equivariant
if the circle acts on $\tsC_{\ft}$ by the action \eqref{eq:S1ZAction} and
on $\tilde N_{\ft,\fs}(\eps)$ by the composition of multiplication with
weight two on the fibers of $\tilde N_{\ft,\fs}(\eps)$ and the circle
action \eqref{eq:S1ZActionOnN}.  The splicing map is $S^1$-equivariant if
the circle acts by scalar multiplication on the sections both before and
after cutting off the sections. Hence, the splicing map $\bgamma'$ is
$S^1$-equivariant if the circle acts on $\sC_{\ft'}$ by the action
\eqref{eq:S1ZAction} and by the action (1) in the lemma statement on the
domain, namely: for $(B,\Psi,\eta)\in\tN_{\ft,\fs}$ and
$\bg\in\Gl_{\ft}(\delta)$, by
$$
[(B,\Psi,\eta),\bg] 
\mapsto 
[\varrho(e^{i\theta})(B,\Psi,e^{i2\theta}\eta),\bg].
$$
Because $\sG_{\fs}$ acts diagonally in the definition of $\tilde
N_{\ft,\fs}(\eps)\times_{\sG_{\fs}}\bar{\Gl}_{\ft}(\delta)$, the action (1)
on the domain is equal to the circle action given by scalar multiplication
with weight two on the fibers of $\tilde N_{\ft,\fs}(\eps)$ and by the
action of the constant $S^1$ in $\sG_{\fs}$ with weight negative one on the
factor $\bar{\Gl}_{\ft}(\delta)$.  This last action is equal to the action
(2) on the domain in the statement of the lemma by the equation
\eqref{eq:GaugeGrpS1OnGluing} comparing the action of the
constant $S^1$ subgroup of $\sG_{\fs}$ on
$\bar{\Gl}_{\ft}(\delta)$ with the standard action:
$$
[\varrho(e^{i\theta})(B,\Psi,e^{i2\theta}\eta),\bg]
=
[(B,\Psi,e^{i2\theta}\eta),\varrho(e^{-i\theta})\bg]
=
[(B,\Psi,e^{i2\theta}\eta),e^{2i\theta}\bg].
$$
This completes the proof.
\end{proof}

We note that there is a projection map,
\begin{equation}
\label{eq:GluingDataProjection}
\tilde N_{\ft,\fs}(\eps)\times_{\sG_{\fs}\times S^1}
\bar{\Gl}_{\ft}(\delta)
\to
M_{\fs}\times X,
\end{equation}
defined via the projection maps $\bar{\Gl}_{\ft}(\delta)\to X$ (given by
the fiber bundle structure in \eqref{eq:CompactifiedGluingBundle}) and
$N_{\ft,\fs}\to M_{\fs}$.

\subsection{The obstruction bundle}
\label{subsec:ObstructionBundle}
The $\SO(3)$-monopole obstruction bundle is a finite-rank,
$S^1$-equivariant subbundle,
\begin{equation}
\label{eq:ObstructionBundleEmbedding}
\begin{CD}
\Upsilon_{\ft',\fs} @>{\bvarphi}>>\fV_{\ft'}
\\
@VVV @VVV
\\
\tilde N_{\ft,\fs}(\eps)\times_{\sG_{\fs}}\Gl_{\ft}(\delta)
@>{\bgamma'}>>\sC_{\ft'}^{*,0}
\end{CD}
\end{equation}
of the $S^1$-equivariant vector bundle
\begin{equation}
\label{eq:InfiniteDimObstruction}
\fV_{\ft'}
=
\tsC^{*,0}_{\ft'}\times_{\sG_{\ft'}}
L^2_{k-1}(\La^+\otimes\fg_{V'} \oplus {V'}^-)
\to \sC^{*,0}_{\ft'}.
\end{equation}
The circle actions in \eqref{eq:InfiniteDimObstruction}
are given by the diagonal circle action, with the action
\eqref{eq:S1ZAction} on $\tsC^{*,0}_{\ft'}$, the
trivial action on $L^2_{k-1}(\La^+\otimes \fg_{V'})$, and by scalar
multiplication on $L^2_{k-1}(V')$.  As described in \cite{FL3}, the fiber
$\Upsilon_{[(B,\Psi,\eta),\bg]}$ over a point 
$[(B,\Psi,\eta),\bg] = (\bgamma')^{-1}([A',\Phi'])$ in the base of the
obstruction bundle, $\Upsilon_{\ft',\fs}$, is given by
\begin{equation}
\label{eq:ObstructionDirecSumDecomposition}
\Upsilon_{[(B,\Psi,\eta),\bg]}
\cong
\Upsilon^s_{[B,\Psi,\eta]}\oplus \Upsilon^i_{[\bA]},
\end{equation}
where $\Upsilon^s_{[B,\Psi,\eta]}$ is a fiber of the background obstruction
bundle (see \S \ref{subsubsec:BackgroundObstruction}), and
$\Upsilon^i_{[\bA]}$ is a fiber of the instanton obstruction bundle (see \S
\ref{subsubsec:InstantonObstruction}), and $[\bA]$ is a point in
$M_k^{s,\natural}(S^4,\delta)$ associated with $[\bg]\in\Gl_{\ft}(\delta)$.

\subsubsection{Background or Seiberg-Witten component of the obstruction
bundle}
\label{subsubsec:BackgroundObstruction}
In \cite{FL3}, the Seiberg-Witten component of the $\SO(3)$-monopole
obstruction bundle was identified with the subbundle constructed in
\cite[\S 3.5.2]{FL2a}, namely $\pi_N^*\Xi_{\ft,\fs}\cong
N_{\ft,\fs}(\eps)\times\CC^{r_\Xi}$ where $\pi_N:N_{\ft,\fs}\to M_{\fs}$ is
the projection, and there is an embedding of
$S^1$-equivariant vector bundles,
\begin{equation}
\label{eq:LowerLevelObstructionEmbedding}
\begin{CD}
\pi_N^*\Xi_{\ft,\fs} @>>>\fV_{\ft}
\\
@VVV @VVV
\\
N_{\ft,\fs}(\eps) @>>>\sC_{\ft}
\end{CD}
\end{equation}
The bundle map $\pi_N^*\Xi_{\ft,\fs}\to\fV_{\ft}$ in
\eqref{eq:LowerLevelObstructionEmbedding}
is $S^1$-equivariant if the circle acts diagonally on
$N_{\ft,\fs}(\eps)\times \CC^{r_{\Xi}}$
by scalar multiplication with weight two on both the fibers of
$N_{\ft,\fs}$ and on $\CC^{r_\Xi}$, while the circle acts
on $\fV_{\ft}$ by the action described following
\eqref{eq:InfiniteDimObstruction}.

The cut-and-paste construction
in \S \ref{subsec:SplicingMap} defines an embedding of vector bundles,
\begin{equation}
\label{eq:DefineBackgroundObstructionBundle}
\begin{CD}
\Upsilon_{\ft',\fs}^{s} @>{\bvarphi_s}>>\fV_{\ft'}
\\
@VVV @VVV
\\
\tilde N_{\ft,\fs}(\eps)\times_{\sG_{\fs}}\Gl_{\ft}(\delta)
@>{\bgamma'}>>\sC_{\ft'}
\end{CD}
\end{equation}
where $\Upsilon_{\ft',\fs}^{s}$ is the pullback
by the projection
$\tilde N_{\ft,\fs}(\eps)\times_{\sG_{\fs}}\Gl_{\ft}(\delta) 
\to N_{\ft,\fs}(\eps)$
of the obstruction bundle $\pi_N^*\Xi_{\ft,\fs}$
appearing in \eqref{eq:PulledBackObstructionBundle}, so that
$$
\Upsilon_{\ft',\fs}^{s} \cong \CC^{r_\Xi}
\times
\tilde N_{\ft,\fs}(\eps)\times_{\sG_{\fs}}\Gl_{\ft}(\delta).
$$
Because this cut-and-paste construction is circle equivariant,
the embedding \eqref{eq:DefineBackgroundObstructionBundle} has the same
circle-equivariance properties as the embedding
\eqref{eq:LowerLevelObstructionEmbedding}, where the circle acts
on $\fV_{\ft'}$ as previously described and diagonally on
$\CC^{r_\Xi}\times
\tilde N_{\ft,\fs}(\eps)\times_{\sG_{\fs}}\Gl_{\ft}(\delta)$,
with weight two on $\CC^{r_\Xi}$ and on $\tilde
N_{\ft,\fs}(\eps)\times_{\sG_{\fs}}\Gl_{\ft}(\delta)$ with the action given
in Lemma \ref{lem:S1EquivarianceSplicing}.

To describe the intersection of the image of the extended gluing map,
associated with the extended splicing map \eqref{eq:ExtendedSplicing} (see
\S \ref{subsec:GluingMap}) with $\bar\sM_{\ft'}$,
we need to extend the background obstruction bundle over
the space of extended gluing data \eqref{eq:CompactifiedGluingBundle}.
Thus, we define the extension by
\begin{equation}
\label{eq:ExtendedBackgroundObstruction}
\bar\Upsilon_{\ft',\fs}^{s} \cong \CC^{r_\Xi}
\times
\tilde N_{\ft,\fs}(\eps)\times_{\sG_{\fs}}\bar{\Gl}_{\ft}(\delta).
\end{equation}
The restriction of $\bar\Upsilon_{\ft',\fs}^{s}$ to the stratum
$N_{\ft,\fs}(\eps)\times X$ of the space of extended gluing data is given
by the pullback of the bundle in the left-hand column of
\eqref{eq:LowerLevelObstructionEmbedding} to $N_{\ft,\fs}(\eps)\times X$
and is embedded into $\fV_{\ft}$ via the diagram
\eqref{eq:LowerLevelObstructionEmbedding}.

\subsubsection{Instanton component of the obstruction bundle}
\label{subsubsec:InstantonObstruction}
In \cite{FL3}, the instanton components of the $\SO(3)$-monopole gluing
obstruction space were identified with cokernels of the Dirac operators,
$$
D_{\bA}:C^\8(S^4,\bV^+)\to C^\8(S^4,\bV^-),
$$
where the induced connections $\hat\bA$ on the $\SO(3)$ bundle $\fg_{\bV}$
over $S^4$ are anti-self-dual. Because the standard round metric on $S^4$
has positive scalar curvature, one has $\Ker D_{\bA} = \{0\}$ and as
$\Ind_\CC D_{\bA} = -1$ (see, for example, \cite[p. 314]{FL1}), the
cokernel bundles,
\begin{equation}
\label{eq:InstantonCokernelBundle}
\Coker\tilde\bD_{\bV} \to \tilde M^{s,\natural}_1(S^4,\delta)
\quad\text{and}\quad
\Coker\bD_{\bV} \to M^{s,\natural}_1(S^4,\delta),
\end{equation}
are complex line bundles with fibers $\Coker\bD_{\bV}|_{\bA} =
\Coker D_{\bA}$. The total space of the bundle $\Coker\bD_{\bV}$ is
identified by the diffeomorphism
\begin{equation}
\label{eq:ExplicitS4CokernelBundle}
\Coker\bD_{\bV}
\cong
\left(\U(2)\times (0,\delta]\right)\times_{S^1}\CC \cong
\SU(2)\times_{\{\pm 1\}}\CC \times (0,\delta],
\end{equation}
where $S^1$ acts {\em diagonally\/} on $\U(2)\times\CC$. The bundle
$\Coker\bD_{\bV}\to M^{s,\natural}_1(S^4,\delta)$ is non-trivial, but
torsion.
There is an action of $\Spinu(4)$ on $\Coker\bD_{\bV}$, covering the
action of $\SO(4)\times\SO(3)$ on $M_1^{s,\natural}(S^4,\delta)$:
this action is most easily understood using the trivializations 
\eqref{eq:TrivSpinRank8Sphere} of
$\bV|_{S^4\less\{n,s\}}$.
The section $\Psi$, in this trivialization, is a map
$$
\RR^4-\{0\} \to \Delta\otimes_\CC\CC^2.
$$
Then $\tilde R\otimes M\in \Spinu(4)$ acts by
$$
\Psi \mapsto (\tilde R\otimes M)\circ\Psi\circ \Ad^c(\tilde R)^{-1}.
$$
The procedure for splicing spinors over $S^4$ onto $X$ defines a map (see
\S \ref{subsubsec:SpliceSpinorsFromSphereToX})
\begin{equation}
\label{eq:UnquotientedDomainOfInstantonObstruction}
\begin{gathered}
\tilde\bvarphi_i:
\tilde N_{\ft,\fs}\times\Fr_{\CCl(T^*X)}(V)
\times\Coker\tilde\bD_{\bV}
\to
\tilde\sC_{\ft'}\times \Gamma(\Lambda^+\otimes\fg_{V'}\oplus {V'}^-).
\end{gathered}
\end{equation}
The map \eqref{eq:UnquotientedDomainOfInstantonObstruction} is
invariant under the diagonal action of $\Spinu(4)$ on
$\Fr_{\CCl(T^*X)}(V)\times\Coker\tilde\bD_{\bV}$.
Note that this $\Spinu(4)$ action covers
the diagonal $\SO(4)$ action on the component $\Fr(T^*X)\times\tilde
M_1^\natural(S^4,\delta)$ of the base and the diagonal $\SO(3)$ action on
$\Fr(\fg_{V})\times\Fr(\fg_{\bV}|_s)$.
In addition,
the map \eqref{eq:UnquotientedDomainOfInstantonObstruction} is equivariant
with respect to the action of $\sG_{\ft'}$ on the image and
\begin{itemize}
\item
The diagonal action of $\sG_{\fs}$ on
$\tilde N_{\ft,\fs}\times \Fr_{\CCl(T^*X)}(V)$,
where $\sG_{\fs}$ acts on $\Fr_{\CCl(T^*X)}(V)$ by the homomorphism
$\sG_{\fs}\embed\sG_{\ft}$ in \eqref{eq:GaugeGroupInclusion},
covering the diagonal action of $\sG_{\fs}$ on
$\tilde N_{\ft,\fs}\times \Fr(\fg_V)$,
\item
The action of $\sG_1$ on $\Coker\tilde\bD_{\bV}$,
covering its diagonal action on the component
$\tilde M_1^{s,\natural}(S^4,\delta)$ of the base.
\end{itemize}
Finally, the map \eqref{eq:UnquotientedDomainOfInstantonObstruction}
is $S^1$-equivariant with respect to the circle action 
$\fV_{\ft'}$ described after
\eqref{eq:InfiniteDimObstruction} and the diagonal circle action
on $\tilde N_{\ft,\fs}\times\Coker\bD_A$, where the circle acts on
$\Coker\bD_A$ by scalar multiplication and on $\tilde N_{\ft,\fs}$ by the
action \eqref{eq:S1ZActionOnN}.  We define a vector bundle over $\tilde
N_{\ft,\fs}\times_{\sG_{\fs}}\Gl_{\ft}(\delta)$ by setting
\begin{equation}
\label{eq:DefnInstantonObstruction}
\Upsilon_{\ft',\fs}^{i}
=
\tilde N_{\ft,\fs}\times_{\sG_{\fs}}
\Fr_{\CCl(T^*X)}(V)\times_{\Spinu(4)}\Coker\bD_{\bV},
\end{equation}
where $\Spinu(4)$ acts diagonally on
$\Fr_{\CCl(T^*X)}(V)\times\Coker\bD_{\bV}$ and trivially on $\tilde
N_{\ft,\fs}$, while $\sG_{\fs}$ acts diagonally on $\tilde
N_{\ft,\fs}\times \Fr_{\CCl(T^*X)}(V)$ as described above and trivially on
$\Coker\bD_{\bV}$.  We then have the

\begin{lem}
\label{lem:InstantonObstructionSplicingDomain}
The map $\tilde\bvarphi_i$ descends to an embedding of $S^1$-equivariant
vector bundles,
\begin{equation}
\label{eq:InstantonObstructionSplicingDomain}
\begin{CD}
\Upsilon_{\ft',\fs}^{i}
@> \bvarphi_i >> \fV_{\ft'}
\\
@VVV @VVV
\\
\tilde N_{\ft,\fs}\times_{\sG_{\fs}}\Gl_{\ft}(\delta)
@> \bgamma' >> \sC_{\ft'}
\end{CD}
\end{equation}
where the group actions are as described in the preceding paragraph.
\end{lem}

\subsection{Construction of the gluing map}
\label{subsec:GluingMap}
The Uhlenbeck compactification 
\eqref{eq:UhlenbeckCompactifiedSplicingDomain} of the virtual moduli space
model $\sM_{\ft',\fs}^{\stab}$ has a smooth, circle invariant stratification
(where the strata are manifolds with boundary),
\begin{equation}
\label{eq:ModelStratification}
\begin{aligned}
\bar\sM_{\ft',\fs}^{\stab}
&=
\tN_{\ft,\fs}(\eps)\times_{\sG_{\fs}}\Gl_{\ft}(\delta)
\\
&\quad\sqcup 
(N_{\ft,\fs}(\eps)-M_{\fs})\times X
\\
&\quad\sqcup 
M_{\fs}\times X.
\end{aligned}
\end{equation}
If $r_N$ is the complex rank of $N_{\ft,\fs}\to M_{\fs}$
and $d_{s}(\fs)=\dim M_{\fs}$, then
because the fiber of $\Gl_{\ft}(\delta)\to X$ is four-dimensional,
the dimensions of the strata in \eqref{eq:ModelStratification} are
given by (see the rank and dimension relations \eqref{eq:SWDimRelations}
and \eqref{eq:SWBundleRankRelations}) 
\begin{align*}
\dim\left(\tN_{\ft,\fs}(\eps)\times_{\sG_{\fs}}\Gl_{\ft}(\delta)\right)
&= 2r_N + d_{s}(\fs) +8,
\\
\dim\left((N_{\ft,\fs}(\eps)-M_{\fs})\times X\right)
&= 2 r_N+ d_{s}(\fs) + 4,
\\
\dim\left(M_{\fs}\times X\right) &= d_{s}(\fs)+4.
\end{align*}
The restriction of the splicing map $\bga'$ to the top stratum
$\sM_{\ft',\fs}^{\stab}$ takes values in $\sC^{*,0}_{\ft'}$; on the middle
stratum, $(N_{\ft,\fs}(\eps)-M_{\fs})\times X$, the splicing map $\bga'$
restricts to the map $\bga_{\ft,\fs}\times\id_X$, which takes values in
$\sC^{*,0}_{\ft}\times X$; on the lowest stratum, $M_{\fs}\times X$, the
splicing map $\bga'$ restricts to the identity map on $M_{\fs}\times X$. The
stratum $M_{\fs}\times X$ is the fixed point set of the circle action
described in Lemma \ref{lem:S1EquivarianceSplicing}.

The splicing map $\bga'$ may be deformed
$S^1$-equivariantly to a `gluing map' $\bga$ with the
properties described in the following 

\begin{thm}
\label{thm:GluingThm}
\cite{FL3, FL4}
For small enough positive $\eps$ and $\delta$, there is a
topological embedding,
\begin{equation}
\label{eq:DomainOfGluingMap}
\bga:
\bar\sM_{\ft',\fs}^{\stab}
\to
\bar{\sC}_{\ft'}
=
\sC_{\ft'} \sqcup (\sC_{\ft}-M_{\fs})\times X \sqcup M_{\fs}\times X,
\end{equation}
restricting to a smooth embedding of the top stratum of
\eqref{eq:ModelStratification} into $\sC_{\ft'}^{*,0}$, the smooth embedding
$\bga_{\ft,\fs}\times\id_X$ of the middle stratum into
$\sC_{\ft}^{*,0}\times X$, and the identity map on the lowest stratum,
$M_{\fs}\times X\subset\sC_{\ft}^0\times X$. There is a smooth,
circle-equivariant section $\bchi_i$ of the instanton obstruction 
bundle \eqref{eq:DefnInstantonObstruction},
$$
\Upsilon_{\ft',\fs}^{i} \to \sM_{\ft',\fs}^{\stab},
$$
and a continuous, circle-equivariant section
$\bchi_s$ of the background obstruction 
bundle \eqref{eq:ExtendedBackgroundObstruction},
$$
\bar\Upsilon_{\ft',\fs}^{s} \to \bar\sM_{\ft',\fs}^{\stab},
$$
which is smooth when restricted to each stratum of $\bar\sM_{\ft',\fs}^{\stab}$
in \eqref{eq:ModelStratification} such that, if
$\bchi=\bchi_s\oplus\bchi_{i}$, then
$$
\bga(\sM_{\ft',\fs}^{\stab}\cap\bchi^{-1}(0))
=
\bga(\sM_{\ft',\fs}^{\stab})\cap\sM_{\ft'},
$$
and
$$
\bga\left(((N_{\ft,\fs}(\eps)-M_{\fs})\times X)\cap \bchi_s^{-1}(0)\right)
=
\bga((N_{\ft,\fs}(\eps)-M_{\fs})\times X)\cap\bar\sM_{\ft'}.
$$
The sections $\bchi_s\oplus \bchi_i$ and $\bchi_s$ of the vector
bundles,
$$
\bar\Upsilon_{\ft',\fs}^{s}\oplus \Upsilon_{\ft',\fs}^{i}
\to \sM_{\ft',\fs}^{\stab}
\quad\text{and}\quad
\bar\Upsilon_{\ft',\fs}^{s}\to (N_{\ft,\fs}-M_{\fs})\times X,
$$
respectively, vanish transversely.
\end{thm}

We shall formally extend the section $\bchi_i$ over the lower strata
of $\bar\sM^{\stab}_{\ft',\fs}$ by setting it equal to zero on the
lower strata.  This is done so that we may discuss $\bchi^{-1}(0)$
as a subspace of $\bar\sM^{\stab}_{\ft',\fs}$; we make no assumptions
about the continuity or transversality of this formal
extension of $\bchi_i$ over the lower strata.

\subsection{Link of a level-one Seiberg-Witten stratum}
\label{subsec:DefnOfLink}
We define the {\em Seiberg-Witten\/} and {\em instanton components\/} of
the {\em virtual link\/} of the stratum $M_{\fs}\times X$ in
$\bar\sM_{\ft',\fs}^{\stab}/S^1$ by
\begin{equation}
\label{eq:DefineSWandInstantonModelThickLink}
\begin{aligned}
\bar\bL_{\ft',\fs}^{\stab,s}
&=
\left(
\rd\tilde N_{\ft,\fs}(\eps)\times_{\sG_{\fs}}\bar{\Gl}_{\ft}(\delta)
\right)/S^1,
\\
\bL_{\ft',\fs}^{\stab,i}
&=
\left(
\tilde N_{\ft,\fs}(\eps)\times_{\sG_{\fs}}\partial\bar{\Gl}_{\ft}(\delta)
\right)/S^1,
\end{aligned}
\end{equation}
with (see \eqref{eq:DefnXGluingParameters} and
\eqref{eq:CompactifiedGluingBundle} for the definitions
of $\Gl_{\ft}(\delta)$ and $\bar{\Gl}_{\ft}(\delta)$, respectively),
\begin{equation}
\label{eq:GlueBoundary}
\partial\bar{\Gl}_{\ft}(\delta)
=
\left(\Fr(\fg_{V})\times_X\Fr(T^*X)\times 
\lambda^{-1}(\delta)\cap M_1^{s,\natural}(S^4)\right)
/(\SO(3)\times\SO(4)),
\end{equation}
and where the circle action on $\bar\sM_{\ft',\fs}^{\stab}$ is induced by the
circle action \eqref{eq:S1ZActionOnN} on $\tilde N_{\ft,\fs}$ and the trivial
action on $\bar{\Gl}_{\ft}(\delta)$ (as described in the first
action in Lemma \ref{lem:S1EquivarianceSplicing}). The map 
$$
\lambda:\sB_1(S^4)\to (0,\8)
$$
appearing in the gluing-data boundary \eqref{eq:GlueBoundary} is defined by 
the scale definition \eqref{eq:DefnCenterScale}. 

We define the {\em virtual link\/} of the stratum $M_{\fs}\times X$ in
$\bar\sM_{\ft',\fs}^{\stab}/S^1$ by setting
\begin{equation}
\label{eq:DefineModelThickLink}
\bar\bL^{\stab}_{\ft',\fs}
=
\bar\bL_{\ft',\fs}^{\stab,s}
\cup\bL_{\ft',\fs}^{\stab,i},
\end{equation}
We let
\begin{equation}
\label{eq:TopLevels}
\bL^{\stab}_{\ft',\fs}
=
\bar\bL^{\stab}_{\ft',\fs}\cap\sM^{\stab}_{\ft',\fs}/S^1
\quad\text{and}\quad
\bL^{\stab,s}_{\ft',\fs}
=
\bar\bL^{\stab,s}_{\ft',\fs}\cap\sM^{\stab}_{\ft',\fs}/S^1
\end{equation}
denote the intersection of these subspaces with the top stratum
$\sM^{\stab}_{\ft',\fs}/S^1$ of $\bar\sM^{\stab}_{\ft',\fs}/S^1$.
Note that the top level, $\bL^{\stab}_{\ft',\fs}$, of the virtual
link is only a topological and not a smooth manifold because of
the `edge':
\begin{equation}
\label{eq:Edge}
\bL^{\stab,i}_{\ft',\fs}\cap\bL^{\stab,s}_{\ft',\fs}.
\end{equation}
Let $\bL_{\ft',\fs}^{\sing}$
be the intersection of $\bar\bL^{\stab}_{\ft',\fs}$ with
the union of the lower levels of $\bar\sM^{\stab}_{\ft',\fs}/S^1$:
\begin{equation}
\label{eq:DefineLSing}
\begin{aligned}
\bL_{\ft',\fs}^{\sing}
&=
\bar\bL^{\stab}_{\ft',\fs} - \bL^{\stab}_{\ft',\fs}
\\
&=
\bar\bL^{\stab}_{\ft',\fs}
\cap\left(\bar\sM^{\stab}_{\ft',\fs}-\sM^{\stab}_{\ft',\fs}\right)/S^1
\\
&=
\rd N_{\ft,\fs}(\eps)/S^1 \times X,
\end{aligned}
\end{equation}
We then have a stratification: 
\begin{equation}
\label{eq:StratificationOfVirtLink}
\bar\bL^{\stab}_{\ft',\fs} 
=
\bL^{\stab}_{\ft',\fs}
\sqcup
\bL_{\ft',\fs}^{\sing}.
\end{equation}
In this article, where $\ell(\ft',\fs)=1$, we see from 
\eqref{eq:DefineLSing} that
$\bL_{\ft',\fs}^{\sing}$ is a closed, smooth manifold.

We define the {\em
link\/} of $(M_{\fs}\times X)\cap\bar\sM_{\ft'}/S^1$ in $\bar\sM_{\ft'}/S^1$
and its top stratum by
\begin{equation}
\label{eq:DefineLink}
\bar\bL_{\ft',\fs}
=
\bga\left(
\bchi^{-1}(0)
\cap
\bar\bL^{\stab}_{\ft',\fs}
\right)
\quad\text{and}\quad
\bL_{\ft',\fs}=\bar\bL_{\ft',\fs}\cap\sM^{*,0}_{\ft'}/S^1.
\end{equation}
For generic choices of the parameters $\eps$ and $\delta$
(which we shall henceforth assume) defining $\bL^{\stab}_{\ft',\fs}$, the
subspaces 
$$
\bL_{\ft',\fs}\cap\bL^{\stab,i}_{\ft',\fs}
\quad\text{and}\quad
\bL_{\ft',\fs}\cap\bL^{\stab,s}_{\ft',\fs}
$$
will be smooth submanifolds of $\bL^{\stab,i}_{\ft',\fs}$
and $\bL^{\stab,s}_{\ft',\fs}$ respectively, transverse to the
common boundary \eqref{eq:Edge}.
The following lemma is the first step in showing that the
intersection numbers on the right-hand-side of
\eqref{eq:RawCobordismSum} are well-defined.

\begin{lem}
\label{lem:GRIntersect0}
Assume $w\in H^2(X;\ZZ)$ is such that $w\pmod{2}$ is good
in the sense of Definition \ref{defn:Good}. Given
a Riemannian metric on $X$ and a pair $(\ft',\fs)$ with $\ell(\ft',\fs)=1$
and $w_2(\ft')\equiv w\pmod{2}$, there are positive constants $\eps_0$ and
$\delta_0$ such that the following hold for all generic choices of 
$\eps\leq \eps_0$ and $\delta\leq\delta_0$ defining
$\bar\bL^{\stab}_{\ft',\fs}$. 
\begin{itemize}
\item
The link $\bar\bL_{\ft',\fs}$ is disjoint from $\bar M^w_\ka$ and
$\bar\sM_{\ft'}^{\red}$ in the stratification
\eqref{eq:StratificationCptPU(2)Space} of $\bar\sM_{\ft'}/S^1$. 
\item 
For all $z\in\AAA(X)$, the intersection,
\begin{equation}
\label{eq:Int}
\bar\sV(z)\cap\bar\sW^{\eta}\cap\bar\bL_{\ft',\fs}
\end{equation}
is contained in the top stratum $\bL_{\ft',\fs}$ of
$\bar\bL_{\ft',\fs}\subset\bar\sM_{\ft'}/S^1$,
and is disjoint from the image of the edge \eqref{eq:Edge} under 
the gluing map $\bga$.
\end{itemize}
\end{lem}

\begin{proof}
Because $w\pmod 2$ is good, the union of strata $\bar M^w_\ka$ is disjoint
from the union of strata $\bar\sM_{\ft'}^{\red}$ (see remarks in \S
\ref{subsec:ASDsingularities}). Hence, for sufficiently
small parameters $\delta$ and $\eps$ in the definition of
$\bar\bL^{\stab}_{\ft',\fs}$, the link $\bar\bL_{\ft',\fs}$ is disjoint
from these strata.

The geometric representatives $\bar\sV(z)\cap\bar\sW^{\eta}$
do not intersect the lower levels $(\bar\sM_{\ft'}-\sM_{\ft'})/S^1$ of
$\bar\sM_{\ft'}/S^1$ except at points in $\bar\sM_{\ft'}^{\red}$
or at points in $\bar M^w_\ka$ by \cite[Corollary 3.18]{FL2b}.
Therefore, the intersection \eqref{eq:Int} is contained in the top
stratum, $\bL_{\ft',\fs}$.

The geometric representatives are transverse to $\sM^{*,0}_{\ft'}/S^1$.
Hence, the intersection
$$
\bga^{-1}(\bar\sV(z)\cap\bar\sW^\eta)\cap \bchi^{-1}(0),
$$
is a smooth submanifold of $\sM^{\stab,*}_{\ft',\fs}/S^1$.
For generic values of $\eps$ and $\delta$, the preceding intersection
will be transverse to the edge and thus, by dimension-counting, disjoint
from the edge.
\end{proof}

Because the intersection \eqref{eq:Int} is contained in the locus of
smooth points of $\bar\bL_{\ft',\fs}$, it will be possible to define
an intersection number once we have discussed the orientation of the link.

We note that the construction of the link in this section applies to links
$\bar\bL_{\ft',\fs}$ of $M_{\fs}\times\Sym^\ell(X)$ when
$\ell(\ft',\fs)>1$.  The main difference is that the additional dilation
parameters needed to describe a neighborhood of the strata
$M_{\fs}\times\Sym^\ell(X)$ result in more complicated `edges' in the
boundaries of this neighborhood.  However, for generic choices of dilation
parameters, these edges will be smooth submanifolds of codimension-one or
greater and the argument in Lemma \ref{lem:GRIntersect0} will show that the
geometric representatives are disjoint from these edges.

\subsection{Orientations}
\label{subsec:Orient}
We now discuss the orientations of the moduli spaces and links appearing in
the $\SO(3)$-monopole cobordism --- in particular links of level-one
Seiberg-Witten moduli spaces --- with respect to which the pairings are
defined. Orientations for moduli spaces and links of top-level
Seiberg-Witten moduli spaces were discussed in \cite[\S 2]{FL2b}.

An orientation for $\sM_{\ft'}$ determines one for $\bL_{\ft',\fs}$ through
the convention introduced in \cite[Equations (2.16), (2.16) \& (2.25)]{FL2b} by
considering $\bL_{\ft',\fs}$ as a boundary of
$(\sM_{\ft'}-\bga(\sM_{\ft',\fs}^{\stab}))/S^1$. Specifically, at a point
$[A,\Phi]\in \bL_{\ft',\fs}$, if 
\begin{itemize}
\item
$\vec r\in
T\sM_{\ft'}^{*,0}$ is an outward-pointing radial vector with respect to the
open neighborhood $\sM_{\ft'}\cap\bga(\sM_{\ft',\fs}^{\stab})$
 and complementary to the tangent space of
$\bL_{\ft',\fs}$, 
\item
$v_{S^1}\in T\sM_{\ft'}^{*,0}$ is
tangent to the orbit of $[A,\Phi]$ under the (free) circle action (where
$S^1\subset\CC$ has its usual orientation), and 
\item
$\la_{\sM}\in\det(T\sM_{\ft'}^{*,0})$ is an orientation for $T\sM_{\ft'}$
at $[A,\Phi]$, 
\end{itemize}
then we define an orientation $\la_{L}$ for
$T\bL_{\ft',\fs}$ at $[A,\Phi]$ by
\begin{equation}
\label{eq:QuotientBoundaryOrientation}
\la_{\sM}
=
-v_{S^1}\wedge\vec r \wedge \tilde{\la}_{L},
\end{equation}
where the lift $\tilde\la_{L}\in \Lambda^{\max-2}(T\sM_{\ft'}^{*,0})$ at
$[A,\Phi]$ of $\la_{L}\in\det(T\bL_{\ft',\fs})\subset
\Lambda^{\max-1}(T(\sM_{\ft'}^{*,0}/S^1))$, 
obeys $\pi_*\tilde\la_{L}=\la_{L}$, if $\pi:\sM_{\ft'}\to\sM_{\ft'}/S^1$ is
the quotient map.

\begin{defn}
\label{defn:BdryOrnLink}
If $O$ is an orientation for
$\sM_{\ft'}$, we call the orientation for $\bL_{\ft',\fs}$ related to
$O$ by equation \eqref{eq:QuotientBoundaryOrientation} the {\em boundary
orientation defined by $O$\/}.  
\end{defn}

We note that the orientation convention
\eqref{eq:QuotientBoundaryOrientation} can be applied more generally to
define an orientation for a quotient $\rd M/S^1$, given an orientation for
a manifold with boundary $M$ with a free circle action.

We now introduce an orientation for $\bL_{\ft',\fs}$, based on 
a choice of orientation for the component
$\bL_{\ft',\fs}^{\stab,i}$ 
(see \eqref{eq:DefineSWandInstantonModelThickLink}), which is useful for
cohomological computations.  First,
observe that the gluing map $\bga$ identifies a relatively open
subspace of $\bL_{\ft',\fs}$ with the zero-locus of the obstruction section
$\bchi$ in $\bL_{\ft',\fs}^{\stab,i}$.
Because this section vanishes
transversely, the normal bundle of the zero locus in
$\bL_{\ft',\fs}^{\stab,i}$ is identified with the obstruction bundle,
$\Upsilon_{\ft',\fs}/S^1$,
which has a complex orientation.  If this normal bundle is given
the complex orientation, orientations of $\bL_{\ft',\fs}$ are thus
determined by orientations for $\bL_{\ft',\fs}^{\stab,i}$. From its definition
\eqref{eq:DefineSWandInstantonModelThickLink}, we see that
$\bL_{\ft',\fs}^{\stab,i}$ can be viewed as a complex disk bundle:
\begin{equation}
\label{eq:IBoundaryDiskBundle}
\bL_{\ft',\fs}^{\stab,i} 
=
\tilde N_{\ft,\fs}(\eps)\times_{\sG_{\fs}\times S^1}
\rd\bar{\Gl}_{\ft}(\delta) 
\to 
\tM_{\fs}\times_{\sG_{\fs}\times S^1}\rd\bar{\Gl}_{\ft}(\delta).
\end{equation}
The circle action in \eqref{eq:IBoundaryDiskBundle},
described as the second circle action in Lemma
\ref{lem:S1EquivarianceSplicing}, is trivial on $\tM_{\fs}$.
Because the circle acts trivially on $\tM_{\fs}$
and $\sG_{\fs}$ acts trivially on the quotient
$\rd\bar{\Gl}_{\ft}(\delta)/S^1$, we can identify the base of
the bundle \eqref{eq:IBoundaryDiskBundle} with the following product,
a compact manifold:
\begin{equation}
\label{eq:ProductDecompositionOfBase}
\begin{aligned}
\bB\bL_{\ft',\fs}^{\stab,i}
:&=
\tM_{\fs}\times_{\sG_{\fs}} \rd\bar{\Gl}_{\ft}(\delta)/S^1
\\
&\cong
M_{\fs}\times \rd\bar{\Gl}_{\ft}(\delta)/S^1.
\end{aligned}
\end{equation}
We may further assume without loss of generality that the point
$[A,\Phi]\in \bL_{\ft',\fs}\subset\sM_{\ft'}/S^1$ --- at which we choose to
compare orientations --- is identified, via the gluing
map $\bgamma:\sM_{\ft',\fs}^{\stab}\to\sC_{\ft'}^{*,0}$, with a point
$$
\bgamma^{-1}([A,\Phi])
=
[[B,\Psi,0],[\bg]],
$$ 
in the base $\bB\bL_{\ft',\fs}^{\stab,i}$ of the bundle
\eqref{eq:IBoundaryDiskBundle}, where $[B,\Psi,0] \in M_{\fs}\subset
N_{\ft,\fs}(\eps)$ and $\bg\in\rd\bar{\Gl}_{\ft}(\delta)$.  The commutative
diagram \eqref{eq:UBundle} shows that the fiber of the disk bundle
\eqref{eq:IBoundaryDiskBundle} over the point $[[B,\Psi,0],[\bg]]$ is
identified with the fiber of $N_{\ft,\fs}(\eps)$ over $[B,\Psi]$.  We
therefore have an isomorphism of tangent spaces at this point:
\begin{equation}
\label{eq:TangentOfIBoundary}
T_{\bgamma^{-1}([A,\Phi])}\bL_{\ft',\fs}^{\stab,i}
\cong
T_0(N_{\ft,\fs}|_{[B,\Psi]})
\oplus 
T_{[B,\Psi]}M_{\fs}\oplus T_{[\bg]}(\rd\bar{\Gl}_{\ft}(\delta)/S^1).
\end{equation}
We define the {\em standard orientation\/} for
$\rd\bar{\Gl}_{\ft}(\delta)/S^1$ by applying the convention
\eqref{eq:QuotientBoundaryOrientation} to an orientation for
$\Gl_{\ft}(2\delta)$, considering $\rd\bar{\Gl}_{\ft}(\delta)/S^1$ as a
boundary of $\Gl_{\ft}(2\delta)/S^1-\Gl_{\ft}(\delta)/S^1$.  
{}From its definition, we see that the space
$\Gl_{\ft}(2\delta)$ is a locally-trivial fiber-bundle over $X$ with fiber
$(0,2\delta)\times\SO(3)$, so the orientation for $X$ and an orientation
for $(0,2\delta)\times\SO(3)$ determine an orientation for
$\Gl_{\ft}(2\delta)$.  The standard orientation for
$\rd\bar{\Gl}_{\ft}(\delta)/S^1$ will then be given through the convention
\eqref{eq:QuotientBoundaryOrientation} by taking the orientation of
$\Gl_{\ft}(2\delta)$ induced by the orientation on
$(0,2\delta)\times\SO(3)$ given by identifying its tangent spaces
with the space of quaternions, $\HH$, with the tangent spaces to $\SO(3)$
being identified with $\Imag(\HH)$.
(This agrees with the convention
given in \cite[\S 3(c) \& \S 3(d)]{DonOrient}.)

\begin{defn}
\label{defn:StdOrnVirtLink}
The {\em standard orientation\/} for $\bL^{\stab,i}_{\ft',\fs}$ is
determined, through the isomorphism \eqref{eq:TangentOfIBoundary}, by the
\begin{itemize}
\item
Standard orientation for $\rd\bar{\Gl}_{\ft}(\delta)/S^1$, 
\item
Complex
orientation of the fibers of the bundle \eqref{eq:IBoundaryDiskBundle}, and the
\item
Orientation of $M_{\fs}$ induced by a homology orientation $\Om$
\cite[\S 6.6]{MorganSWNotes}.
\end{itemize}
\end{defn}

\begin{defn}
\label{defn:StdOrnLink}
The {\em standard orientation\/} for $\bL_{\ft',\fs}$ is determined by the
\begin{itemize}
\item
Standard orientation for $\bL_{\ft',\fs}^{\stab,i}$, 
\item
Complex orientation of the normal bundle of
$\bchi^{-1}(0)\cap\bL^{\stab,i}_{\ft',\fs}$ in 
$\bL^{\stab,i}_{\ft',\fs}$, and the 
\item
Diffeomorphism $\bgamma:\bchi^{-1}(0)\cap\sM^{\stab}_{\ft',\fs}\to 
\sM_{\ft'}\cap\bgamma(\sM^{\stab}_{\ft',\fs})$.
\end{itemize}
\end{defn}

We now relate the two orientations of $\bL_{\ft',\fs}$ which we defined
above: 

\begin{lem}
\label{lem:LinkOrientation}
Suppose that the \spinu structure $\ft$ admits a splitting
$\ft=\fs\oplus \fs\otimes L$. Then the
boundary orientation for $\bL_{\ft',\fs}$ defined by the
orientation $O^{\asd}(\Om,c_1(L))$ for $\sM_{\ft'}$ (see
\cite[Definition 2.3]{FL2b}) agrees with the standard orientation for
$\bL_{\ft',\fs}$.
\end{lem}

\begin{proof}
The standard orientation for $\bL_{\ft',\fs}\subset\sM_{\ft'}/S^1$ is defined,
through the diffeomorphism  
$$
\bgamma:\bchi^{-1}(0)\cap\sM^{\stab}_{\ft',\fs}
\to 
\sM_{\ft'}\cap\bgamma(\sM^{\stab}_{\ft',\fs}),
$$
by applying the convention \eqref{eq:QuotientBoundaryOrientation} to the
submanifold $\bchi^{-1}(0)\cap\sM^{\stab}_{\ft',\fs}$ of
$\sM^{\stab}_{\ft',\fs}$, while the orientation $O^{\asd}(\Om,c_1(L))$ is
defined by applying it to $\sM_{\ft'}$.  Thus, to compare the two
orientations of $\bL_{\ft',\fs}$, it suffices to compare the orientations
of $\bchi^{-1}(0)\cap\sM^{\stab}_{\ft',\fs}$ and $\sM_{\ft'}$ which induce
these orientations of $\bL_{\ft',\fs}$.

An orientation for $\sM_{\ft'}$ is given by an orientation for the
index bundle of the deformation operator $\sD$ defined in
\cite[Equation (2.62)]{FL2a}, whose kernel gives the
tangent spaces of $\sM^{*,0}_{\ft'}$.  Suppose that the point
$(A,\Phi)\in\tilde\sM_{\ft'}$ at which we do the orientation comparisons
is obtained by gluing a {\em framed\/} \spinu connection $\bA\in
\tM^{s,\natural}_1(S^4,\delta)$ onto the background pair $(B,\Psi)\in
\tM_{\fs}$ at a point $x\in X$:
\begin{equation}
\label{eq:SplicedPairForOrns}
(A,\Phi) = (B,\Psi)\#(\bA,0).
\end{equation}
Note that we abuse notation here and omit explicit mention of the frame for
$\fg_{\bV}|_s$.  Then the excision argument \cite[\S 7.1]{DK} yields the
isomorphism
\begin{equation}
\label{eq:DeformExcision}
\det\sD_{[A,\Phi]}\cong\det\sD_{[B,\Psi]}\otimes\det\sD_{[\bA,0]}.
\end{equation}
There are isomorphisms (see \cite[Equations (2.5), (2.11), and
(2.26)]{FL2b}), 
\begin{equation}
\label{eq:DeformComponentIdentity}
\begin{aligned}
\det\sD_{[B,\Psi]} &\cong \det(N_{\ft,\fs}|_{[B,\Psi]})
\otimes \det(T_{[B,\Psi]}M_{\fs}) \otimes
\det(\Xi_{\ft,\fs}|_{[B,\Psi]})^*,
\\
\det\sD_{[\bA,0]} &\cong \det(T_{[\bA]}M^s_1(S^4,\delta)) \otimes
\det(\Coker D_\bA)^*
\\
&\cong \det(T_{[\bA]}M^{s,\natural}_1(S^4,\delta)) \otimes
\det(T_xX) \otimes \det(\Coker D_\bA)^*,
\end{aligned}
\end{equation}
where in the last isomorphism we have used the identification
$M^{s}_1(S^4,\delta)\cong M^{s,\natural}_1(S^4,\delta)\times\RR^4$ and the
identification of a ball in $\RR^4$ with a ball around $x$ in $X$ via the
gluing map. Lemma 2.6 in \cite{FL2b} implies that the orientation
$O^{\asd}(\Om,c_1(L))$ of $\det\sD_{[A,\Phi]}$ is given through the
isomorphisms \eqref{eq:DeformExcision} and
\eqref{eq:DeformComponentIdentity} by the complex orientations of
$N_{\ft,\fs}|_{[B,\Psi]}$, $(\Xi_{\ft,\fs}|_{[B,\Psi]})^*$, and $(\Coker
D_\bA)^*$, the orientation of $T_{[B,\Psi]}M_{\fs}$ determined by the
homology orientation $\Om$ \cite[\S 6.6]{MorganSWNotes}, the standard
orientation of $T_{[\bA]}M^{s,\natural}_1(S^4,\delta)$ defined in
\cite[p. 413]{DonOrient}, and 
the orientation for $T_xX$ given by that of $X$.

We now describe the orientation of $\bchi^{-1}(0)\cap\sM^{\stab}_{\ft',\fs}$
inducing the standard orientation of $\bL_{\ft',\fs}$.  First, the
orientation of $\bchi^{-1}(0)\cap\sM^{\stab}_{\ft',\fs}$ is induced by one
for $\sM^{\stab}_{\ft',\fs}$ by identifying the normal bundle of
$\bchi^{-1}(0)\cap\sM^{\stab}_{\ft',\fs}$ in $\sM^{\stab}_{\ft',\fs}$ with the
obstruction bundle $\Upsilon_{\ft',\fs}$ and using the complex orientation
of this obstruction bundle.
With $[A,\Phi]$ as in \eqref{eq:SplicedPairForOrns}, then by
\eqref{eq:ObstructionDirecSumDecomposition},
\eqref{eq:DefineBackgroundObstructionBundle},
and \eqref{eq:DefnInstantonObstruction}
there is an isomorphism of obstruction
bundle fibers,
\begin{equation}
\label{eq:ObstructionSplitting}
\Upsilon_{\ft',\fs}|_{\bgamma^{-1}([A,\Phi])}
\cong
\Coker D_\bA \oplus \Xi_{\ft,\fs}|_{[B,\Psi]}.
\end{equation}
Hence, the complex orientation of the obstruction bundle matches the
complex orientation of the factors $\Coker D_\bA$ and
$\Xi_{\ft,\fs}|_{[B,\Psi]}$ in \eqref{eq:DeformComponentIdentity}.

We now compare the tangent space 
$T_{\bgamma^{-1}([A,\Phi])}\sM^{\stab}_{\ft',\fs}$ with
the remaining factors on the right-hand-sides of the isomorphisms
\eqref{eq:DeformComponentIdentity}.
The moduli space $\sM_{\ft',\fs}^{\stab}$ can also be written as the disk
bundle appearing on the left-hand-side of the following diagram:
\begin{equation}
\label{eq:UBundle}
\begin{CD}
\tilde N_{\ft,\fs}(\eps)\times_{\sG_{\fs}}\Gl_{\ft}(\delta) @>>>
N_{\ft,\fs}(\eps) \times X
\\
@V\pi_{N,\Gl}VV @V\pi_{N,X} VV
\\
\tM_{\fs}\times_{\sG_{\fs}}\Gl_{\ft}(\delta) @>\pi_{\Gl} >> M_{\fs}\times
X
\end{CD}
\end{equation}
We further suppose, without loss of generality, that
$[A,\Phi]\in\bL_{\ft',\fs}$ corresponds, via the gluing map $\bgamma$, to a
point $\bgamma^{-1}([A,\Phi])$ in the base of the bundle on the
left-hand-side of the diagram \eqref{eq:UBundle}, so
$$
\pi_{\Gl}(\bgamma^{-1}([A,\Phi]))=([B,\Psi],x). 
$$
The fiber of the projection $\pi_{\Gl}$ in diagram \eqref{eq:UBundle} is
$M^{s,\natural}_1(S^4,\delta)$, while the same diagram identifies the fiber
of the projection $\pi_{N,\Gl}$ with the fiber
$N_{\ft,\fs}(\eps)|_{[B,\Psi]}$ of the projection $\pi_{N,X}$. Thus,
\begin{equation}
\label{eq:DirectSumU}
\begin{aligned}
T_{\bgamma^{-1}([A,\Phi])}\sM_{\ft',\fs}^{\stab}
&\cong
    N_{\ft,\fs}|_{[B,\Psi]}
    \oplus T_{\bgamma^{-1}([A,\Phi])}
    (\tM_{\fs}\times_{\sG_{\fs}}\Gl_{\ft}(\delta))
\\
& \cong
    N_{\ft,\fs}|_{[B,\Psi]}
    \oplus
    T_{[B,\Psi]}M_{\fs}
    \oplus
    T_{[\bA]}M^{s,\natural}_1(S^4,\delta)
    \oplus
    T_xX.
\end{aligned}
\end{equation}
If we compare the isomorphisms \eqref{eq:DirectSumU} with
\eqref{eq:TangentOfIBoundary}, we see that the standard
orientation of $\bL_{\ft',\fs}$ is induced, through the convention
\eqref{eq:QuotientBoundaryOrientation}, by the orientation of
$\bchi^{-1}(0)\cap\sM^{\stab}_{\ft',\fs}$ given by the complex orientation of
the normal bundle of $\bchi^{-1}(0)\cap\sM^{\stab}_{\ft',\fs}$ and the
orientation of $\sM^{\stab}_{\ft',\fs}$ defined through the isomorphism
\eqref{eq:DirectSumU} by the complex orientation of
$N_{\ft,\fs}|_{[B,\Psi]}$, the orientation of $T_{[B,\Psi]}M_{\fs}$
determined by the homology orientation $\Om$, the standard orientation of
$T_{[\bA]}M^{s,\natural}_1(S^4,\delta)$, and the orientation of $T_xX$.
This orientation matches the orientation $O^{\asd}(\Om,c_1(L))$ as described
in the paragraph following \eqref{eq:DeformComponentIdentity}. This completes
the proof.
\end{proof}

We shall work with a fixed orientation $O^{\asd}(\Om,w)$
of $\sM_{\ft'}$ in the sum \eqref{eq:CobordismSum}
and thus we include the following lemma on
how the orientations change as the \spinc structure $\fs$
varies.

\begin{lem}
\label{lem:OrientationFactor}
If $\ft$ is a \spinu structure on $X$, let $\ft'$ be the spliced \spinu
structure \eqref{eq:SplicedSpinu}.  Let $\Om$ be a homology orientation and
let $w$ be an integral lift of $w_2(\ft')$.  If $\ft$ admits a splitting
$\ft=\fs\oplus \fs\otimes L$, then the standard orientation for
$\bL_{\ft',\fs}$ and the boundary orientation for $\bL_{\ft',\fs}$ defined
through the orientation $O^{\asd}(\Om,w)$ for $\sM_{\ft'}$ differ by a
factor of
\begin{equation}
\label{eq:OrientChangeFactor}
(-1)^{o_{\ft'}(w,\fs)},
\quad\text{where}\quad o_{\ft'}(w,\fs)
=
\textstyle{\quarter}(w-c_1(L))^2.
\end{equation}
\end{lem}

\begin{proof}
The result follows from Lemma \ref{lem:LinkOrientation} and the identity
$$
O^{\asd}(\Om,w)
=
(-1)^{\textstyle{\frac{1}{4}}(w-c_1(L))^2} O^{\asd}(\Om,c_1(L))
$$
given in \cite[Lemma 2.6]{FL2b}.
\end{proof}


\section{Pullback and extension of cohomology classes. Euler classes
of obstruction bundles}
\label{sec:cohom}
In \S \ref{subsec:ExtensionCohomClass} we compute the pullbacks of the
cohomology classes $\mu_p(\beta)$ and $\mu_c$ to
$\sM^{\stab}_{\ft',\fs}/S^1$ with respect to the gluing map
$\bgamma:\sM^{\stab}_{\ft',\fs}/S^1 \to \sC_{\ft'}^{*,0}/S^1$.  When no
confusion can arise, we denote the pullbacks $\bgamma^*\mu_p(\beta)$ and
$\bgamma^*\mu_c$ by $\mu_p(\beta)$ and $\mu_c$, respectively. We then
observe that these pullbacks are the restrictions of cohomology classes
$\bar\mu_p(\beta)$ and $\bar\mu_c$ on $\bar{\sM}^{\stab,*}_{\ft',\fs}/S^1$,
where $\bar{\sM}^{\stab,*}_{\ft',\fs}$ is defined in
\eqref{eq:DefineComplementOfReduciblesInSplicingDomain}.
In \S \ref{subsec:EulerClassObstruction} we calculate Euler classes of the
obstruction bundle \eqref{eq:ObstructionBundleEmbedding} over
$\sM^{\stab}_{\ft',\fs}/S^1$.  We will see that the background component
\eqref{eq:DefineBackgroundObstructionBundle} of the obstruction bundle
\eqref{eq:ObstructionBundleEmbedding} (and hence its Euler class) extends
over $\bar{\sM}^{\stab,*}_{\ft',\fs}/S^1$ while the Euler class of the
instanton component
\eqref{eq:InstantonObstructionSplicingDomain} extends over
$\bar{\sM}^{\stab,*}_{\ft',\fs}/S^1$ as a rational cohomology class.
Although we only perform the calculations relevant to the stratum
$$
(M_{\fs}\times\Sym^\ell(X))\cap\bar\sM_{\ft'}\subset\bar\sM_{\ft'}
$$
when $\ell(\ft',\fs)=1$ in the present article, we shall usually indicate
the nature of the changes required to address the general case $\ell\geq
1$.

\subsection{Pullbacks of cohomology classes}
\label{subsec:ExtensionCohomClass}
By analogy with our definition (following
\eqref{eq:StratificationCptPU(2)Space}) of $\bar\sM_{\ft}^*$ as the
subspace of $\bar\sM_{\ft}$ represented by ideal $\SO(3)$ monopoles with
reducible associated $\SO(3)$ connections, we define
\begin{equation}
\label{eq:DefineComplementOfReduciblesInSplicingDomain}
\bar\sM_{\ft',\fs}^{\stab,*}
=
\bar\sM_{\ft',\fs}^{\stab} - (M_{\fs}\times X),
\end{equation}
together with an inclusion map
\begin{equation}
\label{eq:AnInclusionUsedToExtendCohomClasses}
\iota:\sM^{\stab}_{\ft',\fs}/S^1\to\bar{\sM}^{\stab,*}_{\ft',\fs}/S^1.
\end{equation}
The circle action on $\bar\sM_{\ft',\fs}^{\stab}$,
given in Lemma \ref{lem:S1EquivarianceSplicing},
is free on $\bar\sM_{\ft',\fs}^{\stab,*}$ and
trivial on $M_{\fs}\times X$.

\begin{defn}
\label{defn:DefnOfNu}
Let $\nu$ be the first Chern class of the circle bundle
\begin{equation}
\label{eq:S1ActionOnStabilized}
\bar\sM_{\ft',\fs}^{\stab,*}\to\bar\sM_{\ft',\fs}^{\stab,*}/S^1,
\end{equation}
where the circle acts diagonally on
$\bar\sM_{\ft',\fs}^{\stab}=\tilde N_{\ft,\fs}(\eps)
\times_{\sG_{\fs}}\bar{\Gl}_{\ft}(\delta)$,
acting on $\tilde N_{\ft,\fs}(\eps)$ by scalar multiplication on
the fibers and on $\bar{\Gl}_{\ft}(\delta)$ by the action described
before equation \eqref{eq:GaugeGrpS1OnGluing}.
\end{defn}

Recall that $\mu_c\in H^2(\sC^{*,0}_{\ft'}/S^1;\ZZ)$ is the first Chern class
\eqref{eq:DefineMuC1} of the line bundle $\LL_{\ft'}\to \sC^{*,0}_{\ft'}/S^1$
defined in \eqref{eq:DefineDetLineBundle};
we now identify its pullback to $\sM^{\stab}_{\ft',\fs}/S^1$.

\begin{lem}
\label{lem:PullbackMuc1}
Let $\mu_c$ be as in the preceding paragraph, let $\nu$ be as in Definition
\ref{defn:DefnOfNu}, and let $\iota$ be the inclusion
\eqref{eq:AnInclusionUsedToExtendCohomClasses}. Then
$$
\bga^*\mu_c=-\iota^*\nu.
$$
\end{lem}

\begin{proof}
{}From \eqref{eq:DefineMuC1} and \eqref{eq:DefineDetLineBundle}, the
cohomology class $\mu_c$ is the first Chern class of
$$
\LL_{\ft'}=\sC^{*,0}_{\ft'}\times_{(S^1,\times -2)}\CC.
$$
By Lemma \ref{lem:S1EquivarianceSplicing}, the gluing map
$\bgamma:\sM_{\ft',\fs}^{\stab}\to\sC^{*,0}_{\ft'}$
is circle-equivariant when the circle acts on $\sC^{*,0}_{\ft'}$ by the action
\eqref{eq:S1ZAction} and on $\sM_{\ft',\fs}^{\stab}$ by the action in
Definition \ref{defn:DefnOfNu}, but with multiplicity two.
Thus,
$$
\bga^*\LL_{\ft'}=\sM_{\ft',\fs}^{\stab}\times_{(S^1,-1)}\CC,
$$
and the conclusion follows.
\end{proof}

\noindent Lemma \ref{lem:PullbackMuc1} and its proof translate easily to
the case $\ell\ge 1$.

Next, we identify the pullbacks to $\sM_{\ft',\fs}^{\stab}$
of the cohomology classes $\mu_p(\beta)$ on $\sC^{*,0}_{\ft'}$.
Let
\begin{equation}
\label{eq:ProjFromCptVirtModSpace}
\pi_{\fs}:\bar\sM_{\ft',\fs}^{\stab,*}\to M_{\fs}
\quad\text{and}\quad
\pi_{X}:\bar\sM_{\ft',\fs}^{\stab,*}\to X
\end{equation}
denote the restrictions of
the components of the projection $\bar\sM_{\ft',\fs}^{\stab}\to
M_{\fs}\times X$ given in \eqref{eq:GluingDataProjection}
to $\bar\sM_{\ft',\fs}^{\stab,*}$.
We define some additional projections: 
\begin{equation}
\label{eq:MoreProjFromCptVirtModSpace}
\begin{aligned}
\pi_{\sM}: \bar\sM_{\ft',\fs}^{\stab,*} \times X
&\to \bar\sM_{\ft',\fs}^{\stab,*},
\\
\pi_{X,2}: \bar\sM_{\ft',\fs}^{\stab,*}\times X
&\to X,
\\
\pi_{\fs,1}=\pi_{\fs}\circ\pi_{\sM}:
\bar\sM_{\ft',\fs}^{\stab,*} \times X
&\to M_{\fs},
\end{aligned}
\end{equation}
We shall use the same notation for the projections when the space
$\bar\sM_{\ft',\fs}^{\stab,*}$ above is replaced, for example, by its circle
quotient, $\bar\sM_{\ft',\fs}^{\stab,*}/S^1$.

Recall from definition \eqref{eq:UniversalU2Bundle} that
$$
\FF_{\ft'}=\tsC^{*}_{\ft'}/S^1\times_{\sG_{\ft'}}\fg_{V'}
$$
is a universal $\SO(3)$ bundle over $\sC^{*}_{\ft'}/S^1\times X$. We now
define an analogous $\SO(3)$ bundle over $\sM_{\ft',\fs}^{\stab}/S^1\times
X$. Using `$\cl$' for convenience here to indicate the `Uhlenbeck
compactification' implicit when ${\Gl}_{\ft}(\delta)$ is replaced by
$\bar{\Gl}_{\ft}(\delta)$, we denote
\begin{equation}
\label{eq:DefnTildeVirtModuli}
\cl\tilde{\sM}_{\ft',\fs}^{\stab,*}
=
\tilde N_{\ft,\fs}(\eps) \times \bar{\Gl}_{\ft}(\delta)
- \left(\tilde M_{\fs}\times X\right),
\end{equation}
so that
$\bar\sM_{\ft',\fs}^{\stab,*}=\cl\tilde{\sM}_{\ft',\fs}^{\stab,*}/\sG_{\fs}$, 
where $\sG_{\fs}$ acts diagonally on
$\tilde N_{\ft,\fs}(\eps)\times\bar{\Gl}_{\ft}(\delta)$. We set
\begin{equation}
\label{eq:UnivOnComplement}
\bar{\FF}_{\ft',\fs}^{\stab,*}
=
\cl\tilde{\sM}_{\ft',\fs}^{\stab,*}
       \times_{\sG_{\fs}\times S^1}(\underline{\RR}\oplus L)
      \to
\bar\sM_{\ft',\fs}^{\stab,*}/S^1 \times X,
\end{equation}
where $\sG_{\fs}$ acts diagonally on
$\cl\tilde{\sM}^{\stab,*}_{\ft',\fs}\times\fg_{V}$ --- acting on
$\fg_{V}\cong\underline{\RR}\oplus L$ by multiplication with weight
negative two on $L$ and on $\cl\tilde{\sM}^{\stab,*}_{\ft',\fs}$ by the
action described above; the circle acts diagonally on
$\cl\tilde{\sM}^{\stab,*}_{\ft',\fs}\times\fg_{V}$ --- by scalar
multiplication on $\tilde N_{\ft,\fs}(\eps)$, by the action on
$\bar{\Gl}_{\ft}(\delta)$ described before equation
\eqref{eq:GaugeGrpS1OnGluing}, and 
on $\fg_{V}$ by scalar multiplication with weight one on $L$.

Lemma \ref{lem:UnivOnComplement} below
compares the restriction of the bundle $\bar{\FF}_{\ft',\fs}^{\stab,*}$ with
the restriction of the pulled-back
$\SO(3)$ bundle $(\bga\times\id_X)^*\FF_{\ft'}$ to the complement of
$(\pi_X\times\id_X)^{-1}(\Delta)$ in
$\sM_{\ft',\fs}^{\stab}/S^1\times X$, where $\Delta\subset X\times X$ is the
diagonal. We restrict the pullback
of the universal bundle to the complement of the
subspace $(\pi_X\times\id_X)^{-1}(\Delta)$ because
the splicing process at a point $x\in X$ only identifies the restrictions
of the bundles $\fg_{V}$ and $\fg_{V'}$ to $X\less\{x\}$.

\begin{lem}
\label{lem:UnivOnComplement}
Suppose the \spinu structure $\ft=(\rho,V)$ admits a splitting
$\ft=\fs\oplus\fs\otimes L$, so that $\fg_{V}\cong
\underline{\RR}\oplus L$, where $\underline{\RR}=X\times\RR$.
If $\sO$ denotes the complement of $(\pi_X\times\id_X)^{-1}(\Delta)$ in
$\sM_{\ft',\fs}^{\stab}/S^1\times X$, then there is an isomorphism of
$\SO(3)$ bundles:
$$
\left.(\bga\times\id_X)^*\FF_{\ft'}\right|_{\sO}
\cong
\left.\bar{\FF}_{\ft',\fs}^{\stab,*}\right|_{\sO}.
$$
\end{lem}

\begin{proof}
Let $\ft'=(\rho,V')$ be the \spinu structure \eqref{eq:SplicedSpinu}
obtained by splicing, over a neighborhood of $x\in X$, the \spinu structure
$\ft=(\rho,V)$ over $X$ with the \spinu structure $(\rho,\bV)$ over $S^4$,
where $-\frac{1}{4}p_1(\fg_{\bV})=1$. We thus obtain an associated $\SO(3)$
bundle $\fg_{V'}$ as in \eqref{eq:SplicedSO(3)Bundle} and a bundle
isomorphism,
$$
\iota_{V,V'}:\fg_{V}|_{X\less \{x\}} \to \fg_{V'}|_{X\less \{x\}}.
$$
Hence, if $\tilde\sO\subset 
\tilde N_{\ft,\fs}(\eps)\times \tilde{\Gl}_{\ft}(\delta)\times X$ 
is the pre-image of $\sO$ under the obvious projection,
we have a bundle map
\begin{equation}
\label{eq:UnivBundleMapOnComplement1}
\tilde\bga'\times \iota_{V,V'}:
\left.\left(
\tilde N_{\ft,\fs}(\eps)\times \tilde{\Gl}_{\ft}(\delta)
\times
\fg_V
\right)\right|_{\tilde\sO}
\to
\tsC^*_{\ft'}\times \fg_{V'},
\end{equation}
where $\tilde\bga'$ is the pre-splicing map \eqref{eq:PreSplicingMap}.  The
map $\tilde\bga'$ is gauge equivariant with respect to the action of
$\sG_{\fs}\times\sG_1$ on the domain and the action of $\sG_{\ft'}$ on the
range.  The group $\sG_1$ of gauge transformations over $S^4$ 
(see \S \ref{subsec:S4ModuliSp}) acts trivially on the restriction
of $\fg_V$ to the complement of the splicing point $x\in X$.
{}From definitions \eqref{eq:DefnXGluingParameters}
and \eqref{eq:DefnTildeVirtModuli} we have
$$
\tilde N_{\ft,\fs}(\eps)\times \tilde{\Gl}_{\ft}(\delta)/\sG_1
=
\iota^*\cl\tilde {\sM}^{\stab,*}_{\ft',\fs},
$$
(where $\iota$ is the inclusion in
\eqref{eq:AnInclusionUsedToExtendCohomClasses})
and so the bundle map
\eqref{eq:UnivBundleMapOnComplement1} descends to a bundle map on
gauge-group quotients,
\begin{equation}
\label{eq:UniversalBundleOnComplement}
\bga'\times \iota_{V,V'}:
\left.\left(\cl\tilde{\sM}^{\stab,*}_{\ft',\fs}\times_{\sG_{\fs}}
\fg_{V}\right)\right|_{\sO}
\to
\tsC^{*}_{\ft'}\times_{\sG_{\ft'}}\fg_{V'}.
\end{equation}
By Lemma \ref{lem:S1EquivarianceSplicing}, the bundle map 
\eqref{eq:UniversalBundleOnComplement} is
circle-equivariant with respect to the circle action on the domain induced
by the action
\eqref{eq:S1ZActionOnN} on the factor $\tilde N_{\ft,\fs}(\eps)$ in 
$\cl\tilde\sM^{\stab,*}_{\ft',\fs}$ 
together with the trivial action on $\Gl_{\ft}(\delta)$
and the circle action on $\tsC^{*,0}_{\ft'}$ induced by
\eqref{eq:S1ZAction}.  The bundle map
\eqref{eq:UniversalBundleOnComplement} thus
descends to a bundle map on circle quotients,
\begin{equation}
\label{eq:UniversalBundleOnComplementS1Quotient}
\bga'\times \iota_{V,V'}:
\left.\left(\cl\tilde{\sM}^{\stab,*}_{\ft',\fs}/S^1\times_{\sG_{\fs}}
\fg_{V}\right)\right|_{\sO}
\to
\FF_{\ft'}.
\end{equation}
The argument in the proof of Lemma \ref{lem:S1EquivarianceSplicing} which
shows that the two circle actions on $\sM^{\stab}_{\ft',\fs}$ (described in
its hypothesis) are equal then implies that the circle action on
$\cl\tilde\sM^{\stab,*}_{\ft',\fs}\times\fg_{V}$ induced by the circle action
\eqref{eq:S1ZActionOnN} on the factor $\tilde N_{\ft,\fs}(\eps)$ in
$\cl\tilde\sM^{\stab,*}_{\ft',\fs}\times\fg_{V}$ and trivial actions on
$\Gl_{\ft}(\delta)$ and on $\fg_V$ is equal to the twice the circle action
described following \eqref{eq:UnivOnComplement}.  The multiplicity of this
action does not affect the quotient, so the bundle given as the domain of
the map \eqref{eq:UniversalBundleOnComplementS1Quotient} is isomorphic to
$\bar{\FF}_{\ft',\fs}^{\stab,*}$ and hence the restrictions of
$(\bga'\times\id_X)^*\FF_{\ft'}$ and $\bar{\FF}_{\ft',\fs}^{\stab,*}$ to
$\sO$ are isomorphic.  Finally, because the gluing map $\bga$ and the
splicing map $\bga'$ are $S^1$-equivariantly homotopic, there is a bundle
isomorphism $(\bga\times\id_X)^*\FF_{\ft'}\cong
(\bga'\times\id_X)^*\FF_{\ft'}$.  This completes the proof of the lemma.
\end{proof}

Since $\mu_p(\beta)=-\frac{1}{4}
p_1(\FF_{\ft'})/\beta$, our next task is to compare the Pontrjagin classes of
the bundles $\iota^*\bar{\FF}_{\ft',\fs}^{\stab,*}$ 
and $(\bga\times\id_X)^*\FF_{\ft'}$ over $\sM^{\stab}_{\ft',\fs}\times X$:

\begin{lem}
\label{lem:CompareCharClassOfUniv}
Continue the hypotheses of Lemma \ref{lem:UnivOnComplement}.
Let $\PD[\Delta]$ be the Poincar\'e dual of the diagonal
$\Delta \subset X\times X$. Then, on $\sM^{\stab}_{\ft',\fs}\times X$, we have
\begin{equation}
\label{eq:CharClassOfUniv1}
(\bga\times\id_X)^*p_1(\FF_{\ft'})
=
(\iota\times\id_X)^*
\left( p_1(\bar{\FF}_{\ft',\fs}^{\stab,*})
-
4(\pi_X\times\id_X)^*\PD[\Delta]\right),
\end{equation}
where $\iota$ is the inclusion
\eqref{eq:AnInclusionUsedToExtendCohomClasses}.
\end{lem}

\begin{proof}
The restriction of the difference,
\begin{equation}
\label{eq:p1Difference}
(\bga\times\id_X)^*p_1(\FF_{\ft'})
-
(\iota\times\id_X)^*p_1(\bar{\FF}_{\ft',\fs}^{\stab,*}),
\end{equation}
to the subspace $\sO$ defined in Lemma \ref{lem:UnivOnComplement}
vanishes by that lemma. 
Therefore, by considering the exact sequence of the pair,
\begin{equation}
\label{eq:DiagonalComplementPair}
\left( \sM^{\stab}_{\ft',\fs}/S^1\times X,\sO
\right),
\end{equation}
we see that the difference \eqref{eq:p1Difference} lies
in the image of the homomorphism
\begin{equation}
\label{eq:CompactCohomHomDiagComplement}
H^4\left( \sM^{\stab}_{\ft',\fs}/S^1\times X,\sO;\RR
\right)
\to
H^4\left(\sM^{\stab}_{\ft',\fs}/S^1\times X;\RR\right)
\end{equation}
appearing in the exact sequence of the pair
\eqref{eq:DiagonalComplementPair}.  Because the map
$\pi_X\times\id_X$ is transverse to $\Delta$, the image
of the homomorphism \eqref{eq:CompactCohomHomDiagComplement}
(see \cite[p. 69]{BT} and \cite[Proposition VIII.11.10]{Dold})
is generated by $(\pi_X\times\id_X)^*\PD[\Delta]$.  The difference
\eqref{eq:p1Difference} is therefore a multiple of
$(\pi_X\times\id_X)^*\PD[\Delta]$.   One can calculate
this multiple by evaluating the difference \eqref{eq:p1Difference}
on a chain which intersects $(\pi_X\times\id_X)^{-1}(\Delta)$
transversely at a single point.
Such computations are carried out in
the proof of \cite[Theorem III.6.1]{FrM} and in \cite[\S 5.3]{DK},
giving equation \eqref{eq:CharClassOfUniv1}
\end{proof}

We see from \eqref{eq:UnivOnComplement} that the cohomology class on the
right-hand side of  \eqref{eq:CharClassOfUniv1} is defined on
$\bar{\sM}^{\stab,*}_{\ft',\fs}/S^1\times X$ and thus extends the class
$(\bga\times\id_X)^*p_1(\FF_{\ft'})$, which is defined on
$\sM^{\stab}_{\ft',\fs}/S^1\times X$.

\begin{defn}
\label{defn:ExtendedCohomClasses}
We define extensions of cohomology classes
from $\sM^{\stab}_{\ft',\fs}/S^1$ to $\bar{\sM}_{\ft',\fs}^{\stab,*}/S^1$:
\begin{itemize}
\item
Let $\barmu_c=-\nu$ be the extension of $\mu_c$
given by Lemma \ref{lem:PullbackMuc1}.  
\item
For $\beta\in H_\bullet(X;\RR)$, let
$\barmu_p(\beta)$ be the extension of (the pullback by $\bgamma$ of)
the cohomology class $\mu_p(\beta)=-\frac{1}{4}p_1(\FF_{\ft'})/\beta$
given by replacing $(\bga\times\id_X)^*p_1(\FF_{\ft'})$ with 
the cohomology class on the right-hand side of \eqref{eq:CharClassOfUniv1}.
\end{itemize}
\end{defn}

In order to use equation \eqref{eq:CharClassOfUniv1} to compute 
$(\bga\times\id_X)^*p_1(\FF_{\ft'})$, we must identify the Pontrjagin class
$p_1(\bar{\FF}_{\ft',\fs}^{\stab,*})$.
The bundle $\bar{\FF}_{\ft',\fs}^{\stab,*}$ admits a reduction
$\bar{\FF}_{\ft',\fs}^{\stab,*}\cong\underline{\RR}\oplus
\bar{\LL}_{\ft',\fs}^{\stab,*}$, 
where $\underline{\RR}=\bar\sM_{\ft',\fs}^{\stab,*}/S^1 \times X\times\RR$ and
\begin{equation}
\label{eq:DefineReducingLineBundle}
\bar{\LL}_{\ft',\fs}^{\stab,*}
=
\cl\tilde\sM_{\ft',\fs}^{\stab,*}
       \times_{\sG_{\fs}\times S^1} L
      \to
\bar{\sM}_{\ft',\fs}^{\stab,*}/S^1 \times X.
\end{equation}
The actions of $S^1$ and $\sG_{\fs}$ in the definition of the line bundle
$\bar{\LL}_{\ft',\fs}^{\stab,*}$ are described in Lemma \ref{lem:UnivOnComplement}.
Since $p_1(\bar{\FF}_{\ft',\fs}^{\stab,*})= c_1(\bar{\LL}_{\ft',\fs}^{\stab,*})^2$, it
suffices to compute $c_1(\bar{\LL}_{\ft',\fs}^{\stab,*})$ and for this purpose we
shall use the following technical tool:

\begin{lem}
\label{lem:DiagonalQuotient}
\cite[Lemma 3.27]{FL2a}
If $Q_i$, for $i=1,2$, are circle bundles over a topological space
$M$ and $L_i=Q_i\times_{S^1}\CC$ are the associated complex line
bundles, and $k\in\ZZ$, then the following hold:
\begin{enumerate}
\item If $V$ is a complex vector bundle over $M$, and
  $e^{i\theta}\in S^1$ acts on the fiber product $Q_1\times_M V$ by
  $e^{i\theta}\cdot(q_1,v)=(e^{i\theta}q_1,e^{ik\theta}v)$, then
$$
(Q_1\times_M V)/S^1\cong L_1^{-k}\otimes V.
$$
\item If $e^{i\theta}\in S^1$ acts on the fiber product $Q_1\times_MQ_2$ by
  $e^{i\theta}\cdot(q_1,q_2)=(e^{i\theta}q_1,e^{ik\theta}q_2)$, then the
  first Chern class of the circle bundle $(Q_1\times_M Q_2)/S^1\rightarrow M$
  is 
$$
c_1((Q_1\times_M Q_2)/S^1) = c_1(Q_2)-kc_1(Q_1), 
$$
where the circle action on $(Q_1\times_M Q_2)/S^1$ is induced by the circle
action on $Q_2$ of weight one.
\end{enumerate}
\end{lem}

Recall that $\mu_{\fs} \in H^2(M_{\fs};\ZZ)$ is the first Chern class of
the base-point fibration over $M_{\fs}$ (see the remark following equation
\eqref{eq:SWClass}); when there is no ambiguity
we may write $\mu_{\fs}$ for the pullback $\pi^*_{\fs}\mu_{\fs}$ to the
total space $\bar\sM_{\ft',\fs}^{\stab}/S^1$.  Recall from \S
\ref{subsubsec:SWMonopoles} that if $b_1(X)=0$ then $\mu_{\fs}\times
1=c_1(\LL_{\fs})$, where $\LL_{\fs}\to
\sC_\fs^0\times X$ is the line bundle \eqref{eq:DefineSWUniversal}.
Moreover, $c_1(L)=c_1(\ft)-c_1(\fs)$ in the definition
\eqref{eq:DefineReducingLineBundle} of the line bundle
$\bar{\LL}_{\ft',\fs}^{\stab,*}$.

In the following we will use $\LL_{\fs}$ to
denote the restriction of the bundle
$\LL_{\fs}\to\sC_\fs^0\times X$ to the subspace $M_{\fs}\times X$.

\begin{lem}
\label{lem:UniversalLineBundle}
Continue the hypotheses of Lemma \ref{lem:UnivOnComplement}. Let
$$
\LL_{\nu}\to\bar{\sM}_{\ft',\fs}^{\stab,*}/S^1
$$
be the complex line bundle associated to the circle bundle
\eqref{eq:S1ActionOnStabilized}.  Then there is an isomorphism of complex
line bundles over $\bar\sM_{\ft',\fs}^{\stab,*}/S^1\times X$:
\begin{equation}
\label{eq:UniversalLineBundleIsom}
\bar{\LL}_{\ft',\fs}^{\stab,*}
\cong
\pi_{\sM}^*\LL_{\nu}^{-1}\otimes
(\pi_{\fs,1}\times\id_X)^*\LL_{\fs}^{\otimes 2}\otimes \pi_{X,2}^*L.
\end{equation}
\end{lem}

\begin{proof}
The argument yielding equation (3.68) in \cite{FL2a} implies that if
$\sG_{\fs}$ acts diagonally on $\tM_{\fs}\times L$ with weight negative two
on $L$, then we have an isomorphism of line bundles over $M_{\fs}\times X$:
\begin{equation}
\label{eq:TwistedSWUniversal1}
\tM_{\fs}\times_{\sG_{\fs}} L
\cong
\LL_{\fs}^{\otimes 2} \otimes \pi_{X,\fs}^* L,
\end{equation}
where $\pi_{X,\fs}:M_{\fs}\times X\to X$ is the projection.
By the preceding isomorphism and the definition 
\eqref{eq:MoreProjFromCptVirtModSpace} of $\pi_{\fs}\times\id_X$,
we have an isomorphism of line bundles over
$\bar\sM_{\ft',\fs}^{\stab,*}/S^1\times X$:
$$
(\pi_{\fs}\times\id_X)^*(\LL_{\fs}^{\otimes 2} \otimes \pi_{X,\fs}^* L)
\cong
\left(\bar\sM_{\ft',\fs}^{\stab,*}/S^1\times X\right)
\times_{M_{\fs}\times X}
\left(\tM_{\fs}\times_{\sG_{\fs}} L\right).
$$
Consequently, from this isomorphism we see that, as bundles over
$\bar\sM_{\ft',\fs}^{\stab,*}/S^1\times X$,
\begin{equation}
\label{eq:IsomorphicReducingBundleNoS1}
\begin{aligned}
{}&\left(\bar\sM_{\ft',\fs}^{\stab,*}\times X\right)
\times_{\bar\sM_{\ft',\fs}^{\stab,*}/S^1\times X}
(\pi_{\fs,1}\times\id_X)^*(\LL_{\fs}^{\otimes 2} \otimes \pi_{X,\fs}^* L)
\\
&\quad\cong
\left(\bar\sM_{\ft',\fs}^{\stab,*}\times X\right)
\times_{M_{\fs}\times X}
\left(\tM_{\fs}\times_{\sG_{\fs}}L\right).
\end{aligned}
\end{equation}
Now, consider the circle quotient of the bundle on the right-hand side of
the preceding isomorphism,
\begin{equation}
\label{eq:IsomorphicReducingBundle2}
\left(\bar\sM_{\ft',\fs}^{\stab,*}\times X
\times_{M_{\fs}\times X}
(\tM_{\fs}\times_{\sG_{\fs}}L)\right)/S^1
\to
\bar\sM_{\ft',\fs}^{\stab,*}/S^1\times X,
\end{equation}
where the circle acts diagonally on the factors $\bar\sM_{\ft',\fs}^{\stab,*}$
and $L$.  The definition of $\LL_{\nu}$ as the line bundle associated to
the circle bundle \eqref{eq:S1ActionOnStabilized},
Lemma \ref{lem:DiagonalQuotient} and
the isomorphism \eqref{eq:IsomorphicReducingBundleNoS1}
then yield the bundle isomorphisms over $\bar\sM_{\ft',\fs}^{\stab,*}\times X$:
\begin{equation}
\label{eq:IsomorphismReducingBundle3}
\begin{aligned}
{}&\left((\bar\sM_{\ft',\fs}^{\stab,*}\times X)
\times_{M_{\fs}\times X}
(\tM_{\fs}\times_{\sG_{\fs}}L)\right)/S^1
\\
&\quad\cong
\pi_{\sM}^*\LL_{\nu}^{-1}\otimes
(\pi_{\fs,1}\times\id_X)^*(\LL_{\fs}^{\otimes 2} \otimes \pi_{X,\fs}^* L)
\\
&\quad\cong
\pi_{\sM}^*\LL_{\nu}^{-1}\otimes
(\pi_{\fs,1}\times\id_X)^*\LL_{\fs}^{\otimes 2} \otimes \pi_{X,2}^* L.
\end{aligned}
\end{equation}
One can check that the map defined below is an isomorphism of line bundles
over $\bar\sM_{\ft',\fs}^{\stab,*}\times X$,
\begin{equation}
\label{eq:BundleIsom1}
\bar{\LL}_{\ft',\fs}^{\stab,*}
=
\cl\tilde\sM^{\stab,*}_{\ft',\fs}\times_{\sG_{\fs}\times S^1}L
\cong
\left((\bar\sM_{\ft',\fs}^{\stab,*}\times X)
\times_{M_{\fs}\times X}
(\tM_{\fs}\times_{\sG_{\fs}}L)\right)/S^1,
\end{equation}
given for $(B,\Psi)\in\tM_{\fs}$,
$(B,\Psi,\eta)\in\tN_{\ft,\fs}(\eps)$, $\bg\in\Gl_{\ft}(\delta)$, and
$\zeta\in L|_x$
by
$$
\left[ (B,\Psi,\eta),\bg,\zeta\right]
\mapsto
\left( ([(B,\Psi,\eta),\bg],x),[B,\Psi,\zeta]\right).
$$
Therefore, the isomorphisms \eqref{eq:BundleIsom1} and
\eqref{eq:IsomorphismReducingBundle3} give the desired isomorphism
\eqref{eq:UniversalLineBundleIsom}.
\end{proof}

\begin{cor}
\label{cor:UniversalLineBundle}
Continue the hypotheses of Lemma
\ref{lem:UniversalLineBundle}. Then, on
$\bar\sM_{\ft',\fs}^{\stab,*}\times X$, we have
\begin{equation}
\label{eq:ChernClassOfReducingLineBundle}
c_1(\bar{\LL}_{\ft',\fs}^{\stab,*})
=
2(\pi_{\fs,1}\times\id_X)^*c_1(\LL_{\fs})
-\pi_{\sM}^*\nu +\pi_{X,2}^*c_1(L).
\end{equation}
\end{cor}

Lemma \ref{lem:CompareCharClassOfUniv}, equation
\eqref{eq:ChernClassOfReducingLineBundle}, and the relation
$p_1(\bar{\FF}_{\ft',\fs}^{\stab,*})= c_1(\bar{\LL}_{\ft',\fs}^{\stab,*})^2$ then yield:

\begin{cor}
\label{cor:CharClassOfUniv}
Continue the hypotheses of Lemma \ref{lem:UnivOnComplement}.
Let $\PD[\Delta]$ be the Poincar\'e dual of the diagonal
$\Delta \subset X\times X$, let
$\nu$ be the class in Definition \ref{defn:DefnOfNu}, and let
$\iota$ be the inclusion \eqref{eq:AnInclusionUsedToExtendCohomClasses}.
Then, on $\sM^{\stab}_{\ft',\fs}/S^1$, we have
\begin{equation}
\label{eq:CharClassOfUniv}
\begin{aligned}
(\bga\times\id_X)^*p_1(\FF_{\ft})
&=
(\iota\times\id_X)^*\left( 2(\pi_{\fs,1}\times\id_X)^*c_1(\LL_{\fs})
-\pi_{\sM}^*\nu +\pi_{X,2}^*c_1(L)\right)^2
\\
&\quad
-
4(\iota\times\id_X)^*(\pi_X\times\id_X)^*\PD[\Delta].
\end{aligned}
\end{equation}
\end{cor}

The equation \eqref{eq:CharClassOfUniv} for $(\bga\times
\id_X)^*p_1(\FF_{\ft'})$ yields the following formulae for the pullback of
$\mu_p(\beta)=-\frac{1}{4}p_1(\FF_{\ft'})/\beta$ to
$\sM^{\stab}_{\ft',\fs}/S^1$: 

\begin{lem}
\label{lem:CohomologyOnReducibleLink}
Continue the hypotheses of Lemma
\ref{lem:UnivOnComplement}.
Let $\iota$ be the inclusion \eqref{eq:AnInclusionUsedToExtendCohomClasses},
let $x\in H_0(X;\ZZ)$ be the
positive generator, and let $h\in H_2(X;\RR)$. Let $\nu,\pi_{\fs}^*\mu_{\fs}\in
H^2(\bar{\sM}^{\stab,*}_{\ft',\fs}/S^1;\ZZ)$ 
be the classes given by Definition \ref{defn:DefnOfNu}
and equation \eqref{eq:SWClass}, respectively, while
$\pi_{\fs}$ is the projection \eqref{eq:ProjFromCptVirtModSpace}.
Then, on $\sM^{\stab}_{\ft',\fs}/S^1$, we have
\begin{equation}
\label{eq:CohomologyOnReducibleLink}
\begin{aligned}
\mu_p(x)
&=  
-\tquarter \iota^*(2\pi_{\fs}^*\mu_{\fs}-\nu)^2+\iota^*\pi_X^*\PD[x]
\in H^4(\sM^{\stab}_{\ft',\fs}/S^1;\RR) ,
\\
\mu_p(h)
&=
 \thalf \iota^*\langle c_1(\fs)-c_1(\ft),h\rangle (2\pi_{\fs}^*\mu_{\fs}-\nu)
+\iota^*\pi_X^*\PD[h]
\in H^2(\sM^{\stab}_{\ft',\fs}/S^1;\RR).
\end{aligned}
\end{equation}
\end{lem}

\begin{proof}
If $\beta\in H_\bullet(X;\RR)$, then $\PD[\Delta]/\beta=\PD[\beta]$ by
\cite[Theorem 30.6]{GreenbergHarper}.  Thus,
$$
(\pi_X\times\id_X)^*\PD[\Delta]/\beta
=
\pi_X^*\PD[\beta].
$$
The assertions now follow from Corollary \ref{cor:CharClassOfUniv} and
standard computations (see, for example, \cite[\S 5.1.4]{DK}).
\end{proof}

Finally, we note that the results of \S \ref{subsec:ExtensionCohomClass},
bearing on cohomology classes on a neighborhood $M_{\fs}\times X$, will
extend (with appropriate modifications, though these do not cause undue
difficulty) to the case of a neighborhood of $M_{\fs}\times\Sym^\ell(X)$
when $\ell>1$.  For example, one must replace the diagonal $\Delta\subset
X\times X$ appearing in Lemma \ref{lem:UnivOnComplement} with an incidence
subset,
$$
\{(\bx,x)\in\Sym^\ell(X)\times X: x\in |\bx|\},
$$
where $|\bx|\subset X$ is the support of a point $\bx\in \Sym^\ell(X)$. 

\subsection{Euler class of the  obstruction bundle}
\label{subsec:EulerClassObstruction}
We first compute the Euler class of the background
component $\bar\Upsilon_{\ft',\fs}^{s}/S^1\to 
\bar{\sM}^{\stab,*}_{\ft',\fs}/S^1$
of the obstruction bundle defined in
\eqref{eq:ExtendedBackgroundObstruction}.

\begin{lem}
\label{lem:EulerOfBackgroundObstruction}
Continue the hypotheses of Lemma \ref{lem:UnivOnComplement}.  Let
$\nu$ be the class in Definition \ref{defn:DefnOfNu}. Let
$\bar\Upsilon_{\ft',\fs}^{s}/S^1\to\bar{\sM}_{\ft',\fs}^{\stab,*}/S^1$
be the extended Seiberg-Witten obstruction bundle
\eqref{eq:ExtendedBackgroundObstruction},
and let $r_\Xi$ denote its complex rank.  Then
\begin{equation}
\label{eq:EulerOfBackgroundObstruction}
e(\bar\Upsilon_{\ft',\fs}^{s}/S^1)=(-\nu)^{r_\Xi}
\in
H^{2r_{\Xi}}(\bar{\sM}_{\ft',\fs}^{\stab,*}/S^1;\ZZ).
\end{equation}
\end{lem}

\begin{proof}
The $S^1$ quotient of the bundle $\bar\Upsilon_{\ft',\fs}^{s}$
by the action described at the end of
\S \ref{subsubsec:BackgroundObstruction}
can be written as
\begin{equation}
\label{eq:BckrndEuler2}
\bar\Upsilon_{\ft',\fs}^{s}/S^1
\cong \bar{\sM}_{\ft',\fs}^{\stab,*}\times_{(S^1,-1)}\CC^{r_\Xi},
\end{equation}
with the factor of negative one appearing because the action is diagonal,
as explained in Lemma \ref{lem:DiagonalQuotient}.  Note that the circle
action described in Lemma \ref{lem:S1EquivarianceSplicing} is twice the
action in the definition of $\nu$ (see Definition \ref{defn:DefnOfNu}), so
the weight two action on $\CC^{r_{\Xi}}$ described at the end of \S
\ref{subsubsec:BackgroundObstruction} becomes a weight one action
here. Equation \eqref{eq:EulerOfBackgroundObstruction} then follows
immediately from equation \eqref{eq:BckrndEuler2} and the definition of
$\nu$.
\end{proof}

We now compute the Euler class of the instanton component 
\eqref{eq:DefnInstantonObstruction} of
the obstruction bundle when $\ell(\ft',\fs)=1$; more work is required to
extend this calculation to the case $\ell(\ft',\fs)>1$. For the proof of
the following lemma, it is convenient to define
\begin{equation}
\label{eq:InstantonBdry}
\rd^i\sM^{\stab}_{\ft',\fs}
=
\tilde N_{\ft,\fs}(\eps)\times_{\sG_{\fs}}\partial\bar{\Gl}_{\ft}(\delta)
\subset
\sM^{\stab}_{\ft',\fs},
\end{equation}
so that $\bL_{\ft',\fs}^{\stab,i}=\rd^i\sM^{\stab}_{\ft',\fs}/S^1$.

\begin{lem}
\label{lem:InstantonEuler} 
Let $\nu$ be the class in Definition \ref{defn:DefnOfNu}. 
Then the Euler class of the
instanton obstruction bundle, $\Upsilon_{\ft',\fs}^{i}/S^1\to
\sM_{\ft',\fs}^{\stab}/S^1$
defined in \eqref{eq:DefnInstantonObstruction},
is given as an element of rational
cohomology by
\begin{equation}
\label{eq:InstantonEuler}
e(\Upsilon_{\ft',\fs}^{i}/S^1) 
=
\thalf\left(\pi_X^*c_1(\ft) -\nu \right)
\in
H^{2}({\sM}_{\ft',\fs}^{\stab}/S^1;\QQ).
\end{equation}
\end{lem}

\begin{proof}
Because $\sM^{\stab}_{\ft',\fs}/S^1$ retracts onto
$\rd^i\sM^{\stab}_{\ft',\fs}/S^1=\bL_{\ft',\fs}^{\stab,i}$, as one can see
from the definition \eqref{eq:DefineSWandInstantonModelThickLink} of
$\bL_{\ft',\fs}^{\stab,i}$, it suffices to compute the Euler class of the
restriction of $e(\Upsilon^i_{\ft',\fs}/S^1)$ to
$\rd^i\sM^{\stab}_{\ft',\fs}/S^1$.  It is easier to compute the Euler class
of the bundle $(\Upsilon_{\ft',\fs}^{i}/S^1)^{\otimes 2}$.  The observation
(following from \eqref{eq:ExplicitS4CokernelBundle}) that
$$
(\Coker\bD)^{\otimes 2}
\cong
\SU(2)\times_{\{\pm\id\}}(\CC\otimes\CC)
=\SO(3)\times\CC,
$$
Lemma \ref{lem:InstantonObstructionSplicingDomain},
and the description of $\Upsilon^i_{\ft',\fs}$ in
\eqref{eq:DefnInstantonObstruction} imply that
$$
\left.(\Upsilon_{\ft',\fs}^{i})^{\otimes 2}
\right|_{\rd^i\sM^{\stab}_{\ft',\fs}}
\cong
\tilde N_{\ft,\fs}(\eps)\times_{\sG_{\fs}}\Fr_{\CCl(T^*X)}(V)
\times_{\Spinu(4)}(\SO(3)\times\CC).
$$
By the description of the circle action on the bundle $\Upsilon_{\ft',\fs}^{i}$
prior to its definition \eqref{eq:DefnInstantonObstruction},
the preceding isomorphism is circle-equivariant if the circle acts on
the factor $\tN_{\ft,\fs}(\eps)$ by the action
\eqref{eq:S1ZActionOnN} and with weight two on the fiber $\CC$
and trivially on the remaining factors; the weight of the
circle action on the fiber $\CC$ is two because $\Upsilon_{\ft',\fs}^{i}$
appears in the tensor-product square on the left-hand side. The action of
$\Spinu(4)$ on the vector bundle $(\Coker\bD)^{\otimes 2}\to
M^{s,\natural}_1(S^4,\delta)$ is given, for $(\hat M,\zeta)\in\SO(3)\times
\CC=(\Coker\bD)^{\otimes 2}$ and $\tilde U\in\Spinu(4)$, by
\begin{equation}
\label{eq:SpinuGrpActionOnCoker}
(\hat M,\zeta) \mapsto (\Ad(\tilde U)\hat M,{\det}^u(\tilde U)^{-1}\zeta),
\end{equation}
where the homomorphisms ${\det}^u:\Spinu(4)\to S^1$ and
$\Ad=(\Ad^u,\Ad^{\fg}):\Spinu(4)\to\SO(4)\times\SO(3)$ are defined in
\eqref{eq:SpinuDeterminant} and \eqref{eq:DefineSpinuHom}, respectively.
To see this, observe that
\begin{enumerate}
\item
The action of $\Spinu(4)$ on
$\Coker\bD$ covers the action of $\SO(4)\times\SO(3)$ on
$M^{s,\natural}_1(S^4,\delta)$,
\item
The central $S^1$
in $\Spinu(4)$ acts on $\Coker\bD$ by scalar multiplication with
weight negative one.
\end{enumerate}
Property (1) implies that
the action of $\Spinu(4)$ on
the component $M^{s,\natural}_1(S^4,\delta)$ of
$$
(\Coker \bD)^{\otimes2}=M^{s,\natural}_1(S^4,\delta)\times\CC
$$ 
is given by the projection
$\Ad:\Spinu\to \SO(4)\times \SO(3)$ and the action of $\SO(4)\times\SO(3)$
on $M^{s,\natural}_1(S^4,\delta)$.  Property (2) implies that the central
$S^1$ acts on $(\Coker\bD)^{\otimes 2}$ by scalar multiplication on the
fibers with weight negative two, just as in the definition of the
homomorphism $\det^u$ before equation \eqref{eq:SpinuDeterminant}.  Any
group action satisfying the above two properties will differ from
the action
\eqref{eq:SpinuGrpActionOnCoker} by a representation of $\Spin(4)\subset
\Spinu(4)$ on $\CC$.
However, by \cite[Proposition 5.1]{BrockertomDieck} there are no such
non-trivial representations. Hence, these two properties characterize the
above group action.  Equations \eqref{eq:AssociatedSpinuBundles}
and \eqref{eq:DetuBundle} yield the bundle
isomorphisms,
\begin{align*}
\Fr_{\CCl(T^*X)}(V)\times_{(\Spinu(4),\Ad)}\SO(3)
&\cong
\Fr(T^*X)\times_X\Fr(\fg_V)\times_{\SO(4)\times\SO(3)}\SO(3)
\\
&=
\rd \bar{\Gl}_{\ft}(\delta),
\\
\Fr_{\CCl(T^*X)}(V)\times_{(\Spinu(4),\det^u)}\CC
&\cong
{\det}^{\frac{1}{2}}(V^+),
\end{align*}
and so 
$$
\Fr_{\CCl(T^*X)}(V)\times_{\Spinu(4)}(\SO(3)\times\CC)
\cong
\rd\bar{\Gl}_{\ft}(\delta)\times_X{\det}^{\frac{1}{2}}(V^+).
$$
The preceding isomorphism is $\sG_{\fs}$-equivariant, where $\sG_{\fs}$
acts on $\Fr_{\CCl(T^*X)}(V)$ as in the description of the instanton
obstruction bundle \eqref{eq:DefnInstantonObstruction}, if $\sG_{\fs}$ acts
trivially on ${\det}^{\frac{1}{2}}(V^+)$ and on
$\rd\bar{\Gl}_{\ft}(\delta)$ by the standard action on $\Fr(\fg_{V})$.
(The action is trivial on ${\det}^{\frac{1}{2}}(V^+)$ because elements of
$\varrho(\sG_{\fs})\subset\sG_{\ft}$ act on the fiber of
$\Fr_{\CCl(T^*X)}(V)$ by elements of $\Spinu(4)$ which are in the kernel of
the homomorphism ${\det}^u:\Spinu(4)\to S^1$.) Therefore,
\begin{equation}
\label{eq:InstantonEulerIsom0}
(\Upsilon_{\ft',\fs}^{i})^{\otimes 2}|_{\rd^i{\sM}^{\stab}_{\ft',\fs}}
\cong
\tilde N_{\ft,\fs}(\eps)\times_{\sG_{\fs}}\rd\bar{\Gl_{\ft}}(\delta)\times_X
{\det}^{\frac{1}{2}}(V^+).
\end{equation}
The bundle isomorphism \eqref{eq:InstantonEulerIsom0} is circle-equivariant
and so, from the definition \eqref{eq:InstantonBdry}
of $\rd^i{\sM}^{\stab}_{\ft',\fs}$,
we see that the isomorphism \eqref{eq:InstantonEulerIsom0}
descends to an isomorphism
\begin{equation}
\label{eq:InstantonEulerIsom1}
(\Upsilon_{\ft',\fs}^{i}/S^1)^{\otimes 2}|_{\rd^i{\sM}^{\stab}_{\ft',\fs}/S^1}
\cong
\left(
\rd^i{\sM}^{\stab}_{\ft',\fs}\times_X{\det}^{\frac{1}{2}}(V^+)
\right)/S^1,
\end{equation}
where the circle acts with weight two on the fibers of
${\det}^{\frac{1}{2}}(V^+)$ and the action \eqref{eq:S1ZActionOnN} on the
factor $\tilde N_{\ft',\fs}(\eps)$ in $\rd^i{\sM}^{\stab}_{\ft',\fs}$.  By
the equivalence of the circle actions on $\sM_{\ft',\fs}^{\stab}$ given by
Lemma \ref{lem:S1EquivarianceSplicing} and because the circle action of
Lemma \ref{lem:S1EquivarianceSplicing} is twice the action in the
Definition \ref{defn:DefnOfNu}, the circle quotient on the right-hand side
of \eqref{eq:InstantonEulerIsom1} is equivalent to one where the circle
acts diagonally by scalar multiplication on the fibers of
${\det}^{\frac{1}{2}}(V^+)$ and by the action of Definition
\ref{defn:DefnOfNu} on $\rd^i{\sM}^{\stab}_{\ft',\fs}
\subset\sM^{\stab}_{\ft',\fs}$.
The desired equation \eqref{eq:InstantonEuler} for the
Euler class of $\Upsilon_{\ft',\fs}^{i}/S^1$ then follows from the
isomorphism \eqref{eq:InstantonEulerIsom1}, the preceding description of
the circle quotient on the right-hand side of
\eqref{eq:InstantonEulerIsom1}, and Lemma \ref{lem:DiagonalQuotient}.
\end{proof}


\section{Intersection numbers and cohomology}
\label{sec:DualCohom}
In this section we prove the topological results necessary
to compute the intersection number 
\begin{equation}
\label{eq:DesiredSWLinkIntersectionNumber}
\#\left(\bar\sV(z)\cap\bar\sW^{\eta} \cap\bar\bL_{\ft',\fs}\right)
\end{equation}
where $\bar\bL_{\ft',\fs}$ is the link of
$M_{\fs}\times\Sym^\ell(X)$ in $\bar\sM_{\ft'}$, with $\ell(\ft',\fs)=1$ and
$\deg(z)+2(\eta+1) = \dim\sM_{\ft'}$.

We solve this problem in two stages. In \S \ref{subsec:Duality}
we prove that the intersection number 
\eqref{eq:DesiredSWLinkIntersectionNumber} is equal to
the pairing of certain cohomology classes with a homology class
$[\bar\bL^{\stab}_{\ft',\fs}]$ which can be understood as a fundamental
class of $\bar\bL^{\stab}_{\ft',\fs}$; see Proposition
\ref{prop:IntersectionNoToCupProduct}.  In \S \ref{subsec:PrelimComp} we
compute enough of the cohomology ring of $\bar\bL^{\stab}_{\ft',\fs}$ to
allow us to compute this pairing.

Since the gluing map $\bgamma:\bar\sM_{\ft',\fs}^{\stab}\to\bar\sC_{\ft'}$
is a circle-equivariant embedding of stratified spaces (a homeomorphism
preserving strata and restricting to a diffeomorphism on smooth strata),
when no confusion can arise we shall follow the convention stated in the
introduction to \S \ref{sec:cohom}. Thus, we shall not explicitly
distinguish in this section between cohomology classes on $\bar\sC_{\ft'}$
and their pullbacks to $\bar\sM_{\ft',\fs}^{\stab}$ or between cycles or
geometric representatives in $\bar\sC_{\ft'}$ and their pre-images in
$\bar\sM_{\ft',\fs}^{\stab}$. For example, we shall simply write
$\sV(\beta)$ and $\mu_p(\beta)$ for the pre-image
$\bgamma^{-1}(\sV(\beta))$ or pullback $\bgamma^*\mu_p(\beta)$ in
$\sM_{\ft',\fs}^{\stab}$ of the corresponding geometric representative or
cohomology class on
$\bgamma(\sM_{\ft',\fs}^{\stab})\subset\sC_{\ft'}^{*,0}$.

It is worth mentioning where we use the various properties of the gluing
maps and obstruction sections described in Theorem \ref{thm:GluingThm}, as
we shall occasionally exploit these properties without comment later this
section. The fact that $\bgamma$ is a homeomorphism from
$\bar\sM_{\ft',\fs}^{\stab}$ into $\bar\sC_{\ft'}$ ensures that the image
of $\bgamma$ contains an open neighborhood of
$\bar\sV(z)\cap\bar\sW\cap\bar\bL_{\ft',\fs}$ in $\bar\sM_{\ft'}$. The
gluing map is a smooth embedding of the top stratum
$\sM_{\ft',\fs}^{\stab}$ into $\sC_{\ft'}^{*,0}$ and, when orientations are
properly taken into account, orientation-preserving; thus, the pre-images
in $\sM_{\ft',\fs}^{\stab}$ of oriented, transverse intersections in
$\sM_{\ft'}$ are again appropriately oriented, and transverse. We use
continuity of $\bgamma$ on $\bar\sM_{\ft',\fs}^{\stab}$
in arguments involving extensions of certain
cocycles on $\bL_{\ft',\fs}$ to $\bar\bL_{\ft',\fs}$; see the introduction
to \S \ref{subsubsec:ExtendingCocycles} for precise statements.

Transversality and smoothness of the obstruction section $\bchi$ is used in
Lemma \ref{lem:RelPairing}, where we use the fact that
$\bgamma(\bchi^{-1}(0))$ is a submanifold of $\sM_{\ft'}$ of the
appropriate `multiplicity'. Continuity of the section $\bchi_s$ on
$\bar\sM_{\ft',\fs}^{\stab}$ is used in the proofs of Lemmas
\ref{lem:Intersection2} and \ref{lem:Intersection3} and in
writing the decomposition \eqref{eq:DecomposeRelEulerClass} of a certain
relative Euler class.

\subsection{Duality}
\label{subsec:Duality}
The goal of this subsection is to prove Proposition
\ref{prop:IntersectionNoToCupProduct}. 

\begin{defn}
\label{defn:ExtendedEulerClasses}
Let $\bar e_s$ and $\bar e_i$ be the extensions
of the Euler classes $e(\Upsilon_{\ft',\fs}^{s}/S^1)$ and
$e(\Upsilon_{\ft',\fs}^{i}/S^1)$ from
$\sM^{\stab}_{\ft',\fs}/S^1$ to $\bar\sM^{\stab,*}_{\ft',\fs}/S^1$ defined by
Lemmas \ref{lem:EulerOfBackgroundObstruction}
and \ref{lem:InstantonEuler}
respectively, and set $\bar e = \bar e_s\smile\bar e_i$.
\end{defn}

\begin{prop}
\label{prop:IntersectionNoToCupProduct}
Assume $w\in H^2(X;\ZZ)$ is such that $w\pmod{2}$ is good
in the sense of Definition \ref{defn:Good}.
Suppose $(\ft',\fs)$ is a pair 
with $\ell(\ft',\fs)=1$ and $w_2(\ft')\equiv w\pmod{2}$.
Let $[\bar\bL^{\stab}_{\ft',\fs}]\in H_{\max}(\bar\bL^{\stab}_{\ft',\fs};\ZZ)$
be the homology class defined in
\eqref{eq:DefineFundClassOfL}.  Let $z\in\AAA(X)$ and $\eta\in\NN$ satisfy
$\deg(z)+2(\eta+1)=\dim\sM_{\ft'}$.
Let $\barmu_p(z)$ and $\barmu_c$ be the extensions
of $\mu_p(z)$ and $\mu_c$ in Definition
\ref{defn:ExtendedCohomClasses} from
$\sM^{\stab}_{\ft',\fs}/S^1$ to $\bar\sM^{\stab,*}_{\ft',\fs}/S^1$.
Give $\bL_{\ft',\fs}$ the standard orientation described in
\S \ref{subsec:Orient}. Then,
\begin{equation}
\begin{aligned}
\#\left( \bar\sV(z)\cap\bar\sW^{\eta} \cap\bar\bL_{\ft',\fs}\right)
&=
\left\langle
    \barmu_p(z)\smile \barmu_c^{\eta}\smile
    \bar e,
    [\bar\bL^{\stab}_{\ft',\fs}]
\right\rangle.
\end{aligned}
\label{eq:CohomFormulation}
\end{equation}
\end{prop}

The corresponding result \cite[Lemma 3.31]{FL2b} for the level-zero
Seiberg-Witten stratum, 
$M_{\fs}\subset \sM_{\ft}/S^1$,
followed trivially from the definition of a geometric representative (see
\cite[Definition 3.4]{FL2b} or \cite[p. 588]{KMStructure}) because the link
of $M_{\fs}\subset \sM_{\ft}/S^1$
is a smooth, compact manifold without boundary
whose fundamental class can be represented by a smooth cycle intersecting
the geometric representatives transversely.

By contrast, the link of the level-one Seiberg-Witten stratum,
$$
(M_{\fs}\times X)\cap\bar\sM_{\ft'}/S^1 \subset \bar\sM_{\ft'}/S^1,
$$
can have non-empty intersection with
a lower level of $\bar\sM_{\ft'}/S^1$. This raises two difficulties which
prevent an immediate translation of the level-zero argument in \cite{FL2b} to
the level-one case here:
\begin{enumerate}
\item
The obstruction section $\bchi$ does not vanish transversely
on the lower level $\bL^{\sing}_{\ft',\fs}$ of $\bar\bL^{\stab}_{\ft',\fs}
=\bL^{\stab}_{\ft',\fs}\sqcup\bL^{\sing}_{\ft',\fs}$
(see \eqref{eq:StratificationOfVirtLink}), so it 
is not immediately obvious that the link 
$\bar\bL_{\ft',\fs}=\bchi^{-1}(0)\cap\bar\bL_{\ft',\fs}^{\stab}$
defines a homology class.  
\item
It is not obvious that the closure $\bar\sV(\beta)$ 
in the compactification $\bar\sM_{\ft',\fs}^{\stab}/S^1$ of the geometric
representative $\sV(\beta)$ in $\sM_{\ft',\fs}^{\stab}/S^1$ 
defines a geometric representative for the cohomology class
$\barmu_p(\beta)$. (See Definition \ref{defn:GeomRepresentative} 
for a review of the concept of a
geometric representative for a cohomology class.)
\end{enumerate}

The intersection number on the left-hand side of \eqref{eq:CohomFormulation} is
a count with sign of the points in the transverse intersection of
the geometric representatives and the zero-locus, $\bL_{\ft',\fs}$, 
of the section $\bchi$ of the obstruction bundle over $\bL_{\ft',\fs}^{\stab}$.
We emphasize that this intersection is contained in the top stratum of
$\bar\bL^{\stab}_{\ft',\fs}$; see Lemma \ref{lem:GRCompactSupport}.
As discussed in the following paragraphs, the  geometric representatives
and the obstruction section define a cohomology class
on $\bL^{\stab}_{\ft',\fs}$
with compact support in the  manifold-with-boundary $\bL^{\stab,i}_{\ft',\fs}$
(defined in \eqref{eq:DefineSWandInstantonModelThickLink}).
This compactly supported cohomology class
is Poincar\'e dual to the intersection
\begin{equation}
\label{eq:TheIntersectionWeKeepReferringTo}
\bar\sV(z)\cap\bar\sW^{\eta} \cap\bar\bL_{\ft',\fs}.
\end{equation}
Therefore, the intersection number
\eqref{eq:CohomFormulation} is equal to the pairing of this compactly
supported cohomology class with the relative fundamental
class $[\bL_{\ft',\fs}^{\stab,i},\rd\bL_{\ft',\fs}^{\stab,i}]$
of $\bL^{\stab,i}_{\ft',\fs}$.

The first step in defining the compactly supported cohomology class
is to define a representative of the cohomology class 
$\mu_p(z)\smile\mu_c^{\eta}$ with specified support.
In \S \ref{subsubsec:DefiningCocycles} we observe that
the geometric representatives $\sV(\beta)$ and $\sW$ define singular
cocycles $c_\beta$ and $c_\sW$, respectively, in the cohomology classes
$\mu_p(\beta)$ and $\mu_c$.
The cocycle $c_\beta$ has support on $\sV(\beta)\cap\sM_{\ft'\fs}^{\stab}/S^1$,
in the sense that $c_\beta$ restricts to zero on
any singular chain not intersecting $\sV(\beta)$.
Therefore, $c_\beta$ defines a relative cohomology class, $[c_\beta]$,
which maps to $\mu_p(\beta)$ under the homomorphism
$$
H^{\deg(\beta)}(\sM_{\ft'\fs}^{\stab}/S^1,
\sM_{\ft'\fs}^{\stab}/S^1-\sV(\beta)\cap\sM_{\ft'\fs}^{\stab}/S^1;\RR)
\to
H^{\deg(\beta)}(\sM_{\ft'\fs}^{\stab}/S^1;\RR),
$$
given by the exact sequence of the pair.  We say a relative cohomology
class $[c]\in H^\bullet(X,A;\RR)$ is a representative of $\mu\in
H^\bullet(X;\RR)$ if $\jmath_A^*[c]=\mu$, where $\jmath_A:(X,\emptyset)\to
(X,A)$ is the inclusion map.  Thus, $[c_\beta]$ is a representative of
$\mu_p(\beta)$.  Similar comments apply to the support of the cocycle
$c_{\sW}$.  We will write $c(z,\eta)$ for the cup-product of the cocycles
defined by the geometric representatives appearing in the left-hand side of
equation \eqref{eq:CohomFormulation}.  The cocycle $c(z,\eta)$ defines a
relative cohomology class, $[c(z,\eta)]$, supported on the intersection of
the geometric representatives.

The {\em relative Euler class\/}
\cite{Kervaire} of the obstruction bundle and section, 
$e(\Upsilon_{\ft',\fs}/S^1,\bchi)$, has support on the zero-locus in
$\sM_{\ft',\fs}^{\stab}$ of $\bchi$.  In \S
\ref{subsubsec:RelEulerClass} (see \eqref{eq:ProductCompactSupport}),
we show that the cup-product $[c(z,\eta)]\smile
e(\Upsilon_{\ft',\fs}/S^1,\bchi)$ has compact support in the top level
$\bL^{\stab}_{\ft',\fs}$ of $\bar\bL^{\stab}_{\ft',\fs}$.  Note that
because the geometric representatives can intersect the lower strata of
$\bar\bL^{\stab}_{\ft',\fs}$, the class $[c(z,\eta)]$ does not have compact
support in $\bL^{\stab,i}_{\ft',\fs}$; it is only the cup-product
$[c(z,\eta)]\smile e(\Upsilon_{\ft',\fs}/S^1,\bchi)$ which has this compact
support.  In Lemma \ref{lem:RelPairing}, we then prove that this
cup-product is Poincar\'e dual to the intersection
\eqref{eq:TheIntersectionWeKeepReferringTo}
and hence the intersection number \eqref{eq:DesiredSWLinkIntersectionNumber}
is equal to the pairing
\begin{equation}
\label{eq:RelPairing0}
\left\langle
    [c(z,\eta)]\smile e(\Upsilon_{\ft',\fs}/S^1,\bchi),
    \left[\bL_{\ft',\fs}^{\stab,i},\rd\bL_{\ft',\fs}^{\stab,i}\right]
\right\rangle.
\end{equation}
Computing the pairing of the relative classes in
\eqref{eq:RelPairing0} directly does not seem practical.
The representatives $[c_\beta]$ of $\mu_p(\beta)$,
$[c_{\sW}]$ of $\mu_c$, and
$e(\Upsilon_{\ft',\fs}/S^1,\bchi)$ of $e(\Upsilon_{\ft',\fs}/S^1)$
in the cup-product in \eqref{eq:RelPairing0}
are elements of different
cohomology rings (compare \eqref{eq:WhereRelEulerObstLives} and
\eqref{eq:WhereRelCupProdCocyleLives}), so it is not possible to compute
their product using only the algebraic structure of the cohomology ring
$H^{\bullet}(\bL_{\ft',\fs}^{\stab,i},\rd\bL_{\ft',\fs}^{\stab,i};\RR)$.
We overcome this difficulty by replacing the pairing
\eqref{eq:RelPairing0}  with
a pairing with relative cohomology classes in
$H^{\max}(\bar\bL^{\stab}_{\ft',\fs},\bar\bL^{\stab,s}_{\ft',\fs};\RR)$.
There is an isomorphism
$$
\jmath^*:H^{\max}(\bar\bL^{\stab}_{\ft',\fs},\bar\bL^{\stab,s}_{\ft',\fs};\RR)
\cong
H^{\max}(\bar\bL^{\stab}_{\ft',\fs};\RR),
$$
because $\bar\bL^{\stab,s}_{\ft',\fs}$ retracts onto the codimension-four
subspace $\bL^{\sing}_{\ft',\fs}$.  Thus, a pairing
with  relative cohomology classes in
$H^\bullet(\bar\bL^{\stab}_{\ft',\fs},\bar\bL^{\stab,s}_{\ft',\fs};\RR)$
will depend only on the image of the
relative cohomology classes under $\jmath^*$, and thus only on the
absolute cohomology classes which these relative cohomology classes
represent.  These absolute cohomology classes are all elements of
the ring $H^\bullet(\bar\bL^{\stab}_{\ft',\fs};\RR)$, allowing us to
compute products in the algebra of that single ring.

To replace the pairing \eqref{eq:RelPairing0} with one involving the
compactification, we construct extensions of the cocycles $c_\beta$ and
$c_{\sW}$ from $\bL^{\stab}_{\ft',\fs}$ to $\bar\bL^{\stab}_{\ft',\fs}$ in
\S \ref{subsubsec:ExtendingCocycles}. (This extension process for cocycles
is simpler than the corresponding process for geometric representatives.)
We find exact cocycles $\delta\theta_\beta$ and $\delta\theta_{\sW}$ such
that $c_{\beta}+\delta\theta_{\beta}=\iota^*\bar c_\beta$ and
$c_{\sW}+\delta\theta_{\sW}=\iota^*\bar c_{\sW}$, where $\bar c_{\beta}$
and $\bar c_{\sW}$ are cocycles on $\bar\bL^{\stab}_{\ft',\fs}$ and $\iota$
is the inclusion map \eqref{eq:AnInclusionUsedToExtendCohomClasses}.  By
keeping track of the support of $\delta\theta_{\beta}$ and
$\delta\theta_{\sW}$, we can prove that replacing $c_\beta$ and $c_{\sW}$
with $\iota^*\bar c_{\beta}$ and $\iota^*\bar c_{\sW}$ does not change the
pairing \eqref{eq:RelPairing0}.  The construction of $\bar c_{\beta}$ and
$\bar c_{\sW}$ implies that they represent the cohomology classes on the
right-hand-side of identity \eqref{eq:CohomFormulation}, namely
$\bar\mu_p(\beta)$ and $\bar\mu_c$.

Having applied this extension method to the cocycle $c(z,\eta)$,
we next apply it to the relative Euler class in the pairing
\eqref{eq:RelPairing0}.
In \S \ref{subsubsec:ElimExcis}, we describe a homotopy of
the obstruction section $\bchi$
which does not change the pairing \eqref{eq:RelPairing0}.
We then recognize the relative Euler class of the obstruction bundle and
this homotoped obstruction section as the restriction of
a relative Euler class on $\bar\bL^{\stab}_{\ft',\fs}$.

These extensions accomplish the goal of replacing the relative cohomology
classes in \eqref{eq:RelPairing0} with relative cohomology classes in
$H^{\bullet}(\bar\bL^{\stab}_{\ft',\fs},\bar\bL^{\stab}_{\ft',\fs};\RR)$.
Proposition \ref{prop:IntersectionNoToCupProduct} then follows, by the
argument sketched above, from Lemma \ref{lem:RelPairing} and our
computations in
\S \ref{subsubsec:ExtendingCocycles} of the cohomology classes
which the extended cocycles represent.

The preceding argument contains, albeit implicitly, solutions to
the two previously mentioned problems preventing a direct translation of
\cite[Corollary 3.11]{FL2b}.  The homotopy of the obstruction section
should define a deformation of $\bar\bL_{\ft',\fs}$ in
$\bar\bL^{\stab}_{\ft',\fs}$ to a subspace whose intersection with the
lower strata of $\bar\bL^{\stab}_{\ft',\fs}$ is sufficiently regular to
define a homology class.  Similarly, changing the cocycle $c_\beta$ to
$c_\beta+\delta \theta_\beta$ should be equivalent to changing the
geometric representative $\sV(\beta)$ (by a cobordism) to a geometric
representative whose closure is a geometric representative for
$\barmu_p(\beta)$.  We believe that the cohomological formulations given
here are easier to construct than these deformations of
$\bar\bL_{\ft',\fs}$ and the geometric representatives.

Employing the techniques of virtual fundamental classes developed in
\cite{Brussee}, \cite{LiTian}, \cite{GraberPand}, \cite{RuanSW}
\cite{RuanTian}, and \cite{SiebertGW}
to construct a homology class representing $\bar\bL_{\ft',\fs}$ --- and thus
eliminating the first of the two problems discussed above --- would not
simplify the proof of Proposition \ref{prop:IntersectionNoToCupProduct}.
While those techniques are well developed for the moduli space of
pseudoholomorphic curves, applying that theory to the moduli spaces in this
article would require additional discussion.  
For example, we would have to construct an extension of the obstruction
bundle $\Upsilon_{\ft',\fs}/S^1$ over an appropriate compactification of
$\sM^{\stab}_{\ft',\fs}/S^1$; this would require further work as the
instanton
component $\Upsilon_{\ft',\fs}^{i}/S^1$ does not extend from
$\sM^{\stab}_{\ft',\fs}/S^1$ to the relatively simple `cone' compactification
$\bar\sM^{\stab}_{\ft',\fs}/S^1$ in 
\eqref{eq:UhlenbeckCompactifiedSplicingDomain} which we currently employ.
In contrast, the method used here --- involving a relative Euler
class and a homotopy of the obstruction section rather than constructing
a homology class for $\bar\bL_{\ft',\fs}$ --- does not require any elaborate
technical apparatus.  Moreover, defining a virtual fundamental class would
only solve the first of the two previously mentioned problems --- defining a
fundamental class for $\bar\bL_{\ft',\fs}$; it would still be necessary to
construct the extensions of the cocycles $c_\beta$ and $c_{\sW}$ and to
prove that they represent the cohomology classes $\barmu_p(\beta)$ and
$\barmu_c$.
As the construction of the extended cocycles $\bar c_\beta$ and $\bar
c_{\sW}$ and the proof that they represent the cohomology classes
$\barmu_p(\beta)$ and $\barmu_c$ make up the bulk of the proof of
Proposition
\ref{prop:IntersectionNoToCupProduct}, a construction of a fundamental
class for $\bar\bL_{\ft',\fs}$ would not shorten this section
appreciably.

\subsubsection{The fundamental class of the virtual link}
\label{subsubsec:FundClass}
We begin by defining the homology class
$[\bar\bL^{\stab}_{\ft',\fs}]\in H_{\max}(\bar\bL^{\stab}_{\ft',\fs};\ZZ)$
referred to in Proposition \ref{prop:IntersectionNoToCupProduct}.
We refer to this homology class as the `fundamental class' of the virtual
link, although a precise definition of a fundamental class for
a stratified space will not be necessary for this article.
The existence and uniqueness of this class is an easy application of the
exact sequence of the pair
$(\bar\bL^{\stab}_{\ft',\fs},\bar\bL^{\stab,s}_{\ft',\fs})$, where
$\bar\bL^{\stab,s}_{\ft',\fs}$ is the neighborhood 
\eqref{eq:DefineSWandInstantonModelThickLink} of the lower
level of $\bar\bL^{\stab}_{\ft',\fs}$. 
The definition of $[\bar\bL^{\stab}_{\ft',\fs}]$
will extend from the case of $\ell(\ft',\fs)=1$ to the
general case of $\ell(\ft',\fs)\geq 1$
as it only requires the existence of a neighborhood
of $\bL^{\sing}_{\ft',\fs}$ in $\bar\bL^{\stab}_{\ft',\fs}$
which retracts onto $\bL^{\sing}_{\ft',\fs}$.
We note, however, that the construction of this neighborhood
is more difficult in the general case.

By definition \eqref{eq:DefineSWandInstantonModelThickLink}
of $\bar\bL^{\stab,s}_{\ft',\fs}$
there is a deformation retraction of $\bar\bL^{\stab,s}_{\ft',\fs}$ onto
$\bL_{\ft',\fs}^{\sing}$ (see \eqref{eq:DefineLSing})
induced by the deformation retraction
$\bar{\Gl}_{\ft'}(\delta)\to X$ (given by shrinking the scale to zero).
We have an inclusion map of pairs:
\begin{equation}
\label{eq:DefineJ}
\jmath: \left( \bar\bL^{\stab}_{\ft',\fs},\emptyset\right)
\to
\left( \bar\bL^{\stab}_{\ft',\fs},\bar\bL^{\stab,s}_{\ft',\fs}\right).
\end{equation}
Because $\bar\bL^{\stab,s}_{\ft',\fs}$ retracts onto
$\bL_{\ft',\fs}^{\sing}$
and because $\bL_{\ft',\fs}^{\sing}$ has codimension
four in $\bar\bL^{\stab}_{\ft',\fs}$,
the inclusion \eqref{eq:DefineJ} induces an isomorphism,
\begin{equation}
\label{eq:JCohomology}
\jmath_*: H_{\max}(\bL^{\stab}_{\ft',\fs};\ZZ)
\cong
H_{\max}(\bL^{\stab}_{\ft',\fs},\bL^{\stab,s}_{\ft',\fs};\ZZ).
\end{equation}
Let 
\begin{equation}
\label{eq:DefineRelFundClass}
\left[\bL_{\ft',\fs}^{\stab,i},\rd\bL_{\ft',\fs}^{\stab,i}\right]
\in
H_{\max}(\bL_{\ft',\fs}^{\stab,i},\rd\bL_{\ft',\fs}^{\stab,i};\ZZ)
\end{equation}
be the relative fundamental class \cite[p. 303]{Spanier}
of the manifold-with-boundary $\bL^{\stab,i}_{\ft',\fs}$.
By the construction \eqref{eq:DefineSWandInstantonModelThickLink}
of $\bL^{\stab,i}_{\ft',\fs}$ and $\bar\bL^{\stab,s}_{\ft',\fs}$, 
the boundary of $\bL^{\stab,i}_{\ft',\fs}$ lies
in $\bL^{\stab,s}_{\ft',\fs}$ and so there is an inclusion map
\begin{equation}
\label{eq:ILinkRelInclusion1}
\bar\jmath_{\bL^i}:
\left(\bL_{\ft',\fs}^{\stab,i},\rd\bL_{\ft',\fs}^{\stab,i}\right)
\to
\left(\bar\bL^{\stab}_{\ft',\fs},\bar\bL^{\stab,s}_{\ft',\fs}\right).
\end{equation}
Using the isomorphism \eqref{eq:JCohomology}, we then define
$[\bar\bL^{\stab}_{\ft',\fs}]$ to be the unique homology class satisfying
\begin{equation}
\label{eq:DefineFundClassOfL}
(\bar\jmath_{\bL^i})_*
\left[\bL_{\ft',\fs}^{\stab,i},\rd\bL_{\ft',\fs}^{\stab,i}\right] 
= 
\jmath_*[\bar\bL^{\stab}_{\ft',\fs}].
\end{equation}

\subsubsection{Deforming the geometric representatives}
\label{subsubsec:DeformGR}
The intersection of the
geometric representatives used to define the Donaldson invariants
in  \cite{KMStructure} with the lower strata
can only be easily understood on the set of triples $[A,\Phi,\bx]$ where
the sets of points representing $\bx\in\Sym^\bullet(X)$ do not intersect the
`suitable neighborhoods' \cite{KMStructure} used to define the geometric
representatives.  In this section we prove
that the intersection number on the right-hand side of
\eqref{eq:CohomFormulation} 
is equal to one defined by geometric representatives where
suitable neighborhoods are replaced by tubular neighborhoods,
simplifying our cohomological computations in
\S \ref{subsubsec:ExtendingCocycles}.

For $\beta\in H_\bullet(X;\RR)$, let $T_\beta$ be a smooth submanifold of
$X$ with fundamental class $[T_\beta]=\beta$.  The geometric representative
$\sV(\beta)$ was defined in \cite[\S 3.2]{FL2b} by pulling back a geometric
representative from the quotient
space of connections over a `suitable' neighborhood
$U_\beta$ of $T_\beta$.  Recall that a {\em suitable neighborhood\/} of
$T_\beta$ was defined in \cite[Definition 3.8]{FL2b} or
\cite[p. 589]{KMStructure} as a smoothing of
the union of a tubular neighborhood of $T_\beta$ and a set of loops
$\{\ga_i\}$, where the loops $\ga_i$ generate $H_1(X;\ZZ/2\ZZ)$, are
mutually disjoint, and are transverse to $T_\beta$. We need in \cite{FL2b}
to use a suitable neighborhood rather than a tubular neighborhood because
if $H^1(X;\ZZ/2\ZZ)\neq 0$, there could be a point $[A,0]\in\sM_{\ft'}$
such that the restriction of the induced $\SO(3)$ connection $\hat A$ on
$\fg_{V}$ to a tubular neighborhood of $T_\beta$ would be reducible
\cite[p. 586]{KMStructure} even if $\hat A$ is not globally reducible.  We
note that the assumption that $w\pmod 2$ is good (as defined prior to
\eqref{eq:DefineASDLink}) implies that there are no reducible, zero-section
pairs in $\bar\sM_{\ft'}$; Lemmas \ref{lem:NoTwistedRed},
\ref{lem:GRCompactSupport}, and \ref{lem:DeformedIntersection}
rely on this constraint on $w$.
The following lemma shows that when defining
geometric representatives near the strata of reducible pairs (but not
zero-section pairs) it suffices to use tubular rather than suitable
neighborhoods of $T_\beta$:

\begin{lem}
\label{lem:NoTwistedRed}
Assume $w\in H^2(X;\ZZ)$ is such that $w\pmod{2}$ is good.
Given a Riemannian metric on $X$ and a pair $(\ft',\fs)$
with $\ell(\ft',\fs)=1$ and $w_2(\ft')\equiv w\pmod{2}$,
there are positive constants $\eps_0$ and $\delta_0$ such that for all
positive $\eps\leq\eps_0$ and $\delta\leq\delta_0$ used in the definition
of ${\sM}^{\stab}_{\ft',\fs}$, the following holds.
For $\beta\in H_\bullet(X;\RR)$, let $T_\beta$ be the submanifold
of $X$ with $[T_\beta]=\beta$ and let
$\nu(\beta)$ be a tubular neighborhood of $T_\beta$.
Then, for any $[A,\Phi,\bx]\in\bga(\bar{\sM}^{\stab,*}_{\ft',\fs})$,
the restriction of $\hat A$ to $\nu(\beta)$ is not reducible.
\end{lem}

\begin{proof}
Because $w_2(\ft')=w_2(\ft)$ is good, there are no zero-section
pairs in $M_{\fs}$.
Thus, if
the parameters $\eps$ and $\delta$ in the definition of
$\sM^{\stab}_{\ft',\fs}$ are sufficiently small, then the neighborhood
$\bga(\bar{\sM}^{\stab,*}_{\ft',\fs})$ is disjoint from $\bar M_{\kappa}^w$.
Hence, for any point $[A,\Phi,\bx]\in \bga(\bar{\sM}^{\stab,*}_{\ft',\fs})$,
the section $\Phi$ is not identically zero and $\hat A$ is
not a reducible connection. By \cite[Theorem 3.10]{FL2b} or \cite[Theorem
5.10]{FL1}, the restriction of $\hat A$ to any open subspace of $X$
cannot be reducible.
\end{proof}

Lemma \ref{lem:NoTwistedRed} implies that there is a geometric
representative, $\sV'(\beta)$ for $\mu_p(\beta)$, pulled back from the
quotient space of irreducible $\SO(3)$ connections, $\sB^*(\nu(\beta))$, over
the tubular neighborhood $\nu(\beta)$.  This geometric representative is
then constructed by the same methods used to construct $\sV(\beta)$ in
\cite{KMStructure}.  The only difference between $\sV(\beta)$ and
$\sV'(\beta)$ is that $\sV(\beta)$ is pulled back from $\sB^*(U_\beta)$
rather than $\sB^*(\nu(\beta))$.

Henceforth, we assume that the $z\in\AAA(X)$ considered in
Proposition \ref{prop:IntersectionNoToCupProduct} is a monomial,
$z=\beta_1\beta_2\cdots\beta_m$.
We can then choose the submanifold $T_{\beta_i}$ defining
the geometric representative $\sV'(\beta_i)$,
the points $x_1,\dots,x_{\eta}$
defining the geometric representatives $\sW$,
and their tubular neighborhoods so that for any point $x\in X$ we have
\begin{equation}
\label{eq:RestrictionCondition}
\sum_{\{\beta_i:\ x\in \nu'(\beta_i)\}} \deg(\beta_i) +
\sum_{\{x_j:\ x\in \nu'(x_j)\}} 4
\le 4,
\end{equation}
where $\nu'(\beta_i)$ is a tubular neighborhood of $T_{\beta_i}$ with
$\nu(\beta_i)\Subset\nu'(\beta_i)$;
compare \cite[Equation (2.7)]{KMStructure}.
Let $\sV'(z)=\cap_i\sV'(\beta_i)$ and let $\bar\sV'(z)$
denote the Uhlenbeck closure of $\sV'(z)$.
The following lemma shows we can define an intersection number
of these new geometric representatives with $\bar\bL_{\ft',\fs}$.

\begin{lem}
\label{lem:GRCompactSupport}
Assume $w\in H^2(X;\ZZ)$ is such that $w\pmod{2}$ is good.
Suppose $(\ft',\fs)$ is a pair 
with $\ell(\ft',\fs)=1$ and $w_2(\ft')\equiv w\pmod{2}$. Then
the following hold:
\begin{enumerate}
\item
The intersection,
$\bar\sV'(z)\cap\bar\sW^{\eta}\cap\bar\bL_{\ft',\fs}$,
is contained in the top stratum $\bL_{\ft',\fs}$ of
$\bar\bL_{\ft',\fs}\subset\bar\sM_{\ft'}/S^1$ and is disjoint
from the edge \eqref{eq:Edge}.
\item
The following intersection is empty,
$$
\bar\sV'(z)\cap\bar\sW^{\eta}
\cap\bgamma\left(\bchi_s^{-1}(0)\cap\bL_{\ft',\fs}^{\sing}\right)
=
\emptyset.
$$
\end{enumerate}
\end{lem}

\begin{proof}
Inequality \eqref{eq:RestrictionCondition} and the
dimension-counting argument used in the proof of \cite[Corollary
3.18]{FL2b} imply that the intersection
$\bar\sV'(z)\cap\bar\sW^{\eta}$ is disjoint from the lower level,
$(\bar\sM_{\ft'}-\sM_{\ft'})\cap\bga(\bar\sM^{\stab,*}_{\ft',\fs})/S^1$.
The first assertion then follows from the proof of
Lemma \ref{lem:GRIntersect0}.
The image $\bga(\bchi_s^{-1}(0)\cap\bL_{\ft',\fs}^{\sing})$ is the
intersection of $\bga(\bar\bL^{\stab}_{\ft',\fs})$ with the union of the
lower levels, $(\bar\sM_{\ft'}-\sM_{\ft'})/S^1$.  Therefore, the second
assertion follows from the first and the fact that $\bgamma$ preserves
strata on $\bar\sM_{\ft',\fs}^{\stab}$.
\end{proof}

We now prove that replacing $\bar\sV(z)$ with $\bar\sV'(z)$ does
not change the intersection number in \eqref{eq:CohomFormulation}.

\begin{lem}
\label{lem:DeformedIntersection}
Assume $w\in H^2(X;\ZZ)$ is such that $w\pmod{2}$ is good.
Suppose $(\ft',\fs)$ is a pair 
with $\ell(\ft',\fs)=1$ and $w_2(\ft')\equiv w\pmod{2}$.
Then, for $z\in\AAA(X)$ and $\eta\in\NN$ as
in Proposition \ref{prop:IntersectionNoToCupProduct}, we have:
\begin{equation}
\label{eq:DeformedIntersection}
\#\left( \bar\sV(z)\cap\bar\sW^{\eta}\cap\bar\bL_{\ft',\fs}\right)
=
\#\left( \bar\sV'(z)\cap\bar\sW^{\eta}\cap\bar\bL_{\ft',\fs}\right).
\end{equation}
\end{lem}

\begin{proof}
Both sides of \eqref{eq:DeformedIntersection} are linear in $z$, so
we may assume without loss of generality that $z=\beta_1\cdots\beta_m$ for
$\beta_i\in H_\bullet(X;\RR)$.

First, we note that the intersection numbers on both sides of
\eqref{eq:DeformedIntersection} do not change if we decrease
the parameter $\delta$ defining $\bL^{\stab}_{\ft',\fs}$
(to another positive generic value) because of the obvious cobordism defined
by this change of parameter.  Thus, in the proof we may assume that the
parameter $\delta$ is as small as desired.

There is a cobordism
$$
\sH \subset \bga(\sM^{\stab}_{\ft',\fs}/S^1) \times [0,1],
$$
with boundaries given by
$$
\left(\sV(\beta)\cap \bga(\sM^{\stab}_{\ft',\fs}/S^1)\right)\times \{0\}
\quad\text{and}\quad
\left(\sV'(\beta)\cap\bga(\sM^{\stab}_{\ft',\fs}/S^1)\right)\times\{1\}.
$$
Such a cobordism exists because $\sV(\beta)$ and $\sV'(\beta)$
are defined by pullbacks (by the appropriate restriction maps)
of zero loci of sections of the same line bundle, for $\beta\in H_2(X;\RR)$,
or by the degeneracy locus of sections of the same vector bundle
for $\beta\in H_0(X)$.
If $z=z_1\beta_1$, where
$z_1\in\AAA(X)$ and $\beta_1\in H_\bullet(X;\RR)$, we replace
$\bar\sV(z)$ with $\bar\sV(z_1)\cap\bar\sV'(\beta_1)$ in the intersection
\eqref{eq:DeformedIntersection} as follows.
By perturbing the  cobordism $\sH$, we can assume that $\sH$
is transverse to
\begin{equation}
\label{eq:Cobordism2}
\left(
\bar\sV(z_1)\cap\bar\sW^{\eta}\cap\bar\bL_{\ft',\fs}\cap\bga(\bL^{\stab,i}_{\ft',\fs})
\right)
\times [0,1].
\end{equation}
The dimension-counting arguments in the proof of Corollary 3.18 in
\cite{FL2b} show that the closure $\bar\sH$ 
in $\bga(\bar\sM^{\stab}_{\ft',\fs}/S^1)\times [0,1]$ of the cobordism $\sH$
will not intersect
$$
\left(
    \bar\sV(z_1)\cap\bar\sW^\eta\cap\left(\bar\bL_{\ft',\fs}-\bL_{\ft',\fs}
    \right)
\right)
\times [0,1].
$$
Then, for a sufficiently small parameter $\delta$, the space
$\bar\sH$ will not intersect
$$
\left(
    \bar\sV(z_1)\cap\bar\sW^{\eta}\cap\bar\bL_{\ft',\fs}
    \cap
    \bga(\bar\bL^{\stab,s}_{\ft',\fs})
\right)
\times [0,1].
$$
Therefore the intersection
$\bar\sH
\cap
\left(
\bar\sV(z_1)\cap\bar\sW^{\eta}\cap\bar\bL_{\ft',\fs}\times [0,1]
\right)$
is contained in the transverse intersection of $\sH$ with the space
\eqref{eq:Cobordism2}. Hence, this last intersection is 
a family of smooth, compact, oriented one-dimensional submanifolds with
one boundary given by the set of points in the intersection on the
left-hand side of \eqref{eq:DeformedIntersection}
and the other boundary given by the set of points in the intersection
$$
\bar\sV(z_1)\cap\bar\sV'(\beta_1)\cap\bar\sW^\eta\cap\bar\bL_{\ft',\fs}.
$$
We now repeat this process with each $\beta_i$
in the product $z$ until $\bar\sV(z)$ has been replaced by $\bar\sV'(z)$.
This completes the proof.
\end{proof}

From this point onwards, we shall work with the geometric representative
$\bar\sV'(z)$ defined by tubular rather than suitable neighborhoods.
By choosing generic parameters $\eps$ and $\delta$ in the definition of
$\bar\bL^{\stab}_{\ft',\fs}$ (see
\eqref{eq:DefineSWandInstantonModelThickLink}), we can also ensure that
these geometric representatives are transverse to
$\bL^{\stab,i}_{\ft',\fs}$.  To simplify notation, we shall omit the primes,
writing $\bar\sV(z)$ for $\bar\sV'(z)$, as the
original geometric representatives defined by suitable neighborhoods will not
reappear in the remainder of this article.

\subsubsection{Geometric representatives and cocycles}
\label{subsubsec:DefiningCocycles}
We now define cocycles dual to the geometric representatives
$\sV(\beta)$ and $\sW$.  We first recall the following definition:

\begin{defn}
\label{defn:GeomRepresentative}
\cite[p 588]{KMStructure}.
Let $Z$ be a smoothly stratified space.  A {\em geometric
representative\/} for a real cohomology class $\mu$ of
dimension $d$ on $Z$ is a closed, smoothly stratified subspace
$\sV$ of $Z$ together with a real coefficient $q$,
the {\em multiplicity\/}, satisfying
\begin{enumerate}
\item
The intersection $Z_0\cap \sV$ of $\sV$ with the top stratum $Z_0$ of
$Z$ has codimension $d$ in $Z_0$ and has an oriented normal
bundle.
\item
The intersection of $\sV$ with all strata of $Z$ other than the top
stratum has codimension $2$ or more in $\sV$.
\item
The pairing of $\mu$ with a homology class $h$ of dimension $d$
is obtained by choosing a smooth singular cycle $\sigma$ representing $h$
whose intersection with all strata of $\sV$ has codimension
$\dim Z_0-d$ in that stratum of $\sV$, and then taking $q$ times
the count (with signs) of the intersection points between the
cycle and the top stratum of $\sV$:
$$
\langle \mu,h\rangle = q\cdot \#((Z_0\cap\sV)\cap\sigma).
$$

\end{enumerate}
\end{defn}

The intersection of the geometric representative with the top
stratum with real coefficient,
$(\sV\cap Z_0,q)$, defines a singular cocycle $c$
representing the restriction of the cohomology class $\mu$
to $Z_0$ with the properties described in the following lemma.

\begin{lem}
\label{lem:DefineSingularCocycle}
Let $Z_0$ be a smooth manifold and let $(\sV,q)$ be a geometric
representative for a real cohomology class $\mu$ on $Z_0$.
Then, there is a singular cocycle $c$ on $Z_0$
which represents the cohomology
class $\mu$ and which vanishes when restricted to $Z_0-\sV$.
\end{lem}

\begin{proof}
For any smooth manifold $W$, let $\Delta^{\8}_\bullet(W)$ denote
the chain complex of smooth singular chains \cite[p. 291]{BredonTopGeom}
and let $S^\bullet_{\8}(W;\RR)=\Hom(\Delta^{\8}_\bullet(W),\RR)$ be
the complex of smooth singular cochains.
For any smooth submanifold $Y\subset W$, we define the complex of 
smooth singular, relative cochains by
$$
S^\bullet_{\8}(W,Y;\RR)
=
\Hom(\Delta^{\8}_{\bullet}(W)/\Delta^{\8}_{\bullet}(Y),\RR).
$$
We will write $H^\bullet_{\8}(W;\RR)$ and $H^\bullet_{\8}(W,Y;\RR)$ for the
homology of the complexes $S^\bullet_{\8}(W;\RR)$ and
$S^\bullet_{\8}(W,Y;\RR)$, respectively. Thus,
$H^\bullet_{\8}(W;\RR)$ is the smooth singular cohomology
of $W$ and $H^\bullet_{\8}(W,Y;\RR)$ is the smooth singular, relative
cohomology of $(W,Y)$.
By the de Rham theorem (see the discussion in \cite[p. 291]{BredonTopGeom}),
there is a functorial isomorphism $H^p(W;\RR)\cong H^\bullet_{\8}(W;\RR)$.
By applying the Five Lemma to the long exact
sequences of the pair $(W,Y)$ in singular and smooth singular cohomology,
we obtain an isomorphism $H^{\bullet}(W,Y;\RR)\cong H^{\bullet}_{\8}(W,Y;\RR)$.

We first define a smooth singular cocycle, that is,
a closed  $c'\in S^\bullet_{\8}(Z_0;\RR)$,
which represents the  cohomology class $\mu$.
For any smooth singular chain $\si$ of dimension equal to that of $\mu$, set
\begin{equation}
\label{eq:CocycleDefinition}
\langle c',\si\rangle = q\cdot \#(\sV\cap \si'),
\end{equation}
where $\si'$ is any smooth singular chain in $Z_0$
which is homologous in $Z_0$ to $\si$ and transverse to $\sV$.
By definition of $\sV$, if $\si$ is a smooth singular cycle in $Z_0$ which
represents a homology class $h$, we then have
$$
\langle c',\si\rangle = \langle \mu,h\rangle,
$$
so $c'$ is a  smooth singular cocycle representing
the cohomology class $\mu$ on $Z_0$.
(By the de Rham isomorphism, it does not matter whether we consider $\mu$ as
a singular cohomology class or a smooth singular cohomology class.)
Because the restriction of $c'$ to $Z_0-\sV$ vanishes, the cocycle $c'$ defines
an element $[c']$ of the smooth singular, relative cohomology
$H^{\bullet}_{\8}(Z_0,Z_0-\sV;\RR)$ which maps to $\mu$
under the homomorphism
$H^{\bullet}_{\8}(Z_0,Z_0-\sV;\RR)\to H^{\bullet}_{\8}(Z_0;\RR)$.
Let $\alpha\in H^{\bullet}(Z_0,Z_0-\sV;\RR)$ be the 
element of relative singular cohomology given by the image of $[c']$
under the isomorphism
$H^{\bullet}_{\8}(Z_0,Z_0-\sV;\RR)\cong H^{\bullet}(Z_0,Z_0-\sV;\RR)$
and let $\hat c\in S^{\bullet}(Z_0,Z_0-\sV;\RR)$ be any representative of
$\alpha$.  Then the element $c\in S^{\bullet}(Z_0;\RR)$ given by the image 
of $\hat c$ under the homomorphism
$S^{\bullet}(Z_0,Z_0-\sV;\RR)\to S^{\bullet}(Z_0;\RR)$
will satisfy the conclusion of the lemma.
\end{proof}

For $\beta\in H_\bullet(X;\RR)$, we let $c_\beta$ be the cocycle on
$\sM_{\ft',\fs}^{\stab}/S^1$ representing $\mu_p(\beta)$,
defined by Lemma \ref{lem:DefineSingularCocycle}
and the geometric representative $\sV(\beta)$;
we let $c_{\sW}$ be the cocycle on $\sM_{\ft',\fs}^{\stab}/S^1$ representing
$\mu_c$, defined as in Lemma \ref{lem:DefineSingularCocycle},
by the geometric representative $\sW$.
Note that by considering the smooth manifold $Z_0$ in
Lemma \ref{lem:DefineSingularCocycle} to be the image of
$\sM^{\stab}_{\ft',\fs}/S^1$ 
under the composition of the gluing map and the restriction map
$r_{\nu(\beta)}$, we can assume the cocycle $c_\beta$ is pulled back
from $\sB^*(\nu(\beta))$.

\subsubsection{Relative Euler classes and intersection numbers}
\label{subsubsec:RelEulerClass}
We next prove that the intersection number in equation
\eqref{eq:CohomFormulation} is equal to a pairing of a
relative cohomology class with the relative fundamental class
$[\bL_{\ft',\fs}^{\stab,i},\rd\bL_{\ft',\fs}^{\stab,i}]$ of
$\bL_{\ft',\fs}^{\stab,i}$. 

\begin{defn}
If $e(V)$ is the Euler class of a real-rank $r$, oriented vector bundle $V$
a CW complex over $M$ and $\chi:U\subset M\to V$ is a continuous,
local section of $V$,
then the {\em relative Euler class\/} \cite{Kervaire} of the pair $(V,\chi)$
is a relative cohomology class,
$$
e(V,\chi)\in H^r(M,U-\chi^{-1}(0);\ZZ),
$$
satisfying the relation
\begin{equation}
\label{eq:RelativeEulerMap}
j^*_{\chi}e(V,\chi)=e(V),
\end{equation}
defined by the inclusion map:
$$
j_{\chi}:(M,\emptyset)\to (M,U-\chi^{-1}(0)).
$$
\end{defn}

The relative Euler class $e(\Upsilon_{\ft',\fs}/S^1,\bchi)$
of the obstruction bundle $\Upsilon_{\ft',\fs}/S^1$
and obstruction section $\bchi$ of
Theorem \ref{thm:GluingThm} thus defines a relative cohomology
class 
\begin{equation}
\label{eq:WhereRelEulerObstLives}
e(\Upsilon_{\ft',\fs}/S^1,\bchi)
\in
H^{2r_\Xi+2}(\bL^{\stab}_{\ft',\fs},\bL^{\stab}_{\ft',\fs}
-\bchi^{-1}(0);\ZZ).
\end{equation}
Recall that $r_\Xi$ is the complex rank of the background component
of the obstruction bundle, $\bar\Upsilon^s_{\ft',\fs}/S^1$, defined
in \S \ref{subsubsec:BackgroundObstruction}.
For $z=\beta_1\cdots\beta_m\in\AAA(X)$, define a subspace
$$
\sK(z,\eta)
=
\sV(z)\cap\sW^{\eta}\cap\bL^{\stab}_{\ft',\fs}
\subset
\bL^{\stab}_{\ft',\fs},
$$
with closure
\begin{equation}
\label{eq:KClosure}
\bar\sK(z,\eta)
=
\bar\sV(z)\cap\bar\sW^{\eta}\cap\bar\bL^{\stab}_{\ft',\fs}
\subset
\bar\bL^{\stab}_{\ft',\fs}.
\end{equation}
The cocycle on $\sM^{\stab}_{\ft',\fs}$,
\begin{equation}
\label{eq:DefineCocycleProduct}
c(z,\eta)
:=
c_{\beta_1}\smile\dots\smile c_{\beta_m}\smile c_{\sW}^{\eta},
\end{equation}
vanishes on the complement of $\sK(z,\eta)$ in
$\sM^{\stab}_{\ft',\fs}$, and
so defines a relative cohomology class
\begin{equation}
\label{eq:WhereRelCupProdCocyleLives}
[c(z,\eta)]
\in
H^{\deg(z)+2\eta}(\bL^{\stab}_{\ft',\fs},
\bL^{\stab}_{\ft',\fs}-\sK(z,\eta);\RR).
\end{equation}
By Lemma \ref{lem:GRCompactSupport} and the continuity
of the gluing map $\bga$, the following intersection is empty:
$$
\bar\sK(z,\eta)\cap \bL_{\ft',\fs}^{\sing}\cap \bchi^{-1}(0)
=
\emptyset.
$$
Therefore, noting that the link definition
\eqref{eq:DefineSWandInstantonModelThickLink} gives
$$
\bar\bL^{\stab,s}_{\ft',\fs}
=
\bL^{\stab,s}_{\ft',\fs}
\sqcup
\left(\rd N_{\ft,\fs}(\eps)/S^1\times X\right),
$$
and that $\bar\bL^{\stab,s}_{\ft',\fs}$ retracts onto the lower
level $\rd N_{\ft,\fs}(\eps)/S^1\times X$ of $\bar\bL^{\stab,s}_{\ft',\fs}$
by shrinking the instanton scale $\delta$, we can also arrange that the
following intersection is empty by taking $\delta$ to
be sufficiently small:
\begin{equation}
\label{eq:SmallEnoughdeltaImpliesEmpty}
\sK(z,\eta)\cap \bL^{\stab,s}_{\ft',\fs}\cap \bchi^{-1}(0)
=
\emptyset.
\end{equation}
Because cup product of relative cohomology is a map \cite[p. 251]{Spanier},
$$
H^i(Y,A)\otimes H^j(Y,B)\to H^{i+j}(Y,A\cup B),
$$
and \eqref{eq:SmallEnoughdeltaImpliesEmpty} yields
$$
\left(\bL^{\stab}_{\ft',\fs}-\bchi^{-1}(0)\right)
\cup
\left(\bL^{\stab}_{\ft',\fs}-\sK(z,\eta)\right)
\supset
\bL^{\stab,s}_{\ft',\fs},
$$
we see from \eqref{eq:WhereRelEulerObstLives} and
\eqref{eq:WhereRelCupProdCocyleLives} that
the cup product $[c(z,\eta)]\smile e(\Upsilon_{\ft',\fs}/S^1,\bchi)$
vanishes on $\bL^{\stab,s}_{\ft',\fs}$ and therefore we have
\begin{equation}
\label{eq:ProductCompactSupport}
[c(z,\eta)]\smile e(\Upsilon_{\ft',\fs}/S^1,\bchi)
\in
H^{\max}\left(\bL^{\stab}_{\ft',\fs},\bL^{\stab,s}_{\ft',\fs};\RR\right).
\end{equation}
The inclusion
\begin{equation}
\label{eq:DefineInclusion}
\iota:\bL^{\stab}_{\ft',\fs}\to\bar\bL^{\stab}_{\ft',\fs},
\end{equation}
induces an isomorphism of relative homology by excision,
\begin{equation}
\label{eq:Excision}
\iota_*:
H_\bullet\left(\bL^{\stab}_{\ft',\fs},\bL^{\stab,s}_{\ft',\fs};\RR\right)
\cong
H_\bullet\left(\bar\bL^{\stab}_{\ft',\fs},
\bar\bL^{\stab,s}_{\ft',\fs};\RR\right).
\end{equation}
Hence, we can express the intersection number in
\eqref{eq:CohomFormulation} as
a pairing of relative homology and cohomology classes: 

\begin{lem}
\label{lem:RelPairing}
If $[\bL_{\ft',\fs}^{\stab,i},\rd\bL_{\ft',\fs}^{\stab,i}]$
is the relative fundamental
class of $\bL^{\stab,i}_{\ft',\fs}$,
if $\bar\jmath_{\bL^i}$ is the inclusion map \eqref{eq:ILinkRelInclusion1}, 
and  
if
$\iota_*$ is the excision isomorphism \eqref{eq:Excision},
then
\begin{equation}
\label{eq:Intersection1}
\begin{aligned}
{}&\#\left( \bar\sV(z)\cap\bar\sW^{\eta} \cap\bar\bL_{\ft',\fs}\right)
\\
&=
\left\langle
[c(z,\eta)]\smile e\left(\Upsilon_{\ft',\fs}/S^1,
    \left.\bchi\right|_{\bL^{\stab,s}_{\ft',\fs}}\right),
\iota_*^{-1}
(\bar\jmath_{\bL^i})_*\left[\bL_{\ft',\fs}^{\stab,i},\rd\bL_{\ft',\fs}^{\stab,i}\right]
\right\rangle.
\end{aligned}
\end{equation}
\end{lem}

\begin{proof}
Let
$$
\jmath_{\bL^i}:
\left(\bL_{\ft',\fs}^{\stab,i},\rd\bL_{\ft',\fs}^{\stab,i}\right)
\to
\left(\bL^{\stab}_{\ft',\fs},\bL^{\stab,s}_{\ft',\fs}\right)
$$
be the inclusion map of pairs.  The equality
$\bar\jmath_{\bL^i}=\iota\circ \jmath_{\bL^i}$  
implies that
$$
\iota_*^{-1}
(\bar\jmath_{\bL^i})_*\left[\bL_{\ft',\fs}^{\stab,i},\rd\bL_{\ft',\fs}^{\stab,i}\right]
=
(\jmath_{\bL^i})_*
\left[\bL_{\ft',\fs}^{\stab,i},\rd\bL_{\ft',\fs}^{\stab,i}\right]
$$
is the image of the relative fundamental class of
$\bL^{\stab,i}_{\ft',\fs}$ under the inclusion $\jmath_{\bL^i}$.  By
\eqref{eq:SmallEnoughdeltaImpliesEmpty}, the intersection on the
left-hand-side of \eqref{eq:Intersection1} is a finite collection of points
in the interior of $\bL^{\stab,i}_{\ft',\fs}$ (with the multiplicities of
the geometric representatives).  The obstruction section $\bchi$ vanishes
transversely on $\bL^{\stab,i}_{\ft',\fs}$ by Theorem \ref{thm:GluingThm}
and the assumption that the parameters $\eps$ and $\delta$ in the
definition of the link are generic, so the relative fundamental class of
the manifold-with-boundary $\bL_{\ft',\fs}-\bL^{\stab,s}_{\ft',\fs}$ is
given by
\begin{equation}
\label{eq:RelHomologyOfLinkInLi}
e\left(\Upsilon_{\ft',\fs}/S^1,
    \left.\bchi\right|_{\bL^{\stab,s}_{\ft',\fs}}\right)
    \cap
(\jmath_{\bL^i})_*
\left[\bL_{\ft',\fs}^{\stab,i},\rd\bL_{\ft',\fs}^{\stab,i}\right].
\end{equation}
The geometric representatives intersect
$\bL_{\ft',\fs}-\bL^{\stab,s}_{\ft',\fs}$ transversely so, by the
definition of a geometric representative and the cocycle $c(z,\eta)$, the
intersection number on the left-hand-side of
\eqref{eq:Intersection1} is given by evaluating $[c(z,\eta)]$
on the relative fundamental class \eqref{eq:RelHomologyOfLinkInLi},
yielding \eqref{eq:Intersection1}.
\end{proof}

\subsubsection{Extending the cocycles}
\label{subsubsec:ExtendingCocycles}
To rewrite the pairing in identity \eqref{eq:Intersection1}
as a pairing of absolute cohomology classes via the
isomorphism $\jmath_*$ in \eqref{eq:JCohomology},
we must first eliminate the excision isomorphism $\iota_*$.
For this reason, we now show how the cocycles
$c_{\beta}$ and $c_{\sW}$ are, modulo exact cocycles,
the restrictions of cocycles on $\bar\bL^{\stab}_{\ft',\fs}$.
Because the cocycle $c_{\beta}$ on $\bL^{\stab}_{\ft',\fs}$
represents the cohomology class 
$\mu_p(\beta)$, which extends to $\barmu_p(\beta)$ on
$\bar\bL^{\stab}_{\ft',\fs}$, by Definition \ref{defn:ExtendedCohomClasses},
we can always find cocycles in the class $\mu_p(\beta)$ which extend
to $\bar\bL^{\stab}_{\ft',\fs}$ and represent
$\barmu_p(\beta)$. However, we emphasize that the
changes of the cocycles $c_\beta$ and $c_{\sW}$ which we construct
here have the property that they do not alter the pairing in
\eqref{eq:Intersection1},  
which is false for arbitrary choices of cocycles in the cohomology
classes $\mu_p(\beta)$ or $\mu_c$.

For any subset $U\subset X$, we define
\begin{equation}
\label{eq:DefineIN}
\sI(U)=
\bar\bL^{\stab}_{\ft',\fs}\cap
\bga^{-1}\left(\{[A,\Phi,\bx]: |\bx|\cap U\neq \emptyset\}\right),
\end{equation}
where $|\bx|\subset X$ denotes the support of $\bx\in\Sym^\ell(X)$; since
$\ell=1$ here, the set on the right-hand side above is simply 
$\{[A,\Phi,x]: x\in U\}$.
Because the cocycles $c_{\beta}$ and $c_{\sW}$ are pulled
back from $\sB^*(\nu(\beta))$ and $\sC^{*,0}_{\ft'}(\nu(x))/S^1$,
respectively, by the composition of the gluing map and the restriction
map, the cocycles $c_{\beta}$ and $c_{\sW}$ naturally extend over
$\bar\bL^{\stab}_{\ft',\fs}-\sI(\nu(\beta))$ and
$\bar\bL^{\stab}_{\ft',\fs}-\sI(\nu(x))$,
respectively.  (Note that we use the continuity of the gluing map
with respect to Uhlenbeck limits here.)
To construct the extension of
$c_{\beta}$ to $\bar\bL^{\stab}_{\ft',\fs}$, we first show
in Lemma \ref{eq:ClassOfCocycle} that the cohomology class
represented by $c_\beta$ on $\bar\bL^{\stab}_{\ft',\fs}-\sI(\nu(\beta))$
is equal to the restriction of $\barmu_p(\beta)$ to this subspace.
The proof of Lemma \ref{eq:ClassOfCocycle} uses
some cohomological computations appearing in
Lemmas \ref{lem:Cohom1} and \ref{lem:UniqueExtension} and 
the continuity of the gluing map mentioned above.
Then, in Lemma \ref{lem:ExtensionPerturbation} we describe the
perturbation, in the form of an exact cocycle,
necessary to extend $c_\beta$ from
$\bar\bL^{\stab}_{\ft',\fs}-\sI(\nu(\beta))$ to a
cocycle $\bar c_\beta$ on $\bar\bL^{\stab}_{\ft',\fs}$.
The construction of $\bar c_\beta$ will show that it
represents the cohomology class $\barmu_p(\beta)$.

We begin by proving the following lemma about the cohomology of the link.
While the proofs of Lemmas \ref{lem:Cohom1} and \ref{lem:UniqueExtension}
apply only to the link of
the stratum $M_{\fs}\times\Sym^\ell(X)$ when
$\ell=1$, the results should be true for any $\ell\geq 1$.
However, the proofs for larger $\ell$ may be difficult so it may
be easier to prove Lemma \ref{eq:ClassOfCocycle} by a direct, if
tedious, analysis of the bundles discussed there.

\begin{lem}
\label{lem:Cohom1}
If $\pi_X:\bL^{\stab}_{\ft',\fs}\to X$ is the projection map,
then the induced map on cohomology is injective:
$$
\pi_X^*: H^\bullet(X;\RR) \to H^\bullet(\bL^{\stab}_{\ft',\fs};\RR).
$$
\end{lem}

\begin{proof}
The definition of $\bL^{\stab}_{\ft',\fs}$ in \eqref{eq:TopLevels}
and the existence of a retraction from
the uncompactified gluing data space $\Gl_{\ft}(\delta)$
to the boundary $\rd\bar{\Gl}_{\ft}(\delta)\subset \Gl_{\ft}(\delta)$
imply that $\bL^{\stab}_{\ft',\fs}$ deformation-retracts
to $\bL^{\stab,i}_{\ft',\fs}$. In turn, $\bL^{\stab,i}_{\ft',\fs}$
deformation-retracts to the space $\bB\bL^{\stab,i}_{\ft',\fs}$
defined in \eqref{eq:ProductDecompositionOfBase}.
These retractions commute with $\pi_X$, so it suffices to
prove the lemma for the restriction of the projection $\pi_X$
to  $\bB\bL^{\stab,i}_{\ft',\fs}$.  This restriction of $\pi_X$
can be written as the composition of the projections
$$
\pi_1:M_{\fs}\times \rd\bar{\Gl}_{\ft}(\delta)/S^1 \to M_{\fs}\times X
\quad\text{and}\quad
M_{\fs}\times X\to X.
$$
The projection $M_{\fs}\times X\to X$ 
induces an injective homomorphism on cohomology, so we need
only verify that $\pi_1$ induces an injective homomorphism on cohomology.
The map $\pi_1$ is  the product of the identity on $M_{\fs}$ and
the projection
$$
\rd\bar{\Gl}_{\ft}(\delta)/S^1\to X.
$$
As will be discussed the proof of Lemma \ref{lem:PairingsWithInstantonLink},
the space $\rd\bar{\Gl}_{\ft}(\delta)/S^1$ can be identified with
the projectivization of a complex rank-two vector bundle $F\to X$.
By \cite[Equation (20.7)]{BT}, the projection map $\PP(F)\to X$ induces
an injective homomorphism on cohomology. This proves the lemma.
\end{proof}

The following lemma will be used to identify the cohomology
classes which the extensions of the cocycles represent;
the result should be true for the links of ideal Seiberg-Witten moduli spaces
$M_{\fs}\times \Sym^\ell(X)$ defined in \cite{FL5} for $\ell\geq 1$.
The difficulties in extending the proof of
Lemma \ref{lem:UniqueExtension} to the case $\ell>1$ lie
first in the use of Lemma \ref{lem:Cohom1} (which is only proven for $\ell=1$)
and secondly in applying the Thom isomorphism theorem
(as is done below, before \eqref{eq:RelCohom}) to the
analogue of $\nu(\beta)$ in $\Sym^\ell(X)$.

\begin{lem}
\label{lem:UniqueExtension}
For $\beta\in H_\bullet(X;\RR)$, let $T_\beta$ be a smooth submanifold of
$X$ with $[T_\beta]=\beta$ and let $\nu(\beta)$ be a tubular neighborhood
of $T_\beta$ in $X$.  Then for $p=\deg(\beta)$, the following direct sum of
restriction maps is injective:
\begin{equation}
\label{eq:MV1}
H^p(\bar\bL^{\stab}_{\ft',\fs}-\sI(\nu(\beta));\RR)
\to
H^p(\bL^{\stab}_{\ft',\fs};\RR)
\oplus
H^p(\bar\bL^{\stab}_{\ft',\fs}-\pi_X^{-1}(\nu(\beta));\RR).
\end{equation}
\end{lem}

\begin{proof}
To prove the homomorphism \eqref{eq:MV1} is injective, we first observe
that it appears in the Meyer-Vietoris sequence for the
open cover
$$
\bar\bL^{\stab}_{\ft',\fs}-\sI(\nu(\beta))
=
\bL^{\stab}_{\ft',\fs}
\cup
\left(
\bar\bL^{\stab}_{\ft',\fs}-\pi_X^{-1}(\nu(\beta))
\right).
$$
The intersection of these two open sets is
$\bL^{\stab}_{\ft',\fs}-\pi_X^{-1}(\nu(\beta))$.
Hence, the homomorphism
\eqref{eq:MV1} will be injective if the restriction map
\begin{equation}
\label{eq:RestrictionMap}
H^{p-1}(\bL^{\stab}_{\ft',\fs};\RR)
\to
H^{p-1}(\bL^{\stab}_{\ft',\fs}-\pi_X^{-1}(\nu(\beta));\RR)
\end{equation}
is surjective.  As noted in the proof of Lemma \ref{lem:Cohom1},
there is a deformation retraction from $\bL^{\stab}_{\ft',\fs}$
to $\bB\bL^{\stab,i}_{\ft',\fs}$ which commutes with $\pi_X$.
Thus, proving the map \eqref{eq:RestrictionMap} is surjective
is equivalent to proving that the following restriction map
is surjective:
\begin{equation}
\label{eq:RestrictionMap2}
H^{p-1}(\bB\bL^{\stab,i}_{\ft',\fs};\RR)
\to
H^{p-1}(\bB\bL^{\stab,i}_{\ft',\fs}-\pi_X^{-1}(\nu(\beta));\RR),
\end{equation}
The surjectivity of the map
\eqref{eq:RestrictionMap2} follows from a discussion
of the following diagram in which the vertical maps
come from the long exact sequences of the pairs:
$$
\begin{CD}
H^{p-1}(X-\nu(\beta);\RR) @> \pi_X^* >>
H^{p-1}(\bB\bL^{\stab,i}_{\ft',\fs}-\pi_X^{-1}(\nu(\beta));\RR)
\\
@V \delta^*_{X,\beta} VV
@V \delta^*_{\bL,\beta} VV
\\
H^p(X,X-\nu(\beta);\RR) @> \pi_X^* >>
H^p(\bB\bL^{\stab,i}_{\ft',\fs},
\bB\bL^{\stab,i}_{\ft',\fs}-\pi_X^{-1}(\nu(\beta));\RR)
\\
@V \jmath_{X,\beta}^* VV @V \jmath_{\bL,\beta}^* VV
\\
H^p(X;\RR) @> \pi_X^* >>
H^p(\bB\bL^{\stab,i}_{\ft',\fs};\RR)
\end{CD}
$$
Because $\pi_X^{-1}(T_\beta)$ is a smooth, codimension-$p$ submanifold of
$\bB\bL^{\stab,i}_{\ft',\fs}$, with the Thom class of its normal bundle
given by the pullback of the Thom class $\pi_X^*\Th(\beta)$ of the normal
bundle of $T_\beta$ in $X$, we see that the relative cohomology
\begin{equation}
\label{eq:RelCohom}
H^p(\bB\bL^{\stab,i}_{\ft',\fs},
\bB\bL^{\stab,i}_{\ft',\fs}-\pi_X^{-1}(\nu(\beta));
\RR),
\end{equation}
is generated by $\pi_X^*\Th(\beta)$.  By Lemma \ref{lem:Cohom1}, the
following class is non-zero:
$$
\jmath_{\bL,\beta}^*\pi_X^*\Th(\beta)
=
\pi_X^*\jmath_{X,\beta}^*\Th(\beta)
=
\pi_X^*\PD[\beta].
$$ 
Therefore, the generator $\pi_X^*\Th(\beta)$ of the relative
cohomology \eqref{eq:RelCohom} is not in the image of
$\delta_{\bL,\beta}^*$, implying that $\delta_{\bL,\beta}^*$ is the zero map.
Hence, the map \eqref{eq:RestrictionMap2} is surjective.  As discussed
previously, this implies that the map \eqref{eq:MV1}
is injective.
\end{proof}

Lemma \ref{lem:UniqueExtension} gives the following
identity between the cohomology class represented by
the cocycle $c_\beta$ and the
restriction of the cohomology class
$\barmu_p(\beta)$:

\begin{lem}
\label{eq:ClassOfCocycle}
For $\beta\in H_\bullet(X;\RR)$, let $c_{\beta}$ be the cocycle
on $\bar\bL^{\stab}_{\ft',\fs}-\sI(\nu(\beta))$ 
defined at the beginning of \S \ref{subsubsec:ExtendingCocycles}
and consider 
the inclusion map
\begin{equation}
\label{eq:InclusionBeta}
\iota_\beta:\bar\bL^{\stab}_{\ft',\fs}-\sI(\nu(\beta))
\to \bar\bL^{\stab}_{\ft',\fs}.
\end{equation}
Let
$[c_\beta]\in H^{\deg(\beta)}(\bar\bL^{\stab}_{\ft',\fs}-\sI(\nu(\beta));\RR)$
denote the cohomology class which $c_\beta$ represents and 
let $\barmu_p(\beta)\in H^{\deg(\beta)}(\bar\bL^{\stab}_{\ft',\fs};\RR)$
be the restriction
from $\bar\sM^{\stab,*}_{\ft',\fs}/S^1$ to $\bar\bL^{\stab}_{\ft',\fs}$
of the class in Definition \ref{defn:ExtendedCohomClasses}.
Then,
\begin{equation}
\label{eq:ExtensionCohomologyClass1}
[c_\beta]=\iota_\beta^*\barmu_p(\beta).
\end{equation}
\end{lem}

\begin{proof}
By the definition of $c_\beta$ in \S \ref{subsubsec:DefiningCocycles},
the restrictions of the cohomology classes in
\eqref{eq:ExtensionCohomologyClass1} to $\bL^{\stab}_{\ft',\fs}$
coincide.  The identity \eqref{eq:ExtensionCohomologyClass1} will then
follow from Lemma \ref{lem:UniqueExtension}
if we can prove that the difference of
the two cohomology classes is in the kernel of
the restriction map given by the second component of the
map \eqref{eq:MV1}:
$$
\iota_1^*:H^{\deg(\beta)}(\bar\bL^{\stab}_{\ft',\fs}-\sI(\nu(\beta));\RR)
\to
H^{\deg(\beta)}(\bar\bL^{\stab}_{\ft',\fs}-\pi_X^{-1}(\nu(\beta));\RR).
$$
First, we describe the image of $\barmu_p(\beta)$ under this map.
Because there is an inclusion
$$
\left(\bar\bL^{\stab}_{\ft',\fs}-\pi_X^{-1}(\nu(\beta))\right)\times T_\beta
\to
\left(\bar\bL^{\stab}_{\ft',\fs}\times X\right)
-
(\pi_X\times\id_X)^{-1}(\Delta),
$$
the restriction of $(\pi_X\times\id_X)^*\PD[\Delta]/\beta$ to
$\bar\bL^{\stab}_{\ft',\fs}-\pi_X^{-1}(\nu(\beta))$ vanishes.
The expression for $\barmu_p(\beta)$ in Definition
\ref{defn:ExtendedCohomClasses} then implies that 
\begin{equation}
\label{eq:Restriction1}
\iota_1^*\iota_\beta^*\barmu_p(\beta)
=
-\tquarter \iota_1^*\iota_\beta^*(p_1(\FF^{\stab}_{\ft',\fs})/\beta).
\end{equation}
We now consider the restriction of $c_\beta$.  Because the geometric
representative $\sV(\beta)$ is pulled back from the quotient space
$\sB^*(\nu(\beta))$ of irreducible $\SO(3)$ connections over $\nu(\beta)$ 
by the map
$$
r_{\nu(\beta)}\circ\bga:
\bar\bL^{\stab}_{\ft',\fs}-\sI(\nu(\beta))
\to
\sB^*(\nu(\beta)),
$$
the cocycle $c_\beta$ is pulled back from a cocycle $c_{T,\beta}$ on
$\sB^*(\nu(\beta))$.  By the construction \cite[pp. 588-592]{KMStructure}
of the geometric representatives, the cocycle $c_{T,\beta}$ represents the
cohomology class $-\tquarter p_1(\FF_\beta)/\beta$ defined by the universal
$\SO(3)$ bundle:
$$
\FF_\beta
:=
\sA^*(\nu(\beta))\times_{\sG} \fg_{V'}|_{T_\beta}
\to
\sB^*(\nu(\beta))\times T_\beta
$$
We conclude that 
\begin{equation}
\label{eq:Restriction2}
[\iota_1^*c_\beta]
=
-\tquarter \iota_1^*(r_{\nu(\beta)}\circ\bga)^*(p_1(\FF_\beta)/\beta).
\end{equation}
To compare the cohomology classes \eqref{eq:Restriction1} and
\eqref{eq:Restriction2}, we compare the bundles $\FF_\beta$
and $\bar\FF^{\stab,*}_{\ft',\fs}$.
If $\iota_T:T_\beta\to X$ is the inclusion map, then there is an
isomorphism
\begin{equation}
\label{eq:PullbackUnivBundleIsom}
((\iota_\beta\circ \iota_1)\times \iota_T)^*\bar\FF^{\stab,*}_{\ft',\fs}
\cong
((r_{\nu(\beta)}\circ\bga\circ\circ \iota_1)\times\id_{T_\beta})^*\FF_\beta.
\end{equation}
(The existence of this isomorphism follows from the method 
used to obtain the
corresponding isomorphism in Lemma \ref{lem:UnivOnComplement}.)
The isomorphism \eqref{eq:PullbackUnivBundleIsom} then implies that
$$
\iota_1^*(r_{\nu(\beta)}\circ\bga)^* (p_1(\FF_\beta)/\beta)
=
\iota_1^*\iota_\beta^*(p_1(\bar\FF^{\stab,*}_{\ft',\fs})/\beta)
$$
which, together with the identities 
\eqref{eq:Restriction1} and \eqref{eq:Restriction2}, yields 
the desired result, \eqref{eq:ExtensionCohomologyClass1}.
\end{proof}

A similar argument gives the

\begin{lem}
\label{lem:ClassOfW}
Let $c_{\sW}$ be the cocycle defined  
on $\bL^{\stab}_{\ft',\fs}-\sI(\nu(x))$ 
which represents $\mu_c$ when restricted to $\bL^{\stab}_{\ft',\fs}$.
Consider the 
inclusion map
\begin{equation}
\label{eq:InclusionW}
\iota_{\sW}: \bar\bL^{\stab}_{\ft',\fs}-\sI(\nu(x))
\to\bar\bL^{\stab}_{\ft',\fs}.
\end{equation}
If $\barmu_c\in H^{2}(\bar\bL^{\stab}_{\ft',\fs};\RR)$
denotes the restriction 
from $\bar\sM^{\stab,*}_{\ft',\fs}$ to $\bar\bL^{\stab}_{\ft',\fs}$ 
of the class in Definition \ref{defn:ExtendedCohomClasses}, then
\begin{equation}
[c_{\sW}]=\iota_{\sW}^*\barmu_c.
\end{equation}
\end{lem}

We now construct the modification of the cocycle $c_\beta$ on
$\bL^{\stab}_{\ft',\fs}$ necessary to extend it to a cocycle on
$\bar\bL^{\stab}_{\ft',\fs}$.  For any space $A$, let $S^p(A;\RR)$ be the
module of real, singular $p$-cochains on $A$.  For any pair $(A,B)$, let
$S^p(A,B;\RR)$ be the module of real, singular $p$-cochains on $A$,
vanishing on $p$-chains in $B$.

\begin{lem}
\label{lem:ExtensionPerturbation}
For $\beta\in H_\bullet(X;\RR)$, let
$\sU(\beta)\subset\bar\bL^{\stab}_{\ft',\fs}$ be any open neighborhood of
$\sI(\nu(\beta)$ satisfying $\sI(\nu(\beta))\Subset \sU(\beta)$.  Let
$\iota_\beta$ be the inclusion map \eqref{eq:InclusionBeta}.  Then, there
is a cochain
\begin{equation}
\label{eq:DefineCompactPerturbationForMuP}
\theta_\beta
\in
S^{\deg(\beta)-1}(\bar\bL^{\stab}_{\ft',\fs}-\sI(\nu(\beta)),
\bar\bL^{\stab}_{\ft',\fs}-\sU(\beta);\RR),
\end{equation}
and a cocycle $\bar c_\beta$ representing the class
$\barmu_p(\beta)\in H^\bullet(\bar\bL^{\stab}_{\ft',\fs};\RR)$ 
in Definition \ref{defn:ExtendedCohomClasses}
such that
$$
\iota_\beta^*\bar c_\beta=c_\beta+\delta\theta_\beta.
$$
\end{lem}

\begin{proof}
Let $\bar c_\beta'$ be any cocycle in the cohomology class
$\barmu_p(\beta)$.
By Lemma \ref{eq:ClassOfCocycle}, there is a cochain
$\theta_0$ of degree $\deg(\beta)-1$
on $\bL^{\stab}_{\ft',\fs}-\sI(\nu(\beta))$ such that
\begin{equation}
\label{eq:BoundaryPerturb1}
\iota_\beta^*\bar c_\beta'=c_{\beta}+\delta\theta_0.
\end{equation}
Because $\sI(\nu(\beta))\Subset \sU(\beta)$, there is an open
subspace $\sU'(\beta)\subset\bar\bL^{\stab}_{\ft',\fs}$ with
$$
\sI(\nu(\beta))\subset \sU'(\beta)\Subset \sU(\beta).
$$
Hence, the following intersection is empty:
$$
\left(
    \bar\bL^{\stab}_{\ft',\fs}-\sU(\beta)
\right) \cap \sU'(\beta)
=
\emptyset.
$$
This implies that the map 
$$
\begin{CD}
S^p
\left(
    \bar\bL^{\stab}_{\ft',\fs}-\sI(\nu(\beta)),
    \bar\bL^{\stab}_{\ft',\fs}-\sU(\beta);
    \RR
\right)
\oplus
S^p\left(
\bar\bL^{\stab}_{\ft',\fs}-\sI(\nu(\beta)),\sU'(\beta);\RR
\right)
\\
@VVV
\\
S^p(\bar\bL^{\stab}_{\ft',\fs}-\sI(\nu(\beta)),
\left(
    \bar\bL^{\stab}_{\ft',\fs}-\sU(\beta)
\right) \cap \sU'(\beta);\RR)
\end{CD}
$$
which is surjective
because the pairs appearing in this diagram are excisive couples; see
\cite[Theorem 4.6.3]{Spanier} and \cite[p. 218]{Spanier},
is actually a map to the space of absolute cochains,
$$
S^p(\bar\bL^{\stab}_{\ft',\fs}-\sI(\nu(\beta)),\emptyset;\RR).
$$
Thus, we can
write $\theta_0=\theta_\beta+\theta_p$, where
\begin{gather*}
\theta_\beta\in S^{\deg(\beta)-1}
\left(
    \bar\bL^{\stab}_{\ft',\fs}-\sI(\nu(\beta)),
    \bar\bL^{\stab}_{\ft',\fs}-\sU(\beta);\RR
\right),
\\
\theta_p\in  S^{\deg(\beta)-1}
\left(
    \bar\bL^{\stab}_{\ft',\fs}-\sI(\nu(\beta)),\sU'(\beta);\RR
\right).
\end{gather*}
Because $\theta_p$ is supported on the complement of
$\sU'(\beta)$ in $\bar\bL^{\stab}_{\ft',\fs}-\sI(\nu(\beta))$, we see that
$\theta_p$ defines a cochain of degree $\deg(\beta)-1$ on
$\bar\bL^{\stab}_{\ft',\fs}$ (by extending $\theta_p$ by zero)
and if we set
$$
\bar c_\beta
=
\bar c_\beta'-\delta\theta_p,
$$
then $\bar c_\beta$ also represents the cohomology class
$\barmu_p(\beta)$ and equation \eqref{eq:BoundaryPerturb1} yields
\begin{align*}
\iota_\beta^*\bar c_\beta
&=
\iota_\beta^*\bar c_\beta'-\delta\theta_p
\\
&=
c_\beta + \delta(\theta_\beta+\theta_p) -\delta\theta_p
\\
&=
c_\beta + \delta\theta_\beta.
\end{align*}
This completes the proof.
\end{proof}

A similar argument yields:

\begin{lem}
\label{lem:ExtensionPerturbationW}
Let $\sU(x)\subset\bar\bL^{\stab}_{\ft',\fs}$ be any open neighborhood of
$\sI(\nu(x)$ satisfying $\sI(\nu(x))\Subset \sU(x)$.  Let $\iota_{\sW}$ be
the inclusion map \eqref{eq:InclusionW}.  Then, there is a cochain
\begin{equation}
\label{eq:CompactPerturbationForW}
\theta_{\sW}
\in
S^1(\bar\bL^{\stab}_{\ft',\fs}-\sI(\nu(x)),
\bar\bL^{\stab}_{\ft',\fs}-\sU(x);\RR)
\end{equation}
and a cocycle $\bar c_{\sW}$ in the cohomology class
$\barmu_c\in H^2(\bar\bL^{\stab}_{\ft',\fs};\ZZ)$ 
in Definition \ref{defn:ExtendedCohomClasses} such that
$$
\iota_{\sW}^*\bar c_{\sW}=c_{\sW}+\delta\theta_{\sW}.
$$
\end{lem}

The construction of $\bar c_\beta$ and $\bar c_{\sW}$ in Lemmas
\ref{lem:ExtensionPerturbation} and \ref{lem:ExtensionPerturbationW} yields
the following useful corollary:

\begin{cor}
\label{cor:CohomClassOfExtensions}
For $\beta\in H_{\bullet}(X;\RR)$, let $\bar c_{\beta}$
and $\bar c_{\sW}$ be the cocycles defined in
Lemmas \ref{lem:ExtensionPerturbation} and
\ref{lem:ExtensionPerturbationW} respectively.
Then, as elements of $H^\bullet(\bar\bL^{\stab}_{\ft',\fs};\RR)$,
\begin{equation}
[\bar c_{\beta}] = \barmu_p(\beta)
\quad\text{and}\quad
[\bar c_{\sW}]= \barmu_c.
\end{equation}
\end{cor}

\subsubsection{Eliminating the excision isomorphism}
\label{subsubsec:ElimExcis}
The next step in the proof of Proposition
\ref{prop:IntersectionNoToCupProduct} is to eliminate the
excision isomorphism $\iota_*$ from  the pairing in \eqref{eq:Intersection1}.
We do this by presenting the relative cohomology classes
as restrictions (that is, the image under $\iota^*$)
of relative cohomology classes on $\bar\bL^{\stab}_{\ft',\fs}$
and using the adjoint relation between $\iota_*$ and $\iota^*$. 

The relative Euler class is not the restriction of a relative
cohomology class on $\bar\bL^{\stab}_{\ft',\fs}$ because the obstruction
bundle $\Upsilon_{\ft',\fs}/S^1$
does not extend over $\bL_{\ft',\fs}^{\sing}$.
We use the decomposition $\Upsilon_{\ft',\fs}
=\Upsilon^s_{\ft',\fs}\oplus\Upsilon^i_{\ft',\fs}$
of the obstruction bundle, with
background component $\Upsilon^s_{\ft',\fs}/S^1$ which
extends over $\bL_{\ft',\fs}^{\sing}$ (yielding 
$\bar\Upsilon^s_{\ft',\fs}/S^1$ in \eqref{eq:ExtendedBackgroundObstruction}), 
and an instanton component (which does not extend),
$\Upsilon^i_{\ft',\fs}/S^1$, to overcome this difficulty.
We observe in Lemma \ref{lem:Intersection2}
that we can take advantage of this decomposition of
$\Upsilon_{\ft',\fs}/S^1$ to perturb the section $\bchi$ 
without changing the intersection number in \eqref{eq:Intersection1}.
Recall from \eqref{eq:DefineRelFundClass} and \eqref{eq:Excision}
that the relative homology class in the pairing \eqref{eq:Intersection2}
below is contained in 
$H_{\max}(\bL^{\stab}_{\ft',\fs},\bL^{\stab,s}_{\ft',\fs};\ZZ)$.

\begin{lem}
\label{lem:Intersection2}
Let $\bchi=\bchi_s\oplus\bchi_i$ be the obstruction section
provided by Theorem \ref{thm:GluingThm}.  Then the cup product
in \eqref{eq:Intersection2} is a cohomology class in
$H^{\max}(\bL^{\stab}_{\ft',\fs},\bL^{\stab,s}_{\ft',\fs};\ZZ)$ and so
the following pairing is well-defined:
\begin{equation}
\label{eq:Intersection2}
\left\langle
[c(z,\eta)]\smile e\left(\Upsilon_{\ft',\fs}/S^1,
\bchi_s|_{\bL^{\stab,s}_{\ft',\fs}}\right),
\iota_*^{-1}
(\bar\jmath_{\bL^i})_*
\left[\bL_{\ft',\fs}^{\stab,i},\rd\bL_{\ft',\fs}^{\stab,i}\right]
\right\rangle.
\end{equation}
Moreover, we can replace the section $\bchi_s$ in the definition of the
relative Euler class in the pairing \eqref{eq:Intersection2} with the
section $\bchi$, so the pairing \eqref{eq:Intersection2} is
equal to the pairing on the right-hand side of equation 
\eqref{eq:Intersection1}:
\begin{equation}
\label{eq:Intersection2-Identity}
\begin{aligned}
{}&\left\langle
[c(z,\eta)]\smile e\left(\Upsilon_{\ft',\fs}/S^1,
\bchi|_{\bL^{\stab,s}_{\ft',\fs}}\right),
\iota_*^{-1}
(\bar\jmath_{\bL^i})_*
\left[\bL_{\ft',\fs}^{\stab,i},\rd\bL_{\ft',\fs}^{\stab,i}\right]
\right\rangle
\\
&=
\left\langle
[c(z,\eta)]\smile e\left(\Upsilon_{\ft',\fs}/S^1,
    \bchi_s|_{\bL^{\stab,s}_{\ft',\fs}}\right),
\iota_*^{-1}
(\bar\jmath_{\bL^i})_*
\left[\bL_{\ft',\fs}^{\stab,i},\rd\bL_{\ft',\fs}^{\stab,i}\right]
\right\rangle.
\end{aligned}
\end{equation}
\end{lem}

\begin{proof}
Because $\bchi^{-1}(0)\subset \bchi_s^{-1}(0)$, there is an inclusion
map of pairs:
$$
\jmath_s:
\left(
    \bL^{\stab}_{\ft',\fs},
    \bL^{\stab,s}_{\ft',\fs}-\bchi_s^{-1}(0)
\right)
\to
\left(
    \bL^{\stab}_{\ft',\fs},
    \bL^{\stab,s}_{\ft',\fs}-\bchi^{-1}(0)
\right).
$$
The sections $\bchi$ and $\bchi_s$ are homotopic through non-vanishing
sections on $\bL^{\stab,s}_{\ft',\fs}-\bchi_s^{-1}(0)$
(by the homotopy $\bchi_s\oplus t\bchi_i$, for $t\in [0,1]$)
and so the relative Euler classes are equal,
\begin{equation}
\label{eq:EqEulerClasses}
\jmath_s^*e\left(\Upsilon_{\ft',\fs}/S^1,\bchi|_{\bL^{\stab,s}_{\ft',\fs}}\right)
=
e\left(\Upsilon_{\ft',\fs}/S^1,\bchi_s|_{\bL^{\stab,s}_{\ft',\fs}}\right)
\end{equation}
as elements of
$$
H^{2r_\Xi+2}
\left(
    \bL^{\stab}_{\ft',\fs},
    \bL^{\stab,s}_{\ft',\fs}-\bchi_s^{-1}(0);\ZZ
\right).
$$
For $\delta$ sufficiently small, the intersection $\bar\sK(z,\eta)\cap
\bar\bL^{\stab,s}_{\ft',\fs}\cap \bchi_s^{-1}(0)$ is empty by Lemma
\ref{lem:GRCompactSupport} and the continuity of the gluing map and
obstruction section $\bchi_s$ with respect to Uhlenbeck limits.  Because
$[c(z,\eta)]$ is supported on $\sK(z,\eta)$ while the relative Euler class
in \eqref{eq:EqEulerClasses} is supported on the complement of
$\bchi_s^{-1}(0)$, the cup-product in \eqref{eq:Intersection2} thus defines
an element of
$$
H^{\max}\left(\bL^{\stab}_{\ft',\fs},\bL^{\stab,s}_{\ft',\fs};\RR\right),
$$
(compare the argument giving \eqref{eq:ProductCompactSupport}),
so the pairing \eqref{eq:Intersection2} is well-defined.
We can then write
\begin{align*}
{}&
\left\langle
[c(z,\eta)]\smile e\left(\Upsilon_{\ft',\fs}/S^1,
    \bchi|_{\bL^{\stab,s}_{\ft',\fs}}\right),
\iota_*^{-1}
(\bar\jmath_{\bL^i})_*
\left[\bL_{\ft',\fs}^{\stab,i},\rd\bL_{\ft',\fs}^{\stab,i}\right]
\right\rangle
\\
&=
\left\langle
[c(z,\eta)]\smile \jmath_s^*e\left(\Upsilon_{\ft',\fs}/S^1,
\bchi|_{\bL^{\stab,s}_{\ft',\fs}}\right),
\iota_*^{-1}
(\bar\jmath_{\bL^i})_*
\left[\bL_{\ft',\fs}^{\stab,i},\rd\bL_{\ft',\fs}^{\stab,i}\right]
\right\rangle
\\
&=
\left\langle
[c(z,\eta)]\smile e\left(\Upsilon_{\ft',\fs}/S^1,
    \bchi_s|_{\bL^{\stab,s}_{\ft',\fs}}\right),
\iota_*^{-1}
(\bar\jmath_{\bL^i})_*
\left[\bL_{\ft',\fs}^{\stab,i},\rd\bL_{\ft',\fs}^{\stab,i}\right]
\right\rangle,
\end{align*}
using equation \eqref{eq:EqEulerClasses} in the last line. 
This proves the desired identity \eqref{eq:Intersection2-Identity}.
\end{proof}

The advantage of the relative Euler class
$e(\Upsilon_{\ft',\fs}/S^1,\bchi_s)$ over $e(\Upsilon_{\ft',\fs}/S^1,\bchi)$
is that it can be written as a cup-product,
\begin{equation}
\label{eq:DecomposeRelEulerClass}
e\left(\Upsilon_{\ft',\fs}/S^1,\bchi_s|_{\bL^{\stab,s}_{\ft',\fs}}\right)
=
\iota^*e\left(\bar\Upsilon^s_{\ft',\fs}/S^1,
\bchi_s|_{\bar\bL^{\stab,s}_{\ft',\fs}}\right)
\smile
e\left(\Upsilon^i_{\ft',\fs}/S^1\right),
\end{equation}
where $\iota$ is the inclusion map \eqref{eq:DefineInclusion}.  Note that
the continuity of the obstruction map $\bchi_s$ with respect to Uhlenbeck
limits is required to define the first relative Euler class appearing on
the right-hand-side of \eqref{eq:DecomposeRelEulerClass} as a continuous
section is needed to define a relative Euler class.  The relative
cohomology class and the absolute cohomology class on the right-hand-side
of equation \eqref{eq:DecomposeRelEulerClass} are both restrictions of
classes on $\bar\bL^{\stab}_{\ft',\fs}$.  We now show that we can also
replace the relative cohomology class $[c(z,\eta)]$ with one which extends
over $\bar\bL^{\stab}_{\ft',\fs}$ without changing the pairing
\eqref{eq:Intersection2}.

Let $z=\beta_1\cdots\beta_m\in\AAA(X)$ be the monomial appearing in
Proposition \ref{prop:IntersectionNoToCupProduct}.  We replace
$c_{\sW}^{\eta}$ with
$$
c_{\sW,1}\smile\dots\smile c_{\sW,\eta},
$$
as it will be necessary to distinguish between the factors of this
cup-product, which are defined via specific choices described below.  Let
$\nu'(\beta_i)$ and $\nu'(x_j)$, $j=1,\dots,\eta$, be tubular neighborhoods
of $T_{\beta_i}$ and the point $x_j$ defining $c_{\sW,j}$, respectively.
Assume these neighborhoods satisfy both $\nu(\beta_i)\Subset\nu'(\beta_i)$
and $\nu(x_j)\Subset\nu'(x_j)$ and the intersection condition
\eqref{eq:RestrictionCondition}.  Therefore, we can find open subspaces
$\sU(\beta_i)$ and $\sU(x_j)$ in $\bar\bL^{\stab}_{\ft',\fs}$ such that
$$
\sI(\nu(\beta))\Subset \sU(\beta)\quad\text{and}\quad
\sI(\nu(x_j))\Subset \sU(x_j),
$$
and such that the intersection
\begin{equation}
\label{eq:IntersectionSupport}
\sU(\beta_{i_1})\cap\dots\cap \sU(\beta_{i_k})
\cap
\sU(x_{j_1})\cap\dots\cap \sU(x_{j_r}),
\end{equation}
is empty unless
\begin{equation}
\label{eq:IntersectionCondition2}
\sum_{p=1}^k\deg(\beta_{i_p}) + 4r\le 4.
\end{equation}
Then let
\begin{equation}
\label{eq:PerturbedCocycles}
\iota_{\beta_i}^*\bar c_{\beta_i}=c_{\beta_i}+\delta\theta_{\beta_i}
\quad
\text{and}
\quad
\iota_{\sW}^*\bar c_{\sW,j}=c_{\sW,j}+\delta\theta_{\sW,j},
\end{equation}
be the cocycles constructed in Lemmas  \ref{lem:ExtensionPerturbation}
and \ref{lem:ExtensionPerturbationW} respectively,
where $\bar c_{\beta_i}$ and $\bar c_{\sW,j}$ are cocycles on
$\bar\bL^{\stab}_{\ft',\fs}$,
$\delta\theta_{\beta_i}$ is supported in
$\sU(\beta_i)$, and $\delta\theta_{\sW,j}$ is
supported in $\sU(x_j)$.
Define
\begin{equation}
\label{eq:DefineExtendedCocycleProduct}
\bar c(z,\eta)
=
\bar c_{\beta_1}\smile\dots\smile \bar c_{\beta_m}\smile 
\bar c_{\sW,1}\smile\dots \smile\bar c_{\sW,\eta},
\end{equation}
so $\iota^*\bar c(z,\eta)=c(z,\eta)$.
We now show that we can replace $c(z,\eta)$ with
$\iota^*\bar c(z,\eta)$  without changing
the pairing \eqref{eq:Intersection2}.

\begin{lem}
\label{lem:Intersection3}
Let $z=\beta_1\cdots\beta_m\in\AAA(X)$ be the
monomial appearing in 
Proposition \ref{prop:IntersectionNoToCupProduct}.
Let $\iota$ be the inclusion map \eqref{eq:DefineInclusion}.
Then we can replace the class
$[c(z,\eta)]$ in the pairing \eqref{eq:Intersection2}
by the class $\iota^*[\bar c(z,\eta)]$:
\begin{equation}
\label{eq:Intersection3}
\begin{aligned}
{}&
\left\langle
[c(z,\eta)]\smile e\left(\Upsilon_{\ft',\fs}/S^1,
\bchi_s|_{\bL^{\stab,s}_{\ft',\fs}}\right),
\iota_*^{-1}
(\bar\jmath_{\bL^i})_*
\left[\bL_{\ft',\fs}^{\stab,i},\rd\bL_{\ft',\fs}^{\stab,i}\right]
\right\rangle
\\
&=
\left\langle
\iota^*[\bar c(z,\eta)]\smile e\left(\Upsilon_{\ft',\fs}/S^1,
\bchi_s|_{\bL^{\stab,s}_{\ft',\fs}}\right),
\iota_*^{-1}
(\bar\jmath_{\bL^i})_*
\left[\bL_{\ft',\fs}^{\stab,i},\rd\bL_{\ft',\fs}^{\stab,i}\right]
\right\rangle.
\end{aligned}
\end{equation}
\end{lem}

\begin{proof}
By construction and \eqref{eq:PerturbedCocycles},
the cocycles $c(z,\eta)$ and $\iota^*\bar c(z,\eta)$ differ by cocycles
of the form (up to a re-ordering which will be seen to be irrelevant)
\begin{equation}
\label{eq:DifferenceTerm}
c(z',\eta-r)\smile \delta\theta_{i_1}\smile\cdots\smile\delta\theta_{i_k}
\smile
\delta\theta_{j_1}\smile \dots\smile\delta\theta_{j_r},
\end{equation}
where $z=z'\beta_{i_1}\cdots\beta_{i_k}$.
Note that we assume $k>0$ or $r>0$,  so
we actually have a difference term in \eqref{eq:DifferenceTerm}
and not just $\iota^*\bar c(z,\eta)$.
Because the intersection \eqref{eq:IntersectionSupport} is empty
unless the condition \eqref{eq:IntersectionCondition2} holds
and because $\delta\theta_{\beta_i}$ and
$\delta\theta_{\sW,j}$ have support in $\sU(\beta_i)$ and $\sU(x_j)$
respectively,
the term \eqref{eq:DifferenceTerm} vanishes unless condition
\eqref{eq:IntersectionCondition2} holds.
Consequently,
\begin{equation}
\label{eq:IntersectionDim}
\begin{aligned}
\deg(z')+2(\eta-r)
&\ge
\deg(z)+2\eta-\left(\sum_{\mu=1}^k\deg(\beta_i)+4r\right)
\\
&\ge
\dim\sM_{\ft'}- 6
\\
&=\dim \sM_{\ft}
\\
&=
\dim\left(\bchi_{\fs}^{-1}(0)\cap\bL_{\ft',\fs}^{\sing}\right)
+2.
\end{aligned}
\end{equation}
The inequality \eqref{eq:IntersectionDim} implies that the following
intersection is empty:
$$
\bar\sK(z',\eta-r)\cap\bL_{\ft',\fs}^{\sing}\cap \bchi_s^{-1}(0)
=
\emptyset.
$$
Hence, because of the continuity of the gluing map and of the obstruction 
section $\bchi_s$ on $\bar\bL^{\stab,s}_{\ft',\fs}$ and by the same argument
used to establish \eqref{eq:SmallEnoughdeltaImpliesEmpty},
the following intersection is also empty
for $\delta$ sufficiently small:
$$
\sK(z',\eta-r)\cap \bL^{\stab,s}_{\ft',\fs}\cap \bchi_s^{-1}(0)
=
\emptyset.
$$
By the reasoning which gave
\eqref{eq:ProductCompactSupport}, the preceding equality implies that 
the cup-product
$$
[c(z',\eta-r)]\smile e\left(\bar\Upsilon^s_{\ft',\fs}/S^1,
\bchi_s|_{\bL^{\stab,s}_{\ft',\fs}}\right)
$$
vanishes on $\bL^{\stab,s}_{\ft',\fs}$ for $\delta$ sufficiently small.  If
we assume that the neighborhoods $\sU(\beta_i)$ and $\sU(x_j)$ are
contained in $\bL^{\stab,s}_{\ft',\fs}$ for all $i$ and $j$, the assumption
that $k>0$ or $r>0$ implies that the difference term
\eqref{eq:DifferenceTerm} is supported in $\bL^{\stab,s}_{\ft',\fs}$
and hence the cup-product of the difference term with
the relative Euler class vanishes.  Therefore
identity \eqref{eq:Intersection3} holds.
\end{proof}

We now have all the ingredients we need to conclude the proof of
Proposition \ref{prop:IntersectionNoToCupProduct} at our disposal:

\begin{proof}[Proof of Proposition \ref{prop:IntersectionNoToCupProduct}]
Equations \eqref{eq:Intersection1},
\eqref{eq:Intersection2-Identity}, 
\eqref{eq:DecomposeRelEulerClass},
and \eqref{eq:Intersection3}
yield the identity:
\begin{equation}
\label{eq:Intersection5}
\begin{aligned}
{}&\#\left( \bar\sV(z)\cap\bar\sW^{\eta} \cap\bar\bL_{\ft',\fs}\right)
\\
&=
\left\langle
\iota^*\left(
[\bar c(z,\eta)]\smile
    e\left(\bar\Upsilon^s_{\ft',\fs}/S^1,
\bchi_s|_{\bar\bL^{\stab,s}_{\ft',\fs}}\right)
\smile \bar e_i\right),
\iota_*^{-1}
(\bar\jmath_{\bL^i})_*
\left[\bL_{\ft',\fs}^{\stab,i},\rd\bL_{\ft',\fs}^{\stab,i}\right]
\right\rangle
\\
&=
\left\langle
[\bar c(z,\eta)]\smile
    e\left(\bar\Upsilon^s_{\ft',\fs}/S^1,
\bchi_s|_{\bar\bL^{\stab,s}_{\ft',\fs}}\right)
\smile \bar e_i,
(\bar\jmath_{\bL^i})_*
\left[\bL_{\ft',\fs}^{\stab,i},\rd\bL_{\ft',\fs}^{\stab,i}\right]
\right\rangle.
\end{aligned}
\end{equation}
Because $\bar c(z,\eta)$ is a cocycle representing
the cohomology class $\barmu_p(z)\smile\barmu_c^{\eta}$
on $\bar\sM_{\ft',\fs}^{\stab,*}/S^1$
by Corollary \ref{cor:CohomClassOfExtensions}
and because of the relation \eqref{eq:RelativeEulerMap} 
between relative and absolute Euler classes 
and the definition of $\bar e_s$ given by
Lemma \ref{lem:EulerOfBackgroundObstruction}, we see that
\begin{equation}
\label{eq:RelativeToAbsoluteClasses}
\begin{aligned}
{}& \jmath^*\left(
[\bar c(z,\eta)]\smile e\left(\bar\Upsilon^s_{\ft',\fs}/S^1,
    \bchi_s|_{\bar\bL^{\stab,s}_{\ft',\fs}}\right)\smile \bar e_i
\right)
\\
&=
\barmu_p(z)\smile\barmu_c^{\eta}
\smile \bar e_s\smile \bar e_i,
\end{aligned}
\end{equation}
where $\jmath$ is the inclusion map \eqref{eq:DefineJ}.  Combining the
expression \eqref{eq:Intersection5} for the intersection number, equation
\eqref{eq:RelativeToAbsoluteClasses}, and the definition of
$[\bar\bL^{\stab}_{\ft',\fs}]$ in \eqref{eq:DefineFundClassOfL} 
then completes the proof of Proposition
\ref{prop:IntersectionNoToCupProduct}.
\end{proof}

\subsection{Cohomological results}
\label{subsec:PrelimComp}
In this section we prove some results about the
cohomology ring of $\bL^{\stab}_{\ft',\fs}$ which
we shall need to compute the right-hand side of
equation \eqref{eq:CohomFormulation} and thus
complete the proof of Theorem \ref{thm:LevelOne}.
In Lemma \ref{lem:PDualOfBaseSpace}, we compute the Poincar\'e dual of
the submanifold of $\bL^{\stab}_{\ft',\fs}$,
\begin{equation}
\label{eq:BaseSpace}
\begin{aligned}
\bB\bL^{\stab,i}_{\ft',\fs}
&=
\tM_{\fs}\times_{\sG_{\fs}\times S^1}\rd\bar{\Gl}_{\ft}(\delta)
\\
&\subset
\tilde N_{\ft,\fs}(\eps)\times_{\sG_{\fs}\times S^1}\rd\bar{\Gl}_{\ft}(\delta)
=
\bL^{\stab,i}_{\ft',\fs}
\subset \bL^{\stab}_{\ft',\fs},
\end{aligned}
\end{equation}
and in Propositions \ref{prop:S1Localization} and
\ref{prop:LinkSegreFormula} we give a formula for division by this Poincar\'e
dual.
Dividing  the cohomology class in the pairing on the right-side of identity
\eqref{eq:CohomFormulation}
by the Poincar\'e dual of the fundamental class of
$\bB\bL_{\ft',\fs}^{\stab,i}$ reduces the pairing with
$[\bar\bL^{\stab}_{\ft',\fs}]$ to one with
$[\bB\bL_{\ft',\fs}^{\stab,i}]$.  Because, as noted in
equation \eqref{eq:ProductDecompositionOfBase}, the space
$\bB\bL_{\ft',\fs}^{\stab,i}$ is a
product
\begin{equation}
\label{eq:ProductDecomp}
\bB\bL_{\ft',\fs}^{\stab,i}
=
M_{\fs}\times \rd\bar{\Gl}_{\ft}(\delta)/S^1,
\end{equation}
the result of
pairing with the fundamental class $[\bB\bL_{\ft',\fs}^{\stab,i}]$
is equal to a product of pairings with $[M_{\fs}]$ and with
$[\rd\bar{\Gl}_{\ft}(\delta)/S^1]$.
The former gives the Seiberg-Witten
invariant, while the latter has been previously computed in \cite{LenessWC}
(see Lemma \ref{lem:PairingsWithInstantonLink}).

A similar technique would work for the link of $M_{\fs}\times\Sym^\ell(X)$,
for $\ell\geq 2$, although the computations of Lemma
\ref{lem:PairingsWithInstantonLink} for the analogue of the space
$\rd\bar{\Gl}_{\ft}(\delta)/S^1$ when $\ell\geq 3$ appear very difficult,
though the case of $\ell=2$ can be deduced from the work of
\cite{LenessWC}.

Both the computation of the Poincar\'e dual of $\bB\bL_{\ft',\fs}^{\stab,i}$
and our formula for division by this Poincar\'e dual use equivariant
cohomology, which we now briefly review.

As customary, $\ES^1\to \BS^1$ denotes the universal circle bundle, where
$\BS^1$ is homotopy equivalent to $\CC\PP^\infty$. For any pair of
topological spaces $(X,Y)$ with a circle action,
\begin{align*}
H^{\bullet}_{S^1}(X;\ZZ)
&=
H^{\bullet}(\ES^1\times_{S^1}X;\ZZ),
\\
H^{\bullet}_{S^1}(X,Y;\ZZ)
&=
H^{\bullet}(\ES^1\times_{S^1}X,\ES^1\times_{S^1}Y;\ZZ),
\end{align*}
are the equivariant cohomology of $X$ and relative equivariant cohomology
of $(X,Y)$, respectively \cite[\S 2]{AtiyahBott}.

\begin{lem}
\label{lem:PDualOfBaseSpace}
There is a continuous map
\begin{equation}
\label{eq:PDMap}
\pi_{N,\nu}:
\bar{\sM}^{\stab,*}_{\ft',\fs}/S^1
\to
N_{\ft,\fs}\times_{S^1}\ES^1,
\end{equation}
such that if $\Th_{S^1}(N_{\ft,\fs})
\in H^\bullet_{S^1}(N_{\ft,\fs},N^0_{\ft,\fs};\ZZ)$ is the Thom
class \cite[Theorem 9.1]{MilnorStasheff}, \cite[\S 6]{BT} of the bundle
\begin{equation}
\label{eq:EquivN}
N_{\ft,\fs}\times_{S^1}\ES^1
\to
M_{\fs}\times \BS^1,
\end{equation}
where $N^0_{\ft,\fs}$ denotes the complement of the zero-section,
and we have an inclusion map,
\begin{equation}
\label{eq:JMapForBL}
j_{\bB\bL}:
\left( \bar\bL^{\stab}_{\ft',\fs},\emptyset\right)
\to
\left( \bar\bL^{\stab}_{\ft',\fs},\bar\bL^{\stab}_{\ft',\fs}
        - \bB\bL_{\ft',\fs}^{\stab,i}
\right),
\end{equation}
then the following hold:
\begin{enumerate}
\item
For any $\om \in H^{d_s(\fs)+6}(\bar\bL^{\stab}_{\ft',\fs};\ZZ)$, where
$\dim\bB\bL^{\stab,i}_{\ft',\fs}=d_s(\fs)+6$, we have
\begin{equation}
\label{eq:PDOfBasePairing}
\left\langle \om\smile j_{\bB\bL}^*\pi_{N,\nu}^*\Th_{S^1}(N_{\ft,\fs}),
[\bar\bL^{\stab}_{\ft',\fs}]\right\rangle 
=
\left\langle \om,[\bB\bL^{\stab,i}_{\ft',\fs}]\right\rangle.
\end{equation}
\item
If $\bc\in H^2(\BS^1;\ZZ)$ is the universal first Chern class and
$$
\pi_B: N_{\ft,\fs}\times_{S^1}\ES^1\to \BS^1,
$$
is the projection and $\nu\in H^2(\bar\sM^{\stab}_{\ft',\fs}/S^1;\ZZ)$
is the cohomology class of Definition \ref{defn:DefnOfNu}, then
$$
\pi_{N,\nu}^*\pi_B^*\bc=\nu.
$$
\item
If $\mu_{\fs}\in H^2(M_{\fs};\ZZ)$ is as defined in \eqref{eq:SWClass} and
$$
\pi_M:N_{\ft,\fs}\times_{S^1}\ES^1\to M_{\fs},
$$
is the projection, then
$$
\pi_{N,\nu}^*\pi_M^*\mu_{\fs}=\pi_{\fs}^*\mu_{\fs}.
$$
\end{enumerate}
\end{lem}

\begin{proof}
We define a bundle map $\pi_N$ by the obvious projection so
that the following diagram commutes:
\begin{equation}
\label{eq:PDBundleMap1}
\begin{CD}
\tilde N_{\ft,\fs}\times_{\sG_{\fs}}\bar{\Gl}_{\ft}(\delta)
@> \pi_N >>N_{\ft,\fs}
\\
@VVV @VVV
\\
\tM_{\fs}\times_{\sG_{\fs}}\bar{\Gl}_{\ft}(\delta) @> \pi_{\fs} >> M_{\fs}
\end{CD}
\end{equation}
By construction, the projection
$\pi_N$ is $S^1$-equivariant where the circle acts
by scalar multiplication on the fibers of $N_{\ft,\fs}$ and
by the action in Definition \ref{defn:DefnOfNu} defining $\nu$
on the fibers of the bundle on the left-hand-side of the diagram.
Let $\iota_\nu$ be the classifying map for the circle
bundle defining $\nu$, covered by the bundle map $\tilde\iota_\nu$,
so that the following diagram commutes:
\begin{equation}
\label{eq:PDBundleMap2}
\begin{CD}
\bar{\sM}^{\stab,*}_{\ft',\fs} @> \tilde\iota_\nu >> \ES^1
\\
@VVV @VVV
\\
\bar{\sM}^{\stab,*}_{\ft',\fs}/S^1
@> \iota_\nu >> \BS^1
\end{CD}
\end{equation}
By construction of the map $\iota_\nu$, we have $\iota_\nu^*\bc=\nu$.
Then, because the maps $\pi_N$
and $\tilde \iota_\nu$ are
$S^1$-equivariant, the product $\pi_N\times\tilde\iota_\nu$
defines the map $\pi_{N,\nu}$ on the circle quotients:
$$
\pi_{N,\nu}:
\bar{\sM}^{\stab,*}_{\ft',\fs}/S^1
\to N_{\ft,\fs}\times_{S^1}\ES^1.
$$
Observe that the intersection of the pre-image under
$\pi_{N,\nu}$ of the zero-section of the bundle \eqref{eq:EquivN}
with $\bar\bL_{\ft',\fs}^{\stab}$ is the base space \eqref{eq:BaseSpace}:
\begin{equation}
\label{eq:IdentifyPreImageOfZeroLocus}
\bar\bL^{\stab}_{\ft',\fs}
\cap
\pi_{N,\nu}^{-1}(M_{\fs}\times\BS^1)
=
\bB\bL_{\ft',\fs}^{\stab,i}.
\end{equation}
Hence, the map $\pi_{N,\nu}$ induces a
map on relative cohomology,
$$
\pi_{N,\nu}^*:
H^\bullet_{S^1}(N_{\ft,\fs},N^{0}_{\ft,\fs};\ZZ)
\to
H^\bullet(
\bar\bL^{\stab}_{\ft',\fs},
\bar\bL^{\stab}_{\ft',\fs}- \bB\bL_{\ft',\fs}^{\stab,i};\ZZ),
$$
where $N^{0}_{\ft,\fs}$ denotes the complement of the
zero-section.  Because
the restriction of $\pi_{N,\nu}$
to $\bL_{\ft',\fs}^{\stab,i}$ defines
a bundle map,
\begin{equation}
\label{eq:BundleDiagram1}
\begin{CD}
\tN_{\ft,\fs}(\eps)\times_{\sG_{\fs}\times S^1}\rd\bar{\Gl}_{\ft}(\delta)
@> \pi_{N,\nu} >>
N_{\ft,\fs}\times_{S^1}\ES^1
\\
@V \pi_{\Gl} VV @V \pi_{N}\times\pi_B VV
\\
\tM_{\fs}\times_{\sG_{\fs}\times S^1}\rd\bar{\Gl}_{\ft}(\delta)
@> \pi_{\fs}\times\iota_{\nu} >>
M_{\fs}\times \BS^1
\end{CD}
\end{equation}
the Thom class of the normal bundle of
$\bB\bL^{\stab,i}_{\ft',\fs}$ in $\bL^{\stab,i}_{\ft',\fs}$
is given by $\pi_{N,\nu}^*\Th_{S^1}(N_{\ft,\fs})$.
The relative fundamental class of the manifold with boundary
$\bL^{\stab,i}_{\ft',\fs}$ is equal to
$[\bL_{\ft',\fs}^{\stab,i},\rd\bL_{\ft',\fs}^{\stab,i}]$, as noted
in \eqref{eq:DefineRelFundClass}.
Hence, because the Thom class of the normal bundle
of a submanifold is equal to
the Poincar\'e dual of that submanifold \cite[p. 371]{BredonTopGeom},
the fundamental class of $\bB\bL^{\stab,i}_{\ft',\fs}$ is given by
\begin{equation}
\label{eq:PDOfBL}
[\bB\bL^{\stab,i}_{\ft',\fs}]=j^*_1\pi_{N,\nu}^*\Th_{S^1}(N_{\ft,\fs}) \cap
\left[\bL_{\ft',\fs}^{\stab,i},\rd\bL_{\ft',\fs}^{\stab,i}\right],
\end{equation}
where $\left[\bL_{\ft',\fs}^{\stab,i},\rd\bL_{\ft',\fs}^{\stab,i}\right]$
is defined following \eqref{eq:JCohomology} and
$$
j_1:(\bar\bL^{\stab}_{\ft',\fs},\bar\bL^{\stab,s}_{\ft',\fs})
\to
(\bar\bL^{\stab}_{\ft',\fs},\bar\bL^{\stab}_{\ft',\fs}-\bB\bL^{\stab,i}_{\ft',\fs})
$$
is the inclusion of pairs.  Observe that
\begin{equation}
\label{eq:InclusionRelation}
j_{\bB\bL}=j_1\circ j,
\end{equation}
where $j:(\bar\bL^{\stab}_{\ft',\fs},\emptyset)\to
(\bar\bL^{\stab}_{\ft',\fs},\bar\bL^{\stab,s}_{\ft',\fs})$
is the inclusion \eqref{eq:DefineJ}.
Thus, for any $\om \in H^{\bullet}(\bar\bL^{\stab}_{\ft',\fs};\ZZ)$, we have
\begin{align*}
{}&
\left\langle
\om\smile j^*_{\bB\bL}\pi_{N,\nu}^*\Th_{S^1}(N_{\ft,\fs}),
[\bar\bL^{\stab}_{\ft',\fs}]
\right\rangle
\\
& =
\left\langle
    j^*\left(\om\smile j_1^*\pi_{N,\nu}^*\Th_{S^1}(N_{\ft,\fs})\right),
    [\bar\bL^{\stab}_{\ft',\fs}]
\right\rangle
\quad\text{(By \eqref{eq:InclusionRelation} 
\& \cite[Statement 5.6.8]{Spanier})}
\\
& =
\left\langle
    \om\smile j_1^*\pi_{N,\nu}^*\Th_{S^1}(N_{\ft,\fs}),
j_*[\bar\bL^{\stab}_{\ft',\fs}]
\right\rangle
\\
& =
\left\langle
    \om,j_1^*\pi_{N,\nu}^*\Th_{S^1}(N_{\ft,\fs})\cap
    \left[\bL_{\ft',\fs}^{\stab,i},\rd\bL_{\ft',\fs}^{\stab,i}\right]
\right\rangle
\quad\text{(By Equation \eqref{eq:DefineFundClassOfL})}
\\
& =
\left\langle \om,[\bB\bL^{\stab,i}_{\ft',\fs}]\right\rangle,
\quad\text{(By Equation \eqref{eq:PDOfBL})}
\end{align*}
which proves equation \eqref{eq:PDOfBasePairing} and thus Assertion (1).

Because the diagram \eqref{eq:PDBundleMap2} commutes,
we observe that $\pi_B\circ\pi_{N,\nu}=\iota_\nu$;
this observation and the fact that
$\iota_\nu^*\bc=\nu$ proves Assertion (2).

Because the diagram \eqref{eq:PDBundleMap1} commutes,
we have $\pi_M\circ\pi_{N,\nu}=\pi_{\fs}$ and this
proves Assertion (3), completing the proof of the lemma.
\end{proof}

We next give an explicit formula for division by an equivariant Thom class
in terms of the Segre classes of the vector
bundle; this division process is also referred to as {\em equivariant
localization\/} --- see \cite[Theorem 7.13]{BerlineGetzlerVergne}. Recall
that the total 
Segre class $s(N)=s_0(N)+s_1(N)+\cdots$, where $s_i(N)\in H^{2i}(M;\ZZ)$,
of a complex vector bundle $N$ over a topological space $M$ is the formal
inverse of the total Chern class $c(N)=1+c_1(N)+c_2(N)+\cdots$, where
$c_i(N)\in H^{2i}(M;\ZZ)$, so that $s(N)c(N)=1$ \cite[p. 69]{Fulton},
\cite[Lemma 4.10]{FL2b}. 

\begin{prop}
\label{prop:S1Localization}
Let $\pi_M:N\to M$ be a complex vector bundle of rank $r$ over a
topological space $M$ with a fundamental class of dimension $d$.  Suppose
$\pi_{B}:\ES^1\times_{S^1}N\to \BS^1$ and $\pi_M:\ES^1\times_{S^1}N\to M$
are the projection maps, where the circle acts diagonally on $N\times
\ES^1$ by scalar
multiplication on the fibers of $N$ and by the standard action on $\ES^1$.
Let $N^0$ denote the complement of the zero-section and let
$\Th_{S^1}(N)\in H^{2r}_{S^1}(N,N^0;\ZZ)$ be
the Thom class of
\begin{equation}
\label{eq:EquivariantN}
\pi_{B}\times\pi_M:\ES^1\times_{S^1}N \to \BS^1\times M.
\end{equation}
If $\bc\in H^2(\BS^1;\ZZ)$ denotes the universal first Chern class, $0\le
k\le [d/2]$ is an integer, and $\alpha\in H^{d-2k}(M;\ZZ)$, then
for integers $m\ge k+r$, one has
\begin{equation}
\label{eq:FormalInverseOfThomClass}
\pi_B^*\bc^{m}\smile \pi_M^*\alpha
=
\iota_n^*\Th_{S^1}(N)\smile
    \left(
        \sum_{j=0}^k (-1)^{r+j}
    \pi_B^*\bc^{m-r-j}\smile \pi_M^*(\alpha\smile s_j(N))
    \right),
\end{equation}
where $\iota_N:(N,\emptyset)\to (N,N^0)$ is the inclusion of pairs,
and the pushforward formula,
\begin{equation}
\label{eq:S1Pushforward}
\pi_{B,*}(\bc^{m}\smile \pi_M^*\alpha)
=(-1)^{k+r}\bc^{m-r-k}\langle \alpha\smile s_k(N),[M]\rangle.
\end{equation}
\end{prop}

\begin{proof}
By the splitting principle we may suppose that $N=\oplus_{i=1}^r L_i$, with
$y_i=c_1(L_i)$, and so Lemma \ref{lem:DiagonalQuotient} implies that
\begin{equation}
\label{eq:EquivariantThomClass}
j^*\Th_{S^1}(N)
=\prod_{i=1}^r (-\pi_B^*\bc+ \pi_M^*y_i)
=\sum_{i=0}^r (-1)^{r-i} \pi_B^*\bc^{r-i} \pi_M^*c_i(N).
\end{equation}
The negative sign above arises because the circle action is diagonal, as
explained in Lemma
\ref{lem:DiagonalQuotient}. The bundle $N$ has total Segre class
$s(N)=\prod_{i=1}^r (1+y_i)^{-1}$, and thus
\begin{align*}
\frac{1}{\Th_{S^1}(N)}
&=
\prod_{i=1}^r \frac{1}{(-\pi_B^* \bc + \pi_M^*y_i)}
\\
&= (-\pi_B^*\bc)^{-r}\prod_{i=1}^r\frac{ 1}{(1-\pi_M^*y_i/\pi_B^*\bc)}
\\
&= \sum_{j=0}^\infty (-\pi_B^*\bc)^{-r-j} \pi_M^*s_j(N).
\end{align*}
This expression for the formal inverse of $\Th_{S^1}(N)$ yields
\begin{equation}
\label{eq:DivisionByThomClass}
\frac{\pi_B^*\bc^{m}\smile \pi_M^*\alpha}{\Th_{S^1}(N)}
=
\sum_{j=0}^k (-1)^{r+j}\pi_B^*\bc^{m-r-j}\smile \pi_M^*(\alpha\smile s_j(N)),
\end{equation}
and equation \eqref{eq:FormalInverseOfThomClass} follows.
Because the pushforward $\pi_{B,*}$ is the composition
of division by the Thom class $\Th_{S^1}(N)$ and integration
over $M$ \cite[\S 6]{BT}, we have
\begin{align*}
\pi_{B,*}(\pi_B^*\bc^{m}\smile \pi_M^*\alpha)
&=
\left(\frac{\pi_B^*\bc^{m}\smile \pi_M^*\alpha}{\Th_{S^1}(N)}\right)/[M]
\\
&=
\left(
\sum_{j=0}^k (-1)^{r+j}\pi_B^*\bc^{m-r-j}\smile \pi_M^*(\alpha\smile s_j(N))
\right)/[M],
\end{align*}
which yields equation \eqref{eq:S1Pushforward}, since $M$ has dimension $d$.
\end{proof}

Lemmas \ref{lem:PullbackMuc1}, \ref{lem:CohomologyOnReducibleLink},
\ref{lem:EulerOfBackgroundObstruction},
and \ref{lem:InstantonEuler} show how to express the cohomology
classes in \eqref{eq:CohomFormulation} in terms of the cohomology
classes $\nu$, $\mu_{\fs}$, and classes pulled back from $X$.
The next proposition reduces the pairing in
\eqref{eq:CohomFormulation} of products of these latter classes
with $[\bL^{\stab}_{\ft',\fs}]$ to
a pairing with the fundamental class
of the submanifold \eqref{eq:BaseSpace}.

\begin{prop}
\label{prop:LinkSegreFormula}
Let $d_s=\dim M_{\fs}$ and let $r_N$ denote the rank of the complex vector
bundle $N_{\ft,\fs}$ over $M_{\fs}$.  Assume $b_1(X)=0$ and $b_2^+(X)$ is
odd.  For integers $0\le i\le d=\half d_s$ and $0\leq k \leq 2$ and any
class $\alpha\in H^{2k}(X;\ZZ)$, we have
\begin{equation}
\label{eq:LinkPairingSegre}
\begin{aligned}
{}&\left\langle
\nu^{d+r_N+3-k-i}\smile \pi_{\fs}^*\mu_{\fs}^i \smile\pi_X^*\alpha,
[\bar\bL^{\stab}_{\ft',\fs}]
\right\rangle
\\
&\quad =
\sum_{j=0}^{d-i} (-1)^{r_N+j}
\left\langle
    \nu^{d+3-i-j-k}
    \smile \pi_{\fs}^*(\mu^i_{\fs}\smile s_j(N_{\ft,\fs}))\smile\pi_X^*\alpha,
[\bB\bL_{\ft',\fs}^{\stab,i}]\right\rangle,
\end{aligned}
\end{equation}
where $\nu\in H^2(\bar\bL^{\stab}_{\ft',\fs};\ZZ)$ is the cohomology class
of Definition \ref{defn:DefnOfNu},
$\mu_{\fs}\in H^2(M_{\fs};\ZZ)$ is the restriction
of the cohomology class defined in \eqref{eq:SWClass},
and
$\pi_X:\bar\bL^{\stab}_{\ft',\fs}\to X$
and
$\pi_{\fs}: \bar\bL^{\stab}_{\ft',\fs}\to M_{\fs}$
are the restrictions of the projections defined in
\eqref{eq:ProjFromCptVirtModSpace}.
\end{prop}

\begin{proof}
Applying the equality $\nu=\iota_\nu^*\bc$ from Lemma
\ref{lem:PDualOfBaseSpace} and using the abbreviations $A=d+r_N+3-k$,
$T=\Th_{S^1}(N_{\ft,\fs})$, and $s_j=s_j(N_{\ft,\fs})$ yields
\begin{equation}
\label{eq:ReduceToPullback}
\begin{aligned}
{}&\left\langle
\nu^{A-i}\smile\pi_{\fs}^*\mu_{\fs}^i \smile\pi_X^*\alpha,
[\bar\bL^{\stab}_{\ft',\fs}]\right\rangle
\\
&=
\left\langle
\pi_{N,\nu}^*
\left(\pi_B^* \bc^{A-i}\smile\pi_{\fs}^*\mu_{\fs}^i
\right)
     \smile\pi_X^*\alpha,
[\bar\bL^{\stab}_{\ft',\fs}]
\right\rangle
\\
&=
    \sum_{j=0}^{d-i} (-1)^{r_N+j}
    \left\langle
        \pi_{N,\nu}^*
        \left(
             \pi_B^* \bc^{A-i-j}
             \smile
              \pi_{\fs}^*( \mu_{\fs}^i\smile s_j)
              \smile
              j^*T
        \right)
        \smile\pi_X^*\alpha,
        [\bar\bL^{\stab}_{\ft',\fs}]
    \right\rangle
\\
&\qquad\text{(By Equation \eqref{eq:FormalInverseOfThomClass})}
\\
&=
\sum_{j=0}^{d-i} (-1)^{r_N+j}
    \left\langle
        (j^*\pi_{N,\nu}^*T)\smile \nu^{A-r_N-i-j}
        \smile \pi_{\fs}^*( \mu_{\fs}^i\smile s_j)
        \smile \pi_X^*\alpha,[\bar\bL^{\stab}_{\ft',\fs}]
    \right\rangle
\\
&=
\sum_{j=0}^{d-i} (-1)^{r_N+j}
    \left\langle
        \nu^{A-r_N-i-j}
        \smile \pi_{\fs}^*( \mu_{\fs}^i\smile s_j)
        \smile \pi_X^*\alpha,[\bB\bL^{\stab,i}_{\ft',\fs}]
    \right\rangle
\end{aligned}
\end{equation}
where the last equality follows from equation \eqref{eq:PDOfBasePairing}.
This proves the desired formula \eqref{eq:LinkPairingSegre}.
\end{proof}

The next step is to use the product decomposition 
\eqref{eq:ProductDecomp} of $\bB\bL_{\ft',\fs}^{\stab,i}$
to reduce the pairing on the right-hand side of equation
\eqref{eq:LinkPairingSegre} to a sum of
products of pairings with $[M_{\fs}]$
and pairings with $[\rd\bar{\Gl}_{\ft}(\delta)/S^1]$.  To accomplish this
we must write the cohomology classes appearing in
\eqref{eq:LinkPairingSegre}, namely $\mu_{\fs}$, $\nu$, and
$\pi_X^*\alpha$, as pullbacks of cohomology classes on the factors
$M_{\fs}$ and $\rd\bar{\Gl}_{\ft}(\delta)/S^1$ of the product; of these
three cohomology classes, only $\nu$ is not obviously such a pullback.

\begin{defn}
\label{defn:NuX}
Let $\nu_{\ft}\in H^2(\rd\bar{\Gl}_{\ft}(\delta)/S^1;\RR)$ be the first
Chern class of the circle bundle
$$
\rd\bar{\Gl}_{\ft}(\delta)\to\rd\bar{\Gl}_{\ft}(\delta)/S^1,
$$
where the circle action is described in the paragraph preceding
equation \eqref{eq:GaugeGrpS1OnGluing}. We shall also let
$\nu_{\ft}$ denote the pullback of this class to the product
$M_{\fs}\times\rd\bar{\Gl}_{\ft}(\delta)/S^1$, where it is the first Chern
class of the circle bundle
$$
M_{\fs}\times\rd\bar{\Gl}_{\ft}(\delta)
\to M_{\fs}\times\rd\bar{\Gl}_{\ft}(\delta)/S^1.
$$
\end{defn}

In the following lemma, we obtain the desired decomposition of
$\nu$ by comparing this class with $\nu_{\ft}$:

\begin{lem}
\label{lem:CompareNu}
Let $\nu$, $\nu_{\ft}$ be the first Chern
classes in Definitions
\ref{defn:DefnOfNu} and \ref{defn:NuX}, respectively. Let $\pi_{\fs}:
M_{\fs}\times
\rd\bar{\Gl}_{\ft}(\delta)/S^1\to M_{\fs}$ be the projection.  Then the
restriction of $\nu$ to $M_{\fs}\times\rd\bar{\Gl}_{\ft}(\delta)/S^1$ is
given by 
\begin{equation}
\label{eq:PiecesOfNu}
\left.\nu\right|_{M_{\fs}\times \rd\bar{\Gl}_{\ft}(\delta)/S^1}
=
\nu_{\ft} +2c_1(\LL_{\fs}),
\end{equation}
and if $b_1(X)=0$, this simplifies to
$$
\left.\nu\right|_{M_{\fs}\times \rd\bar{\Gl}_{\ft}(\delta)/S^1}
=
\nu_{\ft}+2\mu_{\fs}.
$$
\end{lem}

\begin{proof}
To prove the lemma, we must compare the second bundle in Definition
\ref{defn:NuX} with the circle bundle
\begin{equation}
\label{eq:NuBundle}
\tM_{\fs}\times_{\sG_{\fs}}\rd\bar{\Gl}_{\ft}(\delta)
\to
\tM_{\fs}\times_{\sG_{\fs}\times S^1}\rd\bar{\Gl}_{\ft}(\delta).
\end{equation}
The circle action in \eqref{eq:NuBundle} is trivial on $\tilde M_{\fs}$ and
given by the action described prior to \eqref{eq:GaugeGrpS1OnGluing} on
$\rd\bar{\Gl}_{\ft}(\delta)\subset \bar{\Gl}_{\ft}(\delta)$; observe that
the action of the group $\sG_{\fs}$ in the bundle \eqref{eq:NuBundle} is
defined in
\eqref{eq:GaugeGrpS1OnGluing}, for $s\in\sG_{\fs}$, $(B,\Psi)\in\tM_{\fs}$,
and $\bg\in\rd\bar{\Gl}_{\ft}(\delta)$, by
$$
(s, (B,\Psi,\bg))
\mapsto
(s_*(B,\Psi),s^{-2}\bg).
$$
Let $\LL^1_{\fs}$ be the unit sphere 
of the bundle \eqref{eq:DefineSWUniversal}, so
$$
\LL^1_{\fs}= \tM_{\fs}\times_{\sG_{\fs}}(X\times S^1),
$$
with the action \eqref{eq:DefineSWUniversalS1Action} of $\sG_{\fs}$.
Consider the bundle
\begin{equation}
\label{eq:FiberedNuBundle}
\left(\LL^1_{\fs}\times_{X} \rd\bar{\Gl}_{\ft}(\delta)\right)/S^1
\to
M_{\fs}\times \rd\bar{\Gl}_{\ft}(\delta)/S^1,
\end{equation}
where the action of $e^{i\gamma}\in S^1$ 
is given for $[B,\Psi,x,e^{i\theta}] \in
\LL^1_{\fs}|_x$ and $\bg\in \rd\bar{\Gl}_{\ft}(\delta)|_x$  by
$$
\left(
[B,\Psi,x,e^{i\theta}], \bg
\right)
\mapsto
\left(
[B,\Psi,x,e^{i(\theta+\gamma)}],e^{-2i\gamma}\bg
\right).
$$
By Lemma \ref{lem:DiagonalQuotient}, the first Chern class
of the bundle \eqref{eq:FiberedNuBundle} is $\nu_{\ft}+2c_1(\LL_{\fs})$.
To prove \eqref{eq:PiecesOfNu}, it suffices to
show that the bundles \eqref{eq:FiberedNuBundle} and
\eqref{eq:NuBundle} are isomorphic.
Because $\sG_{\fs}$ acts anti-diagonally in the definition of
$\LL_{\fs}$, the map
\begin{gather*}
\left(\LL^1_{\fs}\times_{X} \rd\bar{\Gl}_{\ft}(\delta)\right)/S^1
\to
\tM_{\fs}\times_{\sG_{\fs}}\rd\bar{\Gl}_{\ft}(\delta),
\\
\left[
[B,\Psi,x,e^{i\theta}], \bg
\right]
\mapsto
\left[
(B,\Psi),e^{2i\theta}\bg
\right]
\end{gather*}
is a well-defined bundle isomorphism.
The final statement of the lemma follows from the fact that
$c_1(\LL_{\fs})=\mu$ \cite[Lemma 2.24]{FL2a}.
\end{proof}

As discussed in the paragraph
preceding Definition \ref{defn:NuX}, the results of Lemma
\ref{lem:CompareNu} show how the relevant cohomology classes on
$\bB\bL_{\ft',\fs}^{\stab,i}$ pull back from either $M_{\fs}$ or
$\rd\bar{\Gl}_{\ft}(\delta)/S^1$.  Pairing the appropriate cohomology class
with the fundamental class $[M_{\fs}]$ yields the Seiberg-Witten
invariant.  The required pairings with
$[\rd\bar{\Gl}_{\ft}(\delta)/S^1]$ are given by the following lemma:

\begin{lem}
\label{lem:PairingsWithInstantonLink}
Let $\nu_{\ft}$ be the characteristic class in Definition
\ref{defn:NuX}.
Let $\rd\bar{\Gl}_{\ft}(\delta)/S^1$ have the standard orientation
as defined in the paragraph following \eqref{eq:TangentOfIBoundary}.
If $x\in H_0(X;\ZZ)$ is the positive generator
and $h\in H_2(X;\RR)$, then
\begin{equation}
\label{eq:InstLinkPairing}
\begin{aligned}
\langle \nu_{\ft}\smile \pi_X^*\PD[x],
[\rd\bar{\Gl}_{\ft}(\delta)/S^1]\rangle 
&= 2,
\\
\langle \nu_{\ft}^2\smile \pi_X^*\PD[h],
[\rd\bar{\Gl}_{\ft}(\delta)/S^1]\rangle 
&=
-4\langle (c_1(\fs)-c_1(\ft))\smile\PD[h],[X]\rangle,
\\
\langle \nu_{\ft}^3,[\rd\bar{\Gl}_{\ft}(\delta)/S^1]\rangle 
&=
6(c_1(\fs)-c_1(\ft))^2+2c_1^2(X).
\end{aligned}
\end{equation}
\end{lem}

\begin{proof}
In \cite[\S 3]{LenessWC}, a rank-two, Hermitian vector bundle, $F\to X'$,
where $\pi:X'\to X$ is a finite-degree cover
is constructed with $F\to \pi^*\Gl_{\ft}$
a branched cover of degree negative two.
(Similar bundles are constructed in
\cite{Yang,OzsvathBlowUp,LenessBlowUp}.)
The circle action on $F$ induced by the action of $\CC^*$ and the inclusion
$S^1\subset\CC^*$ maps to twice the circle action on
$\pi^*\Gl_X$.  This implies that the map on quotients, $\PP(F)\to
\pi^*\rd\bar{\Gl}_{\ft}(\delta)/S^1$, has degree negative one
and that the class $\nu_{\ft}$ pulls back
to $2h$ where $h$ is the first Chern class of the circle
bundle $F/\RR^*\to\PP(F)$.
Pulling back the cohomology classes  in
\eqref{eq:InstLinkPairing} to $\PP(F)$  and
applying the computation of the Segre classes
of $F$ in \cite[\S 3]{LenessWC} then yield the formulas in
\eqref{eq:InstLinkPairing}.
\end{proof}

Using Poincar\'e duality, we can use Lemma
\ref{lem:PairingsWithInstantonLink} to compute further pairings 
needed in our proof of Theorem \ref{thm:LevelOne}. For example, if
$\alpha\in H^2(X;\RR)$, then
\begin{equation}
\label{eq:PairingsWithInstantonLink4}
\langle \nu_{\ft}\smile \pi_X^*(\alpha\smile\PD[h]),
[\rd\bar{\Gl}_{\ft}(\delta)/S^1]\rangle
=
\langle \nu_{\ft}\smile\pi_X^*\PD[x],
[\rd\bar{\Gl}_{\ft}(\delta)/S^1]\rangle
\langle\alpha,h\rangle,
\end{equation}
using 
\begin{align*}
\alpha\smile\PD[h]
&=
\langle\alpha\smile\PD[h],[X]\rangle\PD[x]
\\
&=
\langle\alpha,\PD[h]\cap[X]\rangle\PD[x]
\\
&=
\langle\alpha,h\rangle\PD[x],
\end{align*}
by Poincar\'e duality.

\begin{rmk}
In comparing the formulas of this and related articles with computations of
wall-crossing formulas for Donaldson invariants, the cohomology class
$\xi=c_1(Q_{\xi})\in H^2(X;\ZZ)$ of \cite{KotschickMorgan,LenessWC}
defining a reduction $\fg_{V} \cong \underline{\RR}\oplus L$
is equal to the cohomology class
$c_1(L)=c_1(\ft)-c_1(\fs)=c_1(\ft')-c_1(\fs)$
appearing in this article and in \cite{FL2a,FL2b,FKLM}.
\end{rmk}


\section{Intersection with link of level-one Seiberg-Witten moduli space}
\label{sec:comp}
In this section we prove Theorem \ref{thm:WCL1} by computing a formula,
given in Theorem \ref{thm:LevelOne}, for intersection number
\eqref{eq:DesiredSWLinkIntersectionNumber} and applying this formula to the
cobordism sum \eqref{eq:CobordismSum}.  We divide this task into three
parts.  Section \ref{subsec:Algebraic} is devoted to the proof of Theorem
\ref{thm:LevelOne}.  Because the Donaldson invariants are defined using a
moduli space of anti-self-dual connections on $X\#\overline{\CC\PP}^2$ (see
\eqref{eq:DefineDonaldson}), we present a blow-up formula for the link
intersection numbers in \S \ref{subsec:BlowUp}.  Finally, in \S
\ref{subsec:WCL1}, we apply the results of \cite{FKLM} and Theorem
\ref{thm:LevelOne} to the cobordism formula \eqref{eq:CobordismSum} to
prove Theorem \ref{thm:WCL1}.

\subsection{Algebraic computations and completion of proof of Theorem
\protect{\ref{thm:LevelOne}}}
\label{subsec:Algebraic}
With the results of \S \ref{sec:DualCohom}, the proof of the following theorem
is essentially algebraic.
Recall that a \spinu structure over $X$ splits, $\ft=\fs\oplus\fs'$, if and
only if $(c_1(\ft)-c_1(\fs))^2=p_1(\ft)$ \cite[Lemma 3.32]{FL2b}. Hence, a
Seiberg-Witten stratum $M_{\fs}\times\Sym^\ell(X)$ is contained in the
level $\sM_{\ft_\ell}\times\Sym^\ell(X)$ of the space of ideal $\SO(3)$
monopoles enclosing $\bar\sM_{\ft'}$ if and only if
$(c_1(\ft')-c_1(\fs))^2=p_1(\ft')+4\ell$; in this situation,
$\bL_{\ft',\fs}$ 
is the link \eqref{eq:DefineLink} of the stratum
$M_\fs\times\Sym^{\ell}(X)$ in $\bar\sM_{\ft'}/S^1$.

\begin{thm}
\label{thm:LevelOne} Let $X$ be a four-manifold
with $b_1(X)=0$, odd $b_2^+(X)\geq 1$.  Suppose $X$ has a generic
Riemannian metric and a \spinu structure $\ft'$, where $w_2(\ft')$ is {\em
good\/} in the sense of Definition \ref{defn:Good}. Let $\fs$ be a \spinc
structure over $X$ for which $(c_1(\ft')-c_1(\fs))^2=p_1(\ft')+4$.  Let
$\delta$, $m$, and $\eta$ be non-negative integers satisfying
\begin{equation}
\label{eq:DegreeAssumption}
0\leq m\leq [\delta/2]
\quad\text{and}\quad
2(\delta + \eta) = \dim(\sM_{\ft'}/S^1)-1.
\end{equation}
If $x\in H_0(X;\ZZ)$ is the positive generator, $h\in H_2(X;\RR)$,
and $\bL_{\ft',\fs}$ has the standard orientation
defined prior to Lemma \ref{lem:LinkOrientation}, then
\begin{equation}
\label{eq:LevelOne}
\begin{aligned}
{}&\#\left(\bar\sV(h^{\delta-2m}x^m)\cap 
\bar\sW^{\eta}\cap\bL_{\ft',\fs}\right)
\\
&=
(-1)^{m+1+d_s(\fs)/2}2^{-\delta} 2^{d_s(\fs)/2}P^{a,b}_{d_s(\fs)/2}(0)
\langle\mu_{\fs}^{d_s(\fs)/2},[M_{\fs}]\rangle
\left(a_0\langle c_1(\fs)-c_1(\ft'),h\rangle^{\delta-2m}\right.
\\
&\quad+
b_0\langle c_1(\fs)-c_1(\ft'),h\rangle^{\delta-2m-1}\langle c_1(\ft'),h\rangle
\\
&\quad+
\left. a_1\langle c_1(\fs)-c_1(\ft'),h\rangle^{\delta-2m-2}Q_X(h,h)
\right),
\end{aligned}
\end{equation}
where all terms on the right which would have a negative exponent are
omitted,
\begin{align*}
a_0&=
 3(c_1(\fs)-c_1(\ft'))^2+c_1^2(X) +2(c_1(\fs)-c_1(\ft'))\cdot c_1(\ft')
+4(\delta-2m)
-4m,
\\
b_0&= 2(\delta-2m)\frac{P^{a-1,b+1}_{d_s(\fs)/2}(0)}{P^{a,b}_{d_s(\fs)/2}(0)},
\\
a_1&= 4\binom{\delta-2m}{2},
\end{align*}
and $P^{a,b}_{d_s(\fs)/2}(0)$ is the constant coefficient of the Jacobi
polynomial \eqref{eq:DefineJacobiPolynomial} with 
\begin{equation}
\label{eq:Defnab} 
\begin{aligned}
a&=\eta - \textstyle{\frac{1}{2}}d_s(\fs) + 1,
\\
b&= \textstyle{\frac{1}{2}}(2\delta-d_a(\ft')-d_s(\fs))
-\textstyle{\frac{1}{4}}(\chi+\si),
\end{aligned}
\end{equation}
where $d_a(\ft')$ is given by equation \eqref{eq:Transv}.
If $d_s(\fs)=0$, then $P^{a,b}_{d_s(\fs)/2}(0)=1$.
\end{thm}

\begin{proof}
According to Proposition \ref{prop:IntersectionNoToCupProduct}, the
intersection number on the left-hand side of equation \eqref{eq:LevelOne}
is equal to
\begin{equation}
\label{eq:CohomClasses1}
\begin{aligned}
\sP
:&=
\# \left( \bar\sV(h^{\delta-2m}x^m)\cap 
\bar\sW^{\eta} \cap \bL_{\ft',\fs}\right)
\\
&=
\langle
\bga^*\barmu_p(h^{\delta-2m}x^m)\smile\bga^*\barmu_c^{\eta}
\smile \bar e_s \smile \bar e_i,
[\bar\bL^{\stab}_{\ft',\fs}]
\rangle,
\end{aligned}
\end{equation}
where $\bL^{\stab}_{\ft',\fs}$ has the standard orientation defined
in \S \ref{subsec:Orient}.
Writing
\begin{equation}
\label{eq:Defndb-pair}
d = \textstyle{\frac{1}{2}}d_s(\fs)
\quad\text{and}\quad
\fb = \langle c_1(\fs)-c_1(\ft'),h\rangle,
\end{equation}
we apply Lemmas \ref{lem:CohomologyOnReducibleLink},
\ref{lem:EulerOfBackgroundObstruction} and
\ref{lem:InstantonEuler} to expand the cup product
in the pairing \eqref{eq:CohomClasses1} to give
\begin{equation}
\label{eq:CohomClasses1-1}
\begin{aligned}
\sP
&=
\left\langle
\left( \textstyle{\frac{1}{2}} \fb(2\mu_{\fs}-\nu)+\PD[h]\right)^{\delta-2m}
\smile
\left( -\textstyle{\frac{1}{4}} (2\mu_{\fs}-\nu)^2+\PD[x]\right)^{m}
\right.
\\
&\left. \qquad \smile
(-\nu)^{\eta+r_{\Xi}}\left(
\textstyle{\frac{1}{2}}c_1(\ft')-\textstyle{\frac{1}{2}}\nu\right),
[\bar\bL^{\stab}_{\ft',\fs}]\right\rangle,
\end{aligned}
\end{equation}
where, for notational simplicity, we omit explicit mention of pullbacks.
Denoting
\begin{equation}
\label{eq:DefnOfC}
C=(-1)^{m+\eta+r_{\Xi}+1}2^{-\delta-1},
\end{equation}
using the fact that the classes $\PD[x],\PD[h],c_1(\ft')$ are pulled back
from $X$,
and multiplying out, we see that equation \eqref{eq:CohomClasses1-1} becomes
\begin{equation}
\label{eq:CohomClasses1-2}
\begin{aligned}
\sP
&=
C\left\langle \fb^{\delta-2m}(2\mu_{\fs}-\nu)^{\delta}\nu^{\eta+r_{\Xi}+1}
-
\fb^{\delta-2m}(2\mu_{\fs}-\nu)^{\delta}\nu^{\eta+r_{\Xi}}c_1(\ft')
\right.
\\
&\quad
+2(\delta-2m)\fb^{\delta-2m-1}(2\mu_{\fs}-\nu)^{\delta-1}
\nu^{\eta+r_{\Xi}+1}\PD[h]
\\
&\quad
-2(\delta-2m)\fb^{\delta-2m-1}(2\mu_{\fs}-\nu)^{\delta-1}
\nu^{\eta+r_{\Xi}}\PD[h]c_1(\ft')
\\
&\quad
+4\binom{\delta-2m}{2}\fb^{\delta-2m-2}
(2\mu_{\fs}-\nu)^{\delta-2}\nu^{\eta+r_{\Xi}+1}Q_X(h,h)\PD[x]
\\
&\quad
\left.
-
4m\fb^{\delta-2m}(2\mu_{\fs}-\nu)^{\delta-2}\nu^{\eta+r_{\Xi}+1}\PD[x],
[\bar\bL^{\stab}_{\ft',\fs}]\right\rangle.
\end{aligned}
\end{equation}
If $r_N$ denotes the rank of the complex vector bundle $N_{\ft,\fs}\to
M_{\fs}$, then (see equations \eqref{eq:SWDimRelations} and
\eqref{eq:SWBundleRankRelations}) 
\begin{equation}
\label{eq:Dimensions}
\begin{aligned}
\delta+\eta 
&= 
\textstyle{\frac{1}{2}}\dim(\sM_{\ft'}/S^1)-1
\\
&= 
\textstyle{\frac{1}{2}}\dim(\sM_{\ft}/S^1)+2
\quad\text{(By Equation \eqref{eq:Transv})}
\\
&=
d+r_N-r_{\Xi}+2.
\end{aligned}
\end{equation}
The identity \eqref{eq:Dimensions} and the fact that
$\mu_{\fs}^i=0$ when $i>d$, 
because the class $\mu_{\fs}$ is pulled back from $M_{\fs}$, imply that
equation \eqref{eq:CohomClasses1-2} can be rewritten as
\begin{equation}
\begin{aligned}
\label{eq:CohomClass2}
\sP
&=
(-1)^{\delta}C
\sum_{i=0}^{d} \left\langle
\fb^{\delta-2m}
\binom{\delta}{i}2^i(-1)^i
\mu_{\fs}^i\nu^{d+r_N+3-i}
\right.
\\
&\quad
-\fb^{\delta-2m}
 \binom{\delta}{i}2^i(-1)^i
\mu_{\fs}^i\nu^{d+r_N+2-i}c_1(\ft')
\\
&\quad
-2(\delta-2m)\fb^{\delta-2m-1}
\binom{\delta-1}{i}2^i(-1)^i
\mu_{\fs}^i\nu^{d+r_N+2-i}\PD[h]
\\
&\quad
+2(\delta-2m)\fb^{\delta-2m-1}
\binom{\delta-1}{i}2^i(-1)^i
\mu_{\fs}^i\nu^{d+r_N+1-i}\PD[h]c_1(\ft')
\\
&\quad
+4\binom{\delta-2m}{2}\fb^{\delta-2m-2}
\binom{\delta-2}{i}2^i(-1)^i
\mu_{\fs}^i\nu^{d+r_N+1-i}Q_X(h,h)\PD[x]
\\
&\quad
\left.
-4m \fb^{\delta-2m}
\binom{\delta-2}{i}2^i(-1)^i
\mu_{\fs}^i\nu^{d+r_N+1-i}\PD[x],
[\bar\bL^{\stab}_{\ft',\fs}]\right\rangle,
\end{aligned}
\end{equation}
where it is understood that binomial coefficients $\binom{\delta}{i}$
are by definition zero when $i>\delta$.
By applying the formula \eqref{eq:LinkPairingSegre} for division by the
Poincar\'e dual of $[\bB\bL_{\ft',\fs}^{\stab,i}]=[M_{\fs}\times
\rd\bar{\Gl}_{\ft}(\delta)/S^1]$, we see that equation
\eqref{eq:CohomClass2} yields
\begin{equation}
\label{eq:CohomClass2a}
\begin{aligned}
\sP
&=
(-1)^{\delta+r_N}C
\sum_{i=0}^d \sum_{j=0}^{d-i}
\left\langle
\fb^{\delta-2m}
    \binom{\delta}{i}2^i(-1)^{i+j}
     \nu^{d+3-i-j}\mu_{\fs}^is_j(N_{\ft,\fs})
\right.
\\
&\quad
-\fb^{\delta-2m}
    \binom{\delta}{i}2^i(-1)^{i+j}
     \nu^{d+2-i-j}\mu_{\fs}^is_j(N_{\ft,\fs})c_1(\ft')
    \\
&\quad
-2(\delta-2m)\fb^{\delta-2m-1}
    \binom{\delta-1}{i}2^i(-1)^{i+j}
    \nu^{d+2-i-j}\mu_{\fs}^is_j(N_{\ft,\fs})\PD[h]
    \\
&\quad
+2(\delta-2m)\fb^{\delta-2m-1}
    \binom{\delta-1}{i}2^i(-1)^{i+j}
    \nu^{d+1-i-j}\mu_{\fs}^i s_j(N_{\ft,\fs})\PD[h]c_1(\ft')
\\
&\quad
+4\binom{\delta-2m}{2}\fb^{\delta-2m-2}
    \binom{\delta-2}{i}2^i(-1)^{i+j}
     \nu^{d+1-i-j}\mu_{\fs}^i s_j(N_{\ft,\fs}) Q_X(h,h)\PD[x]
\\
&\quad
\left.
-4m \fb^{\delta-2m}
 \binom{\delta-2}{i}2^{i}(-1)^{i+j}
     \nu^{d+1-i-j}\mu_{\fs}^i s_j(N_{\ft,\fs})\PD[x],
[M_{\fs}\times\rd\bar{\Gl}_{\ft}(\delta)/S^1]\right\rangle.
\end{aligned}
\end{equation}
We then make the substitution $\nu=2\mu_{\fs}+\nu_{\ft}$ from equation
\eqref{eq:PiecesOfNu} and expand the powers of $\nu$ in
\eqref{eq:CohomClass2a} as binomial sums.
For $\alpha\in H^{2k}(X;\RR)$, the equality
\begin{equation}
\label{eq:ProductPairing}
\begin{aligned}
{}&\left\langle
    \nu_{\ft}^{d+3-k-i-j-\ell}
    \smile
    \mu_{\fs}^{i+\ell}
    \smile
    s_j(N_{\ft,\fs})
    \smile \pi_X^*\alpha,
    [M_{\fs}\times\rd\bar{\Gl}_{\ft}(\delta)/S^1]
\right\rangle
\\
&\quad =
\begin{cases}
\langle \mu_{\fs}^{d-j} \smile s_j(N_{\ft,\fs}),[M_{\fs}]\rangle
\langle \nu_{\ft}^{3-k}\smile\pi_X^*\alpha,
[\rd\bar{\Gl}_{\ft}(\delta)/S^1] \rangle,
&\text{if $\ell= d-i-j$,}
\\
0, & \text{if $\ell\neq d-i-j$,}
\end{cases}
\end{aligned}
\end{equation}
implies that only one term from each of the binomial expansions of the
powers of $\nu$ in equation \eqref{eq:CohomClass2a} will not vanish.
Therefore, equation \eqref{eq:CohomClass2a} yields
\begin{equation}
\label{eq:CohomClass2b}
\begin{aligned}
\sP
&=
(-1)^{\delta+r_N}C
\sum_{i=0}^d\sum_{j=0}^{d-i}
S_j\langle \mu_{\fs}^d,[M_{\fs}]\rangle
\left\langle
\fb^{\delta-2m}
    \binom{\delta}{i}\binom{d+3-i-j}{d-i-j} (-1)^{i+j}2^{d-j}
    \nu_{\ft}^3
\right.
\\
&
-\fb^{\delta-2m}
    \binom{\delta}{i}\binom{d+2-i-j}{d-i-j}(-1)^{i+j} 2^{d-j}
    \nu_{\ft}^{2}c_1(\ft')
    \\
&
-2(\delta-2m)\fb^{\delta-2m-1}
    \binom{\delta-1}{i}\binom{d+2-i-j}{d-i-j}(-1)^{i+j} 2^{d-j}
    \nu_{\ft}^2\PD[h]
    \\
&
+2(\delta-2m)\fb^{\delta-2m-1}
    \binom{\delta-1}{i}\binom{d+1-i-j}{d-i-j}(-1)^{i+j} 2^{d-j}
    \nu_{\ft}\PD[h]c_1(\ft')
\\
&
+4\binom{\delta-2m}{2}\fb^{\delta-2m-2}
    \binom{\delta-2}{i}\binom{d+1-i-j}{d-i-j}(-1)^{i+j} 2^{d-j}
    Q_X(h,h)\nu_{\ft}\PD[x]
\\
&
\left.
-4m \fb^{\delta-2m}
     \binom{\delta-2}{i}\binom{d+1-i-j}{d-i-j}(-1)^{i+j} 2^{d-j}
    \nu_{\ft}\PD[x],[\rd\bar{\Gl}_{\ft}(\delta)/S^1]
\right\rangle,
\end{aligned}
\end{equation}
where we use Lemma \ref{lem:TopOfReducibleNormal} to see that
$$
\left\langle \mu_{\fs}^{d-j} s_{j}(N_{\ft,\fs}),[M_{\fs}] \right\rangle
=
 S_j\langle \mu_{\fs}^d,[M_{\fs}]\rangle,
$$
and the constants $S_j$ are defined by
\begin{equation}
\label{eq:DefineSegre}
S_j
=
\sum_{k=0}^j 2^k\binom{-n'_s}{k}\binom{-n''_s}{j-k}.
\end{equation}
See \cite[Lemma 4.16]{FL2b} for identities involving binomial coefficients
and extensions of the definition of $\binom{n}{p}$ to the case $n\leq 0$.
Before applying Lemma \ref{lem:PairingsWithInstantonLink} to compute 
each of the pairings with $[\rd\bar{\Gl}_{\ft}(\delta)/S^1]$ on the
right-hand side of equation \eqref{eq:CohomClass2b}, we simplify the
combinatorial factors appearing in each term in \eqref{eq:CohomClass2b}
using the identities, for $v=0,\dots,3$,
\begin{equation}
\label{eq:CombinatorialFactor}
\begin{aligned}
{}&
\sum_{i=0}^d\sum_{j=0}^{d-i}
    (-1)^{i+j}2^{d-j}
    \binom{\delta-v}{i}
    \binom{d+3-v-i-j}{d-i-j} S_j
\\
&=
\sum_{i=0}^d\sum_{j=0}^{d-i}\sum_{k=0}^j
    (-1)^{i+j}2^{d-j+k}
    \binom{\delta-v}{i}
    \binom{d+3-v-i-j}{d-i-j}
    \binom{-n'_s}{k}
    \binom{-n''_s}{j-k}
\\
&\qquad\text{(By Equation \eqref{eq:DefineSegre})}
\\
&=
2^d P^{a,b}_d(0)
\quad\text{(By Equation \eqref{eq:CombinatorialFactor1})},
\end{aligned}
\end{equation}
where $P^{a,b}_d(\zeta)$ is the Jacobi polynomial
\eqref{eq:DefineJacobiPolynomial}, with constants
\begin{equation}
\label{eq:ab1}
a=3+n'_s+n''_s -\delta
\quad\text{and}\quad
b=\delta-n'_s-4-d.
\end{equation}
Because the constant on the left-hand side of equation
\eqref{eq:CombinatorialFactor}
is independent of $v$, we can factor it out of all the terms in equation
\eqref{eq:CohomClass2b} except from the fourth term
(the one containing $\PD[h]c_1(\ft')$), where the
binomial factors can be rewritten as
$$
\binom{(\delta+1)-2}{i}\binom{d+1-i-j}{d-i-j},
$$
and thus, by replacing $\delta$ by $\delta+1$
in the definition of $a$ and $b$ in
\eqref{eq:ab1}, the combinatorial expression
from the fourth term of equation
\eqref{eq:CohomClass2b} corresponding to equation
\eqref{eq:CombinatorialFactor} will give the Jacobi polynomial
$P^{a-1,b+1}_d(0)$. 
Hence, by substituting \eqref{eq:CombinatorialFactor}, equation
\eqref{eq:CohomClass2b} takes the shape
\begin{equation}
\label{eq:Pairing2}
\begin{aligned}
\sP
&=
(-1)^{\delta+r_N}C
2^d P^{a,b}_d(0)
\langle \mu_{\fs}^d,[M_{\fs}]\rangle
\left\langle \fb^{\delta-2m} (\nu_{\ft}^3 - \nu_{\ft}^2c_1(\ft'))
\right.
\\
&\quad
-2(\delta-2m)\fb^{\delta-2m-1}\nu_{\ft}^2\PD[h]
+2(\delta-2m)
    \frac{P^{a-1,b+1}_d(0)}{P^{a,b}_d(0)}
    \fb^{\delta-2m-1}\nu_{\ft}\PD[h]c_1(\ft')
\\
&\quad
\left.
+4\binom{\delta-2m}{2}\fb^{\delta-2m-2}Q_X(h,h)\nu_{\ft}\PD[x]
-4m \fb^{\delta-2m} \nu_{\ft}\PD[x],
[\rd\bar{\Gl}_{\ft}(\delta)/S^1]\right\rangle.
\end{aligned}
\end{equation}
Lemma \ref{lem:PairingsWithInstantonLink} and equation 
\eqref{eq:PairingsWithInstantonLink4} provide formulae for the
pairings:
\begin{gather*}
\langle \nu_{\ft}\PD[x],[\rd\bar{\Gl}_{\ft}(\delta)/S^1]\rangle,
\\
\langle \nu_{\ft}\PD[h]c_1(\ft'),[\rd\bar{\Gl}_{\ft}(\delta)/S^1]\rangle,
\\
\langle \nu_{\ft}^2\PD[h],[\rd\bar{\Gl}_{\ft}(\delta)/S^1]\rangle,
\\
\langle \nu_{\ft}^3,[\rd\bar{\Gl}_{\ft}(\delta)/S^1]\rangle.
\end{gather*}
Hence, Lemma \ref{lem:PairingsWithInstantonLink} and the definition
\eqref{eq:DefnOfC} of $C$ imply that
equation \eqref{eq:Pairing2} can be written as
\begin{equation}
\label{eq:Pairing3}
\begin{aligned}
\sP
&=
(-1)^{\delta+\eta+r_N+r_{\Xi}+m+1}2^{d-\delta-1} P^{a,b}_d(0)
\langle \mu_{\fs}^d,[M_{\fs}]\rangle
\\
&\quad\times
\left(\fb^{\delta-2m}\left(6(c_1(\fs)-c_1(\ft'))^2 + 2c_1^2(X)\right.\right.
\\
&\qquad
\left. +4(c_1(\fs)-c_1(\ft'))\cdot c_1(\ft')
+8(\delta-2m)-8m\right)
\\
&\qquad
+4(\delta-2m)\fb^{\delta-2m-1}
\frac{P^{a-1,b+1}_d(0)}{P^{a,b}_d(0)}
\langle c_1(\ft'),h\rangle
\\
&\qquad
\left.
+8\fb^{\delta-2m-2}\binom{\delta-2m}{2}Q_X(h,h)
\right).
\end{aligned}
\end{equation}
We write the power of $2$ as $2^{d_s(\fs)/2}2^{-\delta}2^{-1}$, since
$d=d_s(\fs)/2$ by definition \eqref{eq:Defndb-pair}, with the factor $2^{-1}$
being absorbed into the coefficients $a_0$, $b_0$, $a_1$, whose resulting
formulae agree by inspection with those following equation
\eqref{eq:LevelOne}.

Equation \eqref{eq:Dimensions} implies that the power of $(-1)$ 
in equation \eqref{eq:Pairing3} simplifies to 
$$
(-1)^{d_s(\fs)/2+m+1},
$$  
which agrees with the power of $(-1)$ in \eqref{eq:LevelOne}. 

Finally, we simplify the combinatorial coefficients $P^{a,b}_d(0)$
in \eqref{eq:Pairing3}. Because 
\begin{align*}
d+n'_s(\ft,\fs)+n''_s(\ft,\fs) 
&=
\textstyle{\frac{1}{2}}\dim \sM_{\ft}
\quad\text{(By Equation \eqref{eq:SWDimRelations})}
\\
&=
\textstyle{\frac{1}{2}}\dim \sM_{\ft'} -3
\quad\text{(By Equation \eqref{eq:Transv})}
\\
&=
\delta+\eta -2
\quad\text{(By Equation \eqref{eq:DegreeAssumption})},
\end{align*}
we have, for $a$ as defined in \eqref{eq:ab1},
\begin{align*}
n''_s+n'_s-\delta+3
&=
\eta -d+1
\\
&=
a,
\end{align*}
agreeing with the claimed formula for $a$ in \eqref{eq:Defnab}.
Then, using
\begin{align*}
n'_s 
&= 
\textstyle{\frac{1}{2}}d_a(\ft) +\textstyle{\frac{1}{4}}(\chi+\si)
\quad\text{(By Equation \eqref{eq:NormalComponentDims})}
\\
&=
\textstyle{\frac{1}{2}}d_a(\ft')-4+\textstyle{\frac{1}{4}}(\chi+\si),
\end{align*}
and $d_s(\fs)=2d$, we have, for $b$ as defined in \eqref{eq:ab1},
\begin{align*}
\delta-4-n'_s-d
&=
\textstyle{\frac{1}{2}}(2\delta-d_a(\ft')-d_s(\fs))
-\textstyle{\frac{1}{4}}(\chi+\si)
\\
&=
b,
\end{align*}
agreeing with the claimed formula for $b$ in \eqref{eq:Defnab}.
This completes the proof of Theorem \ref{thm:LevelOne}.
\end{proof}

It remains to prove the combinatorial identity
\eqref{eq:CombinatorialFactor}.

\begin{lem}
For integers $A,M,N,d$ with $d\geq 0$, and $v=0,\dots,3$, and
$P^{a,b}_d(\zeta)$ the Jacobi polynomial \eqref{eq:DefineJacobiPolynomial},
we have
\begin{equation}
\label{eq:CombinatorialFactor1}
2^d P^{a,b}_d(0)
=
\sum_{i=0}^d\sum_{j=0}^{d-i}\sum_{k=0}^j
(-1)^{i+j}2^{d-j+k}
    \binom{A-v}{i}
    \binom{d+3-v-i-j}{d-i-j}
    \binom{M}{k}
    \binom{N}{j-k},
\end{equation}
where the constants $a$ and $b$ are given by
$$
a=3-N-A-M
\quad\text{and}\quad
b=A+M-4-d.
$$
\end{lem}

\begin{proof}
For $r\in\RR$ and $\ell\in \NN$ it is convenient to define 
$$
(r)_\ell:=r(r+1)\cdots (r+\ell-1).
$$
We then recall the following identities (see \cite[Lemma 4.16]{FL2b}):
\begin{equation}
\label{eq:BinomialToRising}
\binom{r}{\ell} =\frac{(-1)^{\ell} (-r)_{\ell}}{\ell !}
\quad\text{and}\quad
(r)_{\ell} = (-1)^{\ell} (1-r-\ell)_{\ell}.
\end{equation}
The preceding identities yield
\begin{equation}
\label{eq:BinomialReversal}
\begin{aligned}
\binom{d+3-v-i-j}{d-i-j}
&=
\frac{(-1)^{d-i-j}(-d-3+i+j+v)_{d-i-j}}{(d-i-j)!}
\\
&=
\frac{(4-v)_{d-i-j}}{(d-i-j)!}
\\
&=
(-1)^{d-i-j}
\binom{v-4}{d-i-j}.
\end{aligned}
\end{equation}
Let $\sS$ denote right-hand-side of equation
\eqref{eq:CombinatorialFactor1} and observe that, using equation
\eqref{eq:BinomialReversal}, we can write $\sS$ in the form
\begin{equation}
\label{eq:CombinatorialFactor2}
\sS
=
\sum_{i=0}^d
\sum_{j=0}^{d-i}
\sum_{k=0}^j
(-1)^d2^{d-j+k}
\binom{A-v}{i}
\binom{v-4}{d-i-j}
\binom{M}{k}\binom{N}{j-k}.
\end{equation}
With the substitution $u=j-k$, equation
\eqref{eq:CombinatorialFactor2} becomes
\begin{equation}
\label{eq:CombinatorialFactor3}
\begin{aligned}
\sS
&=
(-2)^d
\sum_{i=0}^d \sum_{u=0}^{d-i}\sum_{j=u}^{d-i}
2^{-u}
\binom{A-v}{i}
\binom{v-4}{d-i-j}
\binom{M}{j-u}
\binom{N}{u}
\\
&=
(-2)^d
\sum_{i=0}^d \sum_{u=0}^{d-i}
 2^{-u}
\binom{A-v}{i}
\binom{N}{u}
\sum_{j'=0}^{d-i-u}
\binom{v-4}{d-i-u-j'}
\binom{M}{j'}
\end{aligned}
\end{equation}
where we have set $j'=j-u$.  Applying the Vandermonde convolution identity
\cite[Lemma 4.16 (5)]{FL2b} to the sum over $j'$ in the right-hand side of
\eqref{eq:CombinatorialFactor3}, we see that 
\begin{equation}
\begin{aligned}
\label{eq:CombinatorialFactor4}
\sS
&=
(-2)^d
\sum_{i=0}^d \sum_{u=0}^{d-i}
 2^{-u}
\binom{A-v}{i}
\binom{N}{u}
\binom{v+M-4}{d-i-u}
\\
&=
(-2)^d
\sum_{u=0}^d  2^{-u}\binom{N}{u}
\sum_{i=0}^{d-u}
    \binom{A-v}{i}
    \binom{v+M-4}{d-i-u}
\quad\text{(Reversing summation order)}
\\
&=
(-2)^d
\sum_{u=0}^d  2^{-u}
    \binom{N}{u}
    \binom{A+M-4}{d-u}
    \quad\text{(By \cite[Lemma 4.16(5)]{FL2b})}.
\end{aligned}
\end{equation}
We now express this sum in terms of the hypergeometric
function \cite[\S 9.10]{GradshteynRyzhik6}
$$
\HG(a,b;c;\zeta)
=
\sum_{k=0}^\infty \frac{(a)_k(b)_k}{(c)_k k!}\zeta^k,
\quad
\zeta\in\CC.
$$
Using the first identity in \eqref{eq:BinomialToRising},
equation \eqref{eq:CombinatorialFactor4} yields
\begin{equation}
\label{eq:CombinatorialFactor5}
\sS
=
(-2)^d \sum_{u=0}^d 2^{-u}
\frac{(-1)^{u}(-N)_{u}}{u!} \frac{(-1)^{d-u}(4-A-M)_{d-u}}{(d-u)!}.
\end{equation}
Then applying the identities (see \cite[Lemma 4.16]{FL2b})
\begin{equation}
(a)_{d-u}= \frac{(-1)^u (a)_d}{(1-a-d)_u}
\quad\text{and}\quad
(d-u)!=\frac{d!}{(-1)^u(-d)_u}
\end{equation}
to equation \eqref{eq:CombinatorialFactor5} gives
\begin{equation}
\label{eq:CombinatorialFactor6}
\begin{aligned}
\sS
&=
2^d
\sum_{u=0}^d 2^{-u}
\frac{(-N)_u (-1)^u(4-A-M)_d(-1)^u(-d)_u}{u!(A+M-3-d)_u d!}
\\
&=
\frac{2^d (4-A-M)_d}{d!}
\sum_{r=0}^d  \frac{(-d)_u(-N)_u}{(A+M-3-d)_u u!} 2^{-u}
\\
&=
\frac{2^d (4-A-M)_d}{d!}
\HG(-d,-N;A+M-3-d;\textstyle{\frac{1}{2}})
\\
&=
\frac{(-2)^d(4-A-M)_d}{(A+M-3-d)_d}
P^{3-N-A-M,A+M-4-d}_d(0)
\quad\text{(By \cite[Equation (4.40)]{FL2b})}.
\end{aligned}
\end{equation}
Finally, we note that
$$
\frac{(-2)^d(4-A-M)_d}{(A+M-3-d)_d}=2^d
$$
by the second identity in \eqref{eq:BinomialToRising}.
This completes the proof of the lemma.
\end{proof}

\subsection{The blow-up formula for level-one Seiberg-Witten link pairings}
\label{subsec:BlowUp}
To apply Theorem \ref{thm:LevelOne} and express the link pairing
\eqref{eq:LevelOne} in terms of Seiberg-Witten invariants of $X$ and thus
express, via \eqref{eq:CobordismSum}, the Donaldson invariants of $X$ in
terms of Seiberg-Witten invariants of $X$, we must apply the
$\SO(3)$-monopole cobordism to the blow-up, $\tilde
X=X\#\overline{\CC\PP}^2$, of $X$ in view of our definitions
\eqref{eq:DefineDonaldson} and \eqref{eq:DefSW} of the four-manifold
invariants. Our definitions of the gauge-theoretic four-manifold invariants
incorporate the blow-up formula because the classes $w+\PD[e]\pmod{2}$ and
$c_1(\fs)\pm\PD[e]-\Lambda\pmod{2}$ are always good in the sense of
Definition \ref{defn:Good}, for arbitrary $w$, $c_1(\fs)$, and $\Lambda$.

Following the discussion in \cite[\S 4.5]{FL2b}, every
\spinc structure on $\tilde X$ is given by $\fs\#\fs_{2k-1}$, the result of
splicing a \spinc structure $\fs$ on $X$ and the \spinc structure
$\fs_{2k-1}$ on $\overline{\CC\PP}^2$ with $c_1(\fs_{2k-1})=(2k-1)\PD[e]$,
$k\in \ZZ$. Then,
\begin{equation}
\label{eq:SpincBlowUps}
c_1(\fs\#\fs_{2k-1}) = c_1(\fs) + (2k-1)\PD[e].
\end{equation}
The dimensions of the Seiberg-Witten moduli spaces are related by
\begin{equation}
\label{eq:SWModuliDimensionBlowUp}
d_s(\fs\#\fs_{2k-1}) = d_s(\fs) -k(k-1).
\end{equation}
For the cases $k=1$ or $0$, we denote $\fs^+=\fs\#\fs_1$ and
$\fs^-=\fs\#\fs_{-1}$, so that
\begin{equation}
\label{eq:DefnPMspinc}
c_1(\fs^\pm)=c_1(\fs)\pm\PD[e].
\end{equation}
If $d_s(\fs\#\fs_{2k-1})\geq 0$, then $SW_X(\fs\#\fs_{2k-1})=SW_X(\fs)$
according to the blow-up formula \cite[Theorem 1.4]{FSTurkish},
\cite[Theorem 3.2]{OzsvathSzaboAdjunct} (see the statement given in
\cite[Theorem 4.20]{FL2b}). 

Recall from \cite[Lemma 4.19]{FL2b} that given a
\spinu structure $\ft$ on $X$, there is a \spinu structure $\tilde \ft$ on
$\tilde X$ such that
\begin{equation}
\label{eq:BlowUpSpinuChar}
c_1(\tilde\ft)=c_1(\ft),
\quad
p_1(\tilde \ft)=p_1(\ft) - 1,
\quad
w_2(\tilde \ft)\equiv w_2(\ft)+\PD[e]\pmod 2.
\end{equation}
We then have the following extension of Lemma 4.19 in \cite{FL2b}.

\begin{lem}
\label{lem:BlowUpSpinU}
Continue the notation of the preceding paragraph. Then the following hold
for each integer $\ell\geq 0$:
\begin{enumerate}
\item
There is a bijection between the
\begin{itemize}
\item
Set of strata $\iota(M_{\fs}(X))\times\Sym^\ell(X)$ in $I\sM_{\ft}(X)$,
and the
\item
Set of families of strata
$\iota(M_{\fs_{k^\pm}}(X))\times\Sym^{\ell-j}(X)$ in
$I\sM_{\tilde\ft}(\tilde X)$, where $\fs_{k^\pm}=\fs\#\fs_{2k^\pm-1}$ and
$j,k^\pm\in\ZZ_{\geq 0}$ obey $j=k^\pm(k^\pm-1)\leq \min\{\ell,d_s(\fs)\}$. 
\end{itemize}
\item
The preceding correspondence further restricts to a bijection between the
\begin{itemize}
\item
Set of strata in level $\ell$ with $SW_X(\fs)\neq 0$, and the
\item
Set of families of strata in levels $\ell-j$ with $SW_X(\fs_{k^\pm})\neq 0$,
where $j,k^\pm\in\ZZ_{\geq 0}$ obey $j=k^\pm(k^\pm-1)\leq
\min\{\ell,d_s(\fs)\}$. 
\end{itemize}
\item
If $X$ has SW-simple type then so has $\tilde X$ and the preceding
correspondence further restricts to a bijection between the
\begin{itemize}
\item
Set of strata in level $\ell$ with $SW_X(\fs)\neq 0$, and the
\item
Set of pairs of strata  in level $\ell$ with $SW_X(\fs^\pm)\neq 0$.
\end{itemize}
\end{enumerate}
\end{lem}

\begin{proof}
{}From \cite[Lemma 3.32]{FL2b} there is a splitting
$\ft_{\ell}=\fs\oplus\fs'$ if and only if
$(c_1(\fs)-c_1(\ft))^2=p_1(\ft)+4\ell$, recalling that $\ft_\ell$ is a
\spinu structure with $p_1(\ft_\ell)=p_1(\ft)+4\ell$, 
$c_1(\ft_\ell)=c_1(\ft)$, and $w_2(\ft_\ell)=w_2(\ft)$. Hence, $M_{\fs}$
is contained in the level $\sM_{\ft_\ell}\times\Sym^\ell(X)$ if and only if
\begin{equation}
\label{eq:Level}
\ell 
= 
\ell(\ft,\fs) 
:= 
\textstyle{\frac{1}{4}}((c_1(\fs)-c_1(\ft))^2 - p_1(\ft)).
\end{equation}
Substituting \eqref{eq:SpincBlowUps} and \eqref{eq:BlowUpSpinuChar} into
\eqref{eq:Level}, we see that 
$$
\ell(\tilde \ft,\fs\#\fs_{2k-1}) 
=
\ell(\ft,\fs)-k(k-1)
\leq
\ell(\ft,\fs).
$$
Assertion (1) follows from the preceding formula relating levels
and the relation \eqref{eq:SWModuliDimensionBlowUp} between the dimensions
of the Seiberg-Witten moduli spaces. Assertions (2) and (3) follow from the
blow-up formula \cite{FSTurkish, OzsvathSzaboAdjunct} and the
definition of SW-simple type.
\end{proof}

\begin{rmk}
\label{rmk:BlowUpSpinULevel01}
If $\ell=0$ or $1$ in Cases (1) or (2) of Lemma
\ref{lem:BlowUpSpinU}, then we must have $j=0$, $k^+=1$, $k^-=0$, and
$\ell(\ft,\fs) = \ell(\tilde\ft,\fs^\pm)$. Hence, for $\ell\in\{0,1\}$ and
$\fs\in\Spinc(X)$, there is an injective map from the 
\begin{itemize}
\item
Set of strata, $M_{\fs}\times\Sym^\ell(X)$, in $I\sM_{\ft}$ with
$SW_X(\fs)\neq 0$, to the
\item
Set of pairs of strata, $M_{\fs^\pm}\times\Sym^\ell(\tilde X)$,
in $I\sM_{\tilde\ft}$ with $SW_X(\fs^\pm)\neq 0$.
\end{itemize}
Furthermore, if all $\fs\in\Spinc(X)$ with $\ell(\ft,\fs)\geq 2$ have
$SW_X(\fs)=0$, then the preceding map is bijective and all
$\tilde\fs\in\Spinc(\tilde X)$ with $\ell(\tilde\ft,\tilde\fs)\geq 2$ have
$SW_{\tilde X}(\tilde\fs)=0$.
\end{rmk}

\begin{prop}
\label{prop:LevelOneBlowUp}
Continue the hypotheses of Theorem \ref{thm:LevelOne} leading to equation
\eqref{eq:LevelOne}, except we omit the requirement that $M_{\fs}$ contains
no zero-section pairs.
Let the links $\bL_{\tilde\ft',\fs^\pm}$ have the standard orientation.
Then, for $k$ even,
\begin{equation}
\label{eq:LevelOneBlowUp}
\begin{aligned}
{}&(-1)^{o_{\tilde \ft'}(w+\PD[e],\fs^+)}\
\#
\left(\bar\sV(h^{\delta-2m-k}e^{k+1}x^m)\cap \bar\sW^{\eta}
\cap \bL_{\tilde \ft',\fs^+}\right)
\\
+&(-1)^{o_{\tilde \ft'}(w+\PD[e],\fs^-)}\
\#
\left(\bar\sV(h^{\delta-2m-k}e^{k+1}x^m)\cap \bar\sW^{\eta}
\cap \bL_{\tilde \ft',\fs^-}\right)
\\
&=
(-1)^{o_{\ft}(w,\fs)+m+1-d_s(\fs)/2}
2^{-\delta}2^{d_s(\fs)/2}P^{a,b}_{d_s(\fs)/2}(0) 
SW_X(\fs)
\\
&\quad\times\left(
\tilde a_0\langle c_1(\fs)-c_1(\ft'),h\rangle^{\delta-2m-k}
\right.
\\
&\quad 
+
\tilde b_0
\langle c_1(\fs)-c_1(\ft'),h\rangle^{\delta-2m-k-1}\langle c_1(\ft'),h\rangle
\\
&\quad\left.
+
\tilde a_1\langle c_1(\fs)-c_1(\ft'),h\rangle^{\delta-2m-k-2}Q_X(h,h)
\right),
\end{aligned}
\end{equation}
where the sign change factor $o_{\ft}(w,\fs)$ is defined
by \eqref{eq:OrientChangeFactor}, the coefficients $\tilde a_0,\tilde
b_0,\tilde a_1$ are related to those of Theorem \ref{thm:LevelOne} by
$$
\tilde a_0 = a_0-4\binom{k+1}{2},
\quad
\tilde b_0 = b_0-2k\frac{P^{a-1,b+1}_{d_s(\fs)/2}(0)}{P^{a,b}_{d_s(\fs)/2}(0)},
\quad\text{and}\quad
\tilde a_1 = 4\binom{\delta-2m-k}{2},
$$
and $P^{a,b}_{d_s(\fs)/2}(0)$ is given by
\eqref{eq:DefineJacobiPolynomial}. If $k$ is 
odd, the sum on the left-hand side of \eqref{eq:LevelOneBlowUp} is zero.
\end{prop}

\begin{proof}
{}From \cite[Equation (4.55)]{FL2b} we see that
$o_{\tilde\ft'}(w+\PD[e],\fs^+)=o_{\ft'}(w,\fs)-1$ and that
$o_{\tilde\ft'}(w+\PD[e],\fs^-)=o_{\ft'}(w,\fs)$, so the left-hand side of
\eqref{eq:LevelOneBlowUp} is $(-1)^{o_{\ft'}(w,\fs)-1}$ times
\begin{equation}
\label{eq:BlowUpOrientationSum}
\#\left(\bar\sV(h^{\delta-2m-k}e^{k+1}x^m)\cap \bar\sW^{\eta}
\cap \bL_{\tilde \ft',\fs^+}\right)
-
\#\left(\bar\sV(h^{\delta-2m-k}e^{k+1}x^m)\cap \bar\sW^{\eta}
\cap \bL_{\tilde \ft',\fs^-}\right),
\end{equation}
where the links $\bL_{\tilde\ft',\fs^\pm}$ have the standard orientation.
Suppose $h_i=h$, for $0\leq i\leq \delta -2m-k$
and $h_i=e$ for $\delta-2m-k< i\le \delta+1-2m$.
Applying the polarization identity \cite[p. 396]{FrM} to the formula
\eqref{eq:LevelOne} to compute the pairings
$$
\#
\left(\bar\sV(h^{\delta-2m-k}e^{k+1}x^m)\cap \bar\sW^{\eta}
\cap \bL_{\tilde \ft',\fs^\pm}\right),
$$
and noting that $\dim\sM_{\tilde\ft'}(\tilde X)
=\dim\sM_{\ft'}(X)+2$ by equation \eqref{eq:Transv}, gives
(using $d=d_s(\fs)/2=d_s(\fs^\pm)/2$ for brevity)
\begin{equation}
\label{eq:PolarizedLinkPairing}
\begin{aligned}
{}&
\#
\left(\bar\sV(h^{\delta-2m-k}e^{k+1}x^m)\cap \bar\sW^{\eta}
\cap \bL_{\tilde \ft',\fs^\pm}\right)
\\
&=
(-1)^{m+1+d} 2^{-\delta-1}
2^dP^{a^\pm,b^\pm}_d(0)
\langle\mu_{\fs^\pm}^d,M_{\fs^\pm}(\tilde X)\rangle
\left(
a_0^\pm \prod_{i=1}^{\delta+1-2m}\langle c_1(\fs^\pm)-c_1(\tilde\ft'),h_i\rangle
\right.
\\
&\quad
+
2\frac{P^{a^\pm-1,b^\pm+1}_d(0)}{P^{a^\pm,b^\pm}_d(0)}\sum_{j=1}^{\delta+1-2m}
    \prod_{\substack{i=1\\i\neq j}}^{\delta+1-2m}
\langle c_1(\fs^\pm)-c_1(\tilde\ft'),h_i\rangle
                    \langle c_1(\tilde \ft'),h_j\rangle
\\
&\quad
\left.
+4\sum_{1\leq j<l\leq\delta+1-2m}
\prod_{\substack{i=1\\i\neq j,l}}^{\delta+1-2m}
\langle c_1(\fs^\pm)-c_1(\tilde\ft'),h_i\rangle
                                Q_{\tilde X}(h_j,h_l)
\right),
\end{aligned}
\end{equation}
where, using $c_1(\fs^\pm)=c_1(\fs)\pm\PD[e]$ and
$c_1^2(\tilde X)=c_1^2(X)-1$, and $c_1(\tilde\ft')=c_1(\ft')$,
\begin{equation}
\label{eq:BlowUpCoeff1}
\begin{aligned}
a_0^\pm &= 3\left(c_1(\fs^\pm)-c_1(\tilde \ft')\right)^2 +c_1^2(\tilde X)
        +2\left( c_1(\fs^\pm)-c_1(\tilde \ft')\right)\cdot c_1(\tilde\ft')
\\
&\qquad        +4(\delta+1-2m) -4m
        \\
        &=3( c_1(\fs)-c_1(\ft'))^2 -3 +c_1^2(X)-1
        +2( c_1(\fs)-c_1(\ft'))\cdot c_1(\ft')
\\
&\qquad        +4(\delta-2m) -4m +4
        \\
        &=3( c_1(\fs)-c_1(\ft'))^2  +c_1^2(X)
        +2( c_1(\fs)-c_1(\ft'))\cdot c_1(\ft')
\\
&\qquad        +4(\delta-2m) -4m,
\end{aligned}
\end{equation}
and where $a^\pm,b^\pm$ are given by
\begin{equation}
\label{eq:BlowUpCoeffPab}
\begin{aligned}
a^\pm
&= \eta-d +1
\\
&= a,
\\
b^{\pm}
&= 
\textstyle{\frac{1}{2}}(2(\delta+1) -d_a(\tilde\ft') -d_s(\fs^\pm))
-\textstyle{\frac{1}{4}}(\chi(\tilde X)+\si(\tilde X))
\\
&=
\textstyle{\frac{1}{2}}(2\delta -d_a(\ft')-d_s(\fs))
-\textstyle{\frac{1}{4}}(\chi+\si) 
\\
&= b.
\end{aligned}
\end{equation}
Note that the factor of
$(\delta+1-2m)$ in the
coefficient $b_0/2$ and the
coefficient $a_1/4$ in equation \eqref{eq:LevelOne}
are absorbed by the application of the polarization identity
\cite[p. 396]{FrM}. Using $c_1(\tilde\ft')=c_1(\ft')$ and the identities
\begin{align*}
Q_{\tilde X}(e,e) &= -1, 
\\
Q_{\tilde X}(e,h) &= 0,
\\
\langle c_1(\fs^\pm)-c_1(\tilde \ft'),e\rangle &= \mp 1,
\\
\langle c_1(\fs^\pm)-c_1(\tilde \ft'),h\rangle
&=
\langle c_1(\fs)-c_1(\ft'),h\rangle,
\\
\langle c_1(\tilde \ft'),e\rangle &=0,
\\
\langle c_1(\tilde \ft'),h\rangle  &= \langle c_1(\ft'),h\rangle,
\end{align*}
we simplify the terms on the right-hand side of equation
\eqref{eq:PolarizedLinkPairing}
(recall that $h_i=h\in H_2(X;\RR)$, $1\leq i\leq
\delta-k-2m$, and $h_i=e$ for $\delta-k-2m<i\leq \delta+1-2m$) to give:
\begin{equation}
\label{eq:Prod1}
\prod_{i=1}^{\delta+1-2m}
\langle c_1(\fs^\pm)-c_1(\tilde\ft'),h_i\rangle
=
(\mp 1)^{k+1}\langle c_1(\fs)-c_1(\ft'),h\rangle^{\delta-k-2m},
\end{equation}
and 
\begin{equation}
\label{eq:Prod2}
\begin{aligned}
{}&\sum_{j=1}^{\delta+1-2m}
    \prod_{\substack{i=1\\ i\neq j}}^{\delta+1-2m}
\langle c_1(\fs^\pm)-c_1(\tilde\ft'),h_i\rangle
                    \langle c_1(\tilde \ft'),h_j\rangle
\\
&\quad = (\mp 1)^{k+1}(\delta-2m-k)
    \langle c_1(\fs)-c_1(\ft'),h\rangle^{\delta-k-2m-1}
    \langle c_1(\ft'),h\rangle,
\end{aligned}
\end{equation}
and
\begin{equation}
\label{eq:Prod3}
\begin{aligned}
{}&\sum_{1\leq j<l\leq \delta+1-2m}
\prod_{\substack{i=1\\ i\neq j,l}}^{\delta+1-2m}
\langle c_1(\fs^\pm)-c_1(\tilde\ft'),h_i\rangle Q_{\tilde X}(h_j,h_l)
\\
&=
(\mp 1)^{k+1}\binom{\delta-2m-k}{2}
    \langle c_1(\fs)-c_1(\ft'),h\rangle^{\delta-k-2m-2}Q_X(h,h)
\\
&\quad    +
(\mp 1)^{k-1}\binom{k+1}{2}
    \langle c_1(\fs)-c_1(\ft'),h\rangle^{\delta-k-2m}
Q_{\tilde X}(e,e).
\end{aligned}
\end{equation}
Thus, if $k+1$ is {\em even\/}, the terms
\eqref{eq:Prod1}, \eqref{eq:Prod2}, and \eqref{eq:Prod3}
will {\em cancel\/} when we substitute them into equation
\eqref{eq:PolarizedLinkPairing} and subtract the results in
\eqref{eq:BlowUpOrientationSum}, so in this case the difference 
\eqref{eq:BlowUpOrientationSum} and hence the 
pairing \eqref{eq:LevelOneBlowUp} are zero, as claimed. If $k+1$ is {\em
odd,\/} the terms
\eqref{eq:Prod1}, \eqref{eq:Prod2}, and \eqref{eq:Prod3}
will {\em add\/} when we substitute them into equation
\eqref{eq:PolarizedLinkPairing}  and subtract the results in
\eqref{eq:BlowUpOrientationSum}. Hence, equation \eqref{eq:LevelOneBlowUp}
follows from \eqref{eq:BlowUpOrientationSum} together with the preceding
substitutions and substitution of equation
\eqref{eq:BlowUpCoeff1} for
the coefficient $a_0^\pm$ and equation \eqref{eq:BlowUpCoeffPab} for
$a^\pm,b^\pm$. Note that the factor of $2$ obtained by adding like
terms and the coefficient $2^{-\delta-1}$ in
\eqref{eq:PolarizedLinkPairing} yields the coefficient $2^{-\delta}$ in
\eqref{eq:LevelOneBlowUp} and that the factors $(-1)^{k+1}=-1$
and $(-1)^{o_{\ft'}(w,\fs)-1}$ (mentioned before
\eqref{eq:BlowUpOrientationSum}) yield the factor $(-1)^{o_{\ft'}(w,\fs)}$
appearing on the right-hand-side of \eqref{eq:LevelOneBlowUp}. 
\end{proof}

\subsection{Proofs of main theorems}
\label{subsec:WCL1}
We now apply the computation of Theorem \ref{thm:LevelOne} to the sum
\eqref{eq:CobordismSum} to prove Theorem \ref{thm:SumFormula} and hence
Theorem \ref{thm:WCL1}. 

\begin{proof}[Proof of Theorem \ref{thm:Main}]
We shall derive equation \eqref{eq:Main} from the basic cobordism identity
\eqref{eq:CobordismSum} applied to the blow-up $\tilde
X=X\#\overline{\CC\PP}^2$, the definition \eqref{eq:DefineDonaldson} for the
Donaldson invariant, the blow-up formula
\eqref{eq:LevelOneBlowUp} for level-one Seiberg-Witten link pairings, and
the formula from
\cite[Equation (4.52)]{FL2b} for level-zero Seiberg-Witten link pairings.

Since $\deg(z)=2\delta$ obeys condition
\eqref{mod8}, we can choose an integer $p$ such that $p\equiv w^2\pmod{4}$ and
\begin{equation}
\label{eq:DegreeOfDonaldson}
\delta=-p -\textstyle{\frac{3}{4}}(\chi+\sigma).
\end{equation}
According to the last paragraph of \cite[\S 2.1.3]{FL2a}, there is a
\spinu structure $\ft'$ over $X$ with $c_1(\ft')=\La$, $w_2(\ft')\equiv w\pmod
2$, and $p_1(\ft)=p$. Let $\tilde\ft'$ be the related \spinu structure over
$\tilde X=X\#\overline{\CC\PP}^2$ defined prior to Lemma
\ref{lem:BlowUpSpinU}. {}From equation \eqref{eq:Transv} we see that
\begin{equation}
\label{eq:DiracIndex1}
\begin{aligned}
n_a(\tilde\ft')
&=
\textstyle{\frac{1}{4}}(p_1(\tilde\ft') +\La^2-\si(\tilde X))
\\
&=
\textstyle{\frac{1}{4}}(p_1(\ft') +\La^2-\sigma)
\\
&=
\textstyle{\frac{1}{4}}(i(\La)-\delta)
\\
&=
n_a(\ft').
\end{aligned}
\end{equation}
Therefore, the hypothesis that $\delta<i(\La)$ is equivalent to
$n_a(\tilde \ft')=n_a(\ft')>0$ and hence the moduli space
$M_{\kappa+1/4}^{w+\PD[e]}(\tilde X)$ has real codimension greater than or
equal to two in $\sM_{\tilde\ft'}(\tilde X)$. Moreover,
$w_2(\tilde\ft')=w+\PD[e]\pmod{2}$ is good and so the cobordism
formula \eqref{eq:CobordismSum} applies.

In general, by Remark 3.36 in \cite{FL2b}, the Seiberg-Witten stratum
$M_{\fs}\times\Sym^\ell(X)$ corresponding to a splitting
$\ft_\ell'=\fs\oplus\fs'$ lies in level
$$
\ell(\ft,\fs)=\textstyle{\frac{1}{4}}(\delta-r(\Lambda,c_1(\fs))),
$$
of the space of ideal $\SO(3)$ monopoles $I\sM_{\ft'}(X)$,
since $2\delta=d_a(\ft')=\dim M_\kappa^w(X)$.
By hypothesis, $\delta=r(\Lambda)+4$ and $r(\Lambda)\leq r(\Lambda,c_1(\fs))$
for all $\fs$ with $SW_X(\fs)\neq 0$ (by definition
\eqref{eq:SWInTopLevelFunction} of $r(\Lambda)$), so our choice of $\delta$
implies that
$$
0\leq \ell(\ft,\fs) \leq 1,
$$
if $M_{\fs}\times\Sym^\ell(X)\subset I\sM_{\ft'}$ and $SW_X(\fs)\neq 0$.

By hypothesis, $X$ is `effective' in the sense of Definition
\ref{defn:Effective} and so the only non-zero Seiberg-Witten contributions to
$D_X^w(h^{\delta-2m}x^m)$ arise from moduli spaces $M_{\fs}(X)$ with
$SW_X(\fs)\neq 0$ contained in levels $\ell(\ft,\fs)=0$ or $1$ of
$I\sM_{\ft'}(X)$ or, equivalently, from moduli spaces
$M_{\fs^\pm}(\tilde X)$ contained in levels $\ell(\tilde\ft,\fs^\pm)=0$ or $1$
of $I\sM_{\tilde\ft'}(\tilde X)$, again when $SW_X(\fs)\neq 0$
by Remark \ref{rmk:BlowUpSpinULevel01}.

Therefore, the cobordism identity \eqref{eq:CobordismSum}
applied to the moduli space $\bar\sM_{\tilde\ft'}(\tilde X)/S^1$
and the
definition \eqref{eq:DefineDonaldson} of $D_X^w(h^{\delta-2m}x^m)$ give
\begin{equation}
\label{eq:Sum1}
\begin{aligned}
{}&D^w_X(h^{\delta-2m}x^m)
\\
&=
-2^{1-n_a(\ft')}\sum_{\fs\in\Spinc(X)}
\left(
    (-1)^{o_{\tilde\ft}(w+\PD[e],\fs^+)}
    \#\left(
\bar\sV(h^{\delta-2m}ex^m)\cap\bar\sW^{n_a(\ft)-1}\cap\bL_{\tilde\ft',\fs^+}
\right)
\right.
\\
&\qquad  \left.  +
    (-1)^{o_{\tilde\ft}(w+\PD[e],\fs^-)}
   \#\left(
\bar\sV(h^{\delta-2m}ex^m)\cap\bar\sW^{n_a(\ft)-1}\cap\bL_{\tilde\ft',\fs^-}
\right)
\right).
\end{aligned}
\end{equation}
By the remarks in the preceding paragraphs the sum over \Spinc(X) has
potentially non-zero terms when 
\begin{itemize}
\item
$\ell(\ft,\fs)=0$ or, equivalently,
$r(\Lambda,c_1(\fs))=\delta=r(\Lambda)+4$, and 
\item
$\ell(\ft,\fs)=1$ or, equivalently,
$r(\Lambda,c_1(\fs))=\delta-4=r(\Lambda)$.
\end{itemize}
Recall that the links $\bL_{\tilde\ft,\fs^\pm}$ are empty by definition
when $\ell(\tilde\ft',\fs^\pm)<0$. Substituting
\cite[Equation (4.52)]{FL2b} into equation \eqref{eq:Sum1}
to compute the terms with $\ell(\ft,\fs)=0$ and substituting the blow-up
formulas \eqref{eq:LevelOneBlowUp} into equation \eqref{eq:Sum1} to compute
the terms with $\ell(\ft,\fs)=1$ yields
\begin{equation}
\label{eq:Sum2}
\begin{aligned}
D^w_X(h^{\delta-2m}x^m)
&=
2^{1-n_a(\ft')-\delta}
\sum_{\substack{\fs\in\Spinc(X)\\ r(\Lambda,c_1(\fs))=\delta}}
(-1)^{o_{\ft}(w,c_1(\fs))+m+d_s(\fs)/2}
\\
&\qquad
\times
2^{d_s(\fs)/2}
P^{a-1,b}_{d_s(\fs)/2}(0)SW_X(\fs)
\langle c_1(\fs)-\La,h\rangle^{\delta-2m}
\\
&+
2^{1-n_a(\ft')-\delta}
\sum_{\substack{\fs\in\Spinc(X)\\ r(\Lambda,c_1(\fs))=\delta-4}}
(-1)^{o_{\ft}(w,c_1(\fs))+m+d_s(\fs)/2}
\\
&\qquad
\times
2^{d_s(\fs)/2}
P^{a,b}_{d_s(\fs)/2}(0)SW_X(\fs)
\left(
a_0\langle c_1(\fs)-\La,h\rangle^{\delta-2m}
\right.
\\
&\qquad
+ b_0\langle c_1(\fs)-\La,h\rangle^{\delta-2m-1}\langle \La,h\rangle
\\
&\qquad
\left.
+ a_1\langle c_1(\fs)-\La,h\rangle^{\delta-2m-2}Q_X(h,h)
\right),
\end{aligned}
\end{equation}
where the coefficients $a_0,b_0,a_1$ are as given in the statements
of Theorems \ref{thm:LevelOne} and
\ref{thm:SumFormula} (since they coincide with those of Proposition
\ref{prop:LevelOneBlowUp} when $k=0$).
{}From equation \eqref{eq:DiracIndex1} for $n_a(\ft')$, we see that
\begin{equation}
\label{eq:InvariantFormulaForna}
\begin{aligned}
1-n_a(\ft')-\delta
&=
1-\textstyle{\frac{1}{4}}(i(\La)-\delta)-\delta
\\
&=
1-\textstyle{\frac{1}{4}}i(\La)-\textstyle{\frac{3}{4}}\delta,
\end{aligned}
\end{equation}
and so the power of $2$ in equation \eqref{eq:Sum2} matches that in
equation \eqref{eq:Main}. Finally, via \cite[Equation (4.62)]{FL2b}, we see
that 
$$
m+o_{\ft}(w,\fs)
\equiv
m+\textstyle{\frac{1}{2}}(\si-w^2) 
+ \textstyle{\frac{1}{2}}(w^2+c_1(\fs)\cdot(w-\Lambda))
\pmod{2},
$$
and so the power of $(-1)$ in equation \eqref{eq:Sum2} also matches that in
equation \eqref{eq:Main}. Thus, substituting the preceding formulas for
powers of $2$ and $(-1)$ into equation
\eqref{eq:Sum2} yields equation \eqref{eq:Main}. 

Lastly, we simplify the expressions given in Theorem \ref{thm:LevelOne} for
the constants $a$ and $b$ in $P^{a,b}_{d_s(\fs)/2}(0)$. {}From equation
\eqref{eq:Defnab} and the fact that $d_a(\ft')=2\delta$, we have
$$
b = -\textstyle{\frac{1}{2}}d_s(\fs)-\textstyle{\frac{1}{4}}(\chi+\sigma).
$$
Similarly, equation \eqref{eq:DegreeAssumption} for $\eta$,
the hypothesis that $\delta=r(\Lambda)+4$, and equation
\eqref{eq:Transv} for $\dim\sM_{\ft'}$ yield
\begin{align*}
\eta
&=
\textstyle{\frac{1}{2}}(\dim(\sM_{\ft'}/S^1)-1)-\delta
\\
&=
\textstyle{\frac{1}{2}}(d_a(\ft')+2n_a(\ft')-2)-\delta
\\
&=
n_a(\ft')-1,
\end{align*}
so equation \eqref{eq:DiracIndex1} for $n_a(\ft')$ then implies that
\begin{equation}
\label{eq:EqnForDelta_c}
\eta
= 
\textstyle{\frac{1}{4}}(i(\Lambda)-\delta)-1,
\end{equation}
and thus equation \eqref{eq:Defnab} gives
$$
a 
= 
\textstyle{\frac{1}{4}}(i(\Lambda)-\delta) 
- \textstyle{\frac{1}{2}}d_s(\fs).
$$
This completes the proof of Theorem \ref{thm:Main}.
\end{proof}

\begin{proof}[Proof of Theorem \ref{thm:SumFormula}]
As $\Lambda\in B^\perp$ and $X$ has SW-simple type by
hypothesis, equations \eqref{eq:PositiveDiracIndexFunction} and
\eqref{eq:SimpleTypeBperpformR} imply that 
$$
i(\La)+r(\La) = 2c(X).
$$
Since $\delta=r(\La)+4$, this gives
\begin{equation}
\label{eq:SimpleILambda}
i(\La)=2c(X)-r(\La)=2c(X)-\delta+4.
\end{equation}
In particular, $r(\Lambda,c_1(\fs))=r(\Lambda)$ and thus $\ell(\ft,\fs)=1$
for all $\fs\in\Spinc(X)$ with $SW_X(\fs)\neq 0$.

We simplify the expression for the power of $2$ in equation
\eqref{eq:Main}, noting that we now have $d_s(\fs)=0$ 
whenever $SW_X(\fs)\neq 0$ since $X$ has SW-simple type.  Substituting the
equation \eqref{eq:SimpleILambda} for $i(\Lambda)$ yields
$$
1-\textstyle{\frac{1}{4}}i(\La)-\textstyle{\frac{3}{4}}\delta =
-\textstyle{\frac{1}{2}}(c(X)+\delta),
$$
matching the power of $2$ in equation \eqref{eq:LevelOneSum}.

The Jacobi polynomial constants $P_{d_s(\fs)/2}^{a,b}$ are equal to $1$
when $d_s(\fs)=0$, irrespective of the values of $a$, $b$.  The assertions
concerning the power of $(-1)$ and the coefficients $a_0,b_0,a_1$ follow
immediately from the fact that $c_1(\fs)\in B$ and $\Lambda\in B^\perp$.
This completes the proof of Theorem \ref{thm:SumFormula}.
\end{proof}

In addition to the hypotheses of Theorem \ref{thm:SumFormula}, 
the hypotheses of Theorem \ref{thm:WCL1} requires that $X$ be abundant
in the sense of \S \ref{sec:Introduction}.

\begin{lem}
\label{lem:Abundance}
\cite[Lemma 2.2]{FKLM}
If $X$ is an abundant four-manifold, then the following hold:
\begin{enumerate}
\item
There are classes $\La_0,\La_1\in B^\perp$ with
$\Lambda_0\equiv\Lambda_1\pmod{2}$ such that
$\La_0^2=-(\chi+\sigma)$ and $\La_1^2=4-(\chi+\sigma)$.
\item
        There is a class $\Lambda \in 2 B^{\bot}$ with
    $\Lambda^{2}=-(\chi+\sigma)$ if $-(\chi+\sigma)\equiv
    0\pmod 8$ and $\Lambda^{2}=4-(\chi+\sigma)$ if
    $-(\chi+\sigma)\equiv 4\pmod 8$.
\end{enumerate}
\end{lem}

\begin{proof}
For Assertion (1), the abundance condition implies that there are classes
$e_1,e_2\in B^\perp$ with $e_{1}\cdot e_{1} = e_{2}\cdot e_{2}=0$ and
$e_1\cdot e_2=1$.  Let $t=\quarter(\chi+\sigma)$.  Then the classes
$\La_0=e_1-2te_2$ and $\La_1=e_1+(2-2t)e_2$ will do.

For Assertion (2), the class $\Lambda_{2}= 2e_{1}-te_{2}$ will satisfy the
conclusion when $t\equiv 0 \pmod 2$, so we can take $\Lambda=\Lambda_2$,
while $\Lambda_{3}=2e_{1}+(1-t)e_{2}$ will work when $t\equiv 1 \pmod 2$,
and we can take $\Lambda=\Lambda_3$.
\end{proof}

\begin{rmk}
Although the statement of Lemma \ref{lem:Abundance} is identical to that of
Lemma 2.2 in \cite{FKLM}, the definitions of the classes $\Lambda_0,
\Lambda_1$ here differs slightly from those of \cite{FKLM}. Contrary to the
assertion in the proof of Lemma 2.2 in \cite{FKLM}, the classes
$\Lambda_0,\Lambda_1$ defined in \cite{FKLM} do not satisfy $\La_0\equiv
\La_1\pmod 2$. The proof given here corrects that error in the proof of
\cite[Lemma 2.2]{FKLM}.

Assertion (1) is only used in the proof of
\cite[Theorem 1.1]{FKLM} --- when applying the sign-change formula for
Donaldson invariants --- and not in the proof of \cite[Theorem 1.3]{FKLM},
which is the only part of that article where Assertion (2) is used.

In the present article, Assertion (2) is not required and we only use
Assertion (1) to prove the existence
of the class $\La$ in  Theorem \ref{thm:WCL1}. 
\end{rmk}

\begin{proof}[Proof of Theorem \ref{thm:WCL1}]
The existence of $\La\in B^\perp$ and $w\in H^2(X;\ZZ)$ with
$\La^2=4-(\chi+\si)$ and $w-\La$ characteristic follows immediately
from Assertion (1) of Lemma \ref{lem:Abundance}.  Take $\La=\La_1$
and set $w=\La+v$ where $v$ is characteristic.

We shall apply Theorem \ref{thm:SumFormula} and the vanishing result in
\cite[Theorem 1.1]{FKLM}. By hypothesis, $\La\in B^\perp$ and so
$\Lambda^2\equiv 0\pmod{2}$ (see first paragraph of proof of Theorem 1.1 in
\cite[\S 4.6]{FL2b}). Since $\chi+\sigma$ is even, it is convenient to
write $\La^2=2j-(\chi+\si)$, where $j\in\ZZ$, so equations
\eqref{eq:PositiveDiracIndexFunction} and \eqref{eq:SimpleTypeBperpformR} 
would then take the form
$$
r(\La)=c(X)-2j \quad\text{and}\quad i(\La)=c(X)+2j.  
$$
The constraints $\delta\le r(\La)+4$ and $\delta< i(\La)$ in Theorem
\ref{thm:SumFormula} are thus equivalent to (see
\cite[Figure 1]{FKLM}):
$$
\delta\le c(X)-2j+4\quad{\text{and}}\quad \delta < c(X)+2j.
$$
Therefore, when using Theorem \ref{thm:SumFormula} to compute
$D_X^w(h^{\delta-2m}x^m)$, we can choose $\delta\in\NN$ no larger than
$\delta=c(X)$, with $\La^2=2-(\chi+\si)$ when $j=1$, or
$\La^2=4-(\chi+\si)$ when $j=2$; as we shall see from equation
\eqref{eq:DParityVanishing}, we have $D_X^w(h^{c(X)+1-2m}x^m)=0$, while
\eqref{eq:LowDegreeSWVanishing} shows that the term in $\bS\bW_X^w(h)$ of
degree $c(X)+1$ in $h$ is also zero.

We first observe that both sides of the second identity in
\eqref{eq:WCL1Equation} vanish in
sufficiently low degree. Theorem 1.1 in \cite{FKLM} implies that, for $v$
characteristic, 
\begin{equation}
\label{eq:CharvLowDegreeSWVanishing}
\sum_{\fs\in \Spinc(X)}
(-1)^{\half(v^2+v\cdot c_1(\fs))}SW_X(\fs)
\langle c_1(\fs),h\rangle^d =0,
\end{equation}
if $d<c(X)-2$ or (see \cite[\S 2]{FKLM}) if $d\not\equiv c(X)\pmod 2$,
since the series has parity (using $v^2\equiv \sigma\pmod{8}$)
\begin{equation}
\label{eq:charwparity}
\begin{aligned}
-v^2-\textstyle{\frac{3}{4}}(\chi+\sigma)
&\equiv
-\sigma-\textstyle{\frac{3}{4}}(\chi+\sigma)\pmod{8}
\\
&\equiv
-\textstyle{\frac{1}{4}}(7\chi+11\sigma)\pmod{4}
\\
&=
c(X).
\end{aligned}
\end{equation}
By hypothesis of Theorem \ref{thm:WCL1}, $w-\Lambda$ is characteristic and
so we may write
\begin{equation}
\label{eq:Defnw}
w=v+\Lambda,
\end{equation}
for some $v\equiv w_2(X)\pmod{2}$. Since $\La\in B^\perp$, then
equation \eqref{eq:CharvLowDegreeSWVanishing} implies that
\begin{equation}
\label{eq:LowDegreeSWVanishing}
\sum_{\fs\in \Spinc(X)}
(-1)^{\half(w^2+w\cdot c_1(\fs))}SW_X(\fs)
\langle c_1(\fs),h\rangle^d =0,
\end{equation}
if $d<c(X)-2$ or if $d\not\equiv c(X)\pmod 2$. In particular, aside from
the fact that $w$ need not be characteristic, the vanishing result for the
Seiberg-Witten series in equation \eqref{eq:WCL1Equation} restates Theorem
1.1 in \cite{FKLM}. 

According to the hypothesis of Theorem \ref{thm:WCL1}, we have
\begin{equation}
\label{eq:LambdaSquare}
\Lambda^2 = 4-(\chi+\sigma),
\end{equation}
so $r(\Lambda)=c(X)-4$ and $i(\Lambda)=c(X)+4$.  Because $v$ is
characteristic we have $v\cdot\Lambda\equiv \Lambda^2\pmod{2}$ and as
$\Lambda^2\equiv 0\pmod{4}$, then $w^2\equiv (v+\La)^2\equiv
v^2\pmod{4}$. Therefore, equations
\eqref{eq:charwparity} and condition \eqref{mod8} show that
\begin{equation}
\label{eq:DParityVanishing}
D^{w}_X(h^{\delta-2m}x^m)=0,
\quad\text{when}\quad
\delta\not\equiv c(X) \pmod 4.
\end{equation}
(Note that the alternative solutions to the $r(\Lambda)$ and $i(\Lambda)$
constraints yielding $\delta\leq c(X)$, namely those with
$\Lambda^2=2-(\chi+\sigma)$, would yield $\delta\equiv c(X)+2\pmod{4}$ and
so for that choice of $\Lambda^2$ we could choose $\delta$ no larger than
$c(X)-2$.)  Hence, from definition \eqref{eq:DefineDonaldsonSeries}, the
potentially non-zero terms in the Donaldson series $\bD^{w}_X(h)$ take the
form
$$
\frac{1}{(c(X)-4i)!} D^{w}_X(h^{c(x)-4i})
\quad\text{and}\quad
\frac{1}{2(c(X)-4i-2)!} D^{w}_X(h^{c(x)-4i-2}x).
$$
Because $r(\La)=c(X)-4$, Theorem 1.4(a) in \cite{FL2b} gives
\begin{equation}
\label{eq:LowDegreeVanishing}
D^{w}_X(h^{\delta-2m}x^m)
=0 \quad\text{for $\delta< c(X)-4$},
\end{equation}
while Theorem 1.4(b) in \cite{FL2b} (with $\delta=c(X)-4$) yields
\begin{equation}
\label{eq:LowDegreeEquality}
\begin{aligned}
D^{w}_X(h^{c(X)-4-2m}x^m)
&=
2^{2-c(X)}(-1)^{m+1}\sum_{\fs\in\Spinc(X)}
(-1)^{\half (w^2+c_1(\fs)\cdot w)}
\\
&\quad\times
\SW_X(\fs)\langle c_1(\fs)-\La,h\rangle^{c(X)-4-2m}.
\end{aligned}
\end{equation}
Replacing the terms $\langle c_1(\fs)-\La,h\rangle^{c(X)-4-2m}$ in equation
\eqref{eq:LowDegreeEquality} by their binomial expansions and applying
equation \eqref{eq:LowDegreeSWVanishing} gives
\begin{equation}
\label{eq:HighestDegreeVanDInvts}
D^{w}_X(h^{c(X)-4})
=
0 
=
D^{w}_X(h^{c(X)-6}x).
\end{equation}
Equations \eqref{eq:DParityVanishing}, \eqref{eq:LowDegreeVanishing}, and
\eqref{eq:HighestDegreeVanDInvts} then imply that the
terms of $\bD^{w}_X(h)$ of degree less than $c(X)-2$ in $h$ are zero,
yielding the vanishing result for the Donaldson series stated in
\eqref{eq:WCL1Equation}.
{}From equation \eqref{eq:LowDegreeSWVanishing} the terms of
$\bS\bW^{w}_X(h)$ of degree less than $c(X)-2$ in $h$ are also zero.

Thus, to obtain the second identity in \eqref{eq:WCL1Equation},
it suffices to prove that
\begin{equation}
\label{eq:DInvarPoint}
D^{w}_X(h^{c(X)-2}x)
=
2^{3-c(X)}\sum_{\fs\in \Spin^c(X)}(-1)^{\half (w^2+w\cdot c_1(\fs))}
  SW_X(\fs) \langle c_1(\fs),h\rangle^{c(X)-2},
\end{equation}
and, noting that $\frac{1}{2}\frac{c(X)!}{(c(X)-2)!}=\binom{c(X)}{2}$,
\begin{equation}
\label{eq:DInvar}
\begin{aligned}
D^{w}_X(h^{c(X)})
&=
2^{2-c(X)}\sum_{\fs\in \Spin^c(X)}(-1)^{\half (w^2+w\cdot c_1(\fs))}
      SW_X(\fs)
\\
&\qquad\times      \left(
       \langle c_1(\fs),h\rangle^{c(X)}
       +
       \binom{c(X)}{2}\langle c_1(\fs),h\rangle^{c(X)-2}Q_X(h,h)
      \right).
\end{aligned}
\end{equation}
Since $\La^2=4-(\chi+\sigma)$, $\delta=c(X)$,
$c(X)=-\frac{1}{4}(7\chi+11\sigma)$, and $c_1^2(X)=2\chi+3\sigma$ by
definition \eqref{eq:Definec1Squared}, the coefficients in equation
\eqref{eq:LevelOneSum} for $D_X^{w}(h^{\delta-2m}x^m)$ (with
$w-\Lambda$ characteristic) become 
\begin{equation}
\label{eq:Coefficients2}
\begin{aligned}
a_0&= 4c_1^2(X) + \La^2 +4c(X)-12m
\\
&=4-12m,
\\
b_0&=2(c(X)-2m),
\\
a_1&=4\binom{c(X)-2m}{2}.
\end{aligned}
\end{equation}
According to equation \eqref{eq:LowDegreeSWVanishing}, when we expand the terms
$$
(-1)^{\half (w^2+c_1(\fs)\cdot w)}SW_X(\fs)\langle c_1(\fs)-\La,h\rangle^d
$$
in equation \eqref{eq:LevelOneSum} into a binomial sum 
of terms of the form
$$
(-1)^{\half (w^2+c_1(\fs)\cdot w)}SW_X(\fs)
(-1)^k\binom{d}{k}\langle c_1(\fs),h\rangle^{d-k}\langle \La,h\rangle^k
$$
and sum over $\fs\in \Spinc(X)$, the sums of terms with $d-k<c(X)-2$ or
$d-k\not\equiv c(X)\pmod 2$ will vanish. Therefore, applying equation
\eqref{eq:LowDegreeSWVanishing} to the right-hand side of equation
\eqref{eq:LevelOneSum} with $\delta=c(X)$ and $m=1$, only the terms of
degree $\delta-2$ in $h$ will be non-zero. Hence, only the term with
coefficient $a_0$ will be non-zero and equation
\eqref{eq:LowDegreeSWVanishing} implies that it only contributes
\begin{equation}
\label{eq:a0leftdeltaminus2}
a_0\langle c_1(\fs),h\rangle^{c(X)-2}.
\end{equation}
Since $\delta=c(X)$ and $a_0=-8$ by equation \eqref{eq:Coefficients2}, the
power of $2$ the right-hand side of equation \eqref{eq:LevelOneSum} reduces
to $-c(X)+3$. As $m=1$, and $a_0=-8$, and $w^2\equiv \si\pmod 8$ (as
$w\equiv w_2(X)\pmod{2}$), the power of $(-1)$ in the right-hand side of
\eqref{eq:LevelOneSum} reduces to $\frac{1}{2}(w^2+w\cdot
c_1(\fs))$. Therefore, by these remarks and observation
\eqref{eq:a0leftdeltaminus2} we obtain equation \eqref{eq:DInvarPoint} for
$D^w_X(h^{c(X)-2}x)$.

Next, applying equation \eqref{eq:LowDegreeSWVanishing} to the right-hand
side of equation \eqref{eq:LevelOneSum} with $\delta=c(X)$
and $m=0$, only the terms of degree $\delta$ in $h$ will be non-zero. In
particular, equation \eqref{eq:LowDegreeSWVanishing} implies that the term
in the right-hand side of equation
\eqref{eq:LevelOneSum} with coefficient $a_0$ contributes only
\begin{equation}
\label{eq:a0leftdelta}
a_0\left(
\langle c_1(\fs),h\rangle^{c(X)}
+\binom{c(X)}{2} \langle c_1(\fs),h\rangle^{c(X)-2}\langle\La,h\rangle^2
\right).
\end{equation}
Similarly the term on the right in equation
\eqref{eq:LevelOneSum} with coefficient $b_0$ contributes only
\begin{equation}
\label{eq:b0leftdelta}
-b_0(c(X)-1)\langle c_1(\fs),h\rangle^{c(X)-2}\langle\La,h\rangle^2.
\end{equation}
Finally, the term on the right in equation
\eqref{eq:LevelOneSum} with coefficient $a_1$ contributes only
\begin{equation}
\label{eq:a1leftdelta}
a_1\langle c_1(\fs),h\rangle^{c(X)-2}Q_X(h,h).
\end{equation}
Hence, as $a_0=4$, $b_0=2c(X)$, and $a_1=4\binom{c(X)}{2}$ by equation
\eqref{eq:Coefficients2} when $m=0$, noting that $\delta=c(X)$ and
$w^2\equiv \si\pmod 8$, and combining observations \eqref{eq:a0leftdelta},
\eqref{eq:b0leftdelta}, and \eqref{eq:a1leftdelta},
we see that equation \eqref{eq:LevelOneSum} yields
\begin{align*}
D^w_X(h^{c(X)})
&=
2^{-c(X)}\sum_{\fs\in \Spinc(X)}
(-1)^{\half(w^2+w\cdot c_1(\fs))}
SW_X(\fs)
\\
&\quad\times\left(
4\langle c_1(\fs),h\rangle^{c(X)}
+4\binom{c(X)}{2} \langle c_1(\fs),h\rangle^{c(X)-2}\langle\La,h\rangle^2
\right.
\\
&\quad
-2c(X)(c(X)-1)\langle c_1(\fs),h\rangle^{c(X)-2}\langle\La,h\rangle^2
\\
&\quad
\left.
+4\binom{c(X)}{2}\langle c_1(\fs),h\rangle^{c(X)-2}Q_X(h,h)
\right)
\\
&=
2^{2-c(X)}
\sum_{\fs\in \Spinc(X)}
(-1)^{\half(w^2+w\cdot c_1(\fs))}
SW_X(\fs)
\\
&\quad\times\left(
\langle c_1(\fs),h\rangle^{c(X)}
+
\binom{c(X)}{2}\langle c_1(\fs),h\rangle^{c(X)-2}Q_X(h,h)
\right),
\end{align*}
which proves equation \eqref{eq:DInvar} for $D^w_X(h^{c(X)})$.
Thus, by the remark preceding equation \eqref{eq:DInvarPoint}, this proves
the second identity in \eqref{eq:WCL1Equation} and completes the proof of
Theorem \ref{thm:WCL1}.
\end{proof}

\begin{rmk}
\label{rmk:ProblemWithW}
If one attempted to compute $D^w_X$, where $w$ is characteristic, instead
of $D^{w+\La}_X$, the requirement that $\La-w$ also be characteristic would
imply that $\La\equiv 0\pmod 2$ so $\La=2\La_{2}$ for some $\La_{2}\in
H^2(X;\ZZ)$.  If $B$ is non-empty, then $c_1(\fs)\cdot
\La=2(c_1(\fs)\cdot \La_{2})=0$ implies $\La_{2}^2\equiv 0\pmod 2$ so
$\La^2\equiv 0\pmod 8$.

If $\chi+\sigma\equiv 4\pmod 8$, there exists a class $\La\in 2B^\perp$
with $\La^2=4-(\chi+\sigma)$ by Lemma \ref{lem:Abundance}.
In this case, Theorem \ref{thm:WCL1} would hold for $D^w_X$ with
$w$ characteristic.  In the case $\chi+\sigma\equiv 0\pmod 8$,
the identity
$$
c(X)-r(\La) = \La^2 +(\chi+\sigma)
$$
and the requirement that $\La^2\equiv 0\pmod 8$
would imply that $c(X)\equiv r(\La)\pmod 8$.
Now the Seiberg-Witten moduli spaces appear in level $\ell$, where
$$
\ell=\textstyle{\frac{1}{4}}( \delta - r(\La)).
$$
Thus, if we wished to compute degree $\delta=c(X)$ Donaldson invariants in
the case $\chi+\sigma\equiv 0\pmod 8$, we would need $\ell\equiv 0\pmod 2$.  If
$\ell=0$, the relations $\delta\le r(\La)$ and $\delta < i(\La)$ in
\cite[Theorem 2.1]{FKLM} imply that $\delta < c(X)$.  The cases $\ell\ge 2$
have not yet been computed, though work is in progress
\cite{FL3, FL4, FL5} and 
some partial results have been reported by Kronheimer and Mrowka
\cite{MrowkaPrincetonMorseTalk}. In the meantime, if $\chi+\sigma\equiv 0\pmod
8$, we cannot yet compute Donaldson invariants $D^w_{X}$ of degree $c(X)$
with this methods of this article in full generality.
\end{rmk}


\ifx\undefined\bysame
\newcommand{\bysame}{\leavevmode\hbox to3em{\hrulefill}\,}
\fi


\begin{thebibliography}{10}

\bibitem{AtiyahBott}
M.~F. Atiyah and R.~Bott, {\em The moment map and equivariant cohomology},
  Topology {\bf 23} (1984), 1--28.

\bibitem{AtiyahJones}
M.~F. Atiyah and J.~D.~S. Jones, {\em Topological aspects of {Y}ang-{M}ills
  theory}, Comm. Math. Phys. {\bf 61} (1978), 97--118.

\bibitem{BPV}
W.~Barth, C.~Peters, and A.~Van de~Ven, {\em Compact complex surfaces},
  Springer, New York, 1984.

\bibitem{BerlineGetzlerVergne}
N.~Berline, E.~Getzler, and M.~Vergne, {\em Heat kernels and {D}irac
  operators}, Springer-Verlag, Berlin, 1992.

\bibitem{BT}
R.~Bott and L.~Tu, {\em Differential forms in algebraic topology}, Springer,
  New York, 1982.

\bibitem{BredonTopGeom}
G.~E. Bredon, {\em Topology and geometry}, Springer-Verlag, New York, 1997.

\bibitem{BrockertomDieck}
T.~Br{\"o}cker and T.~tom Dieck, {\em Representations of compact {L}ie groups},
  Springer, New York, 1985.

\bibitem{Brussee}
R.~Brussee, {\em The canonical class and the {$C^\infty$}-properties of
  {K\"a}hler surfaces}, New York J. Math. (electronic) {\bf 2} (1996),
  103--146, \mbox{arXiv:alg-geom/9503004}.

\bibitem{CaoZhouWC}
H.~Cao and J.~Zhou, {\em Equivariant cohomology and wall crossing formulas in
  {S}eiberg-{W}itten theory}, Math. Res. Lett. {\bf 5} (1998), 711--721,
\mbox{arXiv:math.DG/9804134}.

\bibitem{Dold}
A.~Dold, {\em Lectures on algebraic topology}, Springer-Verlag, Berlin, 1995.

\bibitem{DonOrient}
S.~K. Donaldson, {\em The orientation of {Y}ang-{M}ills moduli spaces and
  4-manifold topology}, J. Differential Geom. {\bf 26} (1987), 397--428.

\bibitem{DK}
S.~K. Donaldson and P.~B. Kronheimer, {\em The geometry of four-manifolds},
  Oxford Univ. Press, Oxford, 1990.

\bibitem{FeehanGenericMetric}
P.~M.~N. Feehan, {\em Generic metrics, irreducible rank-one {PU(2)} monopoles,
  and transversality}, Comm. Anal. Geom. {\bf 8} (2000), 905--967,
  \mbox{arXiv:math.DG/9809001}.

\bibitem{FKLM}
P.~M.~N. Feehan, P.~B. Kronheimer, T.~G. Leness, and T.~S. Mrowka, {\em {PU(2)}
  monopoles and a conjecture of {M}ari{\~n}o, {M}oore, and {P}eradze}, Math.
  Res. Lett. {\bf 6} (1999), 169--182, \mbox{arXiv:math.DG/9812125}.

\bibitem{FLGeorgia}
P.~M.~N. Feehan and T.~G. Leness, 
{\em {PU(2)} monopoles and relations between four-manifold
  invariants}, Topology Appl. {\bf 88} (1998), 111--145, 
\mbox{arXiv:dg-ga/9709022}.

\bibitem{FL1}
\bysame, {\em {PU(2)} monopoles. {I}: {R}egularity, {U}hlenbeck compactness,
  and transversality}, J. Differential Geom. {\bf 49} (1998), 265--410,
  \mbox{arXiv:dg-ga/9710032}.

\bibitem{FL2a}
\bysame, {\em {PU(2)} monopoles and links of top-level {S}eiberg-{W}itten
  moduli spaces}, J. Reine Angew. Math. {\bf 538} (2001), to appear,
  \mbox{arXiv:math.DG/0007190}.

\bibitem{FL2b}
\bysame, {\em {PU(2)} monopoles. {II}: {T}op-level {S}eiberg-{W}itten moduli
  spaces and {W}itten's conjecture in low degrees}, J. Reine Angew. Math. {\bf
  538} (2001), to appear, \mbox{arXiv:dg-ga/9712005}.

\bibitem{FL3}
\bysame, {\em {PU(2)} monopoles. {III}: {E}xistence of
  gluing and obstruction maps}, submitted to a print journal, 
\mbox{arXiv:math.DG/9907107}.

\bibitem{FL4}
\bysame, {\em {PU(2)} monopoles. {IV}: {S}urjectivity of gluing maps}, in
  preparation.

\bibitem{FL5}
\bysame, {\em {PU(2)} monopoles. {V}}, in preparation.

\bibitem{FSTurkish}
R.~Fintushel and R.~Stern, {\em Immersed spheres in 4-manifolds and the
  immersed {T}hom conjecture}, Turkish J. Math. {\bf 19} (1995), 145--157,
http://math.uci.edu/$\sim$rstern.

\bibitem{FrM}
R.~Friedman and J.~W. Morgan, {\em Smooth four-manifolds and complex surfaces},
  Springer, Berlin, 1994.

\bibitem{Fulton}
W.~Fulton, {\em Intersection theory}, Springer, Berlin, 1984.

\bibitem{Goettsche}
L.~G{\"o}ttsche, {\em Modular forms and {D}onaldson invariants for 4-manifolds
  with {$b^+=1$}}, J. Amer. Math. Soc. {\bf 9} (1996), 827--843,
  \mbox{arXiv:alg-geom/9506018}.

\bibitem{GraberPand}
T.~Graber and R.~Pandharipande, {\em Localization of virtual classes}, Invent.
  Math. {\bf 135} (1999), 487--518, \mbox{arXiv:alg-geom/9708001}.

\bibitem{GradshteynRyzhik6}
I.~S. Gradshteyn and I.~M. Ryzhik, {\em Table of integrals, series, and
  products}, sixth ed., Academic, San Diego, CA, 2000, Translated from the
  Russian, Translation edited and with a preface by Alan Jeffrey and Daniel
  Zwillinger.

\bibitem{GreenbergHarper}
M.~J. Greenberg and J.~R. Harper, {\em Algebraic topology, a first course},
  Benjamin/Cummings, Reading, Mass., 1981.

\bibitem{Kervaire}
M.~A. Kervaire, {\em Relative characteristic classes}, Amer. J. Math. {\bf 79}
  (1957), 517--558.

\bibitem{KotschickMorgan}
D.~Kotschick and J.~W. Morgan, {\em {SO(3)} invariants for four-manifolds with
  {$b^+=1$}, {II}}, J. Differential Geom. {\bf 39} (1994), 433--456.

\bibitem{KMStructure}
P.~B. Kronheimer and T.~S. Mrowka, {\em Embedded surfaces and the structure of
  {D}onaldson's polynomial invariants}, J. Differential Geom. {\bf 43} (1995),
  573--734, http://www.math.harvard.edu/$\sim$kronheim.

\bibitem{LM}
H.~B. Lawson and M-L. Michelsohn, {\em Spin geometry}, Princeton Univ. Press,
  Princeton, NJ, 1988.

\bibitem{LenessBlowUp}
T.~G. Leness, {\em Blow-up formulae for {SO(3)}-{D}onaldson polynomials}, Math.
  Z. {\bf 227} (1998), 1--26.

\bibitem{LenessWC}
\bysame, {\em Donaldson wall-crossing formulas via topology}, Forum Math. {\bf
  11} (1999), 417--457, \mbox{arXiv:dg-ga/9603016}.

\bibitem{LiTian}
J.~Li and G.~Tian, {\em Virtual moduli cycles and {G}romov-{W}itten invariants
  of general symplectic manifolds}, Topics in symplectic $4$-manifolds (Irvine,
  CA, 1996), Internat. Press, Cambridge, MA, 1998, 
  pp.~47--83, \mbox{arXiv:alg-geom/9608032}.

\bibitem{LiLiu}
T-J. Li and A-K. Liu, {\em General wall-crossing formula}, Math. Res. Lett.
  {\bf 2} (1995), 797--810.

\bibitem{MilnorStasheff}
J.~W. Milnor and J.~D. Stasheff, {\em Characteristic classes}, Princeton Univ.
  Press, Princeton, NJ, 1974.

\bibitem{MooreWitten}
G.~Moore and E.~Witten, {\em Integration over the $u$-plane in {D}onaldson
  theory}, Adv. Theor. Math. Phys. {\bf 1} (1997), 298--387, 
\mbox{arXiv:hep-th/9709193}.

\bibitem{MorganSWNotes}
J.~W. Morgan, {\em The {S}eiberg-{W}itten equations and applications to the
  topology of smooth four-manifolds}, Princeton Univ. Press, Princeton, NJ,
  1996.

\bibitem{MrowkaPrincetonMorseTalk}
T.~S. Mrowka, {\em Marston {M}orse memorial lectures}, Institute for Advanced
  Study, Princeton, NJ, April, 1999.

\bibitem{OTWall}
C.~Okonek and A.~Teleman, {\em Seiberg-{W}itten invariants for manifolds with
  {$b^+=1$} and the universal wall crossing formula}, Internat. J. Math. {\bf
  7} (1996), 811--832, \mbox{arXiv:alg-geom/9603003}.

\bibitem{OzsvathBlowUp}
P.~S. Ozsv{\'a}th, {\em Some blowup formulas for {SU(2)} {D}onaldson
  polynomials}, J. Differential Geom. {\bf 40} (1994), 411--447.

\bibitem{OzsvathSzaboAdjunct}
P.~S. Ozsv{\'a}th and Z.~Szab{\'o}, {\em Higher type adjunction formulas in
  {S}eiberg-{W}itten theory}, \mbox{arXiv:math.DG/0005268}.

\bibitem{PTLectures}
V.~Y. Pidstrigatch, Lectures at the {N}ewton {I}nstitute in {D}ecember 1994,
  {O}berwolfach in {M}ay 1996, and the {N}ewton {I}nstitute in {N}ovember 1996.

\bibitem{PTLocal}
V.~Y. Pidstrigatch and A.~N. Tyurin, {\em Localisation of {D}onaldson
  invariants along the {S}eiberg-{W}itten classes}, 
\mbox{arXiv:dg-ga/9507004}.

\bibitem{RuanSW}
Y.~Ruan, {\em Virtual neighborhoods and the monopole equations}, Topics in
  symplectic $4$-manifolds (Irvine, CA, 1996), Internat. Press, Cambridge, MA,
  1998, pp.~101--116, \mbox{arXiv:alg-geom/9611021}.

\bibitem{RuanTian}
Y.~Ruan and G.~Tian, {\em A mathematical theory of quantum cohomology}, J.
  Differential Geom. {\bf 42} (1995), 259--367.

\bibitem{SalamonSWBook}
D.~Salamon, {\em Spin geometry and {S}eiberg-{W}itten invariants},
  Birkh{\"a}user, Boston, to appear.

\bibitem{SiebertGW}
B.~Siebert, {\em Symplectic {G}romov-{W}itten invariants}, New trends in
  algebraic geometry (Warwick, 1996), Cambridge Univ. Press, Cambridge, 1999,
  pp.~375--424, \mbox{arXiv:alg-geom/9608005}.

\bibitem{Spanier}
E.~H. Spanier, {\em Algebraic topology}, Springer, New York, 1966.

\bibitem{TauFrame}
C.~H. Taubes, {\em A framework for {M}orse theory for the {Y}ang-{M}ills
  functional}, Invent. Math. {\bf 94} (1988), 327--402.

\bibitem{TelemanGenericMetric}
A.~Teleman, {\em Moduli spaces of {PU(2)}-monopoles}, Asian J. Math. {\bf 4}
  (2000), 391--435, \mbox{arXiv:math.DG/9906163}.

\bibitem{Witten}
E.~Witten, {\em Monopoles and four-manifolds}, Math. Res. Lett. {\bf 1} (1994),
  769--796, \mbox{arXiv:hep-th/9411102}.

\bibitem{Yang}
H-J. Yang, {\em Transition functions and a blow-up formula for {D}onaldson
  polynomials}, {Ph.D}. thesis, Columbia University, 1992.

\end{thebibliography}
\end{document}